# W.B.VASANTHA KANDASAMY

# BIALGEBRAIC STRUCTURES AND SMARANDACHE BIALGEBRAIC STRUCTURES

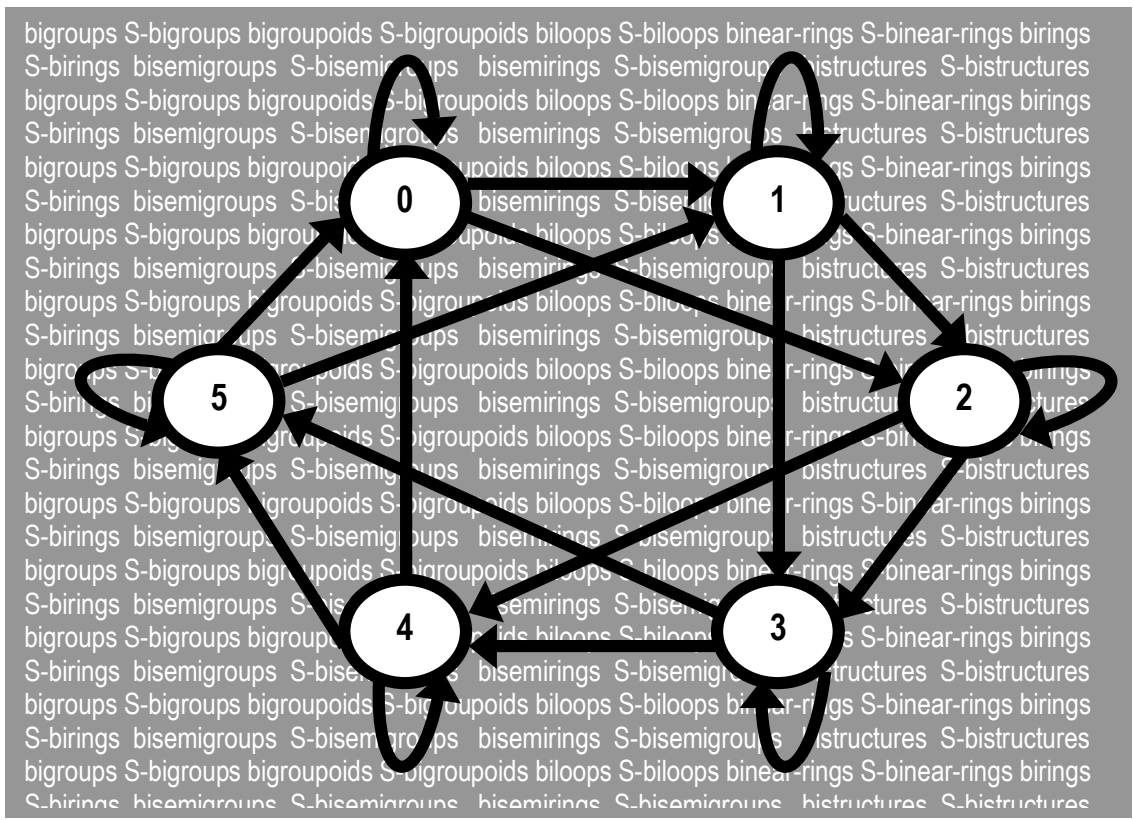

2003

# Bialgebraic Structures and Smarandache Bialgebraic Structures

W. B. Vasantha Kandasamy
Department of Mathematics
Indian Institute of Technology, Madras
Chennai – 600036, India
e-mail: *vasantha@iitm.ac.in*
web: *http://mat.iitm.ac.in/~wbv*

2003

# CONTENTS

**Preface**

**Chapter One: INTRODUCTORY CONCEPTS**

**Chapter Two: BIGROUPS AND SMARANDACHE BIGROUPS**

**Chapter Three: BISEMIGROUPS AND SMARANDACHE BISEMIGROUPS**

**Chapter Four: BILOOPS AND SMARANDACHE BILOOPS**

**Chapter Five: BIGROUPOIDS AND SMARANDACHE BIGROUPOIDS**

**Chapter Six: BIRINGS AND SMARANDACHE BIRINGS**

# Preface

The study of bialgebraic structures started very recently. Till date there are no books solely dealing with bistructures. The study of bigroups was carried out in 1994-1996. Further research on bigroups and fuzzy bigroups was published in 1998. In the year 1999, bivector spaces was introduced. In 2001, concept of free De Morgan bisemigroups and bisemilattices was studied. It is said by Zoltan Esik that these bialgebraic structures like bigroupoids, bisemigroups, binear rings help in the construction of finite machines or finite automaton and semi automaton. The notion of non-associative bialgebraic structures was first introduced in the year 2002. The concept of bialgebraic structures which we define and study are slightly different from the bistructures using category theory of Girard's classical linear logic. We do not approach the bialgebraic structures using category theory or linear logic.

We can broadly classify the study under four heads :

  i.   bialgebraic structures with one binary closed associative operation : *bigroups and bisemigroups*
  ii.  bialgebraic structures with one binary closed non-associative operation: *biloops and bigroupoids*
  iii. bialgebraic structures with two binary operations defined on the biset with both closure and associativity: *birings, binear-rings, bisemirings and biseminear-rings*. If one of the binary operation is non-associative leading to the concept of *non-associative biring, binear-rings, biseminear-rings and bisemirings*.
  iv.  Finally we construct bialgebraic structures using bivector spaces where a bigroup and a field are used simultaneously.

The chief aim of this book is to give Smarandache analogous to all these notions for Smarandache concepts finds themselves accommodative in a better analysis by dissecting the whole structures into specified smaller structure with all properties preserved in them. Such sort of study is in a way not only comprehensive but also more utilitarian and purpose serving. Sometimes several subsets will simultaneously enjoy the same property like in case of defining Smarandache automaton and semi-automaton where, in a single piece of machine, several types of submachines can be made present in them -- thereby making the operation economical and time-saving.

Bistructures are a very nice tool as this answers a major problem faced by all algebraic structures – groups, semigroups, loops, groupoids etc. that is the union of two subgroups, or two subrings, or two subsemigroups etc. do not form any algebraic structure but all of them find a nice bialgebraic structure as bigroups, birings, bisemigroups etc. Except for this bialgebraic structure these would remain only as sets without any nice algebraic structure on them. Further when these bialgebraic structures are defined on them they enjoy not only the inherited qualities of the algebraic structure from which they are taken but also several distinct algebraic properties that are not present in algebraic structures. One such is the reducibility or the irreducibility of a polynomial, or we can say in some cases a polynomial is such that it cannot be reducible or irreducible. Likewise, we see in case of groups an element can be a Cauchy element or a non-Cauchy element or neither.

This book has ten chapters. The first chapter is unusually long for it introduces all concepts of Smarandache notions on rings, groups, loops, groupoids, semigroup, semirings, near ring, vector spaces and their non-associative analogues. The second chapter is devoted to the introduction of bigroups and Smarandache bigroups. The notion of Smarandache bigroups is very new. The introduction of bisemigroups and Smarandache bisemigroups is carried out in chapter three. Here again a new notion called biquasi groups is also introduced. Biloops and Smarandache biloops are introduced and studied in chapter four. In chapter five we define and study the bigroupoid and Smarandache bigroupoid. Its application to Smarandache automaton is also introduced. Chapter six is devoted to the introduction and study of birings and Smarandache birings both associative and non-associative. Several marked differences between birings and rings are brought out. In chapter seven we introduce bisemirings, Smarandache bisemirings, bisemivector spaces, Smarandache bisemivector spaces. Binear rings and Smarandache binear rings are introduced in chapter eight. Chapter nine is devoted to the new notion of bistructures and bivector spaces and their Smarandache analogue. Around 178 problems are suggested for any researcher in chapter ten. Each chapter has an introduction, which brings out clearly what is dealt in that chapter. It is noteworthy to mention in conclusion that this book totally deals with 460 Smarandache algebraic concepts.

I deeply acknowledge the encouragement given by Dr. Minh Perez, editor of the *Smarandache Notions Journal* for writing this book-series. As an algebraist, for the past one-year or so, I have only been involved in the study of the revolutionary and fascinating Smarandache Notions, and I owe my thanks to Dr. Perez for all the intellectual delight and research productivity I have experienced in this span of time.

I also thank my daughters Meena and Kama and my husband Dr. Kandasamy without whose combined help this book would have been impossible. Despite having hundreds of mathematicians as friends, researchers and students I have not sought even a single small help from any of them in the preparation of this book series. I have been overwhelmingly busy because of this self-sufficiency – juggling my teaching and research schedules and having my family working along with me late hours every night in order to complete this – but then, I have had a rare kind of intellectual satisfaction and pleasure.

I humbly dedicate this book to Dr. Babasaheb Ambedkar (1891-1956), the unparalleled leader of India's two hundred million dalits. His life was a saga of struggle against casteist exploitation. In a land where the laws decreed that the "low" caste untouchables must not have access to education, Dr. Ambedkar shocked the system by securing the highest academic honours from the most prestigious universities of the world. After India's independence, he went on to frame the Constitution of India (the longest in the world) – making laws in a country whose bigoted traditional laws were used to stifle the subaltern masses. His motto for emancipation and liberation was "Educate. Organize. Agitate." Education – the first aspect through which Dr. Ambedkar emphasized the key to our improvement – has become the arena where we are breaking the barriers. Through all my years of fighting against prejudice and discrimination, I have always looked up to his life for getting the courage, confidence and motivation to rally on and carry forward the collective struggles.

**Chapter 1**

# INTRODUCTORY CONCEPTS

This chapter has seven sections. The main aim of this chapter is to introduce several of the Smarandache concepts used in this book. If these Smarandache concepts are not introduced, the reader may find it difficult to follow when the corresponding Smarandache bistructures are given. So we have tried to be very brief, only the main definitions and very important results are given. In the first section we just recall definition of groups, loops and S-loops. Section two is devoted to the recollection of notions about semigroups and more about S-semigroups. In section three we introduce the concepts of groupoids and S-groupoids. Section four covers the notions about both rings and non-associative rings and mainly their Smarandache analogue.

In the fifth section we give the notions of semirings and Smarandache semirings. Also in this section we give the concepts of semivector spaces and their Smarandache analogue. In section six concepts on near-rings and Smarandache near-rings are given to the possible extent as the very notion of binear-rings and Smarandache binear-rings are very new. In the final section we give the notions of vector spaces and Smarandache vector spaces.

## 1.1 Groups, loops and S-loops

In this section we just recall the definitions of groups, loops and Smarandache loops (S-loops) for the sake of completeness and also for our notational convenience, as we would be using these notions and notations in the rest of the book. Also we will recall some of the very basic results, which we feel is very essential for our study and future reference.

**DEFINITION 1.1.1:** *A non-empty set G, is said to form a group if in G there is defined a binary operation, called the product and denoted by '•' such that*

   i. *a, b $\in$ G implies a • b $\in$ G.*
   ii. *a, b, c $\in$ G implies (a • b) • c = a • (b • c).*
   iii. *There exists an element e $\in$ G such that a • e = e • a = a for all a $\in$ G.*
   iv. *For every a $\in$ G there exists an element $a^{-1}$ $\in$ G such that a • $a^{-1}$ = $a^{-1}$ • a = e.*

*If a • b = b • a for all a, b $\in$ G we say G is a abelian or a commutative group. If the group G has only a finite number of elements we call G a group of finite order otherwise G is said to be of infinite order. If a • b $\neq$ b • a, for atleast a pair of elements a, b $\in$ G, then we call G a non-commutative group.*

*Notation:* Let X = {$x_1, x_2, \ldots, x_n$}, the set of all one to one mappings of the set X to itself under the product called composition of mappings. Then this is a group. We denote this by $S_n$ called the symmetric group of degree n. We will adhere to this notation and the order of $S_n$ is n!. $D_{2n}$ will denote the dihedral group of order 2n. That is $D_{2n}$ = {a, b | $a^2$ = $b^n$ = 1, bab = a} = {1, a, b, $b^2$, …, $b^{n-1}$, ab, $ab^2$, …, $ab^{n-1}$}. |$D_{2n}$| =

2n. $G = \langle g \mid g^n = 1 \rangle$ is the cyclic group of order n; i.e. $G = \{1, g, g^2, \ldots, g^{n-1}\}$. $A_n$ will denote the alternating subgroup of the symmetric group $S_n$ and

$$|A_n| = \frac{|S_n|}{2} = \frac{n!}{2}.$$

We call a proper subset H of a group G to be a subgroup if H itself is a group under the operations of G. The following classical theorems on group are just recalled.

**LAGRANGE THEOREM:** *If G is a finite group and H is a subgroup of G then o(H) is a divisor of o(G).*

**Note**: o(G) means the number of elements in G it will also be denoted by |G|.

**CAUCHY THEOREM (FOR ABELIAN GROUPS):** *Suppose G is a finite abelian group and p/o(G) where p is a prime number, then there is an element $a \neq e \in G$ such that $a^p = e$.*

**SYLOW'S THEOREM (FOR ABELIAN GROUPS):** *If G is an abelian group of finite order and if p is a prime number, such that $p^\alpha / o(G)$, $p^{\alpha+1} \nmid o(G)$ then G has a subgroup of order $p^\alpha$.*

**CAYLEY'S THEOREM:** *Every group is isomorphic to a subgroup of $S_n$ for some appropriate n.*

**CAUCHY THEOREM:** *If p is a prime number and p / o(G), then G has an element of order p.*

For more results about group theory please refer [23 & 27]. Now we proceed on to recall some basic concepts on loops, a new class of loops using $Z_n$, n prime and n > 3 and about identities on loops and several other properties about them.

**DEFINITION 1.1.2:** *A non-empty set L is said to form a loop if in L is defined a binary operation called product and denoted by '•' such that*

  i. *for all a, b ∈ L we have a • b ∈ L.*
  ii. *there exists an element e ∈ L such that a • e = e • a = a for all a ∈ L.*
  iii. *for every ordered pair (a, b) ∈ L × L there exists a unique pair (x, y) ∈ L such that a • x = b and y • a = b, '•' defined on L need not always be associative.*

*Example 1.1.1*: Let L be a loop given by the following table:

| * | e | a | b | c | d |
|---|---|---|---|---|---|
| e | e | a | b | c | d |
| a | a | e | c | d | b |
| b | b | d | a | e | c |
| c | c | b | d | a | e |
| d | d | c | e | b | a |

Clearly (L, •) is non-associative with respect to '•'. It is important to note that all groups are in general loops but loops in general are not groups.

**MOUFANG LOOP:** *A loop L is said to be a Moufang loop if it satisfies any one of the following identities:*

    i.    *(xy) (zx) = (x(yz))x.*
    ii.   *((xy)z) y = x(y (zy)).*
    iii.  *x (y (xz))= ((xy)x)z*

*for all x, y, z ∈ L.*

**BRUCK LOOP:** *Let L be a loop, L is called a Bruck loop if (x (yx)) z = x (y (xz)) and $(xy)^{-1} = x^{-1} y^{-1}$ for all x, y, z ∈ L.*

**BOL LOOP:** *A loop L is called a Bol loop if ((xy)z) y = x((yz)y) for all x, y, z ∈ L.*

**ALTERNATIVE LOOP:** *A loop L is said to be right alternative if (xy) y = x(yy) for all x, y ∈ L and left alternative if (x x) y = x (xy) for all x, y ∈ L, L is said to be alternative if it is both a right and a left alternative loop.*

**WEAK INVERSE PROPERTY LOOP:** *A loop L is called a weak inverse property loop if (xy) z = e imply x(yz) = e for all x, y, z ∈ L; e is the identity element of L.*

**DEFINITION 1.1.3:** *Let L be a loop. A non-empty subset H of L is called a subloop of L if H itself is a loop under the operations of L. A subloop H of L is said to be a normal subloop of L if*

    i.    *xH = Hx.*
    ii.   *(Hx) y = H (xy).*
    iii.  *y(xH) = (yx)H,*

*for all x, y ∈ L. A loop is simple if it does not contain any non-trivial normal subloop.*

**DEFINITION 1.1.4:** *If x and y are elements of a loop L, the commutator (x, y) is defined by xy = (yx) (x, y). The commutator subloop of a loop L denoted by L' is the subloop generated by all of its commutators that is ⟨ {x ∈ L / x = (y, z) for some y, z ∈ L }⟩ where for A ⊂ L, ⟨A⟩ denotes the subloop generated by A.*

**DEFINITION 1.1.5:** *If x, y, z are elements of a loop L an associator (x, y, z) is defined by (xy)z = (x(yz)) (x, y, z). The associator subloop of a loop L denoted by A(L) is the subloop generated by all of its associators, that is A(L) = ⟨{x ∈ L / x = (a, b, c ) for some a, b, c ∈ L }⟩.*

**SEMIALTERNATIVE LOOP:** *A loop L is said to be semialternative if (x, y, z) = (y, z, x) for all x, y, z ∈ L where (x, y, z) denotes the associator of elements x, y, z ∈ L.*

**DEFINITION 1.1.6:** *Let L be a loop the left nucleus $N_\lambda = \{a \in L \,/\, (a, x, y) = e$ for all for all $x, y \in L\}$ is a subloop of L. The right nucleus $N_\rho = \{a \in L \mid (x, y, a) = e$ for all $x, y \in L\}$ is a subloop of L. The middle nucleus $N_\mu = \{a \in L \,/\, (x, a, y) = e$ for all $x, y \in L\}$ is a subloop of L. The nucleus N(L) of the loop L is the subloop given by $N(L) = N_\lambda \cap N_\mu \cap N_\rho$.*

**THE MOUFANG CENTER:** *The Moufang center C(L) is the set of elements of the loop L which commute with every element of L i.e. $C(L) = \{x \in L \,/\, xy = yx$ for all $y \in L\}$. The center Z(L) of a loop L is the intersection of the nucleus and the Moufang center i.e. $Z(L) = C(L) \cap N(L)$.*

**DEFINITION 1.1.7:** *Let $L_n(m) = \{e, 1, 2, 3, ..., n\}$ be a set where $n > 3$, n is odd and m a positive integer such that $(m, n) = 1$ and $(m - 1, n) = 1$ with $m < n$. Define on $L_n(m)$ a binary operation '•' as follows:*

    i.    $e \bullet i = i \bullet e = i$ *for all* $i \in L_n(m)$.
    ii.    $i \bullet i = i^2 = e$ *for all* $i \in L_n(m)$.
    iii.    $i \bullet j = t$ *where* $t = [mj - (m-1)i] \pmod{n}$

*for all $i, j \in L_n(m)$. $i \neq j$, $i \neq e$ and $j \neq e$.*

*Then $L_n(m)$ is a loop.*

***Example 1.1.2:*** Consider $L_5(2) = \{e, 1, 2, 3, 4, 5\}$.

The composition table for $L_5(2)$ is given below.

| * | e | 1 | 2 | 3 | 4 | 5 |
|---|---|---|---|---|---|---|
| e | e | 1 | 2 | 3 | 4 | 5 |
| 1 | 1 | e | 3 | 5 | 2 | 4 |
| 2 | 2 | 5 | e | 4 | 1 | 3 |
| 3 | 3 | 4 | 1 | e | 5 | 2 |
| 4 | 4 | 3 | 5 | 2 | e | 1 |
| 5 | 5 | 2 | 4 | 1 | 3 | e |

This loop is of order 6, which is non-associative and non-commutative.

$L_n$ denotes the class of all loops $L_n(m)$ for a fixed n and varying m's satisfying the conditions $m < n$, $(m, n) = 1$ and $(m - 1, n) = 1$, that is $L_n = \{L_n(m) \mid n > 3$, n odd, $m < n$, $(m, n) = 1$ and $(m - 1, n) = 1\}$. This class of loops will be known as the new class of loops and all these loops are of even order. Several nice properties are enjoyed by these loops, which is left for the reader to discover.

Now we proceed on to recall the definition of S-loops.

**DEFINITION 1.1.8:** *The Smarandache loop (S-loop) is defined to be a loop L such that a proper subset A of L is a subgroup with respect to the operations of L that is $\phi \neq A \subset L$.*

*Example 1.1.3:* Let L be a loop given by the following table. L is a S-loop as every pair $A_i = \{e, a_i\}$; $i = \{1, 2, 3, 4, 5, 6, 7\}$ are subgroups of L.

| * | e | $a_1$ | $a_2$ | $a_3$ | $a_4$ | $a_5$ | $a_6$ | $a_7$ |
|---|---|---|---|---|---|---|---|---|
| e | e | $a_1$ | $a_2$ | $a_3$ | $a_4$ | $a_5$ | $a_6$ | $a_7$ |
| $a_1$ | $a_1$ | e | $a_5$ | $a_2$ | $a_6$ | $a_3$ | $a_7$ | $a_4$ |
| $a_2$ | $a_2$ | $a_5$ | e | $a_6$ | $a_3$ | $a_7$ | $a_4$ | $a_1$ |
| $a_3$ | $a_3$ | $a_2$ | $a_6$ | e | $a_7$ | $a_4$ | $a_1$ | $a_5$ |
| $a_4$ | $a_4$ | $a_6$ | $a_3$ | $a_7$ | e | $a_1$ | $a_5$ | $a_2$ |
| $a_5$ | $a_5$ | $a_3$ | $a_7$ | $a_4$ | $a_1$ | e | $a_2$ | $a_6$ |
| $a_6$ | $a_6$ | $a_7$ | $a_4$ | $a_1$ | $a_5$ | $a_2$ | e | $a_3$ |
| $a_7$ | $a_7$ | $a_4$ | $a_1$ | $a_5$ | $a_2$ | $a_6$ | $a_3$ | e |

**DEFINITION 1.1.9:** *The Smarandache Bol loop (S-Bol loop) L, is defined to be a S-loop L such that a proper subset A, $A \subset L$, which is a subloop of L (A not a subgroup of L) is a Bol loop.*

Similarly we define S-Bruck loop, S-Moufang loop, S-right(left) alternative loop. Clearly by this definition we may not have every S-loop to be automatically a S-Bol loop or S-Moufang loop or S-Bruck loop or so on.

**THEOREM 1.1.1:** *Every Bol loop is a S-Bol loop but not conversely.*

*Proof*: Left as an exercise for the reader to prove.

The same result holds good in case of Moufang, Bruck and alternative loops.

**DEFINITION 1.1.10:** *Let L and L' be two S-loops with A and A' its subgroups respectively. A map $\phi$ from L to L' is called S-loop homomorphism if $\phi$ restricted to A is mapped to a subgroup A' of L', i.e. $\phi: A \to A'$ is a group homomorphism. It is not essential that $\phi$ be defined on whole of A.*

*The concept of S-loop isomorphism and S-loop automorphism are defined in a similar way.*

Several properties about loops can be had from [5, 6] and that of S-loops can be had from [115, 116].

## 1.2 Semigroups and S-semigroups

In this section we introduce the notions of semigroups and S-semigroups. These notions will be used in building bigroups and other bi structures. As these notions are very recent [121], we felt it essential to introduce these concepts.

Several of their properties given about S-semigroups will be used for the study of further chapters.

**DEFINITION 1.2.1:** *Let S be a non-empty set on which is defined a binary operation '•', (S, •) is a semigroup.*

   i. *If for all a, b ∈ S we have a • b ∈ S.*
   ii. *a • (b • c) = (a • b) • c for all a, b, c ∈ S.*

*A semigroup in which a • b = b • a for all a, b ∈ S, then we call S a commutative semigroup. If S has a unique element e ∈ S such that a • e = e • a = a for all a, b ∈ S then we call the semigroup a monoid or a semigroup with unit. If the number of elements in S is finite we say the semigroup S is of finite order and denote the order of S by o(S) or |S|. If the number of elements in S is not finite i.e. infinite we say S is of infinite order.*

**DEFINITION 1.2.2:** *Let (S, •) be a semigroup, a proper subset P of S is said to be a subsemigroup of S if (P, •) is a semigroup.*

**DEFINITION 1.2.3:** *Let (S, •) be a semigroup. P a proper subset of S. P is called a (right) left ideal of S if (ar), ra ∈ P for all r ∈ S and a ∈ P. P is said to be a two sided ideal if P is simultaneously a left and a right ideal of S.*

The concept of maximal ideal, principal ideal and prime ideal can be had from any textbook on algebra [22, 26, 27, 30, 31, 32].

*Notation:* Let $X = \{a_1, \ldots, a_n\}$ where $a_i$'s are distinct, that is $|X| = n$. Let S(X) denote the set of all mappings of the set X to itself. Then S(X) under the operations of composition of mappings is a semigroup.

This semigroup will be addressed in this book as the symmetric semigroup on n elements. Clearly $o(S(X)) = n^n$ and S(X) is a non-commutative moniod. Further S(X) contains $S_n$ the symmetric group of degree n as a proper subset.

**DEFINITION 1.2.4:** *A Smarandache semigroup (S-semigroup) is defined to be a semigroup A such that a proper subset X of A is a group with respect to the same binary operation on A. (X ≠ ϕ and X ≠ A but X ⊂ A).*

*Example 1.2.1*: Let S(5) be the symmetric semigroup. S(5) is a S-semigroup.

*Example 1.2.2*: Let $Z_{10} = \{0, 1, 2, \ldots, 9\}$ be a semigroup under multiplication modulo 10. Clearly $Z_{10}$ is a S-semigroup for X = {2, 4, 6, 8} is a group.

**DEFINITION 1.2.5:** *Let S be a S-semigroup. If every proper subset A of S which is a group is commutative, then we say the S-semigroup is a Smarandache commutative semigroup (S-commutative semigroup). If atleast one of them is commutative then we say the S-semigroup is Smarandache weakly commutative (S-weakly commutative). If every proper subset which is subgroup is cyclic then we call S a Smarandache cyclic semigroup (S-cyclic semigroup). If atleast one of the subgroup is cyclic we call the semigroup Smarandache weakly cyclic (S-weakly cyclic).*

Several interesting results can be obtained in this direction but the reader is requested to refer [121].

**DEFINITION 1.2.6:** *Let S be a S-semigroup. A proper subset A of S is said to be a Smarandache subsemigroup (S-subsemigroup) of S if A itself is a S-semigroup, that is A has a proper subset B (B $\subset$ A) such that B is a group under the operations of S.*

Several nice characterizations theorems can be had from [121].

**DEFINITION 1.2.7:** *Let S be a S-semigroup. If A $\subset$ S is a proper subset of S and A is a subgroup which cannot be contained in any other proper subsemigroup of S then we say A is the largest subgroup of S.*

*Suppose A is the largest subgroup of S and if A is contained in a proper subsemigroup X of S then we call X the Smarandache hyper subsemigroup (S-hyper subsemigroup) of S. Thus we say the S-semigroup S is Smarandache simple (S-simple) if S has no S-hyper subsemigroup.*

It is interesting to note that $Z_{19}$ = {0, 1, 2, 3, …, 18} is the semigroup under product modulo 19; we see $|Z_{19}|$ = 19. Take the set X = {1, 18}. Clearly X is subgroup as well as a subsemigroup with o(X) = 2; we see 2 $\nmid$ 19.

So we are interested in introducing the concept of Smarandache Lagrange theorem.

**DEFINITION 1.2.8:** *S be a finite S-semigroup. If the order of every proper subset, which is a subgroup of S, divides the order of the S-semigroup S then we say S is a Smarandache Lagrange semigroup (S-Lagrange semigroup). If there exists atleast one subgroup in S which divides the order of the S-semigroup we call S a Smarandache weakly Lagrange semigroup (S-weakly Lagrange semigroup). If the order of no subgroup divides the order of the S-semigroup S then we call S a Smarandache non-Lagrange semigroup (S-non-Lagrange semigroup). The semigroups $Z_p$ (p a prime) and $Z_p$ the semigroup under multiplication modulo p falls under the S-non-Lagrangian semigroups.*

The concepts of S-p-Sylow semigroups Cauchy elements and several other interesting results can be found in the book [121], the reader is expected to refer it for more information.

## 1.3 Groupoids and S-groupoids

In this section we give a brief sketch of the results about groupoids and S-groupoids. Further this section also recalls the new classes of groupoids built using the ring of integers $Z_n$, and Z and Q the field of rationals. We give some important results about S-groupoids and groupoids and S-groupoids, which satisfy special identities like Bol, Bruck, Moufang and alternative.

For more about groupoids and S-groupoids the reader is requested to refer the book [111, 114].

**DEFINITION 1.3.1:** *Given an arbitrary set P a mapping of P $\times$ P into P is called a binary operation on P. Given such a mapping $\sigma$: P $\times$ P $\rightarrow$ P we use it to define a* *product '∗' in P by declaring a ∗ b = c if σ (a, b) = c or equivalently we can define it with more algebraic flavor as:*

*A non-empty set of elements G is said to form a groupoid if in G is defined a binary operation called the product denoted by '∗' such that a ∗ b ∈ G for all a, b ∈ G.*

*It is important to note that the binary operation '∗' defined on G need not be in general be associative i.e. a ∗ (b ∗ c) ≠ (a ∗ b) ∗ c for a, b, c ∈ G: so we can roughly say the groupoid (G, ∗) is a set on which is defined a non-associative binary operation which is closed on G.*

**DEFINITION 1.3.2:** *A groupoid (G, ∗) is said to be a commutative groupoid if for every a, b ∈ G we have a ∗ b = b ∗ a. A groupoid G is said to have an identity element e in G if a ∗ e = e ∗ a = a for all a ∈ G.*

*We call the order of the groupoid G to be the number of elements in G and we denote it by o(G) or |G|. If o(G) = n, n < ∞ we say G is a finite groupoid other wise G is said to be an infinite groupoid.*

**DEFINITION 1.3.3:** *Let (G, ∗) be a groupoid. A proper subset H ⊂ G is a subgroupoid, if (H, ∗) is itself a groupoid. All semigroups are groupoids; so groupoids form a most generalized class of semigroups.*

**DEFINITION 1.3.4:** *A groupoid G is said to be a Moufang groupoid if it satisfies the Moufang identity (xy) (zx) = (x (yz)) x for all x, y, z ∈ G.*

*A groupoid G is said to be a Bol Groupoid if G satisfies the Bol identity; ((xy) z)y = x((yz)y) for all x, y, z ∈ G.*

*A groupoid is said to be a P-groupoid if (xy) x = x (yx) for all x, y ∈ G.*

*A groupoid G is said to be right alternative if it satisfies the identity (xy) y = x (yy) for all x, y ∈ G and G is said to be left alternative if (xx)y = x(xy) for all x, y ∈ G. A groupoid is alternative if it is both right and left alternative simultaneously.*

Several properties about these groupoids with examples can be had from the two textbooks [5, 114].

**DEFINITION 1.3.5:** *Let (G, ∗) be a groupoid. A proper subset H of G is said to be subgroupoid of G if (H, ∗) is itself a groupoid. A non-empty subset P of the groupoid G is said to be a left ideal of the groupoid G if*

      i. *P is a subgroupoid.*
     ii. *For all x ∈ G and a ∈ P.*
*x ∗ a ∈ P.*

*One can similarly define right ideal of a groupoid G. We say P is an ideal of the groupoid G if P is simultaneously a left and a right ideal of G.*

**DEFINITION 1.3.6:** *Let G be a groupoid. A subgroupoid V of G is said to be a normal subgroupoid of G if*

i. *aV = Va.*
ii. *(Vx) y = V(xy).*
iii. *y (xV) = (yx) V*

*for all a, x, y ∈ G. A groupoid is said to be simple if it has no nontrivial normal subgroupoids.*

Now when do we call a groupoid G itself normal.

**DEFINITION 1.3.7:** *A groupoid G is normal if*

i. *xG = Gx.*
ii. *G(xy) = (Gx) y.*
iii. *y(xG) = (yx)G*

*for all x, y ∈ G.*

**DEFINITION 1.3.8:** *Let G be a groupoid, H and K be two proper subgroupoids of G with H ∩ K = ϕ; we say H is conjugate with K if there exists a x ∈ H such that H = xK or Kx (or in the mutually exclusive sense).*

**DEFINITION 1.3.9:** *Let $(G_1, \theta_1), (G_2, \theta_2), \ldots, (G_n, \theta_n)$ be n groupoids with $\theta_i$ binary operations defined on each $G_i$, i = 1, 2, ..., n. The direct product of $G_1, \ldots, G_n$ denoted by $G = G_1 \times \ldots \times G_n = \{g_1, \ldots, g_n) \mid g_i \in G_i\}$ by component wise multiplication of $G_i$; G becomes a groupoid. For if $g = (g_1, \ldots, g_n)$ and $h = (h_1, \ldots, g_n)$ then $g * h = \{(g_1 \theta_1 h_1, g_2 \theta_2 h_2, \ldots, g_n \theta_n h_n)\}$. Clearly $g * h$ belongs to G; so G is a groupoid.*

Unlike in groups in groupoids we can have either left or right identity; we define them.

**DEFINITION 1.3.10:** *Let (G, •) be a groupoid we say an element e ∈ G is a left identity if e • a = a for all a ∈ G. Similarly right identity of the groupoid can be defined; if e ∈ G happens to be simultaneously both right and left identity we say the groupoid G has an identity. Similarly we say an element a ∈ G is a right zero divisor if a • b = 0 for some b ≠ 0 in G and $a_1$ in G has left zero divisor if $b_1 • a_1 = 0$. We say G has a zero divisor if a • b = 0 and b • a = 0 for a, b ∈ G \ {0}.*

**DEFINITION 1.3.11:** *Let G be a groupoid, the center of the groupoid G is C(G) = {a ∈ G | ax = xa for all x ∈ G}.*

**DEFINITION 1.3.12:** *Let (G, •) be a groupoid of order n. We say b, a ∈ G is a conjugate pair if a = b • x (or xb for some x ∈ G) and b = a • y (or ya for some y ∈ G). An element a in G said to be right conjugate with b in G if we can find x, y ∈ G such that a • x = b and b • y = a (x • a = b and y • b = a).*

Similarly we define left conjugate. It is a very well known fact that we do not have many natural examples of groupoids; here we define four new classes of groupoids built using $Z_n$ the set of integers addition and multiplication modulo n.

**DEFINITION 1.3.13:** *Let $Z_n = \{0, 1, 2, ..., n–1\}$, $n \geq 3$. For $a, b \in Z_n$ define a binary operations '*' on $Z_n$ as follows. $a * b = ta + ub \pmod{n}$ where $t, u$ are 2 distinct elements in $Z_n \setminus \{0\}$ such that $(t, u) = 1$ ; '+' here is the usual addition modulo n. Clearly $\{Z_n, *, (t, u)\}$ is a groupoid. Now for varying $t, u \in Z_n \setminus \{0\}$, $t$ and $u$ distinct such that $(t, u) = 1$, we have a class of groupoids; we denote this class of groupoids by $Z(n)$. $Z(n) = \{ Z_n (t, u), *, (t, u) = 1\}$.*

*Let $Z_n = \{0, 1, 2, ..., n-1\}$, $n \geq 3$, $n < \infty$. Define operation '*' on $Z_n$ by $a * b = ta + bu \pmod{n}$ where $t, u \in Z_n \setminus \{0\}$ ($t$ and $u$ need not always be relatively prime but $t \neq u$). Then $\{Z_n, * (t, u)\}$ is a groupoid we denote this class of groupoids by $Z^*(n) = \{Z_n (t, u), *, (t, u)\}$ thus we have $Z^*(n) \supset Z(n)$. Now using $Z_n$ if we select $t, u \in Z_n \setminus \{0\}$ such that $t$ can also be equal to $u$ then we get yet another new class of groupoids which we denote by $Z^{**}(n)$ thus $Z^{**}(n) = \{Z_n, *, (t, u)\}$ This class of groupoids completely contains the class of groupoids $Z^*(n)$ and $Z(n)$. Thus $Z(n) \subset Z^*(n) \subset Z^{**}(n)$.*

*Now we define yet another new class of groupoids using $Z_n$. We define groupoids using $Z_n$, by for $a, b \in Z_n$ choose any pair of element $(t, u)$ in $Z_n$ and define for $a, b \in Z_n$, $a * b = ta + ub \pmod{n}$. Now we denote this class of groupoids by $Z^{***}(n)$. Clearly $Z(n) \subset Z^*(n) \subset Z^{**}(n) \subset Z^{***}(n)$. Further on $Z$ or $Q$ or $R$ we define '*' by $a * b = ta + bu$. This $(Z, *)$ forms a groupoid of infinite order. Similarly $(Q, *)$ and $(R, *)$. Thus on $Z$ we have infinite number of groupoids of infinite order.*

For more about groupoids please refer [5, 11].

We now proceed on to define Smarandache groupoids.

**DEFINITION 1.3.14:** *A Smarandache groupoid (S-groupoid) G is a groupoid which has a proper subset S, $S \subset G$ such that S under the operations of G is a semigroup. If G is a S-groupoid and if the number of elements in G is finite we say G is finite otherwise G is of infinite order.*

**DEFINITION 1.3.15:** *Let $(G, *)$ be a S-groupoid. A non-empty subset H of G is said to be a Smarandache subgroupoid (S-subgroupoid) if H contains a proper subset $K \subset H$ such that K is a semigroup under the operation '*'.*

**THEOREM 1.3.1:** *Every subgroupoid of a S-groupoid need not in general be a S-subgroupoid of S.*

*Proof*: Left for the reader as an exercise.

**THEOREM 1.3.2:** *Let G be a groupoid having a S-subgroupoid then G is a S-groupoid.*

*Proof***:** Straightforward by the very definition.

**DEFINITION 1.3.16:** *Let G be a S-groupoid if every proper subset of G, which is a semigroup, is commutative then we call G a Smarandache commutative groupoid (S-commutative groupoid). (It is to be noted that G need not be a commutative groupoid, it is sufficient if every subset which is a semigroup is commutative). We say G is a S-weakly commutative groupoid if G has atleast one proper subset which is a semigroup is commutative.*

The following theorem is left as an exercise for the reader to prove.

**THEOREM 1.3.3:** *Every S-commutative groupoid is a S-weakly commutative groupoid.*

**DEFINITION 1.3.17:** *A Smarandache left ideal (S-left ideal) A of the S-groupoid G satisfies the following conditions.*

  i. *A is a S-subgroupoid.*
  ii. *$x \in G$ and $a \in A$ then $x a \in A$.*

*Similarly we can define Smarandache right ideal (S-right ideal). If A is both a S-right ideal and S-left ideal simultaneously then we say A is a Smarandache ideal (S-ideal) of G.*

**DEFINITION. 1.3.18:** *Let G be a S-subgroupoid of G. We say V is a Smarandache seminormal groupoid (S-seminormal groupoid) if*

  i. *$aV = X$ for all $a \in G$,*
  ii. *$Va = Y$ for all $a \in G$,*

*where either X or Y is a S-subgroupoid of G but X and Y are both subgroupoids. V is said to be a Smarandache normal groupoid (S-normal groupoid) if $aV = X$ and $Va = Y$ for all $a \in G$ where both X and Y are S-subgroupoids of G.*

**THEOREM 1.3.4:** *Every S-normal groupoid is a S-seminormal groupoid and not conversely.*

*Proof*: Straightforward hence left for the reader as an exercise.

**DEFINITION 1.3.19:** *Let G be a S-groupoid H and P be two subgroupoids of G. We say H and P are Smarandache semiconjugate subgroupoids (S-semiconjugate subgroupoids) of G if*

  i. *H and P are S-subgroupoids of G.*
  ii. *H = xP or Px or*
  iii. *P = xH or Hx for some $x \in G$.*

*We call two subgroupoids H and P of a groupoid G to be Smarandache conjugate subgroupoids (S-conjugate subgroupoids) of G if*

  i. *H and P are S-subgroupoids of G.*
  ii. *H = xP or Px and*
  iii. *P = xH or Hx.*

The following theorem which directly follows from the very definitions is left as an exercise for the reader.

**THEOREM 1.3.5:** *Let G be a S-groupoid. If P and K are two S-subgroupoids of G, which are S-conjugate, then they are S-semiconjugate and the converse in general is not true.*

**DEFINITION 1.3.20:** *Let G be a S-groupoid. We say G is Smarandache inner commutative (S-inner commutative) if every S-subgroupoid of G is inner commutative.*

Several interesting results can be obtained in this direction connecting commutativity and the inner commutativity.

**DEFINITION 1.3.21:** *Let G be a groupoid, G is said to be a Smarandache Moufang groupoid (S-Moufang groupoid) if there exists $H \subset G$ such that H is a S-subgroupoid of G and (xy) (zx) = (x (yz) ) x for all x, y, z $\in$ H.*

*If every S-subgroupoid of a groupoid G satisfies the Moufang identity then we call G a Smarandache strong Moufang groupoid (S-strong Moufang groupoid).*

On similar lines we define Smarandache Bol groupoid, Smarandache strong Bol groupoid, Smarandache alternative groupoid and Smarandache strong alternative groupoid.

**DEFINITION 1.3.22:** *Let G be a S-groupoid we say G is a Smarandache P-groupoid (S-P-groupoid) if G contains a proper S-subgroupoid A such that $(x*y)*x = x*(y*x)$ for all x, y, $\in$ A. We say G is a Smarandache strong P-groupoid (S-strong P-groupoid) if every S-subgroupoid of G is a S-P- groupoid of G.*

Several interesting results in this direction can be had from [114].

**DEFINITION 1.3.23:** *Let $G_1, G_2, ..., G_n$ be n- groupoids. We say $G = G_1 \times G_2 \times ... \times G_n$ is a Smarandache direct product of groupoids (S-direct product of groupoids) if G has a proper subset H of G which is a semigroup under the operations of G. It is important to note that each $G_i$ need not be a S-groupoid for in this case G will obviously be a S-groupopid.*

**DEFINITION 1.3.24:** *Let $(G_1, \bullet)$ and $(G_2, *)$ be any two S-groupoids. A map $\phi$ from $G_1$ to $G_2$ is said to be a Smarandache homomorphism (S-homomorphism) if $\phi: A_1 \to A_2$ where $A_1 \subset G_1$ and $A_2 \subset G_2$ are semigroups of $G_1$ and $G_2$ respectively that is $\phi (a \bullet b) = \phi (a) * \phi (b)$ for all a, b $\in A_1$.*

*We see that $\phi$ need not be even defined on whole of $G_1$. Further if $\phi$ is 1-1 we call $\phi$ a Smarandache groupoid isomorphism (S-groupoid isomorphism).*

For more about groupoids and S-groupoids the reader is requested to refer [114].

## 1.4 Rings, S-rings and SNA-rings.

In this section we introduce the basic notions of rings especially Smarandache rings (S-rings) as we cannot find much about it in literature. Further we define non-associative rings and Smarandache non-associative rings (SNA-rings). The reader is expected to have a good background of algebra and almost all the algebraic structures thoroughly for her/him to feel at home with this book. For S-rings and SNA-rings please refer [119, 120]. We briefly recollect some results about rings and fields and give more about S-rings and SNA-rings.

**DEFINITION 1.4.1:** *Let $(R, +, \bullet)$ be a non-empty set on which is defined two binary operations '+' and '$\bullet$' satisfying the following conditions:*

   i. *$(R, +)$ is a group under '+'.*
   ii. *$(R, \bullet)$ is a semigroup.*
   iii. *$a \bullet (b + c) = a \bullet b + a \bullet c$ and $(b + a) \bullet c = b \bullet c + a \bullet c$ for all $a, b\ c \in R$.*

*Then we call R a ring. If in R we have $a \bullet b = b \bullet a$ for all $a, b \in R$ then R is said to be a commutative ring. If in particular R contains an element 1 such that $a \bullet 1 = 1 \bullet a = a$ for all $a \in R$, we call R a ring with unit. R is said be a division ring if R is a ring such that R has no non-trivial divisors of zero. R is an integral domain if R is a commutative ring and has no non-trivial divisors of zero. $(R, +, \bullet)$ is said to be a field if $(R \setminus \{0\}, \bullet)$ is a commutative group.*

Q the set of rationals is a field under usual addition and multiplication, R the set of reals is also a field; where as Z the set of integers is an integral domain, $Z_p$, p a prime i.e. $Z_p = \{0, 1, 2, \ldots, p – 1\}$ is also a field. A field F is said to be of characteristic 0 if $nx = 0$ for all $x \in F$ forces $n = 0$ (where n is a positive integer). We say $F_p$ is a field of characteristic p if $px = 0$ for all $x \in F_p$ and p a prime number. $Z_2 = \{0, 1\}$, $Z_3 = \{0, 1, 2\}$ and $Z_7 = \{0, 1, 2, 3, \ldots, 6\}$ are fields of characteristic 2, 3 and 7 respectively. Let F be a field, a proper subset A of F is said to be a subfield, if A itself under the operations of F is a field.

For example in the field of reals R we have Q the field of rationals to be a subfield. If a field has no proper subfields other than itself then we say the field is a prime field. Q is a prime field of characteristic 0 and $Z_{11} = \{0, 1, 2, \ldots, 10\}$ is the prime field of characteristic 11.

**DEFINITION 1.4.2:** *Let $(R, +, \bullet)$ be any ring. S a proper subset of R is said to be a subring of R if $(S, +, \bullet)$ itself is a ring. We say a proper subset I of R is an ideal of R if*

   i. *I is a subring.*
   ii. *$r \bullet i$ and $i \bullet r \in I$ for all $i \in I$ and $r \in R$.*

*(The right ideal and left ideal are defined if we have either $i \bullet r$ or $r \bullet i$ to be in I 'or' in the mutually exclusive sense). An ideal I of R is a maximal ideal of R if J is any other ideal of R and $I \subset J \subset R$ then either $I = J$ or $J = R$. An ideal K of R is said to be a minimal ideal of R if $(0) \subset P \subset K$ then $P = (0)$ or $K = P$. We call an ideal I to be*

*principal if I is generated by a single element. We say an ideal X of R is prime if x • y ∈ X implies x ∈ X or y ∈ X.*

Several other notions can be had from any textbook on ring theory. We just define two types of special rings viz group rings and semigroup rings.

**DEFINITION 1.4.3:** *Let R be a commutative associative ring with 1 or a field and let G be any group. The group ring RG of the group G over the ring R consists of all finite formal sums of the form $\sum_i \alpha_i g_i$ (i-runs over a finite number) where $\alpha_i \in R$ and $g_i \in G$ satisfying the following conditions:*

i. $\sum_{i=1}^{n} \alpha_i g_i = \sum_{i=1}^{n} \beta_i g_i \Leftrightarrow \alpha_i = \beta_i$ *for i = 1,2,…,n.*

ii. $\left(\sum_{i=1}^{n} \alpha_i g_i\right) + \left(\sum_{i=1}^{n} \beta_i g_i\right) = \sum_{i=1}^{n} (\alpha_i + \beta_i) g_i$

iii. $\left(\sum_{i=1}^{n} \alpha_i g_i\right)\left(\sum_{j=1}^{n} \beta_j h_j\right) = \sum \gamma_K m_K$ *where* $\gamma_K = \sum \alpha_i \beta_j$, $m_K = g_i h_j$.

iv. $r_i g_i = g_i r_i$ *for* $r_i \in R$ *and* $g_i \in G$.

v. *Clearly* $1.G \subseteq RG$ *and* $R.1 \subseteq RG$.

**DEFINITION 1.4.4:** *The semigroup ring RS of a semigroup S with unit over the ring R is defined as in case of definition 1.4.3 in which G is replaced by S, the semigroup.*

Now we proceed on to define the concept of loop rings and groupoid rings. These give a class of rings which are non-associative.

**DEFINITION 1.4.5:** *Let R be a commutative ring with 1 or a field and let L be a loop. The loop ring of the loop L over the ring R denoted by RL consists of all finite formal sums of the form $\sum_i \alpha_i m_i$ (i runs over a finite number) where $\alpha_i \in R$ and $m_i \in L$ satisfying the following conditions:*

i. $\sum_{i=1}^{n} \alpha_i m_i = \sum_{i=1}^{n} \beta_i m_i \Leftrightarrow \alpha_i = \beta_i$ *for* $i = 1,2,3, \ldots ,n$.

ii. $\left(\sum_{i=1}^{n} \alpha_i m_i\right) + \left(\sum_{i=1}^{n} \beta_i m_i\right) = \sum_{i=1}^{n} (\alpha_i + \beta_i) m_i$.

iii. $\left(\sum_{i=1}^{n} \alpha_i p_i\right)\left(\sum_{j=1}^{n} \beta_j s_j\right) = \sum \gamma_K m_K$ *where* $m_K = p_i s_j$, $\gamma_K = \sum \alpha_i \beta_j$.

iv.  $r_i m_i = m_i r_i$ for all $r_i \in R$ and $m_i \in L$.

v.  $r \sum_{i=1}^{n} r_i m_i = \sum_{i=1}^{n} (rr_i) m_i$

for all $r \in R$ and $\Sigma r_i m_i \in RL$. RL is a non-associative ring with $0 \in R$ as its additive identity. Since $I \in R$ we have $L = 1.L \subseteq RL$ and $R.e = R \subseteq RL$ where e is the identity element of L.

**Note:** If we replace the loop L by a groupoid with 1 in the defintion 1.4.5 we get groupoid rings, which are groupoid over rings.

This will also form a class of non-associative rings. Now we recall, very special properties in rings.

**DEFINITION 1.4.6:** *Let R be a ring. An element $x \in R$ is said to be right quasi regular (r.q.r) if there exists a $y \in R$ such that $x \circ y = 0$; and x is said to be left quasi regular (l.q.r) if there exists a $y' \in R$ such that $y' \circ x = 0$. An element x is quasi regular (q.r) if it is both right and left quasi regular, y is known as the right quasi inverse (r.q.i) of x and y' is the left-quasi inverse (l.q.i) of x. A right ideal or a left ideal in R is said to be right quasi regular (l-q-r or qr respectively) if each of its element is right quasi regular (l-q r or q-r respectively). Let R be a ring. An element $x \in R$ is said to be a regular element if there exists a $y \in R$ such that $x y x = x$. The Jacobson radical J(R) of a ring R is defined as follows:*

$J(R) = \{a \in R \,/\, a R$ *is a right quasi regular ideal of R*$\}$.

In case of non-associative rings we define the concept of regular element in a different way. We roughly say a ring R is non-associative if the operation '•' on R is non-associative.

**DEFINITION 1.4.7:** *Let R be a non-associative ring. An element $x \in R$ is said to be right regular if there exists a $y \in R$ ($y' \in R$) such that $x(yx) = x$ (($xy'$) $x = x$). A ring R is said to be semisimple if $J(R) = \{0\}$ where J(R) is the Jacobson radical of R.*

Several important properties can be obtained in this direction. We now proceed on to prove the concept of Smarandache rings (S-rings) and Smarandache non-associative rings (SNA-rings for short).

**DEFINITION 1.4.8:** *A Smarandache ring (S-ring) is defined to be a ring A, such that a proper subset of A is a field with respect to the operations induced. By proper subset we understand a set included in A, different from the empty set, from the unit element if any and from A.*

*These are S-ring I, but by default of notation we just denote it as S-ring.*

*Example 1.4.1:* Let F[x] be a polynomial ring over a field F. F[x] is a S-ring.

***Example 1.4.2:*** Let $Z_6 = \{0, 1, 2, \ldots, 5\}$, the ring of integers modulo 6. $Z_6$ is a S-ring for take A = $\{0, 2, 4\}$ is a field in $Z_6$.

**DEFINITION 1.4.9:** *Let R be a ring. R is said to be a Smarandache ring of level II (S-ring II) if R contains a proper subset A ($A \neq 0$) such that*

 i. *A is an additive abelian group.*
 ii. *A is a semigroup under multiplication '•'.*
 iii. *For a, b ∈ A; a • b = 0 if and only if a = 0 or b = 0.*

The following theorem is straightforward hence left for the reader to prove.

**THEOREM 1.4.1:** *Let R be a S-ring I then R is a S-ring II.*

**THEOREM 1.4.2:** *Let R be a S-ring II then R need not be a S-ring I.*

*Proof*: By an example, Z[x] is a S-ring II and not a S-ring I.

**DEFINITION 1.4.10:** *Let R be a ring, R is said to be a Smarandache commutative ring II (S-commutative ring II) if R is a S-ring and there exists atleast a proper subset A of R which is a field or an integral domain i.e. for all a, b ∈ A we have ab = ba. If the ring R has no proper subset, which is a field or an integral domain, then we say R is a Smarandache-non-commutative ring II (S-non-commutative ring II).*

Several results can be obtained in this direction, the reader is advised to refer [120].

It is an interesting feature that for S-ring we can associate several characteristics.

**DEFINITION 1.4.11:** *Let R be a S-ring I or II we say the Smarandache characteristic (S-characteristic) of R is the characteristic of the field, which is a proper subset of R (and or) the characteristic of the integral domain which is a proper subset of R or the characteristic of the division ring which is a proper subset of R.*

Now we proceed on to define Smarandache units, Smarandache zero divisors and Smarandache idempotents in a ring. This is defined only for rings we do not assume it to be a S-ring.

**DEFINITION 1.4.12:** *Let R be a ring with unit. We say $x \in R \setminus \{1\}$ is a Smarandache unit (S-unit) if there exists a $y \in R$ with*

 i. *xy = 1.*
 ii. *there exists a, b ∈ R \ {x, y, 1}*
    a. *xa = y or ax = y or*
    b. *yb = x or by = x and*
    c. *ab = 1.*

*ii(a) or ii(b) is satisfied it is enough to make it a S-unit.*

**THEOREM 1.4.3:** *Every S-unit of a ring is a unit.*

*Proof*: Straightforward by the very definition.

Several results in this direction can be had from [120]. Now we proceed on to define Smarandache zero divisors.

**DEFINITION 1.4.13:** *Let R be a ring we say x and y in R is said to be a Smarandache zero divisor (S-zero divisor) if xy = 0 and there exists a, b ∈ R \ {0, x, y} with*

  i.   $xa = 0$ or $ax = 0$.
  ii.  $yb = 0$ or $by = 0$.
  iii. $ab \neq 0$ or $ba \neq 0$.

The following theorem can be easily proved by examples and by the very definition.

**Theorem 1.4.4:** *Let R be ring. Every S-zero divisor is a zero divisor but all zero divisors in general are not S-zero divisors.*

**DEFINITION 1.4.14:** *Let R be a ring with unit if every unit is a S-unit then we call R a Smarandache strong unit field (S-strong unit field).*

It is to be noted that S-strong unit field can have zero divisors.

**DEFINITION 1.4.15:** *Let R be a commutative ring. If R has no S-zero divisors we say R is a Smarandache integral domain (S-integral domain).*

**THEOREM 1.4.5:** *Every integral domain is a S-integral domain.*

*Proof***:** Straightforward; hence left for the reader to prove all S-integral domains are not integral domains.

*Example 1.4.3:* Let $Z_4 = \{0, 1, 2, 3\}$ be the ring of integers; $Z_4$ is an S-integral domain and not an integral domain.

**DEFINITION 1.4.16:** *Let R be a non-commutative ring. If R has no S-zero divisors we call R a Smarandache division ring (S-division ring). It is easily verified that all division rings are trivially S-division rings.*

Now we proceed on to define Smarandache idempotents of a ring R.

**DEFINITION 1.4.17:** *Let R be a ring. An element $0 \neq x \in R$ is a Smarandache idempotent (S-idempotent) of R if*

  ii.  $x^2 = x$.
  iii. There exists $a \in R \setminus \{1, x, 0\}$ such that
       a. $a^2 = x$.
       b. $xa = a$ $(ax = a)$ or $xa = x$ $(or\ ax = x)$

*or in ii(a) is in the mutually exclusive sense. Let $x \in R \setminus \{0, 1\}$ be a S-idempotent of R i.e. $x^2 = x$ and there exists $y \in R \setminus \{0, 1, x\}$ such that $y^2 = x$ and $yx = x$ or $xy = y$. We call y the Smarandache co idempotent (S-co idempotent) and denote the pair by (x, y).*

**DEFINITION 1.4.18:** *Let R be a ring. A proper subset A of R is said to be a Smarandache subring (S-subring) of R if A has a proper subset B which is a field and A is a subring of R. The Smarandache ideal (S-ideal) is defined as an ideal A such that a proper subset of A is a field (with respect to the induced operations).*

***Example 1.4.4:*** Let $Z_6$ = {0, 1, 2, 3, 4, 5}. Clearly I = {0, 3} and J = {0, 2, 4} are not S-ideals only ideals.

***Example 1.4.5:*** Let $Z_{12}$ = {0, 1, 2, ... , 11} be the ring of integers modulo 12. I = {0, 2, 4, 8, 6, 10} is a S-ideal of $Z_{12}$.

**DEFINITION 1.4.19:** *Let R be a S-ring, B a proper subset of R, which is a field. A non-empty subset C of R is said to be a Smarandache pseudo right ideal (S-pseudo right ideal) of R related to A if*

   i.   *(C, +) is an additive abelian group.*
   ii.  *For $b \in B$ and $s \in C$ we have $sb \in C$.*

*On similar lines we define Smarandache pseudo left ideal. A non-empty subset X of R is said to be a Smarandache pseudo ideal if X is both a S-pseudo right ideal and S-pseudo left ideal.*

**DEFINITION 1.4.20:** *Let R be a ring. I a S-ideal of R; we say I is a Smarandache minimal ideal (S-minimal ideal) of R if we have $J \subset I$ where J is another S-ideal of R then J = I is the only ideal. Let R be a S-ring and M be a S-ideal of R, we say M is a Smarandache maximal ideal (S-maximal ideal) of R if we have another S-ideal N such that $M \subset N \subset R$ then the only possibility is M = N or N = R.*

**DEFINITION 1.4.21:** *Let R be a S-ring and I be a S-pseudo ideal related to A, $A \subset R$ (A is a field) I is said to be a Smarandache minimal pseudo ideal (S-minimal pseudo ideal) of R if $I_1$ is another S-pseudo ideal related to A and $\{0\} \subset I_1 \subset I$ implies $I = I_1$ or $I_1 = \{0\}$. Thus minimality may vary with the different related fields. Let R be a S-ring, M is said to be Smarandache maximal pseudo ideal (S-maximal pseudo ideal) of R related to the field A, $A \subset R$ if $M_1$ is another S-pseudo ideal related to A and if $M \subset M_1$ then $M = M_1$.*

**DEFINITION 1.4.22:** *Let R be a S-ring, a S-pseudo ideal I related to a field A, $A \subset R$ is said to be Smarandache cyclic pseudo ideal (S-cyclic pseudo ideal) related to A if I can be generated by a single element. Let R be S-ring, a S-pseudo ideal I of R related to the field A is said to be Smarandache prime pseudo ideal (S-prime pseudo ideal) related to A if $x y \in I$ implies $x \in I$ or $y \in I$.*

Several nice and interesting results about them can be had from [120].

Now we proceed on to define S-subring II, S-ideal II and S-pseudo ideal II.

**DEFINITION 1.4.23:** *Let R be a S-ring II, A is a proper subset of R is a Smarandache subring II (S-subring II) of R if A is a subring and A itself is a S-ring II. A non-empty* *subset I of R is said to be a Smarandache right ideal II (S-right ideal II) of R (S-left ideal II of R) if*

> i. *I is a S-subring II.*
> ii. *Let $A \subset I$ be an integral domain or a division ring in I then $ai \in I$. $(ia \in I)$ for all $a \in A$ and $i \in I$. If I is simultaneously S-right ideal II and S-left ideal II then I is called a Smarandache ideal II (S-ideal II) of R related to A. If R has no S-ideals I (II) we call R a Smarandache simple ring (S-simple ring) I (II), we call R a Smarandache pseudo simple I (II) (S-pseudo simple I (II)) if R has no S-pseudo ideal I (or II).*

**DEFINITION 1.4.24:** *Let R be a S-ring I (II). I a S-ideal I (II) of R, R/I = {a + I | a ∈ R} is a Smarandache quotient ring I (II) (S-quotient ring I (II)) of R related to I.*

**DEFINITION 1.4.25:** *Let R be a S-ring, I a S-pseudo ideal of R. R / I = {a + I / a ∈ R} is a Smarandache pseudo quotient ring (S-pseudo quotient ring) I of R.*

Now we proceed on to define Smarandache modules.

**DEFINITION 1.4.26:** *The Smarandache R-module (S-R-module) is defined to be an R-module (A, +, ×) such that a proper subset of A is a S-algebra (with respect to the same induced operations and another '×' operation internal on A) where R is a commutative unitary S-ring and X its proper subset which is a field.*

Notions of S-right (left) module I and II are defined as follows:

**DEFINITION 1.4.27:** *Let R be a S-ring I (or II). A non-empty set B which is an additive abelian group is said to be a Smarandache right (left) module I (or II) (S-right (left) module I (or II)) relative to the S-subring I (or II) of A if $D \subset A$ where D is a field (Division ring or an integral domain) then $DB \subset B$ and $BD \subset B$ i.e. bd (and db) are in B with $b(d + c) = bd + bc$ for all $d, c \in D$ and $b \in B$.*

*On similar lines we can define S-left module. If I is simultaneously both a S-left and S-right module than we call I a S-module I (II).*

**DEFINITION 1.4.28:** *Let $(R, +, \bullet)$ be a S-ring. B be a proper subset of A $(B \subset A)$ which is a field. A set M is said to be a Smarandache pseudo right (left) module (S-pseudo right (left) module) of R related to B if*

> i. *(M, +) is an additive abelian group.*
> ii. *For $b \in B$ and $m \in M$, $m \bullet b \in M$ $(b \bullet m \in M)$.*
> iii. *$(m_1 + m_2) \bullet b = m_1 \bullet b + m_2 \bullet b$; $(b \bullet (m_1 + m_2) = b \bullet m_1 + b \bullet m_2)$ for $m_1, m_2 \in M$ and $b \in B$.*

*If M is simultaneously a S-pseudo right module and S-pseudo left module then we say M is a Smarandache pseudo module (S- pseudo module) related to B.*

For more literature please refer [120].

**DEFINITION 1.4.29:** *Let R be a ring, we say the ring R satisfies the Smarandache ascending chain conditions (S-A.C.C for brevity) if for every ascending chain of S-ideals $I_j$ of R that is $I_1 \subset I_2 \subset I_3 \subset ...$ is stationary in the sense that for some integer $p \geq 1$, $I_p = I_{p+1} = ....$ Similarly R is said to satisfy Smarandache descending chain conditions (S-D.C.C) for brevity) if every descending chain $N_1 \supset N_2 \supset ... \supset N_k \supset ...$ of S-ideals $N_j$ of R is stationary. On similar lines one can define Smarandache A.C.C and Smarandache D.C.C of S-right ideals and S-left ideals of a ring. Thus we say a ring R is said to be Smarandache left Noetherian (S-left Noetherian) if the S-A.C.C on S-left ideals is satisfied. Similarly we say a ring R is Smarandache right Artinian (S-right Artinian) if the S-left ideals of R satisfies the S-D.C.C. condition.*

**DEFINITION 1.4.30:** *Let R be a ring. A nilpotent element $0 \neq x \in R$ is said to be a Smarandache nilpotent (S-nilpotent) element if $x^n = 0$ and there exists a $y \in R \setminus \{0, x\}$ such that $x^r y = 0$ or $yx^s = 0$; $r, s > 0$ and $y^m \neq 0$ for any integer $m > 1$.*

**DEFINITION 1.4.31:** *An element $\alpha \in R \mid \{0\}$ is said to be Smarandache semi-idempotent I (S-semi-idempotent I).*

  i. *if the ideal generated by $\alpha^2 - \alpha$ that is $R (\alpha^2 - \alpha) R$ is a S-ideal I and $\alpha \notin R (\alpha^2 - \alpha) R$ or $R = R (\alpha^2 - \alpha) R$.*

*We call $\alpha$ a Smarandache semiidempotent III (S-semi-idempotent III).*

  ii. *if the ideal generated by $\alpha^2 - \alpha$ i.e. $R (\alpha^2 - \alpha) R$ is a S-ideal II and $\alpha \notin R (\alpha^2 - \alpha) R$ or $R = R (\alpha^2 - \alpha) R$.*

**DEFINITION 1.4.32:** *Let R be a non-commutative ring. A pair of distinct elements $x, y \in R$ different from the identity of R which are such that $xy = yx$ is said to be a pseudo commutative pair of R if $xay = yax$ for all $a \in R$. If in a ring R every commutative pair happens to be a pseudo commutative pair of R then R is said to be a pseudo commutative ring. A pair of elements $x, y \in A$, $A \subset R$ A a S-subring of R such that $xy = yx$ is said to be Smarandache pseudo commutative pair (S-pseudo commutative pair) of R if $xay = yax$ for all $a \in A \subset R$. If in a S-subring $A \subset R$, every commuting pair happens to be a S-pseudo commutative pair then we call R a Smarandache pseudo commutative ring (S-pseudo commutative ring).*

**DEFINITION 1.4.33:** *Let R be a ring. We say R is a Smarandache strongly regular ring (S-strongly regular ring) if R contains a S-subring B such that for every $x, y$ in B we have $(xy)^n = xy$ for some integer $n = n(x, y) > 1$.*

Several results in this direction can be had form [120].

**DEFINITION 1.4.34:** *Let R be a ring. R is said to be Smarandache quasi commutative (S-quasi commutative) ring if for any S-subring A or R we have $ab = b^\gamma a$ for every a, $b \in A$; $\gamma \geq 1$.*

**DEFINITION 1.4.35:** *Let R be a ring. An element $x \in R$ is said to be Smarandache seminilpotent (S-seminilpotent) if $x^n - x$ is S-nilpotent. We call an element $x \in R$ to be* *seminilpotent if $x^n - x = 0$. A ring R is said to be a Smarandache reduced ring (S-reduced ring) if R has no S-nilpotents.*

**DEFINITION 1.4.36:** *Let R be a ring, R is said to be a Smarandache p-ring (S-p-ring) if R is a S-ring and has a subring P such that $x^p = x$ and $px = 0$ for every $x \in P$.*

*We call R a Smarandache E-ring (S-E-ring) if $x^{2n} = x$ and $2x = 0$ for all $x \in P$, where A is a S-subring of R for which P is a subring of A. It is pertinent to mention here S-E-ring is a partial particularization of a S-p-ring.*

**DEFINITION 1.4.37:** *Let R be a ring and P a subring of a S-subring A of R. We say R is a Smarandache pre J-ring (S-pre J-ring) if for every pair a, b $\in$ P we have $a^n b = ab^n$ for some integer n.*

**DEFINITION 1.4.38:** *Let R be a ring. R is said to be a Smarandache semiprime ring (S-semiprime ring) if and only if R has no zero S-ideal I which has elements $x \in I$ with $x^2 = 0$.*

**DEFINITION 1.4.39:** *A commutative ring with 1 is called Smarandache Marot ring (S-Marot ring) if every S-ideal I (or II) of R is generated by regular elements of I i.e. these S-ideals does not contain zero divisors.*

**DEFINITION 1.4.40:** *Let R be a ring. Let A be a S-subring I (or II) of R. Let $I \subset A$ be an S-ideal I (or II) of the S-subring A. Then I is called the Smarandache subsemiideal I (S-subsemiideal I) (or II). If in the ring R we have a S-subsemiideal then we call R a Smarandache subsemiideal ring (S-subsemiideal ring).*

A more interesting notion about rings is the concept of filial rings and Smarandache filial rings.

**DEFINITION 1.4.41:** *Let R be a ring. We say R is a Smarandache filial ring (S-filial ring) if the relation S-ideal in R is transitive, that is if a S-subring J is an S-ideal in a S-subring I and I is a S-ideal of R then J is an S-ideal of R.*

**DEFINITION 1.4.42:** *Let R be a ring. We call R a Smarandache n-ideal ring (S-n-ideal ring) if for every set of n-distinct S-ideals I (or II), $I_1, I_2, \ldots, I_n$ of R and for every distinct set of n elements $x_1, x_2, \ldots, x_n \in R \setminus (I_1 \cup I_2 \cup \ldots \cup I_n)$ we have $\langle x_1 \cup I_1 \cup \ldots \cup I_n \rangle = \langle x_2 \cup I_1 \cup I_2 \ldots \cup I_n \rangle = \ldots = \langle x_n \cup I_1 \cup I_2 \cup \ldots \cup I_n \rangle$, $\langle \rangle$ denotes the ideal generated by $x_i \cup I_1 \cup \ldots \cup I_n$ for $1 \leq i \leq n$.*

**DEFINITION 1.4.43:** *Let R be a ring. A be a S-subring of R. We say R is Smarandache s-weakly regular ring (S-s-weakly regular ring) if for each $a \in A$; $a \in aAa^2A$.*

**DEFINITION 1.4.44:** *Let R be a ring, R is said to be a Smarandache strongly subcommutative (S-strongly subcommutative) if every S-right ideal I (II) of it is right quasi reflexive (we say a right ideal I of R is said to be quasi reflexive if whenever A and B are two right ideals of R with $AB \subset I$ then $BA \subset I$).*

Now we just introduce the concept of Smarandache Chinese ring.

**DEFINITION 1.4.45:** *Let R be a ring. R is said to be a Smarandache Chinese ring (S-Chinese ring) I (or II) if given elements a, b ∈ R and S-ideal I (or II) in R such that ⟨I ∪ J ∪ a⟩ = ⟨I ∪ J ∪ b⟩ and there exist an element c ∈ R such that ⟨I ∪ a⟩ = ⟨I ∪ c⟩ and ⟨J ∪ b⟩ = ⟨J ∪ c⟩.*

**DEFINITION 1.4.46:** *Let R be a ring. R is said to be a Smarandache J-ring (S-J-ring) if R has a S-subring A such that for all a ∈ A we have $a^n = a$, $n > 1$.*

**DEFINITION 1.4.47:** *Let R be a ring. Let $\{S_i\}$ denote the collection of all S-subrings of R. We say R is a Smarandache strong subring (S-strong subring) if every pair of S-subrings of R generates R. We say the ring is a Smarandache strong ideal ring (S-strong ideal ring) if every pair of S-ideals of R generate R.*

*If $\{S_j\}$ denotes the collection of all S-subrings $\{I_J\}$ denotes the collection of all S-ideals; if every pair $\{S_j, I_i\}$ generate R then we say R is a Smarandache strong subring ideal ring (S-strong subring ideal ring).*

**DEFINITION 1.4.48:** *Let R be ring, we say R is a Smarandache weak ideal ring (S-weak ideal ring) if there exists a pair of distinct S-ideals $I_1, I_2$ in R which generate R i.e. $R = ⟨I_1 ∪ I_2⟩$ similarly we define Smarandache weak subring (S-weak subring) if there exists distinct pair of S-subrings $S_1, S_2$ in R which generate R i.e. $R = ⟨S_1 ∪ S_2⟩$.*

Several interesting results in these directions can be had from [120].

**DEFINITION 1.4.49:** *Let R be a ring, A be a S-subring of R. R is said to be a Smarandache f-ring (S-f-ring) if and only if A is a partially ordered ring without non-zero nilpotents and for any a ∈ A we have $a_1, a_2 ∈ R$ with $a_1 \geq 0$, $a_2 \geq 0$; $a = a_1 - a_2$ and $a_1 a_2 = a_2 a_1 = 0$.*

**DEFINITION 1.4.50:** *Let R be a ring. If the set of S-ideals of R is totally ordered by inclusion then we say R is a Smarandache chain ring (S-chain ring) (if the S-ideals are replaced by S-right ideals or S-left ideals we call R a Smarandache right chain ring (S-right chain ring) or a Smarandache left chain ring (S-left chain ring)).*

**DEFINITION 1.4.51:** *Let R be a ring, I be a S-ideal of X. We say I is a Smarandache obedient ideal (S-obedient ideal) of R if we have two ideals X, Y in R (X ≠ Y) such that ⟨X ∩ I, Y ∩ I⟩ = ⟨X, Y⟩ ∩ I. In case every S-ideal of the ring R happens to be a S-obedient ideal of R then we say R is a Smarandache ideally obedient ring (S-ideally obedient ring).*

**DEFINITION 1.4.52:** *Let R be a ring. R is said to be a Smarandache Lin ring (S-Lin ring) if R contains a S-subring B such that B is a Lin ring i.e. every pair of elements in B satisfies the identities.*

$$(xy - yx)^n = xy - yx \text{ or}$$
$$(xy + yx)^n = xy + yx.$$

For a S-Lin ring we do not demand the equality to be true for all elements of R.

Now we define a new notion called super ore condition for S-ring.

**DEFINITION 1.4.53:** *Let R be a ring. We say R satisfies Smarandache super ore condition (S-super ore condition) if R has a S-subring A and for every pair x, y ∈ A we have r ∈ R such that xy = yr.*

Several result in this direction can be obtained by an innovative reader.

**DEFINITION 1.4.54:** *Let R be a ring. R is said to be an ideally strong ring if every subring of R not containing identity is an ideal of R. We call R a Smarandache ideally strong ring (S-ideally strong ring) if every S-subring of R is an S-ideal of R.*

**DEFINITION 1.4.55:** *Let R be a ring $\{A_i\}$ be the collection of all S-ideals of R. If for every pair of ideals $A_1, A_2 \in \{A_i\}$ we have for every $x \in R \setminus \{A_1 \cup A_2\}$; $\langle A_1 \cup x \rangle = \langle A_2 \cup x \rangle$ and they generate S-ideals of R, then we say R is a Smarandache $I^*$ - ring (S-$I^*$ - ring).*

Smarandache quotient ring isomorphism etc can also be defined as in case of ring, the interested reader is requested to refer [120].

**DEFINITION 1.4.56:** *Let R be a ring. R is called a Smarandache F-ring (S-F-ring) if we have a subset A of R which is a S-subring of R and we have for every subset X in R and a non-zero $b \in R$ such that $bA \cap X \neq \phi$. It is pertinent to mention here that we need not take X as a subset of A but nothing is lost even if we take X to be a subset of A. Similarly b can be in A or in R.*

**DEFINITION 1.4.57:** *Let x be an element of R, x is said to be a Smarandache SS-element (SSS-element) of R if there exists $y \in R \setminus \{x\}$ with $x \bullet y = x + y$.*

*If a ring has at least one nontrivial SSS-element then we call R a SSS-ring.*

Now we proceed on to define two more new notions called Smarandache demiring and Smarandache demi module.

**DEFINITION 1.4.58:** *Let R be a commutative ring with 1. A subset S of R is said to be a Smarandache demi subring (S-demi subring) of R if*

  i. *(S, +) is a S-semigroup.*
  ii. *(S, •) is a S-semigroup.*

**DEFINITION 1.4.59:** *Let R be a commutative ring with 1. P is said to be a Smarandache demi module (S-demi module) over R if*

  i. *P is a S-semigroup under '+' and '•'.*
  ii. *There exists a nontrivial S-demi subring V of R such that for every $p \in P$ and $v \in V$, $v \bullet p$ and $p \bullet v \in P$.*
  iii. *$v \bullet (p_1 + p_2) = v \bullet p_1 + v \bullet p_2$.*
  iv. *$v \bullet (v_1 \bullet p) = (v \bullet v_1) \bullet p$ for all $p, p_1, p_2 \in P$ and $v_1, v_2 \in V$.*

*We call T a subset of P to be a Smarandache subdemi module (S-subdemi module) if T is a S-demi module for the same S-demi subring.*

**DEFINITION 1.4.60:** *Let R be a ring. If for every $x \in P$ there exists a S-semiidempotent s of R such that $x \bullet s = s \bullet x = x$ then we call the ring R a Smarandache locally semiunitary ring (S-locally semiunitary ring).*

Similarly the notion of CN-ring can be defined.

**DEFINITION 1.4.61:** *Let R be a ring, a subset S of R is said to be Smarandache closed net (S-closed net) if*

  i. *S is a semigroup.*
  ii. *S is a S-semigroup.*

*We call R a Smarandache CN-ring (S-CN-ring) if $R = \cup S_j$ where $S_j$'s are S closed nets such that $S_i \cap S_j = A$, $i \neq j$, $A \neq S_i$, $A \neq S_j$ where A is a subgroup of $S_j$. We call R a Smarandache weakly CN-ring (S-weakly CN-ring) if $R \subset \cup S_j$.*

Till now we saw mainly properties related to S-ring and several of the new Smarandache notions for more about these concepts please refer [120]. Now we define Smarandache mixed direct product of rings.

**DEFINITION 1.4.62:** *Let $R = R_1 \times R_2$ where $R_1$ is a ring and $R_2$ is an integral domain or a division ring. Clearly R is a S-ring and we call R a Smarandache mixed direct product (S-mixed direct product). We can extend this to several number of rings say $R_1, R_2, \ldots, R_n$. Let $R = R_1 \times R_2 \times \ldots \times R_n$ is called the Smarandache mixed direct product of n-rings (S-mixed direct product of n-rings) if and only if atleast one or some of the $R_i$'s is an integral domain or a division ring.*

Several identities in them can be studied.

**DEFINITION 1.4.63:** *Let R be a ring, an element x in R is said to be a Smarandache weakly super related (S-weakly super related) in R if there exists $\alpha, \beta, \gamma \in A$ such that $(x + \alpha)(x + \beta)(x + \gamma) = x + \alpha\beta(x + \gamma) + \alpha\gamma(x + \beta) + \beta\gamma(x + \alpha)$ where A is a S-subring of R.*

*We call an element $x \in R$ to be Smarandache super related (S-super related) in R if for all $\alpha, \beta, \gamma \in A$ where A is a S-subring such that $(x + \alpha)(x + \beta)(x + \gamma) = x + \alpha\beta(x + \gamma) + \alpha\gamma(x + \beta) + \beta\gamma(x + \alpha)$. If R has no S-subring but R is a S-ring then we say $x \in R$ is S-super related in R.*

We proceed on to define the notion of bisimple.

**DEFINITION 1.4.64:** *Let R be a ring, we say R is Smarandache bisimple (S-bisimple) if it has more than one element and satisfies the following conditions.*

  i. *For any $a \in A$ we have $a \in aA \cap Aa$ where A is a S-subring.*

ii. For any non-zero $a, b \in A$ there is some $c \in A$ such that $aA = cA$ and $Ac = Ab$.

*We call R Smarandache semibisimple (S-semibisimple) if for any $a, b \in A$ where A is a S-subring, we have $c \in A$ such that $aA = cA$ and $Ac = Ab$. We still define Smarandache weakly bisimple (S-weakly bisimple) if for every $a \in A$; A a S-subring of R we have $a \in aA \cap Aa$ and for every pair of elements $a, b \in A$, $aA \subset cA$ and $Ab \subset Ac$ for some $c \in A$.*

**Note:** We can extend the notion of S-bisimple to be Smarandache trisimple (S-trisimple) if $a \in aA \cap Aa \cap aAa$. For more about these Smarandache notions please refer [120]. Finally we go on to define Smarandache n-like ring.

**DEFINITION 1.4.65:** *Let R be a ring, we say R is a Smarandache n-like ring (S-n-like ring) if R has a proper S-subring A of R such $(xy)^n - xy^n - x^n y - xy = 0$ for all $x, y \in R$.*

This is an identity which is of a special type and hence several interesting results can be determined; the notion of Smarandache power joined and (m, n) power joined are introduced.

**DEFINITION 1.4.66:** *Let R be a ring. If for every $a \in A \subset R$ where A is a S-subring there exists $b \in A$ such that $a^n = b^n$ for some positive integer m and n then we say R is a Smarandache power joined ring (S-power joined ring). If we have $a^n = b^m$, m and n different then we call Ra, the Smarandache (m, n) power joined ring (S-(m, n) power joined ring) if for every $x \in A \subset R$ there exists a $y \in A \subset R$ such that $x^m = y^m$ $(x \neq y)$ and $m \geq 2$, then we say R is a Smarandache uniformly power joined ring (S-uniformly power joined ring).*

Several interesting properties can be had from [120]. We briefly recall the definition of Smarandache strongly right commutative, Smarandache quasi semicommutative and concepts like S-semicommutator.

**DEFINITION 1.4.67:** *Let R be a ring. An element $x \in A \subseteq R$ where A is a S-subring of R is said to be a Smarandache quasi semicommutative (S-quasi semicommutative) if there exists $y \in A$ ($y \neq 0$) such that $xy - yx$ commutes with every element of A.*

*If $x \in A \subset R$ is a Smarandache semicommutator (S-semicommutator) of x denoted by $SQ(x) = \{p \in A \;/\; xp - px$ commutes with every element of $A\}$; R is said to be a Smarandache quasi semicommutative ring (S-quasi semicommutative ring) if for every element in A is S-quasi semicommutative.*

*We define Smarandache quasi semicenter (S-quasi semicenter) $SQ(R)$ of R to be $SQ(R) = \{x \in R \mid xp - px$ is S-quasi semicommutative$\}$.*

The reader is assigned the task of obtaining interesting results about these concepts. Now we proceed on to define Smarandache magnifying and shrinking elements of a ring R.

**DEFINITION 1.4.68:** *Let R be a ring. A ⊂ R be a proper S-subring of R. An element v is called Smarandache left magnifying element (S-left magnifying element) of R if vM = A for some proper subset M of A, we say v is Smarandache right magnifying element (S-right magnifying element) if $M_1v$ = A for some proper subset $M_1$ of A. v is said to be a Smarandache magnifying (S-magnifying) if vM = Mv = A for some M a proper subset of A. If v ∈ A, then v is said to be a Smarandache friendly magnifying element (S-friendly magnifying element). If v ∉ A, we call v a Smarandache non-friendly magnifying element (S-non-friendly magnifying element), even if v ∈ R \ A, we still call v a non-friendly magnifying element.*

*An element x is called Smarandache left shrinking element (S-left shrinking element) of R if for some S-subring A of R we have proper subset M of R such that x A = M (M ≠ A) or M ≠ R. We define similarly Smarandache left shrinking (S-left shrinking) and Smarandache shrinking (S-shrinking) if xA = Ax = M. If x ∈ A, we call x a Smarandache friendly shrinking element (S-friendly shrinking element) if x ∉ A, we call x a Smarandache non-friendly shrinking element (S-non-friendly shrinking element).*

**DEFINITION 1.4.69:** *Let R be a ring if R has a proper S-subring A of R such that the S-subring A has only two S-idempotents then we call R a Smarandache dispotent ring (S-dispotent ring). If every S-subring A of R has exactly two S-idempotents then we say R is a Smarandache strong dispotent ring (S-strong dispotent ring).*

**DEFINITION 1.4.70:** *Let R be a ring. X be a S-subring of R. We say R is a Smarandache normal ring (S-normal ring) if aX = Xa for all a ∈ R. We call R a Smarandache strongly normal ring (S-strongly normal ring) if every S-subring X of R is such that aX = Xa for all a ∈ R.*

Several results can be developed in this direction.

**DEFINITION 1.4.71:** *Let R be a ring. If for every S-semigroup, P under addition we have rP = Pr =P for every r ∈ R (r ≠ 0), then we call R a Smarandache G-ring (S-G-ring). If we have for every S-semigroup P under addition and for every r ∈ R we have rP = Pr, then we call R a Smarandache weakly G-ring (S-weakly-G-ring).*

**DEFINITION 1.4.72:** *Let R be a ring. If R has atleast one S-ideal which contains a non-zero S-idempotent then we say R is Smarandache weakly e-primitive (S-weakly e-primitive); we call R a Smarandache e-primitive (S-e-primitive) if every non-zero S-ideal in R contains a non-zero S-idempotent.*

**DEFINITION 1.4.73:** *Let R be a commutative ring. An additive S-semigroup S of R is said to be a Smarandache radix (S-radix) of R if $x^3t$, $(t^2 - t) x^2 + xt^2$ are in S if for every x ∈ S and t ∈ R. If R is a non-commutative ring then for any S-semigroup S of R we say R has Smarandache left-radix (S-left-radix) if $tx^3$, $(t^2 - t)x^2 + t^2x$ are in S if for every x ∈ S and t ∈ R.*

*Similarly we define Smarandache right radix (S-right radix) of R. If S is a simultaneously a S-left radix and S-right radix of a non-commutative ring then we say R has a S-radix.*

**DEFINITION 1.4.74:** *R be a ring. R is said to be a Smarandache SG-ring (S-SG-ring) if $R = \cup S_i$ where $S_i$ are multiplicative S-semigroups such that $S_i \cap S_J = \phi$ if $i \neq j$. We say R is Smarandache weakly SG -ring (S-weakly SG-ring) if $R = \cup S_i$ where $S_i$'s are S-multiplicative semigroups and $S_i \cap S_J \neq \phi$ even if $i \neq j$.*

**DEFINITION 1.4.75:** *Let R be a ring. We say $0 \neq r \in R$ is called a Smarandache insulator (S-insulator) if for r there exists a non-empty subset X of R where X is a S-semigroup under '+' and the right annihilator $r_s = (\{rx \mid x \in X\}) = \{0\}$. A non-zero ring R is said to be Smarandache strongly prime (S-strongly prime) if every non-zero element of R has a finite S-insulator.*

**DEFINITION 1.4.76:** *Let R be a commutative ring and P an additive S-semigroup of R. P is called a Smarandache n-capacitor group (S-n-capacitor group) of R if $x^n P \subseteq P$ for every $x \in R$ and $n \geq 1$ and n a positive integer.*

**DEFINITION 1.4.77:** *Let R be a ring, a pair $x, y \in R$ is said to have a Smarandache subring right link relation (S-subring right link relation) if there exists a S-subring P in $R \setminus \{x, y\}$ such that $x \in Py$ and $y \in P x$. Similarly Smarandache subring link relation (S-subring link relation) if $x \in y P$ and $y \in x P$. If it has both a Smarandache left and right link relation for the same S-subring P, then we say x and y have a Smarandache subring link (S-subring link).*

*We say $x, y \in R$ is Smarandache weak subring link with a S-subring P (S-weak subring link with a S-subring P) in $R \setminus \{x, y\}$, if either $x \in Py$ or $y \in Px$ (or in a strictly mutually exclusive sense) we have a S-subring $Q \neq P$ such that $y \in Q x$ (or $x \in Q y$). We say a pair $x, y \in R$ is said to be Smarandache one way weakly subring link related (S-one way weakly subring link related) if we have a S-subring $P \subset R \setminus \{x, y\}$ such that $x \in Py$ and for no subring $Q \subset R \setminus \{x, y\}$ we have $y \in Q x$.*

**DEFINITION 1.4.78:** *Let R be a ring, A be a S-subring of R. A is said to be a Smarandache essential subring (S-essential subring) of R if the intersection of every other S-subring is zero. If every S-subring of R is S-essential S-subring then we call R a Smarandache essential ring (S-essential ring).*

**DEFINITION 1.4.79:** *Let R be a ring. If for every pair of S-subrings P and Q of R there exists a S-subring T of R ($T \neq R$) such that the S-subrings generated by PT and TQ are equal; i.e. $\langle PT \rangle = \langle TQ \rangle$ then we say the pair P and Q is a Smarandache stabilized pair (S-stabilized pair) and T is called the Smarandache stabilizer (S-stabilizer) of P and Q. A pair of S-subrings A, B of R is said to be a Smarandache stable pair (S-stable pair) if there exists a S-subring C of R such that $C \cup A = C \cup B$ and $\langle C \cup A \rangle = \langle C \cup B \rangle$ where $\langle \rangle$ means the subring generated by $C \cup A$ and $C \cup B$; C is called Smarandache stability S-subring (S-stability S-subring) for the S-stable pair A and B. If every pair of S-subrings of R is a stable pair then we say R is a Smarandache stable ring (S-stable ring).*

The notion of Smarandache hyper ring will find its application in case of birings, which we choose to call as hyper biring. So we define these notions in case of S-rings.

**DEFINITION 1.4.80:** *Let $Z_n$ be a ring with A to be a S-subring of $Z_n$. Define the Smarandache hyper ring I (or II) (S-hyper ring I or (II)) to be a subring of $A \times A$ given by; for any $q \in A$. $(A, q, +) = \{(a_1 + a_2, a_1 + a_2 + q) \mid a_1, a_2 \in A\}$ and $(A, q, \bullet) = \{(a_1 . a_2, a_1.a_2.q) \mid a_1, a_2 \in A\}$.*

Similarly we define Smarandache hyper ring II for any S-subring II of R.

**DEFINITION 1.4.81:** *Let R ring. We say R is Smarandache semiconnected (S-semiconnected) if the center of R contains a finite number of S-idempotents.*

Study of how the S-ideals of a ring, or the S-subrings of a ring is an important one. We know the set of ideals of a ring forms a modular lattice.

The study in this direction is an interesting and an innovative one. We leave the reader to develop in this direction. We study how the S-ideals and S-subrings of a biring look like.

**DEFINITION 1.4.82:** *Let R be a ring. An element $a \in A$ where A is a S-subring of R is said to be Smarandache clean element (S-clean element) of R if it can be expressed as a sum of an idempotent and a unit in R. A ring R is called a Smarandache clean ring (S-clean ring) if every element of R is S-clean. We call an element $a \in R$ to be Smarandache strongly clean (S-strongly clean) if it can be written as a sum of a S-idempotent and a S-unit. If every element $x \in R$ is S-strongly clean then we call the ring R a Smarandache strongly clean ring (S-strongly clean ring).*

Many more notions about S-rings which has not been carried out in this section can be had from [120].

We now proceed on to define the concept of Smarandache properties in non-associative rings. All Smarandache properties which can be modified to non-associative rings in a analogous way or with some simple modifications are left for the reader to develop. Only those properties, which are enjoyed by SNA rings, are recalled for the reader. Any way the reader is advised to refer [119].

**DEFINITION 1.4.83:** *Let R be a non-associative ring, R is said to be a Smarandache non-associative ring (SNA-ring) if R contains a proper subset P such that P is an associative ring under the operations of R.*

As non-associative rings cannot be constructed using the set of reals or integers or rationals or complex or modulo integers without the aid of another algebraic structure it is really difficult to give natural examples of non-associative rings.

The well-known algebraic classes of non-associative rings are

 i. Loop rings i.e. using a loop and a ring they are built analogous to group rings.
 ii. Groupoid rings, built using groupoids over rings.

These non-associative rings enjoy a varied relation only when they satisfy the special types of identities like Moufang, Bol, Bruck, P-identity and right (left) alternative identities.

We in this book do not study Lie algebras or Jordan algebras as non-associative rings.

**DEFINITION 1.4.84:** *Let R be a SNA-ring we call a non-empty subset I of R to be a Smarandache seminormal subring (S-seminormal subring) if*

    i.   *I is a S-subring of R.*
    ii.  *aI = X for all a ∈ R.*
    iii. *I a = Y for all a ∈ R.*

*where either X or Y is a S-subring of R, and X and Y are just subrings of R. We say W is a Smarandache normal subring (S-normal subring) if aV = X and Va = Y for all a ∈ R, where both X and Y are S-subrings of R (Here V is a S-subring of R).*

Almost all relations would follow in a very natural way with simple modifications.

## 1.5 Semirings, S-semirings and S-semivector spaces

This section is devoted to the introduction of semirings, Smarandache semirings and Smarandache semivector spaces. As the study of bisemirings or bisemivector spaces is very new only being introduced in this book. Here we give special importance to the study of Smarandache bisemirings and Smarandache bisemivector spaces. For more about semirings, Smarandache semirings and Smarandache semivector spaces the reader can refer [122]. First we introduce the concept of semirings and semivector spaces.

**DEFINITION 1.5.1:** *Let S be a non-empty set on which is defined two binary operations '+' addition and '•' multiplication satisfying the following conditions:*

    i.   *(S, +) is a commutative monoid.*
    ii.  *(S, •) is a semigroup.*
    iii. *(a + b) • c = a • c + b • c and a • (b + c) = a • b + a • c for all a, b, c in S.*

*(S, +, •) is called the semiring. If in the semiring (S, +, •) we have (S, •) to be a commutative semigroup then we call (S, +, •) a commutative semiring. If (S, •) is a non-commutative semigroup then we call (S, +, •) a non-commutative semiring. If in a semiring (S, +, •) we have (S, +, •) to be a monoid, we call the semiring to be a semiring with unit.*

***Example 1.5.1:*** Let $Z^o = Z^+ \cup \{0\}$ where $Z^+$ denotes the set of positive integers then $Z^o$ under usual '+' and multiplication is a semiring; in fact a commutative semiring with unit. Similarly $Q^o = Q^+ \cup \{0\}$ and $R^o = R^+ \cup \{0\}$ are also commutative semirings with unit where $Q^+$ is the set of positive rationals and $R^+$ is the set of positive reals.

**DEFINITION 1.5.2:** *Let (S, +, •) be a semiring. We say the semiring is of characteristic m if mx = x + ...+ x (m times) equal to zero for all x ∈ S. If no such m exists we say the characteristic of the semiring is 0.*

*Here it is pertinent to mention that certain semirings will have no characteristic associated with it. Thus this is the main difference between a ring and a semiring.*

*The number of elements in a semiring is denoted by o(S) or |S| i.e. the order of S. If o(S) = n with n < ∞ we say the semiring is finite, if n = ∞ we say the semiring is of infinite order.*

***Example 1.5.2***: Let L be the chain lattice given by the following diagram:

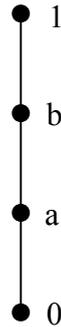

**Figure 1.5.1**

L is semiring with o(L) = 4 and L has no characteristic associated with it.

***Example 1.5.3***: The following distributive lattice S is also a semiring of order 8.

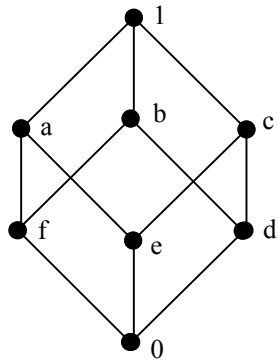

**Figure 1.5.2**

S is commutative and has no characteristic associated with it.

Direct product of semirings can be defined as in case of any algebraic structure.

**DEFINITION 1.5.3**: *Let S be a semiring. P a subset of S. P is said to be a subsemiring if P itself is a semiring. A non-empty subset I of S is called a right (left) ideal of S if*

    i. *I is a subsemiring.*
    ii. *For all $i \in I$ and $s \in S$.*

*we have i s ∈ I (s i ∈ I). I is called an ideal if I is simultaneously a right and a left ideal of S.*

We define the concept of left (right) unit, unit, left (right) zero divisor in a semiring as in case of rings.

**DEFINITION 1.5.4:** *An element $x \in S$, S a semiring is said to be an idempotent if $x^2 = x$.*

**DEFINITION 1.5.5:** *Let S and S' be two semirings. A mapping $\phi : S \to S'$ is called the semiring homomorphism if*

$$\phi(a + b) = \phi(a) + \phi(b) \text{ and}$$
$$\phi(a \bullet b) = \phi(a) \bullet \phi(b) \text{ for all } a, b \in S.$$

*If $\phi$ is one to one and onto we call $\phi$ a semiring isomorphism.*

**DEFINITION 1.5.6:** *Let S be a semiring. We say S is a strict semiring if $a + b = 0$ implies $a = 0$ and $b = 0$. A commutative strict semiring with unit and without divisors of zero is called a semifield.*

***Example 1.5.4:*** $Z^o$ is a semifield.

***Example 1.5.5:*** Any chain lattice L is a semifield.

In example 1.5.4, the semifield is of characteristic zero where as in the semifield given in example 1.5.5, L has no characteristic associated with it. Polynomial semirings are defined analogous to polynomial rings. If S[x] is a polynomial semiring then S[x] is a commutative semiring with 1. To create non-commutative semirings other than matrix semirings we go by the new definitions of the algebraic structures called group semirings and semigroup semirings. By using non-commutative groups and semigroups we can get several new non-commutative semirings.

**DEFINITION 1.5.7:** *Let S be a strict commutative semiring with unit. G any group under multiplication. SG be the group semiring defined analogous to group ring. SG consists of finite formal sums of the form*

$$\sum_{i=1}^{n} s_i g_i$$

*where i runs over a finite number n with $s_i \in S$ and $g_i \in G$ satisfying the following conditions.*

i. $\sum_{i=1}^{n} s_i g_i = \sum_{i=1}^{n} t_i g_i \Leftrightarrow s_i = t_i$ *and $g_i \in G$ for all i.*

ii. $\sum_{i=1}^{n} s_i g_i + \sum_{i=1}^{n} t_i g_i = \sum_{i=1}^{n} (s_i + t_i) g_i$ for all $s_i, t_i \in S$; $g_i \in G$.

iii. $\left(\sum s_i g_i\right)\left(\sum t_j h_j\right) = \sum m_k p_k$ where $m_k = \sum s_i t_j$ and $p_k = g_i h_j$; $g_i h_j \in G$ and $s_i, t_j \in S$.

iv. $s_i g_j = g_j s_i$ for $s_i \in S$ and $g_i \in G$.

v. $s\left(\sum_i s_i g_i\right) = \sum_i s s_i g_i$ for $g_i \in G$ and $s, s_i \in S$.

vi. As $1 \in G$ and $1 \in S$ we have $1.G = G \subseteq SG$ and $S.1 = S \subseteq SG$.

*The group semiring SG will be a commutative semiring when G is commutative and will be a non-commutative semiring when G is non-commutative.*

Now we proceed on to define the concept of semivector spaces using the concept of semifield defined earlier.

**DEFINITION 1.5.8:** *A semivector space V over the semifield S of characteristic zero is the set of elements called vectors with two laws of combination called vector addition and scalar multiplication satisfying the following conditions.*

i. *To every pair of vectors $\alpha$, $\beta$ in V there is associatied a vector in V called the sum which we denote by $\alpha + \beta$.*
ii. *Addition is associative $(\alpha + \beta) + \gamma = \alpha + (\beta + \gamma)$ for all $\alpha, \beta, \gamma \in V$.*
iii. *There exists a vector which we denote by zero such that $0 + \alpha = \alpha + 0 = \alpha$ for all $\alpha \in V$.*
iv. *Addition is commutative i.e $\alpha + \beta = \beta + \alpha$, $\alpha, \beta \in V$.*
v. *For $0 \in S$ and $\alpha \in V$ we have $0 \bullet \alpha = 0$.*
vi. *To every scalar $s \in S$ and every vector $v \in V$ there is associated a unique vector called the product $s \bullet v$, which we denote by sv.*
vii. *Scalar multiplication is associative $(ab)\alpha = a(b\alpha)$ for all $a, b \in S$ and $\alpha \in V$.*
viii. *Scalar multiplication distributes, that is $a(\alpha + \beta) = a\alpha + a\beta$ for all $a \in S$ and $\alpha, \beta \in V$.*
ix. *Scalar multiplication is distributive with respect to scalar addition: $(a + b)\alpha = a\alpha + b\alpha$ for all $a, b \in S$ and for all $\alpha \in V$.*
x. *$1 \bullet \alpha = \alpha$ (where $1 \in S$) and $\alpha \in V$.*

**Example 1.5.6:** $Z^o$ is a semifield and $Z^o[x]$ is a semivector space over $Z^o$.

**Example 1.5.7:** $Q^o$ is a semifield and $R^o$ is a semivector space over $Q^o$.

**Example 1.5.8:** Let $M_{n \times n} = \{(a_{ij}) \mid a_{ij} \in Z^o\}$, the set of all $n \times n$ matrices with entries from $Z^o$. Clearly $M_{n \times n}$ is a semivector space over $Z^o$.

***Example 1.5.9:*** Let $C_n$ be a chain lattice; then $C_n [x]$ is a semivector space over $C_n$ ($C_n$ is a chain lattice with n elements including 0 and 1 i.e. $0 < a_1 < a_2 < \ldots < a_{n-2} < 1$).

Several interesting properties about semivector spaces can be had from [81 and 122]. The main property which we wish to state about semivector spaces is the concept of basis.

**DEFINITION 1.5.9:** *A set of vectors ($v_1, v_2, \ldots, v_n$) in a semivector space V over a semifield S is said to be linearly dependent if there exists a non-trivial relation among them; otherwise the set is said to be linearly independent. (The main difference between vector spaces and semivector spaces is that we do not have negative terms in a semifield over which semivector spaces are built).*

**DEFINITION 1.5.10:** *Let V be a semivector space over the semifield S. For any subset A of V the set of all linear combinations of vectors in A is called the set spanned by A and we shall denote it by ⟨A⟩. Clearly A ⊂ ⟨A⟩. A linearly independent set of a semivector space over the semifield S is called a basis of V if that set can span the semivector space V.*

**DEFINITION 1.5.11:** *A subsemivector space W of a semivector space V over a semifield S is a non-empty subset of V, which is itself a semivector space with respect to the operations of addition and scalar multiplication.*

**DEFINITION 1.5.12:** *Let $V_1$ and $V_2$ be any two semivector spaces over the semifield S. We say a map / function $T : V_1 \to V_2$ is a linear transformation of semivector spaces if $T(\alpha v + u) = \alpha T(v) + T(u)$ for all $u, v \in V_1$ and $\alpha \in S$. A map / function from V to V is called a linear operator of the semivector space V if $T((\alpha v + u) = \alpha T(v) + T(u)$ for all $\alpha \in S$ and $u, v \in V$.*

We proceed on to define notions of Smarandache semirings and semivector spaces. We recall basic properties. The reader can refer [122] for more information.

**DEFINITION 1.5.13:** *The Smarandache semiring S which will be denoted from here onwards by S-semiring is defined to be a semiring S such that a proper subset B of S is a semifield (with respect to the same induced operation). That is $\phi \neq B \subset S$.*

***Example 1.5.10:*** $Z^o \times Z^o$ is a semiring. Clearly B = $Z^o \times \{0\}$ is a semifield. So S is a S-semiring.

If the S-semiring has only a finite number of elements we say the S-semiring is finite other wise infinite.

**DEFINITION 1.5.14:** *Let S be a semiring, A non-empty proper subset A of S is said to be a Smarandache subsemiring (S-subsemiring); if A is a S-semiring; i.e. A has a proper subset P such that P is a semifield under the operations of S.*

**DEFINITION 1.5.15:** *Let S be a S-semiring. We say S is a Smarandache commutative semiring (S-commutative semiring) if S has a S-subsemiring, which is commutative. If the S-semiring has no commutative S-subsemiring then we say S is a Smarandache non-commutative semiring (S-non-commutative semiring).*

**DEFINITION 1.5.16:** *Let S be a semiring. A non-empty subset P of S is said to be Smarandache right (left) ideal (S-right (left) ideal) of S if the following conditions are satisfied.*

   i.   *P is a S-subsemiring.*
   ii.  *For every p ∈ P and A ⊂ P where A is the semifield of P we have for all a ∈ A and p ∈ P ap (pa) is in A.*

*If P is simultaneously both a S-right ideal and a S-left ideal than we say P is a Smarandache ideal (S-ideal).*

**DEFINITION 1.5.17:** *Let S be a semiring. A non-empty proper subset A of S is said to be Smarandache pseudo subsemiring (S-pseudo subsemiring) if the following condition is true. If there exists a subset P of S such that A ⊂ P where P is a S-subsemiring i.e. P has a subset B such that B is a semifield under the operations of S or P itself is a semifield under the operations of S.*

Now we proceed on to define Smarandache pseudo (right) left ideals.

**DEFINITION 1.5.18:** *Let S be a semiring. A non-empty subset P of S is said to be a Smarandache pseudo right (left) ideal (S- pseudo right (left) ideal) of the semiring S if the following conditions are true.*

   i.   *P is a S-pseudo subsemiring i.e. P⊂A, A is a semifield in S.*
   ii.  *For every p ∈ P and every a∈ A; ap ∈ P (pa ∈ P).*

*If P is simultaneously both a S-pseudo right and left ideal we say P is a Smarandache pseudo ideal (S-pseudo ideal).*

*We define Smarandache pseudo dual ideal of S as follows: A non-empty subset P of S is said to be a Smarandache pseudo dual ideal (S-pseudo dual ideal) if the following conditions hold good.*

   i.   *P is a S-subsemiring.*
   ii.  *For every p ∈ P and a ∈ A \ {0}; a + p is in A where A ⊂ P.*

**Example 1.5.11**: Let S be the semiring given by the following diagram:

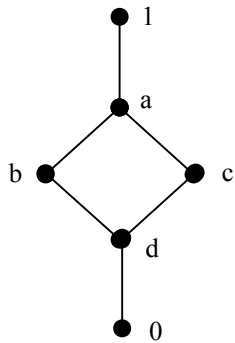

**Figure 1.5.3**

P = {0, a, d, b, 1} ⊂ S, P is a S-dual ideal of S. (For A = {0, 1} is a semifield).

Similarly we can define Smarandache dual ideal of a semiring. This task is left as an exercise for the reader.

**DEFINITION 1.5.19:** *Let S be a S-semiring. S is said to be a Smarandache semidivision ring (S-semidivision ring) if the proper subset A ⊂ S is such that*

iii.  *A is a S-subsemiring which is non-commutative.*
iv.  *A contains a subset P such that P is a semidivision ring that is P is a non-commutative semiring.*

Now the notions of S-zero divisors, S-units and S-idempotents are defined in semirings as in case of rings.

Now we introduce a new notion called Smarandache anti zero divisor.

**DEFINITION 1.5.20:** *Let S be a semiring. An element $x \in S$ is said to be a Smarandache anti zero divisor (S-anti zero divisor) if we can find a y such that $xy \neq 0$ and $a, b \in S \setminus \{0, x, y\}$ such that*

  i.  $ax \neq 0$ or $xa \neq 0$.
  ii. $by \neq 0$ or $yb \neq 0$.
  iii. $ab = 0$ or $ba = 0$.

Several interesting results in this direction can be found in [120].

The concept of S-D.C.C and S-A.C.C conditions can be defined as in case of rings with corresponding minor modifications.

Now we proceed on to define Smarandache semirings of level II.

**DEFINITION 1.5.21:** *A semiring S is said to be a Smarandache semiring of level II (S-semiring of level II) if S contains a proper subset A where A is a field under the operations of S.*

*Example 1.5.12*: Let $S = C_2 \times Q$ where $C_2$ is a chain lattice with two elements and Q is a field; then S is a S-semiring II.

We now define Smarandache mixed direct product of semirings; the main motivation is it gives the mode of construction of S-semiring II.

**DEFINITION 1.5.22:** *Let $S = S_1 \times S_2 \ldots \times S_n$ is called the Smarandache mixed direct product of semirings (S-mixed direct product of semirings) where at least one of the $S_i$'s is a semiring and one of the $S_j$'s is a field, other $S_k$'s can be semifields, rings or* *semirings. So S under component wise operations is a semiring. Clearly it can be proved that S is a S-semiring II.*

All properties like S-ideals S-subsemirings etc can be defined for these S-semirings II and this work is assigned to the reader.

Now we proceed on to define a new notion called Smarandache anti semiring.

**DEFINITION 1.5.23:** *Let R be a ring. We say R is a Smarandache anti semiring (S-anti semiring) if R contains a semiring.*

*Example 1.5.13*: Z be the ring of integers. $S = Z^+ \cup \{0\}$ is a semiring in Z. So Z is a S-anti semiring.

Several examples can be given and it is left as exercise for the reader. Next we proceed on to define Smarandache semifields.

**DEFINITION 1.5.24:** *Let S be a semifield. S is said to be a Smarandache semifield (S-semifield) if a proper subset of S is a K-semialgebra with respect to the same induced operations and an external operator. We call S a Smarandache weak semifield (S-weak semifield) if S contains a proper subset P which is a semifield and P is a S-semifield.*

**DEFINITION 1.5.25:** *Let S be a field or a ring. S is said to be Smarandache anti-semifield (S-anti-semifield) if S has a proper subset which is a semifield. It is easily verified that all fields of characteristic zero are S-anti semifields.*

Now we proceed on to define Smarandache semivector spaces.

**DEFINITION 1.5.26:** *Let G be a semigroup under the operation '+'; S any semifield. Let G be a semivector space over S. G is said to be a Smarandache semivector space (S-semivector space) over S if G is a S-semigroup. We call a subset W of a S-semivector space V to be a Smarandache subsemivector space (S-subsemivector space) if W itself is a S-semivector space.*

**DEFINITION 1.5.27:** *Let V be a semigroup, which is a S-semivector space over a semifield S. A Smarandache basis (S-basis) for V is a set of linearly independent elements which span a S-subsemivector space P of V; that is P is a S-subsemivector space of V, so P is also a S-semigroup.*

We define Smarandache pseudo semivector space.

**DEFINITION 1.5.28:** *Let V be a vector space over S. Let W be a proper subset of V. If W is not a subsemivector space over S but W is a subsemivector space over a proper subset $P \subset S$ then we say W is a Smarandache pseudo semivector space (S-pseudo semivector space) over $P \subset S$.*

The concept of Smarandache anti-semivector space is defined.

**DEFINITION 1.5.29:** *Let V be a vector space over a field F. We say V is a Smarandache anti semivector space (S-anti semivector space) over F if there exists a subspace $W \subset V$ such that W is a semivector space over a semifield $S \subset F$. Here W is just a semigroup under '+'.*

Several interesting properties about these structures can be had from [122].

## 1.6 Near-rings and S-near-rings

In this section we introduce the concept of near-rings and Smarandache near-rings. As it is essential to know these concepts for building up binear-rings and Smarandache binear-rings, which are very new concepts, we have briefly given them. The study of Smarandache near-rings is very recent; introduced only in the year 2002 [43, 118, 126].

**DEFINITION 1.6.1:** *A near-ring is a set N together with two binary operations '+' and '•' such that*

  i. *(N, +) is a group.*
  ii. *(N, •) is a semigroup.*
  iii. *For all $n_1, n_2, n_3 \in N$; $(n_1 + n_2) \cdot n_3 = n_1 \cdot n_3 + n_2 \cdot n_3$ (right distributive law).*

*This near-ring will be termed as right near-ring. If $n_1 \cdot (n_2 + n_3) = n_1 \cdot n_2 + n_1 \cdot n_3$ instead of condition (iii) the set N satisfies, then we call N a left near-ring. Throughout this text we will be using only right near-rings unless otherwise we clearly specify it.*

*Let (N, +, •) be a near-ring |N| or o(N) denotes the order of the near-ring N i.e. the number of elements in N. If $|N| < \infty$ we say the near-ring is of finite order if $|N| = \infty$ we call N an infinite near-ring. We denote by $N_d = \{d \in N \mid d(n_1 + n_2) = dn_1 + dn_2$ for all $n_1, n_2 \in N\}$.*

**DEFINITION 1.6.2:** *Let N be a near-ring; if (N, +) is abelian we call N an abelian near-ring. If (N, •) is commutative we call N itself a commutative near-ring. If $N = N_d$, N is said to be distributive. If all non-zero elements of N are left cancelable we say that N fulfils the left cancellation law. Similarly one defines right cancelable. N is an integral domain if N has no non zero divisors of zero. If $N' = N \setminus \{0\}$ is a group under '•' then N is called a near field. A near-ring which is not a ring will be referred to as non-ring. Similarly a non-field is a near field, which is not a field. A near-ring with the property that $N_d$ generates (N, +) is called a distributively generated near-ring.*

In the theory of near-rings a major role is played by the notion called N-groups.

**DEFINITION 1.6.3:** *Let (P, +) be a group with 0 and let N be a near-ring. Let $\mu: N \times P \rightarrow P$, $(P, \mu)$ is called a N-group if for all $p \in P$ and for all $n, n_1 \in N$ we have $(n + n_1) p = (np + n_1 p)$ and $(nn_1) p = n(n_1 p)$. $N^p$ stands for N-groups.*

A subgroup M of a near-ring N with M.M ⊂ M is called a subnear-ring of N. A subgroup S of $N^p$ with N S ⊂ S is a N-subgroup of P.

**DEFINITION 1.6.4:** *Let N be a near-ring and P a N-group. A normal subgroup I of (N, +) is called an ideal of N if*

   i. IN ⊂ I.
   ii. *for all $n_1$, $n_2$ ∈ N and for all i ∈ I, $n(n_1 + i) - nn_1$ ∈ I.*

*Normal subgroup T of (N, +) with (i) is called right ideal of N while normal subgroup L of (N, +) with (ii) are called left ideals. A normal subgroup S of P is called ideal of $N^p$ if for all s ∈ P; $s_j$ ∈ S and for all n ∈ N; n(s + n) – ns ∈ S.*

*Factor near-ring N / I and factor N-subgroups P/S are defined as in case of rings.*

**DEFINITION 1.6.5:** *A subnear-ring M of N is called invariant if MN ⊂ M and NM ⊂ M. Near-rings are called simple if it has no ideals. $N^p$ is called N-simple if it has no N-subgroups except 0 and P.*

The notions of minimal and maximal ideals in a near-ring are defined as in case of rings.

**DEFINITION 1.6.6:** *Let X, Y be subsets of $N^p$ (X; Y) = {n ∈ N | nY ⊆ X} (0: X) is called the annihilator of X. We denote it by $(X: Y)_N$ as N is the related near-ring. $N^p$ is faithful if (0:P) = {0}.*

The concept of products, direct products and subdirect products in near-rings can be defined in an analogous way as in case of rings with very simple modifications.

The concept of ore conditions in near-rings is little different hence we define them.

**DEFINITION 1.6.7:** *The near-ring N is said to fulfill the left (right) ore condition with respect to a given subsemigroup S of (N, •) if for (s, n) ∈ S × N there exists $n_1$ ∈ N and $s_1$ ∈ S such that n • $s_1$ = s • $n_1$ ($s_1$ • n = $n_1$ • s).*

Now we proceed on to define free near-ring.

**DEFINITION 1.6.8:** *A near-ring $F_x$ ∈ V (V denotes the set of all near-rings) is called a free near-ring in V over N if there exists f: X → F (where X is any non empty set) for all N ∈ V and for all g: X → N there exists a homomorphism h ∈ Hom ($F_x$, N):*
h o f = g

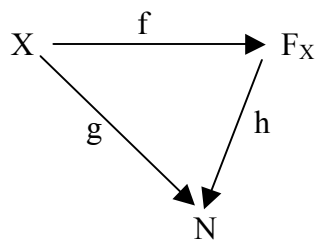

Thus if V = {set of all near-rings}, we simply speak about free near-ring on X. A near-ring is called free if it is free over some set X.

**DEFINITION 1.6.9:** *Let I be an ideal of the near-ring N. I is called a direct summand of N if there exists an ideal J in N such that N = I + J.*

**DEFINITION 1.6.10:** *A finite sequence $N = N_0 \supset N_1 \supset N_2 \supset ... \supset N_n = \{0\}$ of subnear-rings $N_i$ of N is called a normal sequence of N if and only if for all $i \in \{1, 2, ..., n\}$; $N_i$ is an ideal of $N_{i-1}$.*

The concept of A.C.C and D.C.C. can be defined in case of near-rings as in case of rings.

**DEFINITION 1.6.11:** *Let P be an ideal of the near-ring N. P is called a prime ideal if for all ideals I and J of N, $I J \subset P$ implies $I \subset P$ or $J \subset P$. N is called a prime near-ring if {0} is a prime ideal.*

*An ideal S of a near-ring N is called semiprime if and only if for all ideals I of N, $I^2 \subset S$ implies $I \subset S$. N is called a semiprime near-ring if {0} is a semiprime ideal. Let I be an ideal of N; The prime radical of I denoted by $P(I) = \cap P, I \subseteq P$.*

**DEFINITION 1.6.12:** *Let L be a left ideal in a near-ring N. L is called modular if and only if there exists $z \in N$ and for all $n \in N$, $n - nz \in L_z$. In this case we also say that $L_z$ is modular by z and that z is the right identity modulo L since for all $n \in N$; $nz = n$ (mod L). $z \in N$ is called quasi regular if $z \in L_z$, $S \subset N$ is called quasi regular if and only if for all $s \in S$; s is quasi regular.*

The concept of near polynomial rings is defined analogous to polynomial rings similarly matrix near-rings is also built.

For more about these concepts please refer [118, 126].

**DEFINITION 1.6.13:** *Let N be a near-ring. N is said to fulfill the insertion of factors property (IFP) provided that for all a, b, n $\in$ N we have ab = 0 implies anb = 0. N has strong IFP property if every homomorphic image of N has the IFP.*

*N has the strong IFP if and only if for all I in N and for all a, b, n $\in$ N. ab $\in$ I implies anb $\in$ I.*

**DEFINITION 1.6.14:** *Let N be a near-ring. A right ideal I of N is called right quasi reflexive if whenever A and B are ideals of N with $A.B \subset I$ then $b(b' + a) - bb' \in I$ for all $a \in A$ and for all $b, b' \in B$.*

*A near-ring N is strongly subcommutative if every right ideal of it is right quasi reflexive; any near-ring in this class of near-rings is therefore right quasi reflexive.*

Several interesting results using these concepts can be had from the book [43, 118, 126].

**DEFINITION [126]:** *Let N be a near-ring, S a subnormal subgroup of (N, +), S is called a quasi ideal of N if $SN \subset NS \subset S$ where by NS we mean elements of the form $\{n (n' + s) - nn'$ for all $s \in S$ and for all $n, n' \in N\} = NS$.*

**DEFINITION 1.6.15:** *Let N be a left near-ring; N is said to satisfy left self distributive identity if $abc = abac$ for all $a, b, c \in N$. A near-ring N is left permutable if $abc = bac$, right permutable if $abc = acb$, medial if $abcd = acbd$ and right self distributive if $abc = acbc$; for all $a, b, c \in N$.*

**DEFINITION 1.6.16:** *A near-ring R is defined to be equiprime if for all $0 \neq a \in R$ and $x, y \in R$, $arx = ary$ for all $r \in R$ implies $x = y$. If P is an ideal of R, P is called an equiprime ideal of R if R/P is an equiprime near-ring. It can be proved P is an ideal of R if and only if for all $a \in R \setminus P$, $x, y \in P$, $arx - ary \in P$ for all $r \in R$ implies $x - y \in P$.*

Now we define the concept of infra near-ring.

**DEFINITION 1.6.17:** *An infra near-ring (INR) is a triple $(N, +, \bullet)$ where*

  i. *$(N, +)$ is a group.*
  ii. *$(N, \bullet)$ is a semigroup.*
  iii. *$(x + y) \bullet z = x \bullet z - 0 \bullet z + y \bullet z$*

*for all $x, y, z \in N$.*

Several interesting results can be had from [43, 126].

We as in case of rings define n ideal near-ring. Further the marot near-ring N if each regular ideal of N is generated by a regular element of N. Properties of marot near-rings have been studied in [43, 118, 126].

The notions of right loop near-ring was introduced by [45-47].

**DEFINITION 1.6.18:** *The system $N = (N, +, \bullet, 0)$ is called a right loop half groupoid near-ring provided.*

  i. *$(N, +, 0)$ is a loop.*
  ii. *$(N, \bullet)$ is a half groupoid.*
  iii. *$(n_1 \bullet n_2) \bullet n_3 = n_1 \bullet (n_2 \bullet n_3)$ for all $n_1, n_2, n_3 \in N$ for which $n_1 \bullet n_2$, $n_2 \bullet n_3$, $(n_1 \bullet n_2) \bullet n_3$ and $n_1 \bullet (n_2 \bullet n_3) \in N$.*
  iv. *$(n_1 + n_2) \bullet n_3 = n_1 \bullet n_3 + n_2 \bullet n_3$ for all $n_1, n_2, n_3 \in N$ for which $(n_1 + n_2) \bullet n_3, n_2 \bullet n_3 \in N$.*

*(From now onwards we denote $n_1 \bullet n_2$ by $n_1 n_2$). If instead of (iv) N satisfies $n_1 \bullet (n_2 + n_3) = n_1 \bullet n_2 + n_1 \bullet n_3$ for every $n_1, n_2, n_3 \in N$ then we say N is a left loop half groupoid near-ring. We shall denote in any loop L the unique left and right inverse of a by $a_l$ and $a_r$ respectively.*

Now we recall the definition of right near-ring.

**DEFINITION [46]:** *The right loop near-ring N is a system, (N, +, •) of double composition '+' and '•' such that*

i. *(N, +) is a loop.*
ii. *(N, •) is a semigroup.*
iii. *The multiplication '•' is right distributive over addiction + i.e. for all $n_1, n_2, n_3 \in N$; $(n_1 + n_2) \bullet n_3 = n_1 \bullet n_3 + n_2 \bullet n_3$; here after by convention by a loop near-ring we mean only a right loop near-ring.*

*A non empty subset M of a loop near-ring (N, +, •, 0) is said to be a subloop near-ring of N if and only if (M, +, •, 0) is a loop near-ring. Further a loop near-ring N is said to be zero symmetric if and only if $n \bullet 0 = 0$ for every $n \in N$, where '0' is the additive identity. Thus a loop near-ring is said to be a near-ring if the additive loop of the loop near-ring is an additive group.*

The notions of zero divisor, left zero divisor, right zero divisor for loop near-rings are defined as in case of near-rings. The concept of loop near-ring homomorphism is also defined as in case of near-rings. Now we proceed on to define N- subloop.

**DEFINITION [46]:** *Let N be a loop near-ring. An additive subloop A of N is a N-subloop (right N-subloop) if $NA \subset A$ ($AN \subset A$) where $NA = \{na \mid a \in A, n \in N\}$. A non empty subset I of N is called a left ideal in N if*

i. *(L, +) is a normal subloop of (N, +).*
ii. *$n(n_1 + i) + n_r n_1 \in I$ for each $i \in I$ and $n, n_1 \in N$, where $n_r$ denotes the unique right inverse of n.*

*A loop near-ring N is said to be left bipotent if $Na = Na^2$ for every $a \in N$.*

The concepts of center, essential, simple, strictly prime etc can be defined in an analogous way.

As the loop near-rings are defined by [46] in this way when we define loops over near-ring we do not call it as loop near-rings but we call them as near loop rings.

**DEFINITION 1.6.19:** *Let (L, •, 1) be a loop and N a near-ring with a multiplicative identity. The loop over near-ring that is the near loop ring denoted by NL with identity is a non-associative near-ring consisting of all finite formal sums of the form $\alpha = \Sigma \alpha(m) m$; $m \in L$ and $\alpha(m) \in N$. (where supp $\alpha = \{m \mid \alpha(m) \neq 0$, the support of $\alpha$ is finite} satisfying the following operational rules.*

i. *$\Sigma \alpha(m) m = \Sigma \beta(m) m \Leftrightarrow \alpha(m) = \beta(m)$ for all $m \in L$.*
ii. *$\Sigma \alpha(m) m + \Sigma \mu(m) m = \Sigma (\alpha(m) + \mu(m)) m$; $m \in L$ $\alpha(m), \mu(m) \in N$.*
iii. *$(\Sigma \alpha(m)m)(\Sigma \mu(n)n) = \Sigma \gamma(k)k$; $m, n, k \in L$ where $\gamma(k) = \Sigma \alpha(m)\mu(n)$; $mn = k \in L$ (whenever they are distributive other wise it is assumed to be in NL).*
iv. *$n(\Sigma \alpha(m)n) = \Sigma n\alpha(m)m$ for all $n \in N$ and $m \in L$.*
v. *$\alpha(m) m = m \alpha(m)$ for all $\alpha(m) \in N$ and $m \in L$.*

*Dropping the zero components of the formal sum we may write $\alpha = \Sigma \alpha_i m_i$, $i = 1, 2,..., n$. Thus $m \to m.1$, is an embedding of N in NL. After the identification of N with N.1 we shall assume that N is contained NL. It is easily verified NL is a near-ring and it is a non-associative near-ring.*

Now we define NL by taking a loop L and a loop near-ring N, which is taken in this definition as loop near field.

**DEFINITION 1.6.20:** *Let L be a loop. N a loop near field. The loop-loop near-ring NL of the loop L over the loop near-ring N with identity is not even an associative loop but a doubly non-associative near-ring as (NL, +, •) is non-associative with respect to both '+' and '•'. NL consists of all finite formal sums of the form $\alpha = \Sigma\alpha(m)m$; $\alpha(m) \in N$, $m \in L$ such that support of $\alpha = \{m \mid \alpha(m) \neq 0\}$ is a finite set satisfying the following operational rules.*

  i. *$\Sigma \alpha(m) m = \Sigma \beta(m) m \Leftrightarrow \alpha(m) = \beta(m)$ for all $m \in L$. ( we assume in the finite formal sums the coefficient of each element in L occurs only once in the summation. (For otherwise $(\alpha g + \beta g) + \gamma g \neq \alpha g + (\beta g + \gamma g)$ if $(\alpha + \beta) + \gamma \neq \alpha + (\beta + \gamma)$).*

  ii. $\sum_{m \in L}\alpha(m)m + \sum_{m \in L}\beta(m)m = \sum_{m \in L}(\alpha(m) + \beta(m))m$.

  iii. *$\left(\sum \alpha(m)m\right)\left(\sum \beta(n)n\right) = \sum \gamma(k)k$; where $\gamma(k) = \Sigma \alpha(m) \beta(n)$ and $mn = k$; only when $\alpha(m)$'s are distributive in N if $\alpha(m)$'s are not distributive in N we assume the product is in NL.*

  iv. *$n\Sigma\alpha(m)m = \Sigma n\alpha(m)m$ if and only if $n \in N_d$. Dropping the zero components of the formal sum we may write $\alpha = \sum_i \alpha_i m_i$. Thus $n \to n. I_m$; $I_m$ is the identity in the loop L. So $N.1 = N \subset NL$. The element $1.1_m$ acts as the identity element of NL.*

Now we define still a new notion called Lield.

**DEFINITION [118]:** *A non empty set (L, 0, +, •) is called a Lield if*

  i. *L is an abelian loop under '+' with 0 acting as the identity.*
  ii. *L is a commutative loop with respect to '•'. 1 acts as the multiplicative identity.*
  iii. *The distributive laws $a • (b + c) = a • b + a • c$ and $(a + b) • c = a • c + b • c$*

*for all $a, b, c \in L$ holds goods.*

Still a new type of non-associative loop near domains was introduced by [126].

**DEFINITION [126]:** *The algebraic structure (D, +, •, 0, 1) is called a near domain if it satisfies the following axioms:*

    i.   $(D, +, 0)$ is a loop under '+'.
    ii.   $a + b = 0$ implies $b + a = 0$ for all $a, b \in D$.
    iii.   $(D^*, \bullet, 1)$ is a group where $D^* = D \setminus \{0\}$.
    iv.   $0 \bullet a = a \bullet 0 = 0$ for all $a \in D$.
    v.   $a \bullet (b + c) = a \bullet b + a \bullet c$ for all $a, b, c \in D$.
    vi.   For every pair $a, b \in D$ there exists $d_{a,b} \in D^*$ such that for every $x \in D$; $a + (b + x) = (a + b) + d_{a,b} x$.

Now [118 and 126] has defined loop near domains analogous to group rings.

**DEFINITION [118]:** *Let L be a finite loop under '+' and D be a near domain, the loop near domain DL contains elements generated by $d_i m_i$ where $d_i \in D$ and $m_i \in L$ where we admit only finite formal sums satisfying the following:*

    i.   $(DL, +, 0)$ is a loop under '+'.
    ii.   $d_i m_i d_j m_j = d_i d_j m_i m_j$ where we assume $m_i m_j \in DL$ as an element and $d_i d_j$, $d_i, d_j, d_k, d_k d_i \in D$.
    iii.   We assume $(m_i.m_j).m_k = m_i(m_j.m_k)$ where just the elements $m_i m_j m_k$ are juxtaposed and $m_i, m_j, m_k \in L$.
    iv.   We assume only finite sequence of elements of the form $m_1 m_2 \ldots m_k \in DL$.
    v.   $m_i m_j = m_j m_i$ for all $m_i m_j \in L$.
    vi.   $m_1 m_2 \ldots m_k = m'_1 m'_2 \ldots m'_k \Leftrightarrow m_i = m'_i$ for $i = 1, 2, 3, \ldots, k$.
    vii.   Since the loop L under '+' and D is also a loop under '+'; DL inherits the operation of '+' and '0' serves as the additive identity of the loop DL.
    viii.   Since $1 \in D$ we have only $1.L \subset DL$ so L is a subset of DL and $1.m_i = m_i.1 = m_i$ for all $i \in L$.

**DEFINITION 1.6.21:** *Let L be a finite loop and $Z_p$ be the near-ring. $Z_pL$ the loop near-ring of the loop L over the near-ring $Z_p$. The mod p-envelop of the loop $L^* = 1 + U$ where $U = \{\Sigma \alpha_i m_i \mid \Sigma \alpha_i = 0, \alpha_i \in Z_p\}$.*

When the near-rings are non-associative we have several natural identities to be true.

**DEFINITION 1.6.22:** *Let N be a non-associative near-ring we call N a Moufang near-ring if it satisfies any one of the following identities:*

    i.   $(xy)(zx) = (x(yz))x$.
    ii.   $((xy)z)y = x(y(zy))$.
    iii.   $x(y(xz)) = ((xy)x)z$

*for all x, y z in N. We call a near field N a Bruck near-ring if $(x(yx)z = x(y(xz))$ and $(xy)^{-1} = x^{-1} y^{-1}$ for all $x, y, z \in N$. N a non-associative near-ring N is said to be a Bol near-ring if $((xy)z)y = x((yz)y)$ for all $x, y, z \in N$. N is said to be right alternative if $(xy)y = x(yy)$ for all $x, y, \in N$; and left alternative if $(xx)y = x(xy)$. N is said to be alternative if N is simultaneously right and left alternative. A near-ring N is said to have weak-inverse property (WIP) if $(xy)z = 1$ implies $x(yz) = 1$ for all $x, y, z$ in N. N is a P-near-ring if $(xy)x = x(yx)$ for all $x, y \in N$.*

Groupoid near-rings is defined analogous to loop near-rings, where loops are replaced by groupoids. Several properties studied for loop near-ring can also be studied and extended to the case of groupoid near-ring. Several interesting properties about near-rings and seminear-rings can be had from [126]. The concept of seminear-rings can be studied and analogous properties about them can be obtained. For near-rings and seminear-rings please refer [43 and 126].

Now we proceed on to define the new notions of Smarandache near-rings [118]. Several new notions and those that of Smarandache notions can be applied to near-rings. Some of the concepts are recalled in this section.

**DEFINITION 1.6.23:** *N is said to be Smarandache near-ring (S-near-ring) if (N, +, •) is a near-ring and has a proper subset A such that (A, +, •) is a near field.*

*Example 1.6.1:* Let NG be any group near-ring where N is a near field. Then, N is a S-near-ring.

**DEFINITION 1.6.24:** *Let (N, +, •) be a S-near-ring. A non empty proper subset T of N is said to be a Smarandache subnear-ring (S-subnear-ring) if (T, +, •) is a S-near-ring; i.e. T has a proper subset which is a near field.*

**DEFINITION 1.6.25:** *Let (P, +) be a S-semigroup with 0 and let N be a S-near-ring. Let $\mu: N \times Y \to Y$ where Y is a proper subset of P which is a group under the operation of P; (P, $\mu$) is called the Smarandache N-group (S-N-group) if for all $y \in Y$ and for all $n, n_1 \in N$ we have $(n + n_1) y = ny + n_1 y$ and $(nn_1) y = n(n_1 y)$. $S(N^p)$ stands for the S- N-groups.*

**DEFINITION 1.6.26:** *A semigroup M of a near-ring N is called a Smarandache quasi subnear-ring (S-quasi subnear-ring) of N if $X \subset M$ where X is a subgroup of M which is such that $X.X \subset X$. A S-subsemigroup Y of $S(N^p)$ with $NY \subset Y$ is said to be a Smarandache N-subgroup (S-N-subgroup) of P. For any S-near-ring N, a normal subgroup I of (N, +) is called a Smarandache ideal (S-ideal) of N related to X if*

  i. $IX \subset I$.
  ii. $\forall x, y \in X$ and for all $i \in I$, $x(y + i) - xy \in I$.

*where X is the near field contained in N.*

Now we proceed on to define Smarandache homomorphism.

**DEFINITION 1.6.27:** *Let N and $N_1$ be two S-near-rings P and $P_1$ be S-N –subgroups.*

  i. *h: $N \to N_1$ is called a Smarandache near-ring homomorphism (S-near-ring homomorphism) if for all, $m, n \in M$, ($M \subset N$, $M_1 \subset N_1$) we have $h(m + n) = h(m) + h(n)$, $h(mn) = h(m)h(n)$, where $h(n)$ and $h(m) \in M_1$. It is to be noted that h need not even be defined on whole of N.*

ii.  $h: P \rightarrow P_1$ is called Smarandache N subgroup-homomorphism (S-N subgroup-homomorphism) if for all p, q in S ($S \subset P$; $S_1 \subset P_1$) and for all $m \in M \subset N$, $h(p + q) = h(p) + h(q)$ and $h(mp) = mh(p)$; $h(q), h(p)$ and $m\,h(p) \in S_1$.

*Here also we do not demand h to be defined on whole of P.*

**DEFINITION 1.6.28:** *A proper subset S of P is called a Smarandache ideal of $S(N^P)$ (S-ideal of $S(N^P)$) related to Y if*

  i.  *S is a S-normal subgroup of the S-semigroup P.*
  ii. *For all $s_1 \in S$ and $s \in Y$ (Y is the subgroup of P) and for $m \in M$ (M the near field of N) $n(s + s_1) - n\,s \in S$.*

*A S-near-ring is Smarandache simple (S-simple) if it has no S-ideals. $S(N^p)$ is called Smarandache N simple (S-N-simple) if it has no S-normal subgroups except 0 and P.*

**DEFINITION 1.6.29:** *A S-subnear-ring M of a near-ring N is called Smarandache invariant (S-invariant) related to the near field X in N if $MX \subset M$ and $XM \subset M$ where X is a S-near field of N.*

Thus in case of Smarandache invariance it is only a relative concept as a S-subnear-ring M may not be invariant related to every near field in the near-ring N. Let X and Y be S-semigroup of $S(N^P)$ $(X: Y) = \{n \in M \mid nY \subset X\}$ where M is a near field contained in N. $(0, x)$ is called the Smarandache annihilator (S-annihilator) of X.

Now we proceed on to define the notion of Smarandache direct product in near-rings and Smarandache free near-rings.

**DEFINITION 1.6.30:** *Let $\{N_i\}$ be a family of near-rings, which has atleast one S-near-ring. Then this direct product $N_1 \times ... \times N_r = \underset{i \in I}{\times} N_i$ with component wise defined operations '+' and '•' is called Smarandache direct product of near-rings (S-direct product of near-rings).*

**DEFINITION 1.6.31:** *Let N be a near-ring and S a S-subsemigroup of $(N, +)$. The near-ring $N_s$ is called the Smarandache near-ring of left (right) quotients of N (S-near-ring of left (right) quotients of N) with respect to S if*

  i.   *$N_s$ has identity.*
  ii.  *N is embeddable in $N_s$ by a homomorphism h.*
  iii. *For all $s \in S$; $h(s)$ is inveritible in $(N_s, •)$.*
  iv.  *For all $q \in N_s$, there exists $s \in S$ and there exists $n \in N$ such that $q = h(n)h(s)^{-1}$ ($q = h(s)^{-1}h(n)$).*

Now we proceed on to define Smarandache ore conditions.

**DEFINITION 1.6.32:** *The near-ring N is said to fulfill the Smarandache left (right) ore conditions (S-left (right) ore conditions) with respect to a given S-subsemigroup P of $(N, •)$ if for $(s, n) \in S \times N$ there exists $(s_1, n_1) \in S \times N$ such that $n • s_1 = s • n_1$ ($s_1 • n = n_1 • s$).*

**DEFINITION 1.6.33:** *A S-near-ring $F_X \in S(V)$ (X any non empty subset which is S-semigroup under '+' or '•' is called a Smarandache free near-ring (S-free near-ring) in V over X if there exists $f: X \to F$ for all $N \in V$ and for all $g: X \to N$ there exists a S-near-ring homomorphism $h \in S(\text{Hom}(F_X, N)))$, [Here $S(\text{Hom}(F_X, N))$ denotes the collection of all S-homomorphism from $F_X$ to N] such that $h \circ f = g$.*

Smarandache internal direct sum and Smarandache direct summand are defined as in case of near-rings with very simple modifications that ideals in a near-ring are replaced by S-ideals of a near-ring.

Smarandache normal sequence and Smarandache invariant sequence are defined by replacing the subnear-rings by S-subnear-rings and ideals by S-ideals.

**DEFINITION 1.6.34:** *Let P be a S-ideal of the near-ring N. P is called Smarandache prime ideal (S-prime ideal) if for all S-ideals I and J of N, $IJ \subset P$ implies $I \subset P$ or $J \subset P$.*

*We say the S-ideal P is Smarandache semiprime (S-semiprime) if and only if for all S-ideals I of N $I^2 \subset P$ implies $I \subset P$. Let I be a S-ideal of N; the Smarandache prime radical (S-prime radical) of I denoted by $S(P(I)) = \bigcap_{I \subset P} P$. Thus if $S(P(I))$ is a S-prime radical it has prime radical.*

The concept of S-units, S-idempotents and S-zero divisors in case of near-rings can be defined as in case of rings as it is immaterial whether the ring is S-near-ring or not.

**DEFINITION 1.6.35:** *Let P be a S-left ideal of a near-ring N. P is called Smarandache modular (S-modular) if and only if there exists a S-idempotent $e \in N$ and for all $n \in N$; $n - n e \in P$.*

Now we proceed on to define Smarandache quasi regular elements in N.

**DEFINITION 1.6.36:** *$z \in N$ is called Smarandache quasi regular (S-quasi regular) if $z \in L_z$. An S-ideal $P \subset N$ is called S-quasi regular if and only if for all $s \in S$, s is S-quasi regular.*

Concepts like Smarandache principal ideals, Smarandache biregular and Smarandache prime radical are defined by replacing ideals by S-ideals.

Next we define Smarandache mixed direct product of near-rings and seminear-rings.

**DEFINITION 1.6.37:** *A non empty set S together with two binary operations '+' and '•' is said to be a seminear-ring if $(S, +)$ and $(S, •)$ are semi groups and for all s, s', s'' in S we have $(s + s'') • s'' = s • s'' + s' • s''$.*

**DEFINITION 1.6.38:** *$N = N_1 \times ... \times N_n$ is said to be the Smarandache mixed direct product I (S-mixed direct product I) of seminear-rings and near-rings, where $N_i$'s are from the classes of near-rings and seminear-rings. We call $N = N_1 \times N_2 \times ... \times N_n$* *where some of the $N_i$'s are near-rings and some of them are near fields as the Smarandache mixed direct product II (S-mixed direct product II). For $N = N_1 \times ... \times N_n$ where some of the $N_i$'s are rings and some of the $N_i$'s are near-rings as Smarandache mixed direct product III (S-mixed direct product III).*

*$N_1 \times ... \times N_n$ will be called the Smarandache mixed direct product IV (S-mixed direct product IV) if some $N_i$'s are seminear-rings and some of them are semirings.*

Using the concept of S-mixed direct product III we define S-near-ring of level II.

**DEFINITION 1.6.39:** *Let N be near-ring. N is said to be a S-near-ring of level II if N has a proper subset P, $P \subset N$ where P is a ring. We call N a S-seminear-ring of level II if N contains a proper subset P which is a semiring.*

**DEFINITION 1.6.40:** *$N_1 \times ... \times N_n$ is called the Smarandache mixed direct product V (S-mixed direct product V) if some of the $N_i$'s are non-associative seminear-rings, rings and non-associative near-rings. It is easily verified N is a S-quasi near-ring.*

**DEFINITION 1.6.41:** *Let N be a S-near-ring. N is said to fulfill the Smarandache insertion of factors property (S-insertion of factors property) if for all a, b $\in$ N we have ab = 0 implies anb = 0 for all n $\in$ P, P $\subset$ N where P is a near field. We say N, the S-near-ring has Smarandache strong IFP- property (S-strong IFP- property) if and only if for all a, b $\in$ N ab$\in$I implies anb $\in$ I where n $\in$ P, P $\subset$ N and P is a near-field.*

**DEFINITION 1.6.42:** *Let p be a prime. A S-near-ring N is called a Smarandache p-near-ring (S-p-near-ring) provided for all x $\in$ P, $x^p$ = x and px = 0 where P $\subset$ N and P is a semifield.*

**DEFINITION 1.6.43:** *Let N be a near-ring. A S-right ideal I of N is called Smarandache right quasi reflexive (S- right quasi reflexive) if whenever A and B are S-ideals of N with AB $\subset$ I then b (b' + a) – bb' $\in$ I for all a $\in$ A and for all b, b' $\in$ B. N is said to be Smarandache strongly subcommutative (S-strongly subcommutative) if every S-right ideal of it is S-right quasi reflexive.*

**DEFINITION 1.6.44:** *Let N be a near-ring. S a S-subnormal subgroup of (N, +). S is called a Smarandache quasi ideal (S-quasi ideal) of N if SN $\subset$ S and NS $\subset$ S where by NS we mean elements of the form {n (n' +s) – nn' / for all s $\in$ S and for all n, n' $\in$ N} = NS.*

The concept of Smarandache right/left self distributive, Smarandache left permutable, right permutable, Smarandache medial and Smarandache equiprime ideal can be defined in near-rings with suitable modifications.

We just define Smarandache infra near-ring.

**DEFINITION 1.6.45:** *Let (N, +, •) be a triple. N is said to be Smarandache infra near-ring (S-INR) if*

i. (S, +) is a S-semigroup
ii. (N, •) is a semigroup.
iii. (x + y) • z = xz – 0z + yz

for all x, y, z ∈ N.

Several properties of S-INR near-ring can be defined and analyzed. We build more and more near-rings using groups and semigroups which will be called as group near-rings and semigroup near-rings analogous to group rings and semigroup rings. Such study is very meager in near-ring literature.

We call a group near-ring NG of a semigroup S over the near-ring N to be a Smarandache group near-ring (S-group near-ring) if N is a S-near-ring. Thus we see all group near-rings need not in general be a S-group near-ring. Further a group near-ring can be S-near-ring, yet fail to be a S-group near-ring. Several analogous properties true in case of S-group rings can also be studied in case of S-group near-rings. The challenging isomorphism problem, zero divisor conjecture and semisimple problem in case of group near-rings and S-group near-rings remains at a very dormant state.

Now the concept of Smarandache semigroup near-ring is defined in the following manner:

The semigroup near-ring NS is a Smarandache semigroup near-ring if and only if the semigroup S is a S-semigroup. Thus all semigroup near-rings even if they are S-near-rings need not be a S-semigroup near-rings.

All results, which are studied in case of S-group near-rings, can be analyzed in case of S-semigroup near-rings and semigroup near-rings. The three classical problems viz. zero divisor conjecture, semisimplicity problems and the isomorphism problems remains open in case of both semigroup near-rings and S-semigroup near-rings.

Also instead of using near-rings we can use seminear-rings and construct group seminear-rings, semigroup seminear-rings and also their Smarandache analogue. The concept of Smarandache commutative, cyclic, weakly commutative and weakly cyclic are defined in the near-ring. Every one of its S-subnear-ring is S-commutative, S-cyclic or at least one of the S-subnear-ring is cyclic or commutative then we call the near-ring to be S-weakly cyclic or S-weakly commutative. Several interesting results can be determined.

**DEFINITION 1.6.46:** *Let N be any seminear-ring. We say N is a Smarandache seminear-ring II (S-seminear-ring II) if N has a proper subset P where P is a semiring. If the S-seminear-ring has a strict semiring then we call N a Smarandache strict seminear-ring (S-strict seminear-ring). The S-seminear-ring N is Smarandache commutative (S-commutative) if every semiring in N is commutative. We say N is Smarandache weakly commutative (S-weakly commutative) if atleast one semiring in N is commutative.*

**DEFINITION 1.6.47:** *Let $N_1$ and $N_2$ be any two S-seminear-ring II we say a map $\phi: N_1 \to N_2$ is a Smarandache seminear-ring homomorphism II (S-seminear-ring* *homomorphism II) if φ: $A_1 \to A_2$ where $A_1$ and $A_2$ are semirings contained in $N_1$ and $N_2$ respectively and φ is a semiring homomorphism; φ need not be defined on whole of $N_1$.*

Now we proceed onto define Smarandache pseudo seminear-ring.

**DEFINITION 1.6.48:** *N is said to be Smarandache pseudo seminear-ring (S-pseudo seminear-ring) if N is a near-ring and has a proper subset A of N, which is a seminear-ring under the operations of N.*

*The Smarandache pseudo seminear-ring homomorphism (S-pseudo seminear-ring homomorphism) from two S-pseudo seminear-rings $S_1$ to $S_2$ is defined as a near-ring homomorphism from $N_1$ to $N_2$ where $N_1$ and $N_2$ are proper subsets of $S_1$ and $S_2$ respectively and they are near-rings. φ : $N_1 \to N_2$ is a near-ring homomorphism; φ need not be even defined on whole of $S_1$.*

For more about these please refer [118, 126].

Now the concept of Smarandache integral, equiprime and left infra near-rings are defined in an analogous and appropriate ways. These studies are solely assigned to the reader.

**DEFINITION 1.6.49:** *A Smarandache composition near-ring (S-composition near-ring) is a quadruple (C, +, o) where (C, +, •) and C, +, •) are S-near-rings such that (a • b) o c = (a o b) • c for all a, b, c ∈ C.*

**DEFINITION 1.6.50:** *A non zero S-ideal H of G is said to be Smarandache uniform (S-uniform) if for each pair of S-ideals $K_1$ and $K_2$ of G such that $K_1 \cap K_2$ = (0); $K_1 \subset H$, $K_2 \subset H$ implies $K_1$ = (0) or $K_2$ = (0).*

**DEFINITION 1.6.51:** *Let N be a S-near-ring an element a∈ P is a called Smarandache normal element (S-normal element) of N if aN = Na. If aN = Na for every a ∈P then we call N a Smarandache normal near-ring (S-normal near-ring).*

*$S_n(N)$ denotes the set of all S-normal elements of N. $S_n(N)$ = P if and only if N is S-normal near-ring.*

The concept of Smarandache normal subnear-ring can also be defined for any S-subnear-ring A of N. Further if N is a S-near-ring, A be a S-subsemigroup S of N is called a Smarandache normal subsemigroup (S-normal subsemigroup) of S if n S = Sn for all n ∈ P ⊂ N.

**DEFINITION 1.6.52:** *Let N be a S-near-ring. A S-subnear-ring M of N is Smarandache invariant (S-invariant) with respect to an element a ∈ P \ M if a M ⊂ M and Ma ⊂ M. If aM = M and Ma = M we say M is Smarandache strictly invariant (S-strictly invariant) with respect to a. (P is a near field).*

**DEFINITION 1.6.53:** *A commutative S-near-ring N with identity is called a Smarandache Marot near-ring (S-Marot near-ring) if each regular S-ideal of N is* *generated by regular elements where by regular S-ideal we mean an ideal with regular elements that is only non zero divisors of the near-ring N.*

**DEFINITION 1.6.54:** *Let K and I be S-ideals of the S-N-subgroup G. K is said to be the Smarandache complement (S-complement) of I if the following two conditions hold good.*

  i. $K \cap I = (0)$ and
  ii. $K_1$ is an S-ideal of G

*such that $K \subset K_1$ then $K_1 \neq K$ imply $K_1 \cap I \neq (0)$.*

**DEFINITION 1.6.55:** *A subset S of the S-N-subsemigroup G is said to be Smarandache small (S-small) in G if*

  i. *$S + K = G$ where K is an S-ideal of G imply $K = G$.*
  ii. *G is said to be Smarandache hollow (S-hollow) if every proper S-ideal of G is S-small in G and*
  iii. *G is said to have Smarandache finite spanning dimension (S-finite spanning dimension) if for any decreasing sequence of Smarandache NS-subsemigroups (S-NS-subsemigroups);*

*$X \supset X_1 \supset ...$ of G such that $X_i$ is an S-ideal of $X_{i-1}$ there exists an integer such that $X_j$ is S-small in G for all $j \geq k$. Suppose H and K are two S-N-subsemigroup of G. K is said to be Smarandache supplement (S-supplement) for H if $H + K = G$ and $H + K_i \neq G$ for any S-ideal $K_i$ of K. We call N a Smarandache ideally strong (S-ideally strong) if for every S-subnear-ring of N is a S-ideal of N.*

**DEFINITION 1.6.56:** *Let N be a near-ring $\{I_k\}$ be the collection of all S-ideals of N. N is said to be a Smarandache $I^*$ - near-ring (S-$I^*$-near-ring) if for every pair of ideals $I_1, I_2 \in \{I_j\}$ we have for every $x \in N \setminus (I_1 \cup I_2)$ the S-ideal generated by x and $I_1$ and x and $I_2$ are equal, that is $\langle x \cup I_1 \rangle = \langle I_2 \cup x \rangle$.*

**DEFINITION 1.6.57:** *Let S be a set with two binary operations '+' and '•'. (S, +, •) is called a S-seminear-ring if (S, +) and (S, •) are S-semigroups and for all s, s', s'' $\in$ S we have $(s + s') \cdot s'' = s \cdot s' + s' \cdot s''$. A subset P of a S-seminear-ring S is said to be a Smarandache subseminear-ring (S-subseminear-ring) if (P, +, •) is itself a S-seminear-ring.*

*If a S-seminear-ring S has no proper S-subseminear-ring then we call S, Smarandache simple seminear-ring (S-simple seminear-ring). A seminear-ring is commutative if every S-subseminear-ring is commutative. If atleast one S-subseminear-ring is commutative then we call the seminear-ring to be Smarandache weakly commutative seminear-ring (S-weakly commutative seminear-ring).*

**DEFINITION 1.6.58:** *Let (N, +, •) be a seminear-ring. We say N has a Smarandache hyper subseminear-ring (S-hyper subseminear-ring).*

  i. *If (S, +) is a S-semigroup.*

ii. *If A is a proper subset S which is a subsemigroup of S and A contains the largest group of (S, +).*
  iii. *(A, •) is a S-subsemigroup or a group.*

*The Smarandache dual hyper subsemigroup (S-dual hyper subsemigroup) is defined as a proper subset P of a seminear-ring (S, +, •) such that*

  i. *(P, •) is a S-semigroup.*
  ii. *If A is a proper subset of P which is a subsemigroup of P and A contains the largest group of (P, •).*
  iii. *(A, +) is a S-subsemigroup or a group.*

**DEFINITION 1.6.59:** *Let (N, +, •) be a seminear-ring. A proper subset A of N which is a group under the operations of both '+' and '•' is called a Smarandache normal subseminear (S-normal subseminear) if $xA \subseteq A$ and $Ax \subset A$ or $xA = \{0\}$ and $Ax = \{0\}$ for all $x \in N$.*

*A seminear-ring, which has no S-normal subseminear-ring, is called Smarandache pseudo simple (S-pseudo simple). We define Smarandache quotient near-ring (S-quotient near-ring) $N/A = \{Ax \mid x \in N\}$, where A is a S-normal subseminear-ring.*

The concept of Smarandache maximal or minimal can be defined analogous to near-rings.

The notions of Smarandache seminear-ring homomorphism can be defined as in case of near-rings. Further the reader is left with the task of defining S-left ideals, S-right ideals and S-ideals in a S-seminear-ring.

We recall the notion of bipotent seminear-ring and their Smarandache analogues.

**DEFINITION 1.6.60:** *Let N be a S-seminear-ring. N is said to be Smarandache left bipotent (S-left bipotent) if $Pa = Pa^2$ for all $a \in N$ ($P \subset N$ is subgroup of (N, +)). N is said to be Smarandache dually left bipotent (S-dually left bipotent) if $B + a = B + a^2$ for all $a \in N$. ($B \subset N$ is a subgroup of (N, •)). N is said to be Smarandache strongly bipotent (S-strongly bipotent) if*

  i. $Pa = Pa^2$ *and*
  ii. $P + a = P + a^2$

*where (P, +) is subgroup of (N, +) and (P, •) is a subgroup of (N, •).*

We proceed on to define Smarandache seminear-ring II.

**DEFINITION 1.6.61:** *A non empty set N is said to be a Smarandache seminear-ring II (S-seminear-ring II) if (N, +, •) is a seminear-ring having a proper subset A such that $(A \subset N)$; A under the same binary operations of N is a near-ring that is (A, +, •) is a near-ring.*

**DEFINITION 1.6.62:** *N is said to be Smarandache pseudo seminear-ring (S-pseudo seminear-ring) if N is a near-ring and has a proper subset A of N such that A is a seminear-ring under the operations of N.*

**DEFINITION 1.6.63:** *Let N and $N_1$ be two S-pseudo seminear-rings; h: $N_1 \to N_2$ is a Smarandache pseudo seminear-ring homomorphism (S-pseudo seminear-ring homomorphism) if h restricted from A to $A_1$ is a seminear-ring homomorphism.*

**DEFINITION 1.6.64:** *Let N be a S-seminear-ring II. An additive subgroup A of N is called a Smarandache N-subgroup II (S-N-subgroup II) if $NA \subset N$ ($AN \subset A$) where $NA = \{na \mid a \in A; n \in N\}$.*

Concepts of S-ideal, S-left bipotent and S-regular can also be defined in an analogous way.

**DEFINITION 1.6.65:** *A seminear-ring N is said to be Smarandache strictly duo (S-strictly duo) if every SN-subgroup is also a right SN-subgroup. A S-seminear-ring N is called Smarandache irreducible (S-irreducible) i.e. S-simple if it contains only trivial SN-subgroups. We call an S-ideal $P \subset N$ to be Smarandache strictly prime (S-strictly prime) if for any two SN-subgroups A and B of $A_1$ ($A_1 \subset N$, $A_1$ a near-ring of the S-seminear-ring) $AB \subset P$ then $A \subset P$ or $B \subset P$. A S-left ideal B of a S-seminear-ring, N is called Smarandache strictly essential (S-strictly essential) if $B \cap K \neq (0)$ for every non-zero SN-subgroup K of N.*

**DEFINITION 1.6.66:** *A seminear-ring N is strongly subcommutative if every left ideal of it is right quasi reflexive. We call N is Smarandache strongly subcommutative (S-strongly subcommutative) seminear-ring if every S-right ideal of it is right quasi reflexive.*

**DEFINITION 1.6.67:** *Let R be a S-seminear-ring. We call a non empty subset I of R to be a Smarandache common ideal related to A (S-common ideal related to A) ($A \subset R$) and A a near-ring if*

    i.   $x + y \in I$ for all $x, y \in I$.
    ii.  *ax, anx, xa $\in I$.*

*for all $x \in I$ and $a \in A$.*

The concept of Smarandache prime and Smarandache S-ideals in seminear-rings can be defined as in case of near-rings.

**DEFINITION 1.6.68:** *Let $S = [Q, X, \delta]$ be a S-S-semigroup semiautomaton where Q is a S-semigroup under '+' i.e. $G \subset Q$ is such that G is a S-semigroup under '+' and G a proper subset of Q. Let $\delta: G \times X \to G$ is well defined; that is restriction of $\delta$ to G is well defined so $(G, X, \delta) = C$ becomes a semiautomaton; C is called additive if there is some $x_0 \in X$ with $\delta(q, x) = \delta(q, x_0) + \delta(0, x)$ and $\delta(q - q', x_q) = \delta(q, x) - \delta(q', x_0)$, for all $q, q' \in G$ and $x \in X$ then there is some homomorphism*

    *$\psi: G \to G$ and some map $\alpha: X \to G$ with $\psi(x_o) = 0$ and $\delta((q, x) = \psi(q) + \alpha(x)$.*

*Let $\delta_x$ for a fixed $x \in X$ be the map from $G \to G$ ; $q \to \delta_x(q, x)$. Then $\{\delta_x \mid x \in X\}$ generates a subnear-ring N(P) of the near-ring (M(G), +, •) of all mappings on G; this near-ring N(P) is called the Smarandache syntactic near-ring (S-syntactic near-ring).*

For notions and concepts about syntactic near-rings please refer [14, 126]. Thus the S-syntactic near-ring for varying S-semigroups enjoy varied properties.

The most interesting feature about these S-syntactic near-ring is that for one S-S-semigroup semiautomaton we have several group semiautomaton depending on the number of valid groups in the S-semigroup. This is the vital benefit in defining S-S-semigroup automaton and S-syntactic near-ring. So far a given S-S-semigroup semiautomaton we can have several S-syntactic near-ring.

A planar near-ring can be used to construct balanced incomplete block design [BIBD] of high efficiency, we just state here how they are used in developing error correcting codes, as codes which can correct errors is always desirable than codes which can only detect errors.

**DEFINITION 1.6.69:** *A balanced incomplete block design (BIBD) with parameters (v, b, r, k, $\lambda$) is a pair (P, B) with the following properties:*

  i.  *P is a set with v elements.*
  ii. *B = ($B_1$, ..., $B_b$) is a subset of p(P) with b elements.*
  iii. *Each $B_i$ has exactly k elements where k < v each unordered pair (p, q) with p, q, $\in$ P, p $\neq$ q occurs in exactly $\lambda$ elements in B.*

*The set $B_1$, ..., $B_b$ are called the blocks of BIBD. Each a $\in$ P occurs in exactly r sets of B. Such a BIBD is also called a (v, b, r, k, $\lambda$) configuration or 2-(v, k, $\lambda$) tactical configuration or design. The term balance indicates that each pair of elements occurs in exactly the same number of block, the term incomplete means that each block contains less than v - elements. A BIBD is symmetric if v = b.*

*The incidence matrix of a (v, b, r, k, $\lambda$) configuration is the v $\times$ b matrix A = ($a_{ij}$) where*

$$a_{ij} = \begin{cases} 1 \text{ if } i \in B_j \\ 0 \text{ otherwise} \end{cases}$$

*here i denotes the $i^{th}$ element of the configuration. The following conditions are necessary for the existence of a BIBD with parameters v, b, r, k, $\lambda$*

  i.  *bk = rv.*
  ii. *r (k – 1) = $\lambda$ (v – 1).*
  iii. *b $\geq$ v.*

*Recall a near-ring N is called planar or Clay if for all equation x o a = x o b + c, (a, b, c, $\in$ N, a $\neq$ b) have exactly one solution x $\in$ N.*

The study and construction of BIBD from a planar near-ring is given in [118].

The non-associative near-ring and its Smarandache equivalent is given in the following:

**DEFINITION 1.6.70:** *Let (N, +, •) be a non empty set endowed with two binary operations '+' and '•' satisfying the following:*

   i.   *(N, +) is a semigroup.*
   ii.  *(N, •) is a groupoid.*
   iii. *(a + b) • c = a • c + b • c for all a, b, c ∈ N.*

*(N, +, •) is called the right seminear-ring which is non-associative. If (iii) is replaced by a • (b + c) = a • b + a • c for all a, b, c ∈ N then we call (N, +, •) a left seminear-ring.*

Subseminear-ring, ideals and N-subsemigroups are defined analogous to seminear-ring. Now we proceed on to define non-associative Smarandache seminear-ring.

**DEFINITION 1.6.71:** *Let (N, +, •) be a non-associative seminear-ring. N is said to be a Smarandache non-associative seminear-ring of Level I (S-non-associative seminear-ring of Level I) if*

   i.   *(N, +) is a S-semigroup.*
   ii.  *(N, •) is a groupoid.*
   iii. *(a + b) • c = a • c + b • c*

*for all a, b, c ∈ N.*

*N is said to have a Smarandache subseminear-ring (S-subseminear-ring) P ⊂ N, if P is itself a S-seminear-ring. An additive S-subsemigroup A of N is called the Smarandache N-subsemigroup (S-N-subsemigroup) NA ⊂ A (AN ⊂A) where NA={na /n ∈ N, a ∈ A}. A non empty subset I of N is called Smarandache left ideal (S-left ideal) in N if*

   i.   *(I, +) is a normal S-subsemigroup of (N, +).*
   ii.  *n(n₁ + i) +nn₁ ∈ I for each i∈I, n, n₁ ∈ N.*

*A non-empty subset I of N is called a Smarandache ideal (S-ideal) if*

   i.   *I is a S-left ideal.*
   ii.  *IN ⊂ I.*

*The Smarandache seminear-ring which is non-associative is denoted by S-NA seminear-ring.*

**DEFINITION 1.6.72:** *Let (N, +, •) be a non-associative seminear-ring. N is said to be a Smarandache seminear-ring I of type A (S-seminear-ring I of type A) if N has a proper subset P such that (P, +, •) is an associative seminear-ring.*

*N is said to be Smarandache NA-seminear-ring I of type B (SNA-seminear-ring I of type B) if N has a proper subset P where P is a near-ring.*

**DEFINITION 1.6.73:** *Let $N_1$ and $N_2$ be two S-seminear-rings, we say a map $\phi$ from $N_1$ to $N_2$ is a Smarandache non-associative seminear-ring homomorphism (SNA-seminear-ring homomorphism) if*

$$\phi(x + y) = \phi(x) + \phi(y)$$
$$\phi(xy) = \phi(x)\,\phi(y).$$

Now we proceed on to define Smarandache non-associative seminear-ring II.

**DEFINITION 1.6.74:** *Let $(N, +, \bullet)$ be a NA-seminear-ring we say N is a Smarandache NA seminear-ring II (SNA-seminear-ring II) if N has a proper subset P which is a associative seminear-ring.*

Now how to get examples of such S-NA seminear-rings; this is mainly done by defining the concept of Smarandache mixed direct products.

**DEFINITION 1.6.75:** *Let $N_1, N_2, \ldots, N_k$ be k-seminear-rings where atleast one of the seminear-rings is non-associative and atleast one is associative. The direct product of these seminear-rings, $N = N_1 \times \ldots \times N_k$ is called the S-mixed direct product of seminear-rings. This S-mixed direct product gives S-seminear-ring II and S-seminear-ring I of type A.*

Now we define yet another notion called Smarandache strong mixed direct product in the following.

**DEFINITION 1.6.76:** *Let $N_1, N_2, \ldots, N_k$ be a collection of near-rings, non-associative seminear-rings and seminear-rings. The product $N = N_1 \times N_2 \times \ldots \times N_k$ is called the Smarandache strong mixed direct product of NA seminear-rings (S-strong mixed direct product of NA seminear-rings) is a S-NA seminear-ring I of type A and B and also S-NA-seminear-ring II; under component wise operations on N.*

Several interesting results can be had about these new structures. All properties like S-NA-subseminear-ring, S-NA-ideals, S-NA seminear-ring homomorphism etc can be defined in an analogous way as in case of S-seminear-rings.

More properties about these structures can be had from [118].

Now we proceed on to define Smarandache non-associative near-rings.

**DEFINITION 1.6.77:** *Let $(N, +, \bullet)$ be a non empty set. N is said to be a Smarandache non-associative near-ring (SNA-near-ring) if the following conditions are satisfied.*

    i. *$(N, +)$ is a S-semigroup.*
    ii. *$(N, \bullet)$ is a S-groupoid.*
    iii. *$(a + b) \bullet c = a \bullet c + b \bullet c$*

*for all a, b, c ∈ N. If (N, +, •) has an element i ∈ N such that 1 • a = a • 1 = a for all a ∈ N then we all N a Smarandache non-associative near-ring with unit. If (N, •) is a S-loop barring 0 we call N a Smarandache non-associative division near-ring (SNA-division near-ring).*

*If (N \ {0}, •) is a commutative S-loop then we call (N, +, •) a Smarandache non-associative near field (SNA-near-field). We call the S-non-associative near-ring (S-NA-near-ring) to be a Moufang or Bruck or Bol or WIP or alternative as in case of any non-associative ring.*

**DEFINITION 1.6.78:** *Let N be a non empty set. Define two binary operation '+' and '•' satisfying the following conditions:*

   i.   *(N, +) is a S-semigroup*
   ii.  *(N, •) is a groupoid*
   iii. *(a + b) • c = a • c + b • c*

*for all a, b, c ∈ N then we known N is a Smarandache NA-seminear-ring of level I (SNA-seminear-ring of level I).*

*We call (N, +, •) a non-associative seminear-ring N to be a Smarandache NA-seminear-ring of level II (SNA-seminear-ring of level II) if N has a proper subset P, P ⊂ N such that (P, +, •) is a non-associative near-ring.*

All notions like S-pseudo NA-seminear-rings, ideals subsemi-near-rings etc, can be defined in an analogous way.

**DEFINITION 1.6.79:** *Let N be a NA-seminear-ring. We say N is a Smarandache NA-seminear-ring III (SNA-seminear-ring III) if N has a proper subset P where P is an associative seminear-ring.*

*We say N is a Smarandache NA-seminearing level IV (SNA-seminearing level IV) if N has a proper subset P such that (P, +, •) is an associative near-ring.*

Several interesting properties about these S-NA-seminear-rings III and IV are dealt in [118, 126]. We have several classes of S-NA seminear-rings viz. groupoid near-rings. Based on these types of S-NA-seminear-rings we can define Smarandache homomorphism II, III and IV.

Many related concepts on S-NA-seminear-rings can be built. Several properties like Smarandache equiprime, Smarandache strongly semiprime, Smarandache pseudo seminear-ring (SNP ring), Smarandache quasi regular and concept of Smarandache Jacobson radical etc. can be defined and these notions can be had from the book [118, 126].

## 1. 7 Vector Spaces and S-vector spaces

In this section we recall and introduce the interesting concept of Smarandache vector spaces and just recall the basic concepts of vector spaces so that when we define bivector spaces and Smarandache bivector spaces, the reader has a good background of them. For more literature about vector spaces please refer [27, 44, 49, 119] and for the concept of Smarandache vector space which is very new and can be found only in the book [44, 49, 119].

**DEFINITION 1.7.1:** *A vector space consists of the following:*

  i. *a field of scalars.*
 ii. *a set V of objects called vectors.*
iii. *a rule called vector addition which assigns to every pair of vectors $\alpha$, $\beta$ in V, a vector $\alpha + \beta$ in V called the sum of $\alpha$ and $\beta$ in such a way that*

   a. *addition is commutative $\alpha + \beta = \beta + \alpha$.*
   b. *$\alpha + (\beta + \gamma) = (\alpha + \beta) + \gamma$ for all $\alpha, \beta, \gamma \in V$.*
   c. *There is a unique vector 0 in V called the zero vector such that $\alpha + 0 = \alpha$ for all $\alpha$ in V.*
   d. *A rule called scalar multiplication which associates with each scalar c in F and vector $\alpha$ in V a vector $c\alpha$ in V, called the product of c and $\alpha$ in such a way that*

   1. *$1\alpha = \alpha$ for every $\alpha$ in V.*
   2. *$(c_1 c_2)\alpha = c_1(c_2\alpha)$.*
   3. *$c(\alpha + \beta) = c\alpha + c\beta$.*
   4. *$(c_1 + c_2)\alpha = c_1\alpha + c_2\alpha$*

*where c, $c_1$, $c_2$ $\in$ F and $\alpha$, $\beta$ $\in$ V.*

*V is called the vector space over the field F. The term 'vector' is applied to the elements of the set V largely as a matter of convenience.*

***Example 1.7.1***: F be a field, F[x] the polynomial ring; F[x] is a vector space over F.

***Example 1.7.2***: Let V = F × F × F where F is a field. V is a vector space over F.

***Example 1.7.3***: Let $M_{n \times m}$ = {$(a_{ij})$ / $a_{ij}$ $\in$ Q}, set of n × m martrices with entries from Q; Q the field of rationals, $M_{n \times m}$ is a vector space over Q.

**DEFINITION 1.7.2:** *Let V be a vector-space over the field F. A subspace of V is a subset W of V which is itself a vector space over F with the operations of vector addition and scalar multiplication on V.*

**DEFINITION 1.7.3:** *Let S be a set of vectors in a vector space V. The subspace spanned by S is defined to be the intersection W of all subspaces of V which contains S. When S is a finite set say S = {$v_1$, ..., $v_n$} we shall call W the subspace spanned by the vectors $v_1$, $v_2$, ..., $v_n$.*

**DEFINITION 1.7.4:** *Let V be a vector space over F. A subset S of V is said to be linearly dependent (or simply dependent) if there exists distinct vectors $v_1$, $v_2$, ..., $v_n$ in S and scalars $c_1$, $c_2$, ..., $c_n$ in F not all of which are 0 such that $c_1\alpha_1 + ... + c_n\alpha_n = 0$. A* *set which is not linearly dependent is called linearly independent. If the set S contains only finitely many vectors $v_1, v_2, ..., v_n$ we sometimes say that $v_1, v_2, ..., v_n$ are dependent (or independent) instead of saying S is dependent (or independent).*

**DEFINITION 1.7.5:** *Let V be a vector space. A basis for V is a linearly independent set of vectors in V, which spans the space V. The space V is finite dimensional if it has a finite basis. If V does not contain a finite basis then we say V is infinite dimensional. We say a set of vectors $v_1, ..., v_t$ are linearly dependent if we can find non-zero scalars $c_1, ..., c_t$ not all zero such that $c_1v_1 + ... + c_tv_t = 0$.*

**DEFINITION 1.7.6:** *Let V and W be vector spaces over the field F. A linear transformation from V into W is a function T from V into W such that $T(c\alpha + \beta) = cT(\alpha) + T(\beta)$ for all $\alpha, \beta \in V$ and for all scalars $c \in F$. The null-space of T is the set of all vectors $\alpha$ in V such that $T\alpha = 0$. If V is finite dimensional, the rank of T is the dimension of the range of T and nullity of T is the dimension of the nullspace of T.*

$$rank\ T + nullity\ T = dim\ V.$$

*The set of all linear transformations of the vector space V into the vector space W over F denoted by $L_F(V, W)$ is also a vector space over F. Let V be a vector space over F. A linear operator $T : V \to V$; where T is a linear transformation from V to V. The set of all linear operators from V to V is denoted by $L_F(V, V)$. $L_F(V, V)$ is also a vector space over V.*

Now we proceed on to define the notion of characteristic value and characteristic vectors.

**DEFINITION 1.7.7:** *Let V be a vector space over the field F and let T be a linear operator on V. A characteristic value of T is a scalar c in F such that there is a non-zero vector $\alpha$ in V with $T\alpha = c\alpha$. If c is a characteristic value of T then*

  i.   *any $\alpha$ such that $T\alpha = c\alpha$ is called a characteristic vector of T associated with the characteristic value c;*
  ii.  *the collection of all $\alpha$ such that $T\alpha = c\alpha$ is called the characteristic space associated with c.*

Several nice and interesting results can be had from any book on linear algebra.

Now we proceed onto define the notion of Smarandache vector spaces.

**DEFINITION [44, 49]:** *The Smarandache K-vectorial space (S-K-vectorial space] is defined to be K-vectorial space $(A, +, \bullet)$ such that a proper subset A is a K-algebra (with respect to the same induced operation and another '$\times$' operation internal on A where K is a commutative field.).*

***Example 1.7.4:*** Let $V_1 = F_{m \times n}$ be a vector space over Q and $V_2 = Q[x]$ be a vector space over Q. Then $V = V_1 \times V_2$ is the direct product of two vector spaces; hence a vector space over Q. It is easily verified V is a S-K-vectorial space over Q.

**DEFINITION 1.7.8:** *Let A be a K-vectorial space. A proper subset X of A is said to be a Smarandache K-vectorial subspace (S-K-vectorial subspace) of A if X itself is a S-K-vectorial space.*

**DEFINITION 1.7.9**: *Let V be a finite dimensional vector space over a field K. Let B = {$v_1$, ..., $v_n$} be a basis of V. We say B is a Smarandache basis (S-basis) of V if B has a proper subset say A, A $\subset$ B and A $\neq \phi$; A $\neq$ B such that A generates a subspace which is linear algebra over K that is W is the subspace generated by A then W must be a K-algebra with the same operation of V.*

**THEOREM 1.7.1**: *Let A be a K-vectorial space. If A has a S-K-vectorial subspace then A is a S-K-vectorial space.*

*Proof*: Straightforward by the very definitions.

**THEOREM 1.7.2**: *Let V be a vector space over the field K. If B is a S-basis then B is a basis of V.*

*Proof*: Left for the reader as an exercise.

**DEFINITION 1.7.10**: *Let V be a finite dimensional vector space over a field K. Let B = {$v_1$, ..., $v_n$} be a basis of V. If every proper subset of B generates a linear algebra over K then we call B a Smarandache strong basis (S-strong basis) for V.*

**DEFINITION 1.7.11**: *Let V be any vector space over the field K. We say L is a Smarandache finite dimensional vector space (S-finite dimensional vector space) of K if every S-basis has only finite number of elements in it.*

*It is interesting to note that if L is a finite dimensional vector space then L is a Smarandache finite dimensional space (S-finite dimensional space) provided L has a finite S-basis.*

**THEOREM 1.7.3**: *Let V be a vector space over the field K. If A = {$v_1$, ..., $v_n$} is a S-strong basis of V then A is a S-basis of V.*

*Proof*: Direct and the proof is left for the reader as an exercise.

**DEFINITION 1.7.12**: *Let V be a vector space over the field F and let T be a linear operator from V to V. T is said to be a Smarandache linear operator (S-linear operator) on V if V has a S-basis which is mapped by T onto another S-basis of V.*

**DEFINITION 1.7.13**: *Let T be a S-linear operator defined on the space V. A characteristic value c in F associated with T is said to be a Smarandache characteristic value (S-characteristic value) of T if the characteristic vector of T associated with c is a linear algebra. So the eigen vector associated with the S-characteristic values will be called as Smarandache eigen vectors (S-eigen vectors) or Smarandache characteristic vectors (S-characteristic vectors).*

Results on the concepts are taken from [44, 49, 119].

We further suggest all notions pertaining to Smarandache vector spaces can be developed in an usual way as the very S-vector spaces properties are new [44, 49, 119] the reader is advised to define several interesting properties about them.

**Chapter 2**

# BIGROUPS AND SMARANDACHE BIGROUPS

This chapter has two sections. In the first section we recall the definition of bigroups and illustrate it with examples and give some of its properties. In section two we define Smarandache bigroups and define several new notions on them.

## 2.1 Bigroups and its properties

This section is devoted to the recollection of bigroups, sub-bigroups and illustrate it with examples. [34] was the first one to introduce the notion of bigroups in the year 1994. As there is no book on bigroups we give all algebraic aspects of it.

**DEFINITION [34]:** *A set $(G, +, \bullet)$ with two binary operation '+' and '$\bullet$' is called a bigroup if there exist two proper subsets $G_1$ and $G_2$ of $G$ such that*

    i.      $G = G_1 \cup G_2$.
    ii.     $(G_1, +)$ *is a group.*
    iii.    $(G_2, \bullet)$ *is a group.*

*A subset $H (\neq \phi)$ of a bigroup $(G, +, \bullet)$ is called a sub-bigroup, if $H$ itself is a bigroup under '+' and '$\bullet$' operations defined on $G$.*

**THEOREM 2.1.1:** *Let $(G, +, \bullet)$ be a bigroup. The subset $H \neq \phi$ of a bigroup $G$ is a sub-bigroup, then $(H, +)$ and $(H, \bullet)$ in general are not groups.*

*Proof*: Given $(G, +, \bullet)$ is a bigroup and $H \neq \phi$ of $G$ is a sub-bigroup of $G$. To show $(H, +)$ and $(H, \bullet)$ are not bigroup.

We give an example to prove this [101]. Consider the bigroup $G = \{\ldots, -2, -1, 0, 1, 2, \ldots\} \cup \{i, -j\}$ under the operations '+' and '$\bullet$'. $G = G_1 \cup G_2$ where $(G_1, \bullet) = \{-1, 1, i, -i\}$ under product and $G_2 = \{\ldots, -2, -1, 0, 1, 2, \ldots\}$ under '+' are groups.

Take $H = \{-1, 0, 1\}$. $H$ is a proper subset of $G$ and $H = H_1 \cup H_2$ where $H_1 = \{0\}$ and $H_2 = \{-1, 1\}$; $H_1$ is a group under '+' and $H_2$ is a group under product i.e. multiplication.

Thus $H$ is a sub-bigroup of $G$ but $H$ is not a group under '+' or '$\bullet$'.

Now we get a characterization theorem about sub-bigroup in what follows:

**THEOREM [101]:** *Let $(G, +, \bullet)$ be a bigroup. Then the subset $H (\neq \phi)$ of $G$ is a sub-bigroup of $G$ if and only if there exists two proper subsets $G_1, G_2$ of $G$ such that*

i. $G = G_1 \cup G_2$ where $(G_1, +)$ and $(G_2, \bullet)$ are groups.
ii. $(H \cap G_1, +)$ is subgroup of $(G_1, +)$.
iii. $(H \cap G_2, \bullet)$ is a subgroup of $(G_2, \bullet)$.

*Proof*: Let $H (\neq \phi)$ be a sub-bigroup of G, then $(H, +, \bullet)$ is a bigroup. Therefore there exists two proper subsets $H_1, H_2$ of H such that

i. $H = H_1 \cup H_2$.
ii. $(H_1, +)$ is a group.
iii. $(H_2, \bullet)$ is a group.

Now we choose $H_1$ as $H \cap G_1$ then we see that $H_1$ is a subset of $G_1$ and by (ii) $(H_1, +)$ is itself a group. Hence $(H_1 = H \cap G_1, +)$ is a subgroup of $(G_1, +)$. Similarly $(H_2 = H \cap G_2, \bullet)$ is a subgroup of $(G_2, \bullet)$. Conversely let (i), (ii) and (iii) of the statement of theorem be true. To prove $(H, +, \bullet)$ is a bigroup it is enough to prove $(H \cap G_1) \cup (H \cap G_2) = H$.

Now, $(H \cap G_1) \cup (H \cap G_2) = [(H \cap G_1) \cup H] \cap [(H \cap G_1) \cup G_2]$

$= [(H \cup H) \cap (G_1 \cup H)] \cap [(H \cup G_2) \cap (G_1 \cup G_2)]$
$= [H \cap (G_1 \cup H)] \cap [(H \cup G_2) \cap G]$
$= H \cap (H \cap G_2)$ (since $H \subseteq G_1 \cup H$ and $H \cup G_2 \subseteq G$)
$= H$ (since $H \subset H \cup G_2$).

Hence the theorem is proved.

It is important to note that in the above theorem just proved the condition (i) can be removed. We include this condition only for clarification or simplicity of understanding.

Another natural question would be can we have at least some of the classical theorems and some more classical concepts to be true in them.

**DEFINITION 2.1.1:** *Let $(G, +, \bullet)$ be a bigroup where $G = G_1 \cup G_2$; bigroup G is said to be commutative if both $(G_1, +)$ and $(G_2, \bullet)$ are commutative.*

*Example 2.1.1*: Let $G = G_1 \cup G_2$ where $G_1 = Q \setminus \{0\}$ with usual multiplication and $G_2 = \langle g \mid g^2 = I \rangle$ be a cyclic group of order two. Clearly G is a commutative bigroup.

We say the order of the bigroup $G = G_1 \cup G_2$ is finite if the number of elements in them is finite; otherwise we say the bigroup G to be of infinite order.

The bigroup given in example 2.1.1 is a bigroup of infinite order.

*Example 2.1.2*: Let $G = G_1 \cup G_2$ where $G_1 = Z_{10}$ group under addition modulo 10 and $G_2 = S_3$ the symmetric group of degree3. Clearly G is a non-commutative bigroup of finite order. $|G| = 16$.

***Example 2.1.3***: Let $G = G_1 \cup G_2$ be a bigroup with $G_1$ = {set of all n × n matrices under '+' over the field of reals} and $G_2$ = { set of all n × n matrices A with $|A| \neq \{0\}$ with entries from Q}, $(G_1, +)$ and $(G_2, \times)$ are groups and $G = G_1 \cup G_2$ is a non-commutative bigroup of infinite order.

**THEOREM 2.1.2**: *Let $(G, +, \bullet)$ be a bigroup. Let $(H, +, \bullet)$ be a proper sub-bigroup of $(G, +, \bullet)$. Then the order of $(H, +, \bullet)$ in general does not divide the order of G.*

*Proof*: This is evident from the example 2.1.2. $G = G_1 \cup G_2$ where $|G| = 16$, $G_1 = \{Z_{10}, +\}$ and $G_2 = \{S_3\}$. Clearly $H = H_1 \cup H_2$ where $H_1 = \{0, 2, 4, 6, 8\}$ under '+' and

$$H_2 = \left\{ e = \begin{pmatrix} 1 & 2 & 3 \\ 1 & 2 & 3 \end{pmatrix}, p_3 = \begin{pmatrix} 1 & 2 & 3 \\ 2 & 1 & 3 \end{pmatrix} \right\} \text{ is a sub-bigroup of } G_2.$$

Clearly $|H| = 7$ and $7 \nmid 16$. Hence the claim.

Thus we see unlike in a finite group where the order of the subgroup divides the order of the group in finite bigroups we see that in general the order of the sub-bigroups need not divide the order of the bigroup. Thus it is important to note that the classical Lagrange theorem for finite groups does not hold good in case of bigroups.

Now we proceed on to define the concept of normal sub-bigroup.

**DEFINITION 2.1.2**: *Let $(G, +, \bullet)$ be a bigroup with $G = G_1 \cup G_2$. Let $(H, +, \bullet)$ be a sub-bigroup of $(G, +, \bullet)$; where $H = H_1 \cup H_2$ we say $(H, +, \bullet)$ is a normal sub-bigroup of $(G, +, \bullet)$, if $H_1$ is a normal subgroup in $G_1$ and $H_2$ is a normal subgroup of $G_2$.*

Even if one of the subgroups happens to be normal still we do not call H to be a normal sub-bigroup of G.

***Example 2.1.4***: Let $G = \{e, p_1, p_2, p_3, p_4, p_5, 1, g, g^2, g^3, g^4, g^5, g^6, g^7, g^8\} = G_1 \cup G_2$ where $G_1 = S_3$ and $G_2$ is the cyclic group of order 9. $H = \{e, p_4, p_5, g^3, g^6, 1\}$ is a normal sub-bigroup of G; as $H = H_1 \cup H_2$ where $H_1 = \{e, p_4, p_5\}$ which is a normal subgroup of $S_3$ and $H_2 = \{1, g^3, g^6\}$ is a normal subgroup of $G_2$ as $G_2$ is abelian. We see $|H| = 6$ and $|G| = 15$. Clearly $6 \nmid 15$. Thus we see H is a normal sub-bigroup of G.

In view of this we have the following theorem:

**THEOREM 2.1.3**: *Let $(G, +, \bullet)$ be a finite bigroup. If $(H, +, \bullet)$ is a normal sub-bigroup of G. In general o(H) does not divide o(G).*

*Proof*: By an example. In Example 2.1.4 we see the group $G = G_1 \cup G_2$ where $G_1 = S_3$ and $G_2 = \langle g \mid g^9 = 1 \rangle$; G is of order 15. H given in example is a normal sub-bigroup and order of H is 6. But $6 \nmid 15$. Thus we see even in case of normal sub-bigroup H the order of H need not divide order of G.

Now we proceed on to define homomorphism of bigroups.

**DEFINITION 2.1.3:** *Let $(G, +, \bullet)$ and $(K, \oplus, o)$ be any two bigroups where $G = G_1 \cup G_2$ and $K = K_1 \cup K_2$. We say the map $\phi : G \to K$ is said to be a bigroup homomorphism if $\phi$ restricted to $G_1$ (i.e. $\phi / G_1$) is a group homomorphism from $G_1$ to $K_1$ and $\phi$ restricted to $G_2$ (i.e. $\phi / G_2$) is a group homomorphism from $G_2$ to $K_2$.*

*Example 2.1.5*: Let $(G, +, \bullet)$ be a bigroup, $G = G_1 \cup G_2$ where $G_1$ = {group of positive rationals under multiplication} and $G_2$ = {$2 \times 2$ matrices A such that $|A| \neq 0$ under multiplications entries of A are taken from Q}.

Let $K = K_1 \cup K_2$ where $K_2$ = {group of reals under the operation '+'} and $K_1$ = {group of positive reals under '$\times$'}. Define a map $\phi : G \to K$, $\phi/G_1 : G_1 \to K_1$ by $\phi(x) = x$ clearly $\phi$ is a group homomorphism which is well defined

$$\phi/G_2: G_2 \to H_2$$
$$\phi(A) \to (ad - bc)$$

where

$$A = \begin{bmatrix} a & b \\ c & d \end{bmatrix} \in G_2.$$

Thus $\phi/G_2: G_2 \to H_2$ is a group homomorphism. Hence $\phi$ is a bigroup homomorphism.

It is an important feature to be observed that we can have sub-bigroups when the bigroups are of prime order. But in case of group of prime order we do not have proper subgroups.

*Example 2.1.6*: Let $G = G_1 \cup G_2$ be a bigroup of prime order 17; where $G_1 = \langle g \mid g^9 = 1 \rangle$ and $G_2 = D_{2.4} = \{a, b \mid a^2 = b^4 = 1, ab = b^{-1}a\} = \{1, a, b, b^2, b^3, ab, ab^2, ab^3\}$.

Now $H = H_1 \cup H_2 = \{1, g^3, g^6\} \cup \{1, b, b^2, b^3\}$ is a proper sub-bigroup of G. Thus we have a bigroup of prime order can have proper sub-bigroups.

Another main feature about bigroups is that if n is the order of the bigroup n a composite number and if p is a prime such that p/n, then in general the bigroup G may not have a sub-bigroup of order p. Thus the classical Cauchy theorem for bigroups be it commutative or other wise does not in general hold good. We illustrate this in the following:

**THEOREM 2.1.4:** *Let $(G, +, \bullet)$ be a bigroup of order n, n a composite number. p be a prime such that p/n. The bigroup $(G, +, \bullet)$ may not in general have a sub-bigroup of order p.*

*Proof*: Let $G = G_1 \cup G_2$ where $G_1 = \langle g \mid g^9 = 1 \rangle$ and $G_2 = Z_5 = \{\overline{0}, \overline{1}, \overline{2}, \overline{3}, \overline{4}\}$ a group under addition modulo 5.

Clearly $|G| = 14$. 7 is a prime such that 7/14 the other even prime which divides 14 is 2. But it is very easily verified that G has no sub-bigroup of order 7. Hence the claim.

The next natural question would be Sylow theorems for bigroups.

*Example 2.1.7*: Let $G = G_1 \cup G_2$ be a bigroup of finite order; where $G_1 = S_3$ is the symmetric group of order 6 and $G_2$ is a cyclic group of order 6.

Clearly order of the bigroup is 12. We see 2/12, $2^2$/12 and 3/12. So the bigroup can have p-Sylow sub-bigroup of order $4 = 2^2$ and order 3.

Take $H = H_1 \cup H_2$ where

$$H_1 = \left\{ \begin{pmatrix} 1 & 2 & 3 \\ 1 & 2 & 3 \end{pmatrix}, \begin{pmatrix} 1 & 2 & 3 \\ 2 & 1 & 3 \end{pmatrix} \right\}$$

and $H_2 = \{1\}$. $|H| = 3$ and H is a 3-Sylow sub-bigroup of G. Also $H' = H_1 \cup H_2$ where

$$H_1' = \left\{ \begin{pmatrix} 1 & 2 & 3 \\ 1 & 2 & 3 \end{pmatrix}, \begin{pmatrix} 1 & 2 & 3 \\ 2 & 1 & 3 \end{pmatrix} \right\}$$

and $H'_2 = \{1, g^3\}$. We see $|H'| = 4$ and G has a 2 Sylow sub-bigroup. Thus we see that we can define p-Sylow sub-bigroups whenever it exists.

**DEFINITION 2.1.4:** *Let $G = G_1 \cup G_2$ be a finite bigroup of order n and p be a prime. If G has a sub-bigroup of order $p^r$ i.e. $H = H_1 \cup H_2$ where $|H| = p^r$ such that $p^r/n$ but $p^{r+1} \nmid n$ then we say G has a p-Sylow sub-bigroup.*

The above example is an illustration of the above definition. But as in case of groups we cannot always say that Sylow theorems are true in case of bigroups.

**THEOREM 2.1.5:** *Let $(G, +, \bullet)$ be a bigroup of finite order say n. If p is a prime such that p/n; then in general we need not have a p-Sylow sub-bigroup H in G.*

*Proof*: We illustrate this by an example. In example 2.1.7 we see $G = G_1 \cup G_2$ is a bigroup of order 14. 7/14. But G has no sub-bigroup of order 7. Hence the claim. But G has sub-bigroup of order 2 as $K = K_1 \cup K_2$ where $K_1 = \{1\}$ and $K_2 = \{0\}$ i.e. G has a 2-Sylow sub-bigroup of order 2 and no 7-Sylow sub-bigroup.

Now we proceed on to define symmetric bigroup and biorder of a bigroup.

**DEFINITION 2.1.5:** *Let $(G, +, \bullet)$ be a bigroup where $G = G_1 \cup G_2$. We say G is a symmetric bigroup if both $G_1$ and $G_2$ are symmetric groups.*

*Example 2.1.8*: Let $G = G_1 \cup G_2$ where $G_1 = S_6$ be the symmetric group of degree 6 and $G_2 = S_9$ be the symmetric group of degree 9. Clearly G is a symmetric bigroup of degree (6, 9).

**DEFINITION 2.1.6:** *Let $G = G_1 \cup G_2$ be a bigroup such that $G_1 = S_n$ the symmetric group of degree n and $G_2 = S_m$ be the symmetric group of degree m; $m \neq n$: then the bigroup $(G, +, \bullet)$ is symmetric bigroup of degree $(n, m)$.*

Now we proceed on to obtain the analogue of Cayley's theorem for bigroups.

**THEOREM 2.1.6:** *Every bigroup is isomorphic to a sub-bigroup of the symmetric bigroup.*

*Proof*: Let $(G, +, \bullet)$ be a symmetric bigroup with $G = G_1 \cup G_2$ where $G_1 = S_n$ and $G_2 = S_m$.

Suppose $(K, +, \bullet)$ be any bigroup not necessarily symmetric. To show the bigroup K is isomorphic to a sub-bigroup of the symmetric bigroup G.

Given $K = K_1 \cup K_2$. Now we can always find a subgroup in the symmetric group $H_1$ of $G_1$ such that $K_1$ is isomorphic to $H_1$. Similarly $K_2$ isomorphic to a subgroup $H_2$ of $G_2$. Therefore $H = H_1 \cup H_2$ is the sub-bigroup of the symmetric bigroup $G = G_1 \cup G_2 = S_n \cup S_m$. Therefore $K = K_1 \cup K_2$ is isomorphic to $H = H_1 \cup H_2$ which is a sub-bigroup of the symmetric bigroup.

Now we to overcome the problems of obtaining analogous results in case of bigroups like Lagrange or Cauchy we define a new notion in bigroups called biorder.

**DEFINITION 2.1.7:** *Let $(G, +, \bullet)$ be the bigroup with $G = G_1 \cup G_2$. Let $(H, +, \bullet)$ be a finite subset of G which is a sub-bigroup of G. We define biorder of H to be equal to $o(H_1) + o(H_2)$ where $H = H \cap G_1$ and $H_2 = H \cap G_2$ and we denote the biorder of $(H, +, \bullet)$ by B(H).*

**Note:** If the subset $H = G = G_1 \cup G_2$ then biorder of $(H, +, \bullet)$ is the same as the biorder of $(G, +, \bullet)$ that is $B(H) = o(G_1) + o(G_2) = B(G)$ we are forced to get this definition for 2 reasons.

  i.   As we did not assume $(H, +, \bullet)$ to be a proper subset of $(G, +, \bullet)$.
  ii.  To find analogous proofs for Cauchy theorems we need this concept.

**DEFINITION 2.1.8:** *Let $(G, +, \bullet)$ be a finite bigroup; $(H, +, \bullet)$ is a sub-bigroup of the bigroup $(G, +, \bullet)$ we say the biorder of $(H, +, \bullet)$ pseudo divides $o(G)$ if $o(H_1) / o(G_1)$ and $o(H_2) / o(G_2)$ and we denote it by $B(H) /_p o(G)$. If $o(H_1) \nmid o(G_1)$ or $o(H_2) \nmid o(G_2)$ we say biorder of $(H, +, \bullet)$ does not pseudo divide $o(G)$ and is denoted by $B(H) \nmid_p o(G)$, but in reality this will never occur for bigroups except when H is a set.*

We illustrate this by some examples.

*Example 2.1.9*: Let $(G, +, \bullet)$ be a bigroup where $G = G_1 \cup G_2$ with $G_1 = \{1, g, g^2\}$ a cyclic group of degree 3 and $G_2$ be the symmetric group of degree 3. Let $H = \{1, e, p_4, p_5\}$ be a sub-bigroup of $(G, +, \bullet)$ where $H_1 = \{1\}$ and $H_2 = \{e, p_4, p_5\}$.

Clearly $|H_1| = 1$ so $|H_1| / |G|$ and $|H_2| = 3$ and $|H_2| / |G_2|$; so $B(H) /_p o(G)$ that is the order of $(G, +, \bullet)$ is pseudo divisible by the order of the sub-bigroup $(H, +, \bullet)$.

**THEOREM 2.1.7:** *Let (G, +, •) be a bigroup; with G = $G_1 \cup G_2$. (H, +, •) be a subset of (G, +, •) in general o(H) ≠ B(H).*

*Proof*: We prove this statement by the following counter example. Let G = {1, a, b, $b^2$, ab, $ab^2$, g, $g^2$, $g^3$, $g^4$, $g^5$} with G = $G_1 \cup G_2$ where $G_1$ = {1, a, b, $b^2$, ab, $ab^2$} = $D_{23}$ = {a, b | $a^2 = b^3 = 1$, bab = a} and $G_2$ = {g | $g^6$=1}. Let H = {1, b, $b^2$, $g^3$} where H = $H_1 \cup H_2$ with $H_1$ = {1, b, $b^2$} and $H_2$ = {1, $g^3$}. Clearly o(H) = 4, B(H) = o($H_1$) + o($H_2$) = 5. Therefore o(H) ≠ B(H). Also we note that in general the biorder of every subset (H, +, •) of a bigroup (G, +, •) need not divide |G|.

**THEOREM 2.1.8:** *Let (G, +, •) be a bigroup and B(H) be the biorder of the subset (H, +, •) of (G, +, •); then in general B(H) $\nmid_p$ o(G).*

*Proof*: We prove this result by a counter example. Let G = {e, $p_1$, $p_2$, $p_3$, $p_4$, $p_5$, 1, g, $g^2$, ..., $g^7$} be bigroup with G = $G_1 \cup G_2$ where $G_1$ = $S_3$ be the symmetric group of degree three and $G_2$ = ⟨g | $g^8$ = 1⟩ the cyclic group of degree 8.

Clearly order of the bigroup (G, +, •) is 14. Let H = {$p_1$, $p_2$, $p_4$, $p_5$, g, $g^2$, $g^4$} where H = $H_1 \cup H_2$ with $H_1$ = {$p_1$, $p_2$, $p_4$, $p_5$} and $H_2$ = {g, $g^2$, $g^5$}. Clearly o($G_1$) = 6 and o($G_2$) = 8. Now o($H_1$) = 4 and o($H_2$) = 3. We have B(H) = 7. So B(H)/ o(G) that is 7/14. But 4 ∤ 6 and 3 ∤ 8 so it is not pseudo divisible. Therefore B(H)/o(G). But B(H) $\nmid_p$ o(G).

**THEOREM 2.1.9:** *Let (G, +, •) be a bigroup where G = $G_1 \cup G_2$. Let (H, +, •) be a sub-bigroup of (G, +, •), then o (G) is pseudo divisible by biorder B(H) of (H, +, •).*

*Proof*: Let (G, +, •), with G = $G_1 \cup G_2$ be a finite bigroup and H = $H_1 \cup H_2$ be a sub-bigroup of (G, +, •), then B(H) = o($H_1$) + o($H_2$); o($H_1$) / o($G_1$) and o($H_2$) / o($G_2$) so that B(H) $/_p$ o(G).

Since we have that in general Lagrange theorem is not true in case of bigroups. Thus now by using the tool of biorder and pseudo divisibility we obtain a modified Lagrange theorem for finite bigroups.

**THEOREM 2.1.10:** *Let (G, +, •) with G = $G_1 \cup G_2$ be a finite bigroup. Let (H, +, •) be a sub-bigroup of (G, +, •), then B(H) $/_p$ o(G).*

*Proof*: Left for the reader to prove using earlier results.

**THEOREM 2.1.11:** *Let (G, +, •) be a finite bigroup where G = $G_1 \cup G_2$. If a ∈ G then we have three possibilities.*

   i. *a ∈ $G_1$ and a ∉ $G_2$.*
   ii. *a ∈ $G_2$ and a ∉ $G_1$.*
   iii. *a ∈ $G_2 \cap G_2$.*

*In case of (i) and (ii) we have a number p such that $a^p$=e, e identity element of $G_1$ or $G_2$. In case of (iii) we may have two positive integers $p_1$ and $p_2$ such that $a^{p_1} = e_1$ and $a^{p_1} = e_2$ where $e_1$ and $e_2$ are identity elements of $G_1$ and $G_2$ respectively.*

*Proof*: Given (G, +, •) is a bigroup with G = $G_1 \cup G_2$ then we have if a ∈ $G_1$ and a ∉ $G_2$ then by classical Cauchy theorem for finite groups there exists a positive integer $p_1$ such that $a^{p_1} = e_1$ where $e_1$ is the identity element of $G_1$. In case a ∈ $G_2$ and a ∉ $G_1$ we have as before a positive integer $p_2$ such that $a^{p_2} = e_2$ where $e_2$ is the identity element of $G_2$.

Suppose a ∈ $G_1 \cap G_2$ then it implies a ∈ $G_1$ and a ∈ $G_2$ but when a ∈ $G_1$ by the classical Cauchy theorem there exists a number $p_1'$ such that $a^{p_1'} = e_1$ and similarly for a ∈ $G_2$ by the classical Cauchy theorem there exists a number $p_2'$ such that $a^{p_2'} = e_2$, $p_1'$ may not be equal to $p_2'$ or may not be equal to $p_2$. Hence the claim.

*Example 2.1.10*: Let G = $G_1 \cup G_2$. be groups with

$$G_1 = \left\{ \begin{pmatrix} a & b \\ c & d \end{pmatrix} \mid a,b,c,d \in Z_2 \text{ with } ad - bc \neq 0 \pmod 2 \right\}$$

and

$$G_2 = \left\{ \begin{pmatrix} a & b \\ c & d \end{pmatrix} \mid a,b,c,d \in Z_5 \text{ with } ad - bc \neq 0 \pmod 5 \right\}.$$

Clearly both $G_1$ and $G_2$ are non-commutative groups under matrix multiplication. Consider the element

$$\begin{pmatrix} 1 & 0 \\ 1 & 1 \end{pmatrix} \in G$$

we have

$$\begin{pmatrix} 1 & 0 \\ 1 & 1 \end{pmatrix} \in G_1 \text{ and } \begin{pmatrix} 1 & 0 \\ 1 & 1 \end{pmatrix} \in G_2.$$

Clearly if

$$\begin{pmatrix} 1 & 0 \\ 1 & 1 \end{pmatrix} \in G_1$$

then we have

$$\begin{pmatrix} 1 & 0 \\ 1 & 1 \end{pmatrix}^2 = \begin{pmatrix} 1 & 0 \\ 0 & 1 \end{pmatrix}$$

that is order of

$$\begin{pmatrix} 1 & 0 \\ 1 & 1 \end{pmatrix}$$

in $G_1$ is 2. If

$$\begin{pmatrix} 1 & 0 \\ 1 & 1 \end{pmatrix} \in G_2$$

then

$$\begin{pmatrix} 1 & 0 \\ 1 & 1 \end{pmatrix}\begin{pmatrix} 1 & 0 \\ 1 & 1 \end{pmatrix} = \begin{pmatrix} 1 & 0 \\ 2 & 1 \end{pmatrix}$$

and

$$\begin{pmatrix} 1 & 0 \\ 2 & 1 \end{pmatrix}\begin{pmatrix} 1 & 0 \\ 1 & 1 \end{pmatrix} = \begin{pmatrix} 1 & 0 \\ 3 & 1 \end{pmatrix}$$

and

$$\begin{pmatrix} 1 & 0 \\ 3 & 1 \end{pmatrix}\begin{pmatrix} 1 & 0 \\ 1 & 1 \end{pmatrix} = \begin{pmatrix} 1 & 0 \\ 4 & 1 \end{pmatrix}$$

$$\begin{pmatrix} 1 & 0 \\ 4 & 1 \end{pmatrix}\begin{pmatrix} 1 & 0 \\ 1 & 1 \end{pmatrix} = \begin{pmatrix} 1 & 0 \\ 0 & 1 \end{pmatrix}.$$

This in $G_2$ we have

$$\begin{pmatrix} 1 & 0 \\ 1 & 1 \end{pmatrix}^5 = \begin{pmatrix} 1 & 0 \\ 0 & 1 \end{pmatrix}.$$

So for the same

$$x = \begin{pmatrix} 1 & 0 \\ 1 & 1 \end{pmatrix} \in G_1 \cap G_2.$$

We have

$$x^2 = \begin{pmatrix} 1 & 0 \\ 0 & 1 \end{pmatrix} \text{ in } G_1$$

and

$$x^5 = \begin{pmatrix} 1 & 0 \\ 0 & 1 \end{pmatrix} \text{ in } G_2.$$

Hence we make an important observation that an element $x \in (G, +, \bullet)$ may have different prime powers as their order. This is a marked difference between a bigroup and a group.

**DEFINITION 2.1.9:** *Let $(G, +, \bullet)$ be a bigroup with $G = G_1 \cup G_2$. Let $H = H_1 \cup H_2$ be a sub-bigroup of G. We say H is a $(p_1, p_2)$-Sylow sub-bigroup of G if $H_1$ is a $p_1$- Sylow subgroup of $G_1$ and $H_2$ is a $p_2$-Sylow subgroup of $G_2$.*

*Example 2.1.11*: Let G = {e, $p_1$, $p_2$, $p_3$, $p_4$, $p_5$, 1, g, $g^2$, $g^3$, $g^4$, $g^5$, $g^6$=1} be a bigroup, G = $G_1 \cup G_2$ where $G_1$ = {e, $p_1$, $p_2$, $p_3$, $p_4$, $p_5$} and $G_2$ = {1, g, $g^2$, $g^3$, $g^4$, $g^5$}. Clearly |G| = 12. Let H = $H_1 \cup H_2$ where $H_1$ = {1, $p_5$, $p_4$} and $H_2$ = {1, $g^3$}. We see H is a (2,3)-Sylow sub-bigroup of G. Also, G has (3, 2) – Sylow sub-bigroups. Thus we see the Sylow sub-bigroups (2, 3) and (3, 2) are not one and the same, they are distinct.

**THEOREM 2.1.12:** *Let (G, +, •) be a finite bigroup with G = $G_1 \cup G_2$. Let ($p_1$, $p_2$) be a pair of primes such that*

$$p_1^\alpha \,/\, |G_1| \text{ and } p_2^\beta \,/\, |G_2|$$

*then (G, +, •) has ($p_1$, $p_2$)- Sylow subgroup with biorder $p_1^\alpha + p_2^\beta$.*

*Proof*: Let (G, +, •) be a finite bigroup with G = $G_1 \cup G_2$. Let ($p_1$, $p_2$) be a pair of primes such that $p_1^\alpha \,\big|\, |G_1|$ then there exists a $p_1$ –Sylow subgroup $H_1$ of $G_1$ (This is by the classical Sylow theorem for groups). Similarly $p_2^\beta \,\big|\, |G_2|$, thus there exists a $p_2$ – Sylow subgroup $H_2$ of $G_2$ of order $p_2^\beta$. Now H = $H_1 \cup H_2$ is the required ($p_1$, $p_2$)-Sylow sub-bigroup of (G, +, •). We see the biorder of H is $p_1^\alpha + p_2^\beta$.

Now we proceed on to define the notion of conjugate sub-bigroups.

**DEFINITION 2.1.10:** *Let G = $G_1 \cup G_2$ be a bigroup. We say two sub-bigroups H = $H_1 \cup H_2$ and K = $K_1 \cup K_2$ are conjugate sub-bigroups of G if $H_1 \sim K_1$ and $H_2 \sim K_2$ as subgroups of $G_1$ and $G_2$ respectively.*

So that easily the extension of second Sylow theorem can be made.

**THEOREM 2.1.13:** *Let (G, +, •) be a bigroup with G = $G_1 \cup G_2$. The number of ($p_1$, $p_2$)-Sylow sub-bigroups of (G, +, •) for the pair of primes ($p_1$, $p_2$) is of the form $(1+k_1 p_1)(1+k_2 p_2)$.*

*Proof*: The reader is expected to prove this using the classical Sylow theorem for groups.

**DEFINITION 2.1.11:** *Let G = $G_1 \cup G_2$ be a bigroup, the normalizer of a in G is the set N (a) = {x ∈ G | xa = ax}.*

*N(a) is a subgroup of $G_1$ or $G_2$ depending on the fact a ∈ $G_1$ or a ∈ $G_2$. If a ∈ $G_1 \cap G_2$ then one is not in a position to say the structure of N(a).*

The study if a ∈ $G_1 \cap G_2$; when is N(a) a sub-bigroup of G is an open problem.

*We define only when a ∈ $G_1$ or a ∈ $G_2$ 'or' in the mutually exclusive sense about the nature of normalizer.*

$$C_a = \frac{o(G_1)}{o(N(a))} \text{ if } a \in G_1$$

$$C_a = \frac{o(G_1)}{o(N(a))} \text{ if } a \in G_2.$$

What is the structure of

$$\frac{o(G)}{o(N(a))}$$

if $a \in G_1$ or $a \in G_2$ is an important issue?

Now another natural question would be what is the structure of cosets of a bigroup $G = G_1 \cup G_2$?

**DEFINITION 2.1.12:** *Let $G = G_1 \cup G_2$ be a bigroup. Let $H = H_1 \cup H_2$ be the sub-bigroup of G. The right bicoset of H in G for some $a \in G$ is defined to be $H_a = \{h_1 a \mid h_1 \in H_1 \text{ and } a \in G_1 \cap G_2\} \cup \{h_2 a \mid h_2 \in H_2 \text{ and } a \in G_1 \cap G_2\}$ if $a \in G_1$ and $a \notin G_2$, then $H_a = \{h_1 a \mid h_1 \in H_1\} \cup H_2$. If $a \in G_2$ and $a \notin G_1$ then $H_a = \{h_2 a \mid h_2 \in H_2\} \cup H_1$. Thus we have bicosets depends mainly on the way we choose a. If $G_1 \cap G_2 = \phi$ then the bicoset is either $H_1 a \cup H_2$ or $H_1 \cup H_2 a$.*

*Similarly we define left bicoset of a sub-bigroup H of G.*

*The next natural question would be if $H = H_1 \cup H_2$ and $K = K_1 \cup K_2$ be any two sub-bigroups of $G = G_1 \cup G_2$; how to define HK. We define HK as follows;*

$$HK = \left\{ h_1 k_1, h_2 k_2 \;\middle|\; \begin{array}{l} h_1 \in H_1 \\ k_1 \in K_1 \end{array} \text{ and } \begin{array}{l} h_2 \in H_2 \\ k_2 \in K_2 \end{array} \right\}.$$

The study when is HK a bigroup is as follows:

**THEOREM 2.1.14:** *Let $G = G_1 \cup G_2$ be a bigroup. $H = H_1 \cup H_2$ and $K = K_1 \cup K_2$ be any two sub-bigroups of G. HK will be a sub-bigroup if and only if $H_1 K_1$ and $H_2 K_2$ are subgroups of $G_1$ and $G_2$ respectively.*

*Proof*: Straightforward hence left for the reader to prove.

**DEFINITION 2.1.13:** *Let $G = G_1 \cup G_2$ be a bigroup we say $N = N_1 \cup N_2$ is a normal sub-bigroup of G if and only if; $N_1$ is normal in $G_1$ and $N_2$ is normal in $G_2$. We define the quotient bigroup G/N as*

$$\left( \frac{G_1}{N_1} \cup \frac{G_2}{N_2} \right)$$

*which is also a sub-bigroup. Hence the claim.*

**PROBLEMS:**

1. Find the smallest bigroup (i.e. what is the order of the smallest bigroup which is nontrivial).

2. Find for the bigroup $G = G_1 \cup G_2$ (where $G_1 = Z_{12}$ under '+' and $G_2 = S_7$) the sub-bigroups and normal bisubgroups of G.

3. Find for problem 2, $N\left(\begin{pmatrix} 1 & 2 & 3 & 4 & 5 & 6 & 7 \\ 2 & 1 & 4 & 3 & 6 & 7 & 5 \end{pmatrix}\right)$.

4. Let $G = G_1 \cup G_2$ where $G_1 = \langle g \mid g^{16} = 1 \rangle$ and $G_2 = Z_{21}$, the group under '+' modulo 21. Let $H = H_1 \cup H_2$, where $H_1 = \{1, g^4, g^8, g^{12}\}$ and $H_2 = \{0, 7, 14\}$. Find the bicosets of H using all elements in G.

5. Does the bigroup given in problem 4 have Sylow sub-bigroups?

6. Verify, Cauchy theorem for $G = G_1 \cup G_2$, where $G_1 = S_5$ (the symmetric group of degree 5) and

$$G_2 = \left\{ \begin{pmatrix} a & b \\ c & d \end{pmatrix} \middle| ab - bc \neq 0; a, b, c, d \in Z_2 \right\},$$

$G_2$ is a group under matrix multiplication.

7. For the bigroup given in problem 6 find how many sub-bigroups exist?

## 2.2 Smarandache bigroups and its properties

In this section we introduce the concept of Smarandache bigroups. Bigroups were defined and studied in the year 1994 [34]. But till date Smarandache bigroups have not been defined. Here we define Smarandache bigroups and try to obtain several of the classical results enjoyed by groups. Further the study of Smarandache bigroups will throw several interesting features about bigroups in general.

**DEFINITION 2.2.1:** *Let $(G, \bullet, *)$ be a non-empty set such that $G = G_1 \cup G_2$ where $G_1$ and $G_2$ are proper subsets of G. $(G, \bullet, *)$ is called a Smarandache bigroup (S-bigroup) if the following conditions are true.*

   i. *$(G, \bullet)$ is a group.*
   ii. *$(G, *)$ is a S-semigroup.*

***Example 2.2.1*:** Let $G = \{g^2, g^4, g^6, g^8, g^{10}, g^{12} = 1\} \cup S(3)$ where $S(3)$ is the symmetric semigroup. Clearly $G = G_1 \cup G_2$ where $G_1 = \{1, g^2, g^4, g^6, g^8, g^{10}\}$, group under '×' and $G_2 = S(3)$; S-semigroup under composition of mappings. G is a S-bigroup.

**THEOREM 2.2.1:** *Let G be a S-bigroup, then G need not be a bigroup.*

*Proof*: By an example consider the bigroup given in example 2.2.1. G is not a bigroup only a S-bigroup.

*Example 2.2.2*: Let $G = Z_{20} \cup S_5$; $Z_{20}$ is a S-semigroup under multiplication modulo 20 and $S_3$ is the symmetric group of degree 3. Clearly G is a S-bigroup; further G is not a bigroup.

**DEFINITION 2.2.2:** *Let $G = G_1 \cup G_2$ be a S-bigroup, a proper subset $P \subset G$ is said to be a Smarandache sub-bigroup of G if $P = P_1 \cup P_2$ where $P_1 \subset G_1$ and $P_2 \subset G_2$ and $P_1$ is a group or a S-semigroup under the operations of $G_1$ and $P_2$ is a group or S-semigroup under the operations of $G_2$ i.e. either $P_1$ or $P_2$ is a S-semigroup i.e. one of $P_1$ or $P_2$ is a group, or in short P is a S-bigroup under the operation of $G_1$ and $G_2$.*

**THEOREM 2.2.2:** *Let $G = G_1 \cup G_2$ be a S-bigroup. Then G has a proper subset H such that H is a bigroup.*

*Proof*: Given $G = G_1 \cup G_2$ is a S-bigroup. $H = H_1 \cup H_2$ where if we assume $G_1$ is a group then $H_1$ is a subgroup of $G_1$ and if we have assumed $G_2$ is a S-semigroup then $H_2$ is proper subset of $G_2$ and $H_2$ which is a subgroup of $G_2$.

Thus $H = H_1 \cup H_2$ is a bigroup.

**COROLLARY:** *If $G = G_1 \cup G_2$ a S-bigroup then G contains a bigroup.*

Study of S-bigroups is very new as in literature we do not have the concept of Smarandache groups we have only the concept of S-semigroups.

**DEFINITION 2.2.3:** *Let $G = G_1 \cup G_2$ be a S-bigroup we say G is a Smarandache commutative bigroup (S-commutative bigroup) if $G_1$ is a commutative group and every proper subset S of $G_2$ which is a group is a commutative group.*

*If both $G_1$ and $G_2$ happens to be commutative trivially G becomes a S-commutative bigroup.*

*Example 2.2.3:* Let $G = G \cup S(3)$ where $G = \langle g \mid g^9 = 1 \rangle$ Clearly G is a S-bigroup. In fact G is not a S-commutative bigroup.

*Example 2.2.4:* Let $G = G \cup S(4)$ where $G = \langle g \mid g^{27} = 1 \rangle$ and $S(4)$ the symmetric semigroup which is a S-semigroup this S-bigroup is also non-commutative.

**DEFINITION 2.2.4:** *Let $G = G_1 \cup G_2$ be a S-bigroup where $G_1$ is a group and $G_2$ a S-semigroup we say G is a S-weakly commutative bigroup if the S-semigroup $G_2$ has atleast one proper subset which is a commutative group.*

*Example 2.2.5:* The S-bigroups given in example 2.2.3 and 2.2.4 are S-weakly commutative but not a S-commutative bigroup.

**THEOREM 2.2.3:** *Let $G = G_1 \cup G_2$ be a S-commutative bigroup then G is a S-weakly commutative bigroup and not conversely.*

*Proof*: Follows from the fact that if G is a S-commutative bigroup then G is clearly a S-weakly commutative bigroup.

To prove the converse we see the S-bigroups given in example 2.2.3 and 2.2.4 are S-weakly commutative bigroups, which is clearly not a S-commutative bigroup.

**DEFINITION 2.2.5:** *Let $G = G_1 \cup G_2$ be a S-bigroup. We say G is a Smarandache cyclic bigroup (S-cyclic bigroup) if $G_1$ is a cyclic group and $G_2$ is a S-cyclic semigroup.*

**DEFINITION 2.2.6:** *Let $G = G_1 \cup G_2$ be a S-bigroup. We say G is a S-weakly cyclic bigroup if every subgroups of the group $G_1$ is cyclic and $G_2$ is a S-weakly cyclic semigroup.*

**THEOREM 2.2.4:** *Let $G = G_1 \cup G_2$ be a S-cyclic bigroup then G is a S-weakly cyclic bigroup.*

*Proof*: Straightforward by the very definition.

**DEFINITION 2.2.7:** *Let $G = G_1 \cup G_2$ be a S-bigroup. If A be a proper subset of $G_2$ which is the largest group contained in $G_2$; and $G_1$ has no subgroups; then we call $P = G_1 \cup A$ the Smarandache hyper bigroup (S-hyper bigroup).*

*Example 2.2.6:* Let $G = G_1 \cup G_2$ be a S-bigroup, where $G_1 = \langle g \mid g^3 = 1 \rangle$ and $G_2 = S(3)$, then $P = G_1 \cup S_3$ is the S-hyper bigroup of G.

**THEOREM 2.2.5:** *Let $G = G_1 \cup G_2$ be a S-bigroup. The S-hyper bigroup of G if it exists, is a bigroup.*

*Proof*: Direct by the very definition.

**DEFINITION 2.2.8:** *Let $G = G_1 \cup G_2$ be a S-bigroup. We say G is Smarandache simple bigroup (S-simple bigroup) if G has no S-hyper bigroup.*

**DEFINITION 2.2.9:** *Let $G = G_1 \cup G_2$ be a S-bigroup. If the number of elements in G is finite then we say the S-bigroup is of finite order, otherwise we say the S-bigroup is of infinite order.*

*Example 2.2.7:* Let $G = Z \cup S(n)$ be a S-bigroup where Z is the group under '+' and S(n) is the symmetric S-semigroup then G is a S-bigroup of infinite order. The S-bigroups given in examples are 2.2.1 to 2.2.7 are of finite order.

**DEFINITION 2.2.10:** *Let $G = G_1 \cup G_2$ be a finite S-bigroup. If the order of every S-sub-bigroup H of G divides the order of G then we call G a Smarandache Lagrange bigroup (S-Lagrange bigroup). If the order of at least one of the S-sub-bigroups of G* *say H divides the order of G then we call G a Smarandache weakly Lagrange bigroup (S-weakly Lagrange bigroup).*

The following theorem is direct hence left for the reader to prove.

**THEOREM 2.2.6:** *Every S-Lagrange bigroup is a S-weakly Lagrange bigroup.*

**DEFINITION 2.2.11:** *Let $G = G_1 \cup G_2$ be a finite S-bigroup. Let p be a prime such that p divides the order of G. If G has a S-sub-bigroup H of order p or $p^t$ (t >1) then we say G has a Smarandache p-Sylow bigroup (S-p-Sylow bigroup). It is very important to note that $p \,/\, o(G)$ but $p^t \nmid o(G)$ for any t > 1, still we may have S-p-Sylow sub-bigroups having $p^t$ elements in them.*

*Even if $p \,/\, o(G)$ still in case of S-bigroup we may or may not have S-p-Sylow sub-bigroups.*

We illustrate it by the following example.

***Example 2.2.8:*** Let $G = G_1 \cup G_2$ where $G_1 = \langle g \mid g^5 = 1 \rangle$ and $G_2 = S(3)$. G is a S-bigroup. |G| = 33. The primes which divide |G| are 3/33 and 11/33.

Let $H = H_1 \cup H_2$ where $H_1 = \{1, g^3\}$ and $H_2 = \{1, p_2\}$. So H is a sub-bigroup, which is not a S-sub-bigroup. Take $H_1 = \{1, g^2, g^4\} \cup S_3$ then $H_1$ is a S-sub-bigroup and $|H_1|$ = 9. We see the order of S-sub-bigroup need not divide order of G.

***Example 2.2.9:*** $G = G_1 \cup G_2$ where $G_1 = \{g \mid g^{12} = 1\}$ and $G_2 = S(5)$. G is a S-bigroup. $H = H_1 \cup H_2$ where $H_1 = \{g^4, g^8, 1\}$ and $H_2 = S_1(3)$ (By $S_1(3)$ we fix two elements in the mappings and consider only mappings of 3 elements, so $S_1(3)$ is a symmetric semigroup of order 27). Clearly H is a S-sub-bigroup. Now |G| = 3137 and |H| = 30 we see $|H| \nmid |G|$ thus here we wish to make some observations.

We see in case of classical theory if G is a finite group and H a subgroup of G then |H| / |G|. The same result does not hold good in case finite bigroups or even finite S-semigroups and more so in the case of S-bigroups. Thus S-bigroup has a S-hyper subgroup also as from [121]; all S(n) have S-hyper subgroups.

**DEFINITION 2.2.12:** *Let $G = G_1 \cup G_2$ be a S-bigroup of finite order. Let $a \in G$, a is said to be a S-special Cauchy element of G if n / |G| where n > 1 such that $a^n = 1$ (n is the least one).*

*We call $a \in G$ to be a S-Cauchy element of G if $a \in G_1$ (or $a \in G_2$) and $n \,/\, |G_1|$ (or $n \,/\, |G_2|$) if n > 1 such that $a^n = 1$ (n is a least). We see that by no means we can associate at all times the S-Cauchy element and the S-special Cauchy element of a S-bigroup.*

***Example 2.2.10:*** Let $G = G_1 \cup G_2$ where $G_1 = \{g \mid g^5 = 1\}$ and $G_2 = S(3)$; G is a S-bigroup. |G| = 32. Take

$$a = \begin{pmatrix} 1 & 2 & 3 \\ 1 & 3 & 2 \end{pmatrix} \in G$$

a is a S-special Cauchy element of G but a is not a S-Cauchy element of $G_2$; further g $\in$ G is such that $g^5 = 1$, g is not a S-special Cauchy element of G.

It is pertinent to make the following observations:

Every element as in case of groups will not be a Cauchy element. Further an element in a S-bigroup need not in general be a S-Cauchy element or a S-special Cauchy element.

**THEOREM 2.2.7:** *Let $G = G_1 \cup G_2$ be a S-bigroup. Every element in G need not in general be a S-special Cauchy element or a S-Cauchy element of G.*

*Proof*: This is to be proved only by using examples. Consider the example 2.2.10 we can get the proof of the theorem.

Now we proceed on to define Smarandache cosets for S-bigroups.

**DEFINITION 2.2.13:** *Let $G = G_1 \cup G_2$ be a S-bigroup, let $H = H_1 \cup H_2$ be a S-sub-bigroup of G. For any $g \in G$ we define the Smarandache right coset (S-right coset) $Ha = \{H_1a \cup H_2$ if $a \in G_1\}$ or $= \{H_1 \cup H'_2 a$ if $a \in G_2\}$ (if in special case $a \in H_1 \cap H_2$) then $Ha = \{H_1a \cup H'_2a / a \in G_1 \cap G_2)$ where $H'_2 \subset H_2$ is a subgroup of the S-subsemigroup $H_2$. Similarly we define S- left coset a H. We say H is a Smarandache coset (S-coset) if aH = Ha.*

**THEOREM 2.2.8:** *Let $G = G_1 \cup G_2$ be a finite bigroup. If m is a positive integer such that m/ |G| then*

      *i. G need not have a subgroup of order m where m / |G|.*
      *ii. Even if H is a S-sub-bigroup of G then $|H| \nmid |G|$ in general.*

*Proof*: To show if G is a bigroup G can have a S- sub-bigroup of order m such that m / |G|. Let $G = G_1 \cup G_2$, where $G_1 = \langle g \mid g^9 = 1 \rangle$ and $G_2 = S(4)$. Clearly |G| = 265. We see 5 /|G| but G has S-sub-bigroup of order 5; for $H = H_1 \cup H_2$, where $H_1 = \{1\}$ and $H_2 = S(2)$. |H| = 5. Hence the claim.

To show for a S-bigroup $G = G_1 \cup G_2$, H is a S-sub-bigroup of G then in general $|H| \nmid |G|$. Let $G = G_1 \cup G_2$ where $G_1 = \langle g \mid g^5 = 1 \rangle$ and $G_2 = S(3)$. Take $H = H_1 \cup H_2$ where $H_1 = \{1\}$ and $H_2 = S(2)$. Clearly H is a S-sub-bigroup of G and |H| = 5. Clearly 5 $\nmid$ 32. Further we see 2/32 but G has no S-sub-bigroup of order 2. Hence the claim.

**DEFINITION 2.2.14:** *$G = G_1 \cup G_2$ be a finite bigroup. If the order of none of the S-sub-bigroups of G divides the order of G then we call G a Smarandache non-Lagrange bigroup (S-non-Lagrange bigroup).*

***Example 2.2.11:*** Let $G = G_1 \cup G_2$ where $G_1 = \langle g \mid g^4 = 1 \rangle$ and $G_2 = S(3)$. Clearly |G| = 31 a prime. So none of the S-sub-bigroup of G divides the order of G. For instance $H = H_1 \cup H_2$ where $H_2 = S(2)$ and $H_1 = \{1, g^2\}$. Clearly |H| = 6 but 6 $\nmid$ 31 as 31 is a prime.

So what ever be the S-sub-bigroup of G still G will be S-non-Lagrange bigroup.

**DEFINITION 2.2.15:** *Let $G = S(m) \cup S_n$ where $S(m)$ is the symmetric semigroup and $S_n$ is the symmetric group of degree n. G is a S-bigroup if and only if $n > m$. These types of S-bigroups are defined as Smarandache symmetric bigroups (S-symmetric bigroup). Clearly if $n < m$ then G is not a bigroup just a S-semigroup as $S_n \subseteq S(m)$.*

*Example 2.2.12:* Let $G = S_7 \cup S(6)$. Clearly G is a S-bigroup of order $(7! + 6^6)$. Further G is a S-symmetric bigroup. The smallest S-symmetric bigroup is $S_3 \cup S(2)$ and its order is $3! + 2^2 = 6 + 4 = 10$. So all S-symmetric bigroups are of order greater than 10. The concept of S-symmetric bigroup leads to the formulations of Smarandache Cayley theorem for S-bigroups.

**THEOREM (SMARANDACHE CAYLEY'S THEOREM FOR S-BIGROUPS):** *Let $G = G_1 \cup G_2$ be a S-bigroup. Every S-bigroup is embeddable in a S-symmetric bigroup for a suitable m and n.*

*Proof*: We know by the classical Cayley's theorem for groups that every group can be realized as a subgroup of a the symmetric group $S_n$ for a suitable n. By Cayley's theorem for S-semigroups [121] we know every S-semigroup is isomorphic to a S-semigroup $S(m)$ for some m, where m denotes the number of elements in X and $S(m)$ denotes the set of all mappings of the set X to X. Using these theorems we see every S-bigroup is isomorphic to a S-sub-bigroup of $G(m, n)$ where $G(m, n) = S_m \cup S(n)$ with $m > n$. Hence the claim.

**DEFINITION 2.2.16:** *Let $G = G_1 \cup G_2$ be a bigroup. An element $x \in G \setminus \{1\}$ is said to have a Smarandache inverse (S-inverse) $y \in G$ if $xy = 1$ and for $a, b \in G \setminus \{1, x, y\}$ we have*

$$xa = y \ ( or \ ax = y)$$
$$yb = x \ ( or \ by = x)$$

*with $ab = 1$. The pair (x, y) is called the Smarandache inverse pair.*

*If x be the S-inverse of y and (x, y) is a S-inverse pair with the related pair (a, b). If the pair (a, b) happens to be a S-inverse pair not necessarily with (x, y) as related pair then we say (a, b) is the Smarandache co-inverse pair (S-co-inverse pair).*

**THEOREM 2.2.9:** *Let $G = G_1 \cup G_2$, $G_1$ a group. Every S-inverse in $G_1$ or $H \subset G_2$ (H a subgroup of $G_2$) has an inverse in the subset P of G ($P = G_1 \cup H$) but every inverse in G need not have a S-inverse.*

*Proof*: Left for the reader as an exercise.

*Example 2.2.13:* Let $G = G_1 \cup G_2$, where $G_1 = \langle g \mid g^7 = 1 \rangle$ and $G_2 = S(3)$. If $x \in G$ has S-inverse than $x \in G_1$ and $x \notin G_2$.

The study of finding inverses is S-bigroups is left for the reader.

Now we proceed on to define the concept of S-conjugate bigroups.

**DEFINITION 2.2.17:** *Let $G = G_1 \cup G_2$ be a S-bigroup. We say an element $x \in G$ has a Smarandache conjugate (S-conjugate) y in G if*

    i.   *$xa = ay$ for some $a \in G$.*
    ii.  *$ab = bx$ and $ac = cy$ for some b, c in G.*

**THEOREM 2.2.10:** *Let $G = S_n \cup S(m)$ be the S-symmetric bigroup. G has elements, which are S-conjugates.*

*Proof*: Refer [121] for proof.

Now once the concept of cosets are defined in S-bigroups one of the natural question would be how to define Smarandache double cosets.

**DEFINITION 2.2.18:** *Let $G = G_1 \cup G_2$ be a S-bigroup. Let $A = A_1 \cup A_2$ and $B = B_1 \cup B_2$ be any two S-sub-bigroups of G. We define the Smarandache double coset (S-double coset) as $AxB = \{A_1x_1B_1 \cup A_2x_2B_2 / A_1x_1B_1$ provided $x = x_1 \in G_1$ if $x \neq x_1 \notin G_1$ then we just take $A_1B_1$ if $x = x_2 \in G_2$ then take $A_2x_2B_2$ otherwise take $A_2B_2$ for every $x \in G\}$. $AxB = \{A_1B_1 \cup A_2xB_2$ if $x \in G_2\}$ or $= \{A_1xB_1 \cup A_2B_2$ if $x \in G_1\}$ or $= \{A_1xB_1 \cup A_2xB_2$ if $x \in G_1 \cap G_2\}$.*

Finally we proceed on to define S-normal sub-bigroups of a S-bigroup.

**DEFINITION 2.2.19:** *Let $G = G_1 \cup G_2$ be a S-bigroup. We call a S-sub-bigroup A of G to be a Smarandache normal sub-bigroup (S-normal sub-bigroup) of G where $A = A_1 \cup A_2$ if $xA_1 \subset A_1$ and $A_1x \subset A_1$ if $x \in G_1$ and $xA_2 \subset A_2$ and $A_2x \subset A_2$ or $xA_2 = \{0\}$ or $A_2x = \{0\}$ if 0 is an element in $G_2$ and $x \in G_2$ for all $x \in G_1 \cup G_2$.*

**DEFINITION 2.2.20:** *Let $G = G_1 \cup G_2$ be a S-bigroup. We say G is Smarandache pseudo simple bigroup (S-pseudo simple bigroup) if G has no S-normal sub-bigroup.*

Now we proceed on to define Smarandache direct product of S-bigroups (S-direct product of S-bigroups).

**DEFINITION 2.2.21:** *Let $G = G_1 \cup G_2$ be a S-bigroup. We say a proper subset $A = A_1 \cup A_2$ to be a Smarandache maximal S-bigroup (S-maximal S-bigroup) of G if $A_1$ is the largest normal subgroup in G and $A_2$ is the largest proper subset in $G_2$ which is the subgroup of $G_2$.*

It is important to note that for a given S-bigroup one may have several Smarandache maximal sub-bigroups.

*Example 2.2.14:* Let $G = G_1 \cup G_2$ where $G_1 = \langle g \mid g^6 = 1 \rangle$ and $G_2 = S(3)$. $A_1 = \{1, g^3\} \cup S_3$ and $A_2 = \{1, g^2, g^4\} \cup S_3$ are the S-maximal sub-bigroups of G.

*Example 2.2.15:* Let $G = S_4 \cup S(3)$. Then G has only one S-maximal sub-bigroup of G viz. $A = A_4 \cup S_3$. Thus it is very important to note that A is not a S-sub-bigroup of G certainly A is a bigroup of G.

*Example 2.2.16:* Let $G = G_1 \cup G_2$, $G_1 = \langle g \mid g^3 = 1 \rangle$ and $G_2 = S_2$. Then G has no S-maximal sub-bigroup.

**DEFINITION 2.2.22:** *Let $G_1, G_2, \ldots, G_n$ be n- S-bigroups. $G = G_1 \times G_2 \times \ldots \times G_n = \{(s_{11} \cup s_{12}), (s_{21} \cup s_{22}), \ldots, (s_{n1} \cup s_{n2}) \mid s_{j1} \in G_1, s_{j2} \in G_2, \ldots, s_{jn} \in G_n, j = 1, 2\}$ is called the Smarandache direct product of S-bigroup (S-direct product of S-bigroup) $G_1, \ldots, G_n$. If G has a S-maximal sub-bigroup and if $H = H_1 \times \ldots \times H_n$ then each $H_i$ is a S-maximal sub-bigroup of $G_i = 1, 2, \ldots, n$.*

**DEFINITION 2.2.23:** *Let G be a S-bigroup, $A_1 \cup B_1, A_2 \cup B_2, \ldots, A_n \cup B_n$, be non-empty subsets of G i.e. $A_i \neq \phi \subset G_1$ and $B_i \neq \phi \subset G_2$ we say $(A_1 \cup B_1), (A_2 \cup B_2), \ldots, (A_n \cup B_n)$ is the Smarandache internal direct product (S-internal direct product) of G if $G = \{A_1 A_2 \ldots A_n\} \cup \{B_1 B_2 \ldots B_n\} = \{a_1 \ldots a_n \cup b_1 \ldots b_n / a_i \in A_i \text{ and } b_i \in B_i, 1 \leq i \leq n\}$ i.e. $G_1 = A_1 \ldots A_n$ and $G_2 = B_1 \ldots B_n$. If in the internal direct product one of $(A_i \cup B_i)$ happen to be S-maximal sub-bigroup of G then we call G the Smarandache strong internal direct product (S-strong internal direct product).*

**THEOREM 2.2.11:** *If $G = G_1 \cup G_2$ is a S-bigroup and $G = \{(A_1 \ldots A_n) \cup (B_1 \ldots B_n)\}$, $A_1, A_2, \ldots, A_n$ are proper subsets of $G_1$ and $B_1, B_2, \ldots, B_n$ are proper subsets of $G_2$ with $A_i \cup B_i$ to be a S-maximal sub-bigroup i.e. G is the S-strong internal direct product then G is a S-internal direct product.*

*Proof*: Follows directly by the definitions.

**PROBLEMS:**

1. Give an example of a S-commutative bigroup.
2. Give an example of S-weakly cyclic bigroup, which is not a S-cyclic bigroup.
3. Give an example of S-weakly Lagrange bigroup which is not a S- Lagrange bigroup.
4. Does the S-bigroup $G = G_1 \cup G_2$ where $G_1 = \langle g \mid g^{21} = 1 \rangle$, and $G_2 = S(7)$ has a
   i. S-sub-bigroup.
   ii. S-cyclic bigroup.
   iii. S-hyper bigroup.
5. Give an example of a S-Cauchy element in a S-bigroup G, which is not a S-special Cauchy element of G.
6. Give an example of a S-bigroup, which does not satisfy Lagrange theorem for any S-sub-bigroup.
7. Illustrate by an example that there can be S-bigroups in which no element is a S-special Cauchy element.
8. For the S-bigroup $G = G_1 \cup G_2$ where $G_1$ is the dihedral group and $G_2$ is the S-semigroup $S(4)$; let $H = H_1 \cup H_2$ where $H_1 = \{a, 1\}$ and $H_2 = S(3) \subset S(4)$ the S-right coset of H. The S-left coset of H and the S-coset of H.
9. Find all S-sub-bigroups of $G = S(7) \cup Z_{20}$. $Z_{20}$ is the group under '+'.
10. Find an example of a S-Lagrange bigroup.

11. Is $Z_{10} \cup S(3) = G$ a S-weakly Lagrange bigroup? Justify your claim.
12. Let $G = G_1 \cup G_2$ where $G_1 = \langle g \mid g^{10} = 1 \rangle$ and $G_2 = S(3)$. Prove G is a S-non-Lagrange bigroup.
13. Find all the S-inverses in the S-bigroup $G = S(5) \cup D_{2.7}$.
14. Find the element in $G = S(7) \cup S_9$ which has S-conjugates.
15. Let $G = S_5 \cup S(4)$ be a S-bigroup. For the S-sub-bigroups $A = G_1 \cup S(2)$ where

$$G'_1 = \left\langle \left\{ \begin{pmatrix} 1 & 2 & 3 & 4 & 5 \\ 2 & 3 & 4 & 5 & 1 \end{pmatrix} \right\} \right\rangle$$

the cyclic group of order 5. S(2) the symmetric semigroup and $B = G'_1 \cup S(3)$ where

$$G'_1 = \left\langle \left\{ \begin{pmatrix} 1 & 2 & 3 & 4 & 5 \\ 2 & 3 & 4 & 1 & 5 \end{pmatrix} \right\} \right\rangle,$$

the cyclic group of order 4 and S(3) the symmetric semigroup. Find

i. AxB where $x = \begin{pmatrix} 1 & 2 & 3 & 4 & 5 \\ 2 & 3 & 1 & 4 & 5 \end{pmatrix}$.

ii. AyB where $y = \begin{pmatrix} 1 & 2 & 3 & 4 \\ 2 & 3 & 1 & 2 \end{pmatrix}$.

iii. AzB where $z = \begin{pmatrix} 1 & 2 & 3 & 4 & 5 \\ 2 & 1 & 4 & 3 & 5 \end{pmatrix} \in S_5 \cap S(4)$.

16. Does $G = D_{2.9} \cup S(5)$ have S-normal bisubgroup?
17. Give an example of a S-bigroup, which has no S-normal bisubgroup.
18. Give an example of S-pseudo simple bigroup.
19. Is $G = G_1 \cup G_2$ where $G = \langle g \mid g^{11} = 1 \rangle$ and $G_2 = S(3)$ a S-pseudo simple bigroup?
20. Find the S-maximal sub-bigroup of $G = D_{2.7} \cup S(7)$.
21. Give an example of a S-bigroup, which has no S-maximal subgroup.
22. Let $G_1 = D_{2.6} \cup S(3)$ and $G_2 = S_4 \cup S(2)$. Is $G_1 \times G_2 = G$, a S-direct product of S-bigroup?
23. Give an example of a S-strong internal direct product.
24. Illustrate by an example a S-internal direct product of a S-bigroup, which is not a S-strong internal direct product.

**Chapter 3**

# BISEMIGROUPS AND SMARANDACHE BISEMIGROUPS

This chapter has three sections. In the first section we just recall the concept of bisemigroups and illustrate with examples. Then in section two we give the special cases of semilattice and some of its related properties. In section three we introduce and study the concept of Smarandache bisemigroups and derive several interesting results about them.

## 3.1 Bisemigroups and its applications

This section is completely devoted to the recollection of bisemigroup and give some example and properties. For properties about bisemigroups refer [125].

**DEFINITION [125]:** *A bisemigroup is an algebra S = (S, +, •) equipped with associative binary operations + and '•'. Bisemigroups with a common neutral element for the two associative binary operations are called double monoids or bimonoids.*

Some authors term it as binoid. A commutative bisemigroup is a bisemigroup in which both operations are commutative.

*Example 3.1.1:* $(Z, +, •)$ is a bisemigroup which is commutative.

*Example 3.1.2:* $(Q, \times, \div)$ is a bimonoid.

We give a different version of bisemigroup which we will call as bisemigroups of type II. So bisemigroups defined by [125] will be called as bisemigroups of type I.

**DEFINITION 3.1.1:** *A set (S, +, •) with two binary operations '+' and '•' is called a bisemigroup of type II if there exists two proper subsets $S_1$ and $S_2$ of S such that*

  i.   $S = S_1 \cup S_2$.
  ii.  *(S, +) is a semigroup.*
  iii. *(S, •) is a semigroup.*

*Thus by this definition we will always have all bigroups to be bisemigroups but bisemigroups in general are not bigroups.*

*Example 3.1.3:* Let $S = \{Z^+, i, -i, -1\}$, $(S, +, •)$ is a bisemigroup. For we see $S_1 = (Z^+, +)$ is a semigroup and $S_2 = \{1, -1, i, -i\}$ is a semigroup under $\times$. But we see S is not a bigroup only a bisemigroup.

Further S is only a bisemigroup of type II and clearly not a bisemigroup of type I.

**THEOREM 3.1.1:** *All bigroups are bisemigroups but all bisemigroups in general are not bigroups.*

*Proof*: Straightforward by the very definitions.

**DEFINITION 3.1.2:** *A subset $H \neq \phi$ of $(S, +, \bullet)$ with two binary operations '+' and '$\bullet$' is called a sub-bisemigroup if H itself is a bisemigroup i.e. $H = H_1 \cup H_2$ where $H_1 \subset S_1$ and $H_2 \subset S_2$ and $(H_1, +)$ and $(H_2, \bullet)$ are semigroups.*

**THEOREM 3.1.2:** *Let $(S, +, \bullet)$ be bisemigroup. The subset $H \neq \phi$ of S be a sub-bisemigroup. Clearly $(H, +)$ and $(H, \bullet)$ in general need be semigroups.*

*Proof*: The reader is expected to derive the proof as it is straightforward.

*Example 3.1.4:* Let S = {set of all 2 × 3 matrices and 4 × 4 matrices with entries from $Z_6$}. $(S, +, \bullet)$ is a bisemigroup. For take $S_1$ = {set of 2 × 3 matrices with entries from $Z_6$} and $S_2$ = {set of 4 × 4 matrices with entries from $Z_6$}. Clearly $S = S_1 \cup S_2$ and $(S_1, +)$ and $(S_2, +)$ are semigroups so S is a bisemigroup. H = {set of all 2 × 3 and 4 × 4 matrices with entries from P = {0, 2, 4} $\subset Z_6$}. H is a bisemigroup. Clearly $(H, +)$ is not even closed under '+'.

Similarly $(H, \times)$ is not even closed under ×. Thus this example is also one of the examples to counter the above theorem. We obtain a necessary and sufficient condition for the subsemigroup to be a semigroup.

**THEOREM 3.1.3:** *Let $(S, +, \bullet)$ be a bisemigroup. The subset $H (\neq \phi)$ of S is a sub-bisemigroup of S if and only if there exists two proper subsets $S_1, S_2$ of S such that*

    i.     *$S = S_1 \cup S_2$ where $(S_1, +)$ and $(S_2, \bullet)$ are semigroups.*
    ii.    *$(H \cap S_1, +)$ is a subsemigroup of $(S_1, +)$.*
    iii.   *$(H \cap S_2, \bullet)$ is a subsemigroup of $(S_2, \bullet)$.*

*Proof*: Let $H (\neq \phi)$ be a sub-bisemigroup. Therefore there exists two proper subsets $H_1, H_2$ of H such that

    i.     $H = H_1 \cup H_2$.
    ii.    $(H_1, +)$ is a semigroup.
    iii.   $(H, \bullet)$ is a semigroup.

Now we choose $H_1$ as $H \cap S_1$, then we see that $H_1$ is a subset of $S_1$ and by (ii) $(H_1, +)$ is itself a semigroup. Hence $(H_1 = H \cap S_1, +)$ is a subsemigroup of $S_1$. Similarly $(H_2 = H \cap S_2, \bullet)$ is a subsemigroup of $(S_2, \bullet)$. Conversely let (i), (ii) and (iii) of the statement of theorem be satisfied.

To prove $(H, +, \bullet)$ is a bisemigroup, it is enough to prove $(H \cap S_1) \cup (H \cap S_2) = H$.

$(H \cap S_1) \cup (H \cap S_2)$     =     $[(H \cap S_1) \cup H] \cap [(H \cap S_1) \cup S_2]$
    =     $[(H \cup H) \cap (H \cup S_1)] \cap [(H \cup S_2) \cap (S_1 \cup S_2)]$
    =     $[H \cap (S_1 \cup H)] \cap [(H \cup S_2) \cap S]$

$$= H \cap (H \cap S_2)$$
$$= H.$$

Hence the result.

It is important to note that in the above theorem just proved condition (i) can be removed. We include this condition only for clarification or simplicity of understanding.

**DEFINITION 3.1.3:** *Let $(S, +, \bullet)$ be a bisemigroup. A non-empty subset H of S is called the bi-ideal of the bisemigroup S if*

   i. $H = H_1 \cup H_2$; $(H_1, +)$ *and* $(H_2, \bullet)$ *are semigroups i.e. H is a sub-bisemigroup.*

   ii. *If for all $x \in S$ we have $h_1 x, x h_1 \in H$, whenever $h_1 \in H_1$ and $h_2 + x, x + h_2 \in H$ whenever $h_2 \in H_2$ if $x \in S_1 \cap S_2$ then we have $x h_1, h_1 x$ and $x + h_2, h_2 + x \in H$.*

*The notion of right bi-ideal and left bi-ideal is defined as that of left ideals and right ideals in a semigroup. If the bi-ideal H happens to be simultaneously both right and left then we call H a bi-ideal of the bisemigroup.*

***Example 3.1.5:*** Suppose $S = Z \cup Z_{10}$. Z is the semigroup under $\times$ and $Z_{10}$ the semigroup under $\times$ modulo 10, then $H = H_1 \cup H_2$ where $H_1 = \{2, 4, 6, 8, \ldots\}$ and $H_2 = \{0, 5\}$ is a bi-ideal of the bisemigroup S.

Thus we have an interesting theorem.

**THEOREM 3.1.4:** *Let $S = S_1 \cup S_2$ be a bisemigroup. Let $H_1$ and $H_2$ be ideals of $S_1$ and $S_2$ respectively. If $S_1 \cap S_2 = \phi$ then $H = H_1 \cup H_2$ is a bi-ideal of S.*

*Proof*: Given $S = S_1 \cup S_2$ is a bisemigroup. $H_1$ and $H_2$ are ideals of $S_1$ and $S_2$ respectively. If $S_1 \cap S_2 = \phi$ then clearly $H = H_1 \cup H_2$ is a bi-ideal of S.

**DEFINITION 3.1.4:** *Let $(S, +, \bullet)$ and $(L, +, \bullet)$ be any two bisemigroups; where $S = S_1 \cup S_2$ and $L = L_1 \cup L_2$. A map $\phi : S \rightarrow L$ is said to be a bisemigroup homomorphism if $\phi / S_1$ is a semigroup homomorphism from $S_1$ to $L_1$ and $\phi$ restricted to $S_2$ i.e. $\phi / S_2$ is a semigroup homomorphism from $S_2$ to $L_2$.*

**DEFINITION 3.1.5:** *Let $S = S_1 \cup S_2$ be a bisemigroup if $S_1$ and $S_2$ are free semigroups then we call S a free bisemigroup. If $S_1$ is free on the basis $B_1$ and $S_2$ is free on the basis $B_2$ then $B = B_1 \cup B_2$ is called free bi basis of $S = S_1 \cup S_2$, i.e. $S = S_1 \cup S_2$ is called free on the bi basis $B = B_1 \cup B_2$.*

**THEOREM 3.1.5:** *If $S = S_1 \cup S_2$ is a free bisemigroup with free bi basis $B = B_1 \cup B_2$ then $\langle B \rangle = S = \langle B_1 \rangle \cup \langle B_2 \rangle = S_1 \cup S_2$.*

*Proof*: Straightforward by the very concepts of free semigroups.

**THEOREM 3.1.6:** *For any set $B \neq \phi$; $B = B_1 \cup B_2$ there exists a bisemigroup $S = S_1 \cup S_2$ which is a free bisemigroup on B.*

*Proof*: As in case of free semigroups.

Now we know the concept of free semigroup are well exploited in theory of automaton and semi-automaton. Here we define bisemi-automaton and bi-automaton.

**DEFINITION 3.1.6:** *A bisemi-automaton is a triple $X = (Z, A, \delta)$ consisting of two non-empty sets Z and $A = A_1 \cup A_2$ and a function $\delta: Z \times (A_1 \cup A_2) \to Z$, Z is called the set of states and $A = A_1 \cup A_2$ the input alphabet and $\delta$ the next state function of X. We insist that if $X_1 = (Z, A_1, \delta)$ and $X_2 = (Z, A_2, \delta)$ then $X_1$ is not a subsemi-automaton of $X_2$ or $X_2$ is not a subsemi-automaton of $X_1$ but both $X_1$ and $X_2$ are semi-automatons. Thus in case of bisemi-automaton we consider the set A i.e. input alphabets as a union of two sets which may enjoy varied properties.*

*Hence we see in $A = A_1 \cup A_2$, $A_1$ may take states in a varied way and $A_2$ in some other distinct way. So when we take a set A and not the set $A_1 \cup A_2 = A$ there is difference for*

$$Z_r \xrightarrow{a_1} Z_t \quad a_1 \in A_1$$

$$Z_r \xrightarrow{a_1} Z_p \quad a_1 \in A_2$$

*Hence the same $a_1 \in A$ may act in two different ways depending on, to which $A_i$ it belongs to, or it may happen in some cases in $A = A_1 \cup A_2$ with $A_1 \cap A_2 = \phi$, so that they will have a unique work assigned by the function $\delta$.*

**DEFINITION 3.1.7:** *A biautomaton is a quintuple $Y = (Z, A, B, \delta, \lambda)$ where $(Z, A, \delta)$ is a bisemi-automaton and $A = A_1 \cup A_2$, $B = B_1 \cup B_2$ and $\lambda : Z \times A \to B$ is the output function. Here also $Y_1 = (Z, A_1, B_1, \delta, \lambda)$ and $Y_2 = (Z, A_2, B_2, \delta, \lambda)$ then $Y_1$ is not a subautomaton of $Y_2$ and vice versa but both $Y_1$ and $Y_2$ are automatons. The bisemi-automaton and biautomaton will still find its application in a very interesting manner when we replace the sets by free bisemigroups. We denote by $\overline{A} = \overline{A}_1 \cup \overline{A}_2$, $\overline{A}_1$ and $\overline{A}_2$ are free monoids generated by $A_1$ and $A_2$ so $\overline{A}$ denotes the free bimonoid as both have the common identity, $\Lambda$ the empty sequence.*

$$\overline{\delta}(z, \Lambda) = z.$$
$$\overline{\delta}(z, a_1) = \delta(z, a_1)$$
$$\overline{\delta}(z, a_2) = \delta(z, a_2)$$
$$\vdots$$
$$\overline{\delta}(z, a_n) = \delta(z, a_n)$$

*where $a_1, a_2, ..., a_n \in A$. Here $\overline{\delta}(z, a_i)$ may have more than one $z_i$ associated with it.*

$$\overline{\delta}(z, a_1 a_2) = \delta(\overline{\delta}(z, a_1), a_2)$$
$$\vdots$$
$$\overline{\delta}(z, a_1 a_2 \cdots a_r) = \delta(\overline{\delta}(z, a_1 a_2 \cdots a_{r-1}), a_r)$$

*Thus the bisemi-automaton may do several work than the semi-automaton and can also do more intricate work.*

*Now we define biautomaton as follows:*

$$\overline{Y} = (Z, \overline{A}, \overline{B}, \overline{\delta}, \overline{\lambda})$$

where $(Z, \overline{A}, \overline{\delta})$ is a bisemi-automaton. Here $B = B_1 \cup B_2$ and $\overline{B}$ is a free bimonoid.

$$\overline{\lambda}: Z \times \overline{A} \to \overline{B}$$
$$\overline{\lambda}(z, \Lambda) = \Lambda$$
$$\overline{\lambda}(z, a_1) = \lambda(z, a_1)$$
$$\overline{\lambda}(z, a_1 a_2) = \lambda(z, a_1) \overline{\lambda}(\delta(z, a_1), a_2)$$
$$\vdots$$
$$\overline{\lambda}(z, a_1 a_2 \cdots a_n) = \lambda(z, a_1) \overline{\lambda}(\delta(z, a_1), a_2 a_3 \cdots a_n).$$

Thus we see the new biautomaton $\overline{Y}$ can perform multifold and multichannel operations than the usual automaton.

**DEFINITION 3.1.8:** $Y_1 = (Z_1, A, B, \delta_1, \lambda_1)$ is called a sub-biautomaton of $Y_2 = (Z, A, B, \delta, \lambda)$ (in symbols $Y_1 \subset Y_2$ if $Z_1 \subset Z$ and $\delta_1$ and $\lambda_1$ are the restrictions of $\delta$ and $\lambda$ respectively on $Z_1 \times A$ (Here $A = A_1 \cup A_2$ and $B = B_1 \cup B_2$ are subsets).

**DEFINITION 3.1.9:** Let $Y_1 = (Z_1, A_1, B_1, \delta_1, \lambda_1)$ and $Y_2 = (Z_2, A_2, B_2, \delta_2, \lambda_2)$. A biautomaton homomorphism $\phi : Y_1 \to Y_2$ is a triple $\phi = (\xi, \alpha, \beta) \in Z_2^{Z_1} \times A_2^{A_1} \times B_2^{B_1}$ with the property (where $A_1 = A_{11} \cup A_{12}$, $A_2 = A_{21} \cup A_{22}$, $B_1 = B_{11} \cup B_{12}$ and $B_2 = B_{21} \cup B_{22}$).

$$\xi(\delta_1(z, a)) = \delta_2(\xi(z), \alpha(a))$$
$$\beta(\lambda_1(z, a)) = \lambda_2(\xi(z), \alpha(a))$$
$$(z \in Z, a \in A_1 = A_{11} \cup A_{12}).$$

$\phi$ is called a monomorphism (epimorphism, isomorphism) if all functions $\xi$, $\alpha$, $\beta$ are injective (surjective or bijective) where we make relevant adjustments when for a single $a \in A_1$ there are two state changes.

Similarly for a single $a \in A_1$ we may have outputs depending on the section where we want to send the message. Thus we see the biautomaton can simultaneously do the work of two automatons, which is very rare in case of automatons. So in a single machine we have double operations to be performed by it.

***Example 3.1.6:*** Let $Z = \{0, 1, 2, 3\}$, $A = A_1 \cup A_2$, where $A_1 = Z_2 = \{0, 1\}$ and $A_2 = \{0, 1, 2\} = Z_3$. $\delta : (z, a) = z.a$. Is $Y = (Z, A, \delta)$ a bisemi-automaton? The table for Y is

|   | δ | 0 | 1 | 2 | 3 |
|---|---|---|---|---|---|
| $A_1$ | 0 | 0 | 0 | 0 | 0 |
|       | 1 | 0 | 1 | 2 | 3 |
| $A_2$ | 0 | 0 | 0 | 0 | 0 |
|       | 1 | 0 | 1 | 2 | 3 |
|       | 2 | 0 | 2 | 0 | 2 |

The revised table is

| δ | 0 | 1 | 2 | 3 |
|---|---|---|---|---|
| 0 | 0 | 0 | 0 | 0 |
| 1 | 0 | 1 | 2 | 3 |
| 2 | 0 | 2 | 0 | 2 |

But if we consider the semi-automaton $Y_1 = \{Z, A_1, \delta\}$ and $(Z, A_2, \delta) = Y_2$ we get the following graphs:

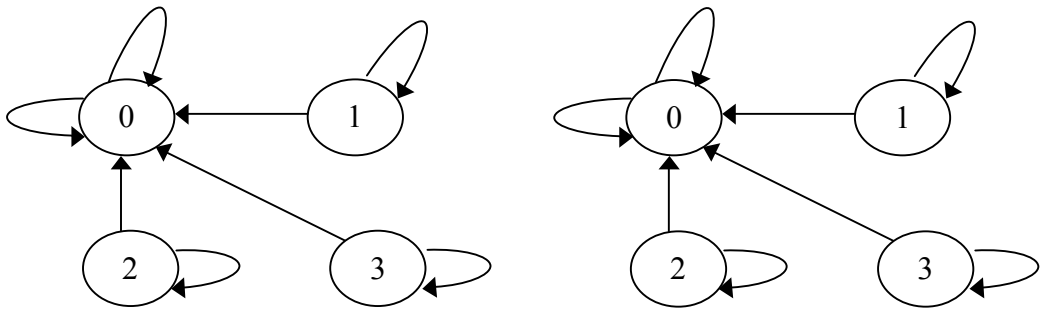

**Figures 3.1.1**

Suppose we take $A = \{0, 1, 2\}$ and divide it as $A'_1 = \{0, 1\}$ and $A'_2 = \{2\}$ then we get the graphs as

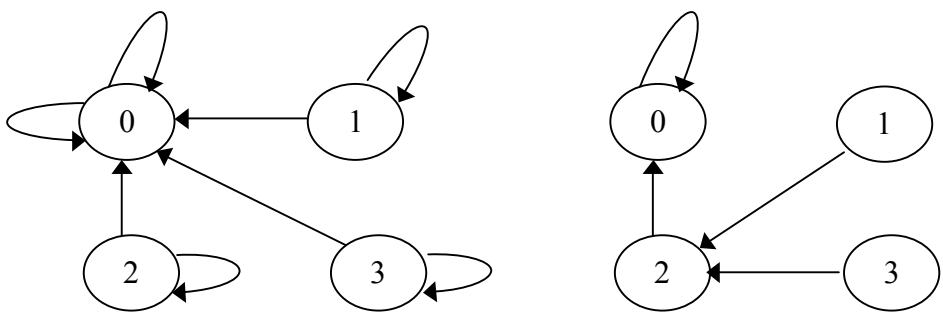

**Figures 3.1.2**

But $X'_1 = \{Z, A'_1, \delta\}$ and $X'_2 = \{Z, A'_2, \delta\}$ is a bisemi-automaton.

**Remark:** It is very important to note that how we choose the subsets $A_1$ and $A_2$ from A if $A_1 \subset A_2$ or $A_2 \subset A_1$. The semi-automaton will cease to be a bisemi-automaton.

Thus the very notion of bisets should have been defined even it may look trivial. So we at this juncture define bisets.

**DEFINITION 3.1.10:** *Let A be a set we say $A = A_1 \cup A_2$ is a biset if $A_1$ and $A_2$ are proper subsets of A such that $A_1 \not\subset A_2$ and $A_2 \not\subset A_1$.*

Unless this practical problem was not illustrated, the very notion of bisets would have been taken up researchers as a trivial concept. Only when bisets plays a vital role in determining whether a semi-automaton or an automaton is a bisemi-automaton or biautomaton.

Thus we give another formulation of the bisemi-automaton and biautomaton. Thus one of the necessary conditions for a semi-automaton to be a bisemi-automaton is that $A = A_1 \cup A_2$ is a biset. If $A = A_1 \cup A_2$ is not a biset then we do not have the semi-automaton to be a bisemi-automaton. But when we construct a bisemi-automaton $A = A_1 \cup A_2$ where $\delta$ is a special mapping then the map gives different positions to the same input symbol depending on which subset the input symbol is.

To this end we give the following example:

*Example 3.1.7*: Let $X = \{Z, A = A_1 \cup A_2, \delta\}$ be the bisemi-automaton, here $Z = \{0, 1, 2\}$, $A_1 = \{0, 1\}$ and $A_2 = \{0, 1, 2, 3\}$. $\delta(z, 0) = z + 0$ if $0 \in A_1$, $\delta(z, 0) = z.0$ if $0 \in A_2$, $\delta(z, 1) = z + 1$ if $1 \in A_1$, $\delta(z, 1) = z.1$ if $1 \in A_2$, $\delta(z, 2) = z.2$ (mod 3), $\delta(z, 3) = z.3$ (mod 3).

Now the diagram for the bisemi-automaton is as follows:

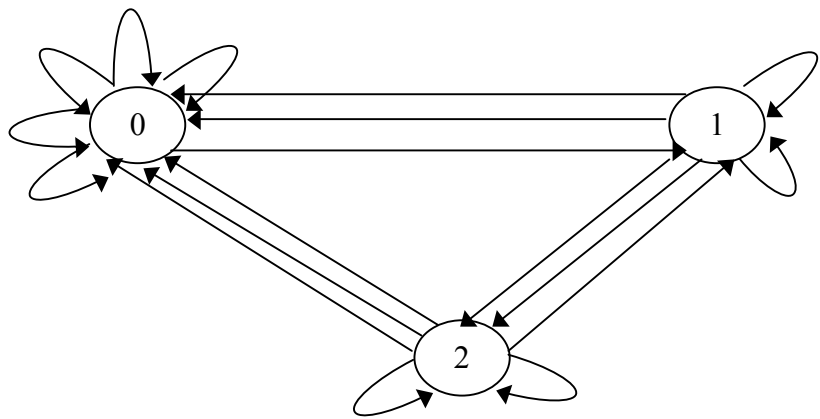

**Figure 3.1.3**

Thus this is a bisemi-automaton but not a semi-automaton as the function is not well defined. Thus we see all semiautomaton cannot be made into bisemiautomaton and in general a bisemiautomaton is not a semiautomaton. What best can be observed is that one can see that by building in a finite machine a bisemiautomaton structure one has the opportunity to make the machine perform more actions at a least cost and according to demand the bisets can be used with minute changes in the function. This bisemi-automaton will serve a greater purpose than that of the semi-automaton, secondly when we generate the free semigroups using the bisets or the sets in a bisemi-automaton we get nice results. Several applications can be found in this case.

Now at this juncture we wish to state that for the automaton also to make it a biautomatons both the input alphabets A and the output alphabets B must be divided into bisets to be the basic criteria to be a bi-automaton. The theory for this is similar to that of bisemi-automaton where in case of bisemi-automaton we used only the input alphabet A but in this case we will use both the input alphabets A and the output alphabets B.

Thus while converting a automaton to biautomaton we basically demand the sets A and B be divided into bisets such that the bisets give way to two distinct subautomaton. But if we are constructing biautomaton we define on $A = A_1 \cup A_2$ and $B = B_1 \cup B_2$ such that $\delta$ and $\lambda$ are not identical on $A_1 \cap A_2$ and $B_1 \cap B_2$ if they are non-empty.

The reader is given the task of finding results in this direction for this will have a lot of application in constructing finite machines.

Similarly when they are made into new biautomaton we take A to the free bisemigroup generated by $A_1$ and $A_2$ and B also to be the free bisemigroup generated by $B_1$ and $B_2$. We give two examples, one a automaton made into a biautomaton another is a constructed biautomaton.

*Example 3.1.8:* Let $Y = (Z, A, B, \delta, \lambda)$ where $Z = \{0, 1, 2, 3\}$, $A = \{0, 1, 2, 3, 4, 5\} = \{0, 2, 4\} \cup \{1, 3, 5\} = A_1 \cup A_2$, $B = \{0, 1, 2, 3\} = \{0, 2\} \cup \{1, 3\} = B_1 \cup B_2$.

$$\delta(z, a) = z + a \pmod 4 \text{ if } a \in A_1,$$
$$\delta(z, a) = z.a \pmod 4 \text{ if } a \in A_2.$$
$$\lambda(z, a) = z.a \pmod 4 \text{ if } a \in B_1,$$
$$\lambda(z, a) = z + a \pmod 4 \text{ if } a \in B_2.$$

Now we find the graphical representation of them.

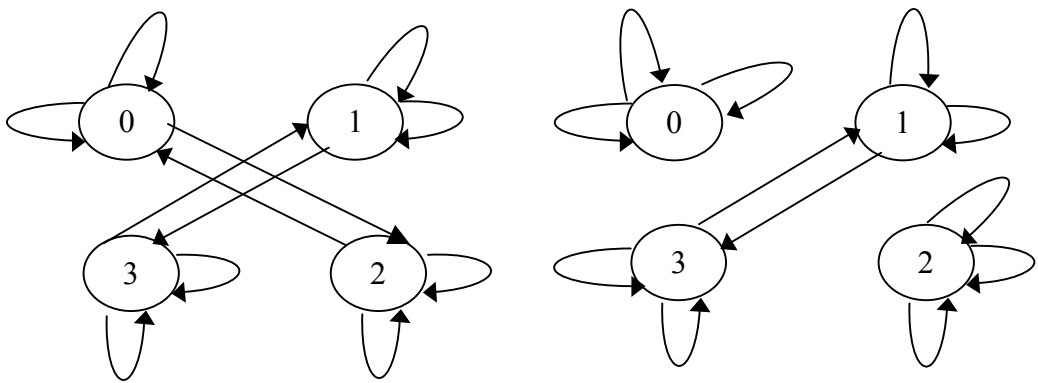

**Figures 3.1.4**

Output table:

| λ | 0 | 1 | 2 | 3 |
|---|---|---|---|---|
| 0 | 0 | 0 | 0 | 0 |
| 2 | 0 | 2 | 0 | 2 |

| λ | 0 | 1 | 2 | 3 |
|---|---|---|---|---|
| 1 | 0 | 1 | 2 | 3 |
| 3 | 0 | 3 | 2 | 1 |

We claim Y is the biautomaton as the output functions do not tally in both the tables. As the main motivation is study of Smarandache bialgebraic structure we give one more example of an automaton, which is a biautomaton.

*Example 3.1.9:* Let Y = (Z, A, B, δ, λ) where A = {0, 1, 2} and B = {0, 1, 2, 3}. Z = {0, 1, 2, 3, 4}. δ(z, a) = z.a (mod 5), λ(z, a) = 2za (mod 4).

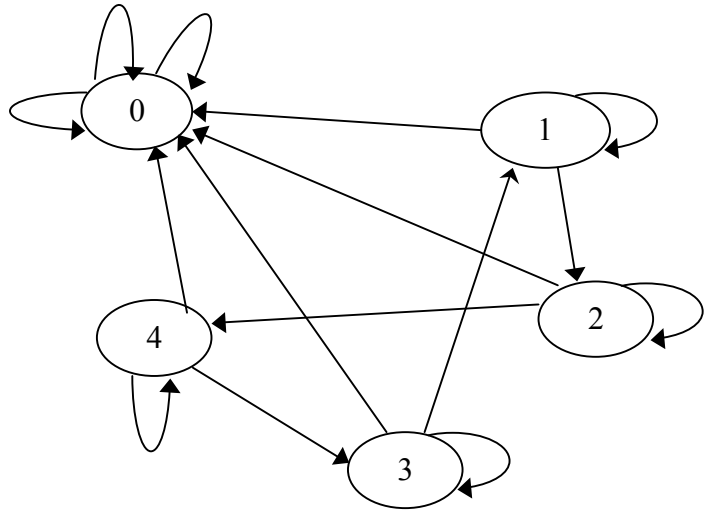

**Figure 3.1.5**

Next state function table:

| δ | 0 | 1 | 2 | 3 | 4 |
|---|---|---|---|---|---|
| 0 | 0 | 0 | 0 | 0 | 0 |
| 1 | 0 | 1 | 2 | 3 | 4 |
| 2 | 0 | 2 | 4 | 1 | 3 |

*Output table:*

| δ | 0 | 1 | 2 | 3 | 4 |
|---|---|---|---|---|---|
| 0 | 0 | 0 | 0 | 0 | 0 |
| 1 | 0 | 2 | 0 | 2 | 0 |
| 2 | 0 | 0 | 0 | 0 | 0 |
| 3 | 0 | 2 | 0 | 2 | 0 |

If we take B = B$_1$ ∪ B$_2$ where B$_1$ = {0, 1} and B$_2$ = {2, 3} and A = A$_1$ ∪ A$_2$ = {0, 1} ∪ {2} we get a biautomaton.

**PROBLEMS**

1. Give an example of a non-commutative bisemigroup.
2. Is S = S$_1$ ∪ S$_2$ where S$_1$ = Z$_{15}$ semigroup under multiplication modulo 15 and S$_2$ = S(5) a commutative bisemigroup?
    i. Find the order of S.
    ii. Does S contain sub-bisemigroup?
    iii. Can S have S-bi-ideal?
3. Give an example of a bisemiautomaton with 5 sets of states.
4. Find a biautomaton with 3 states i.e., (Z = {Z$_1$, Z$_2$, Z$_3$}).
5. Is Y = (Z, A, B, δ, λ) where A = A$_1$ ∪ A$_2$, Z = {0, 1, 2, … , 6} = Z$_7$, B = B$_1$ ∪ B$_2$, δ (z, a) = z + a (mod 7), δ (z, a) = 3za (mod 4) for suitable A and B be a biautomaton?

## 3.2 Biquasi groups and its properties

In this section we still define a new notion called biquasi groups and obtain its properties. All bigroups are trivially biquasi groups i.e. biquasi groups are a generalization of bigroups and all biquasi groups are bisemigroup. Thus bisemigroups form the most generalized class in case of associative structures.

Thus {Class of bigroups} ⊂ {Class of biquasi groups} ⊂ {Class of bisemigroups}. We indicate some of its basic properties and leave the rest for the reader.

**DEFINITION 3.2.1:** *Let (X, +, •) be a set with X = X$_1$ ∪ X$_2$ we call X a biquasi group if (X$_1$, +) is a group and (X, •) is a semigroup.*

*Example 3.2.1:* Take X = {1, 2, 3, 4, 0, 6, 8} under '×' - modulo 5 and × modulo 10. Then X$_1$ = {1, 2, 3, 4} is a group under '×' modulo 5 and X$_2$ = {0, 2, 4, 6, 8} is semigroup under '×' modulo 10. Thus X is a biquasi group.

**THEOREM 3.2.1:** *Every bigroup is a biquasi group and not conversely.*

*Proof:* By the very definition of bigroup and biquasi group we see every bigroup is a biquasi group.

To prove the converse we see from example 3.2.1 the set X is a biquasi group. Clearly not a bigroup as X$_2$ = {0, 2, 4, 6, 8} is not a group under '×' modulo 10. Hence the claim.

**THEOREM 3.2.2:** *Every biquasi group is a bisemigroup. But a bisemigroup in general is not a biquasi group.*

*Proof:* By the very definitions of bisemigroups and biquasi groups we see every biquasi group is a bisemigroup. To show every bisemigroup in general is not a biquasi group consider $B = B_1 \cup B_2$ where $B_1 = \{0, 1, 2, \ldots, 8\}$ under multiplication and $B_2 = S(3)$ the symmetric semigroup. Clearly B is a bisemigroup and it is not a biquasi group.

**DEFINITION 3.2.2:** *Let $X = X_1 \cup X_2$ be a biquasi group. Let $P \subset X$ be a proper subset of X we say P is a sub-biquasi group if P itself under the operations of X is a biquasi group.*

*Example 3.2.2:* Let $X = Z_{15} \cup S_3$. $(Z_{15}, \times)$ is a semigroup and $S_3$ is the group of degree 3. Clearly

$$P = \{0, 5, 10\} \cup \left\{1, \begin{pmatrix} 1 & 2 & 3 \\ 1 & 3 & 2 \end{pmatrix}\right\}$$

is a sub-biquasi group of X.

**DEFINITION 3.2.3:** *Let $(X, +, \bullet)$ be a biquasi group; $X = X_1 \cup X_2$; $(X_1, +)$ is a group and $(X_2, \bullet)$ is a semigroup. $P = P_1 \cup P_2$ be a sub-biquasi group of X. We say P is a normal sub-biquasi group of X if $P_1$ is a normal subgroup of $(X_1, +)$ and $(P_2, \bullet)$ is an ideal of $(X_2, \bullet)$.*

**DEFINITION 3.2.4:** *Let $(X, +, \bullet)$ be a biquasi group. We say X is of finite order if the number of elements in X is finite. We denote the order of X by $o(X)$ or $|X|$. If the number of elements in X is infinite we say X is of infinite order.*

**DEFINITION 3.2.5:** *Let $(X, +, \bullet)$ be a biquasi group of finite order. If the order of every proper sub-biquasi group P of X divides the order of X then we call the biquasi group to be a Lagrange biquasi group.*

It is to be noted that unlike finite groups where the order of every subgroup divides the order of the group we see in case of biquasi group of finite order every finite sub-biquasi group need not divide the order of the biquasi group.

We illustrate it by the following example.

*Example 3.2.3:* Let $(X, +, \bullet)$ be a biquasi group; with $X = X_1 \cup X_2$ where $X_1 = \{g \mid g^9 = 1\}$ and $X_2 = S(2)$ and $|X| = 13$. We see

$$P = \left\{1, g^3, g^6, \begin{pmatrix} 1 & 2 \\ 1 & 2 \end{pmatrix}\right\}$$

is a sub-biquasi group.

$$P = P_1 \cup P_2 = \{1, g^3, g^6\} \cup \left\{\begin{pmatrix} 1 & 2 \\ 1 & 2 \end{pmatrix}\right\}$$

is of order four, but 4 ∤ 13. Thus (X, +, •) is not a Lagrange biquasi group.

**DEFINITION 3.2.6:** *Let (X, +, •) be a biquasi group of finite order. An element $x \in X$ will be called as Cauchy element if there exists a number n such that $n \mid |X|$ and $x^n = e$.*

We may have biquasi groups, which may not have Cauchy element. We illustrate this by an example.

*Example 3.2.4:* Let (X, +, •) be a biquasi group. $X = G_1 \cup S(3)$ where $G_1 = \langle g \mid g^4 = 1 \rangle$ and S(3) is the symmetric semigroup of order 3. Clearly |X| = 31. We have no element in X to be a Cauchy element of X.

*Example 3.2.5:* Let (X, +, •) be a biquasi group. $X = X_1 \cup X_2$; $X_1 = \langle g \mid g^{10} = 1 \rangle$ and $X_2 = S(2)$. |X| = 14. Take $g^5 \in X$ we see $(g^5)^2 = 1$ so 2/14 i.e. $g^5$ is a Cauchy element. Take $g^2 \in X$; $(g^2)^5 = 1$ so 5 ∤ 14 thus $g^2$ is not a Cauchy element.

Hence in a same biquasi group we may have some Cauchy elements and some elements, which are not Cauchy elements.

**DEFINITION 3.2.7:** *Let (X, +, •) be a biquasi group of finite order say n i.e. o(X) = n. If p is a prime such that $p^\alpha \mid n$ and $p^{\alpha+1} \nmid n$ but X has a sub-biquasi group P of order $p^\alpha$ then we call P a p-Sylow sub-biquasi group. Unlike finite group we see in case of biquasi group we may have a sub-biquasi group of order $p^\alpha$ but $p^\alpha \nmid$ (order of X) or o(P)/ |X| but still X has no sub-biquasi group of order $P^\alpha$.*

*Example 3.2.6:* Let (X, +, •) be a biquasi group, where $X = Z_{12} \cup Z_{10}$ with $Z_{12}$ is a semigroup under '×' modulo 12 and $Z_{10}$ is the group under '+' modulo 10. |X| = 22, 11/22 take P = {0, 2, 4, 6, 8, 10} ∪ {0, 2, 4, 6, 8}. P is a biquasi group and |P| / 22 Also $P_1$ = {1} ∪ {0}; 2/22 |$P_1$| = 2. Thus X is a Sylow biquasi group.

**DEFINITION 3.2.8:** *Let $X = X_1 \cup X_2$ be a biquasi group of finite order n. If for every prime p, $p^\alpha / n$ and $p^{\alpha+1} \nmid n$ we have a sub-biquasi group of order p we call X a Sylow biquasi group.*

The above example is a Sylow biquasi group. Now we define commutative bisemigroups of a special type.

**DEFINITION 3.2.9:** *Let $X = X_1 \cup X_2$ be bisemigroup which is commutative and in which every x in X is an idempotent. We call this class of semigroups as bisemilattices.*

The definition of bisemilattice given above by [125]. Several interesting results in this direction can be had from [125].

Thus bisemilattice which are bisemigroups in which every element is an idempotent cannot yield nontrivial subgroup in them; for us to define Smarandache analogue so we have left them without deeper discussions.

**PROBLEMS:**

1. What is the order of the smallest biquasi group?
2. Give an example of a biquasi group of order 29.
3. Find all sub-biquasi groups of $(X, +, \bullet)$ where $X = X_1 \cup X_2$; $X_1 = Z_{12}$ is a group under '+' modulo 12 and $X_2 = S(2)$ is a biquasi group of order 16.
4. Is the biquasi group given in problem 3 a Sylow biquasi group? Justify your answer.
5. Does the biquasi group $X = X_1 \cup X_2$ where $X_1 = S_7$ and $X_2 = S(4)$ have normal sub-biquasi group? If so find them. What is the order of X?
6. Find all sub-biquasi groups for X given in problem 5. Does X have Cauchy element?

## 3.3 S-bisemigroups and S-biquasi groups and its properties

In this section we define the notion of Smarandache bisemigroups and Smarandache biquasi groups and discuss and introduce several interesting properties. As this concept is very new the reader is advised to assimilate the definitions of bigroups, Smarandache bigroups, bisemigroups and biquasi groups.

**DEFINITION 3.3.1:** *Let $(S, +, \bullet)$ be a bisemigroup. We call $(S, +, \bullet)$ a Smarandache bisemigroup (S-bisemigroup) if S has a proper subset P such that $(P, +, \bullet)$ is a bigroup under the operations of S.*

*Example 3.3.1:* Let $(S, +, \bullet)$ be a bisemigroup where $S = S_1 \cup S_2$ with $S_1 = S(4)$ and $S_2 = Z_{12} = \{0, 1, 2, \ldots, 11\}$ under multiplication modulo 12. $(S, +, \bullet)$ is a S-bisemigroup as $P = S_4 \cup \{4, 8\}$ is a bigroup. Hence the claim.

All bisemigroups in general are not S-bisemigroups.

*Example 3.3.2:* Let $(S, +, \bullet)$ be a bisemigroup. $S = S_1 \cup S_2$ where $S_1$ and $S_2$ are given by the following tables:

Table of $S_1$

| * | 0 | a | b | c |
|---|---|---|---|---|
| 0 | 0 | 0 | 0 | 0 |
| a | 0 | a | 0 | 0 |
| b | 0 | 0 | b | 0 |
| c | 0 | 0 | 0 | c |

Table of $S_2$

| * | 0 | $x_1$ | $x_2$ | $x_3$ | $x_4$ |
|---|---|---|---|---|---|
| 0 | 0 | 0 | 0 | 0 | 0 |
| $x_1$ | 0 | $x_1$ | 0 | 0 | $x_1$ |
| $x_2$ | 0 | 0 | $x_2$ | $x_2$ | 0 |
| $x_3$ | 0 | 0 | $x_2$ | $x_3$ | 0 |
| $x_4$ | 0 | $x_1$ | 0 | 0 | $x_4$ |

Cleary every element in both $S_1$ and $S_2$ are idempotents we see S is not S-bisemigroup. Thus we have the class of bisemilattices to be the set of bisemigroups, which are never S-bisemigroups.

**DEFINITION 3.3.2:** *Let (S, +, •) be a bisemigroup. A proper sub-bisemigroup P of S is said to be a Smarandache sub-bisemigroup (S-sub-bisemigroup) if P itself has a proper subset Q such that Q under the operations of S is a bigroup.*

**THEOREM 3.3.1:** *Let (S, +, •) be a bisemigroup; if S has a proper S-sub-bisemigroup then S is a S-bisemigroup.*

*Proof*: Direct by the very definition, hence left for the reader to prove.

**DEFINITION 3.3.3:** *Let (S, +, •) be a S-bisemigroup. We say (S, +, •) is Smarandache commutative bisemigroup (S-commutative bisemigroup) if every proper subset P of S which is a bigroup happens to be a commutative bigroup.*

*We say S is a Smarandache weakly commutative bisemigroup (S-weakly commutative bisemigroup) if atleast one of the proper subset, which is a bigroup, is a commutative bigroup.*

**DEFINITION 3.3.4:** *Let (S, +, •) be a S-bisemigroup. We say S is a Smarandache cyclic bisemigroup (S-cyclic bisemigroup) if every proper subset P of S that is a bigroup is a cyclic bigroup. If atleast one of the proper subset P of S that is a bigroup is cyclic then we say S is a Smarandache weakly cyclic bisemigroup (S-weakly cyclic bisemigroup).*

**Note**: Let S be a S-bisemigroup, the number of elements in S is called the order of S denoted by o(S) or |S| if the number of elements in S is finite we call S a finite order S-bisemigroup, if number of elements in S is infinite, we call S a S-bisemigroup of infinite order.

*Example 3.3.3:* Let (S, +, •) be a S-bisemigroup; where $S = Z_6 \cup S(2)$. Clearly (S, +, •) is a S-cyclic bisemigroup.

*Example 3.3.4:* Let (S, +, •) be a S-bisemigroup where $S = S(3) \cup S_4$. S is S-weakly commutative but not S-commutative for $P = \{S_3 \cup A_4\}$ is a proper subset which is not a commutative bigroup as $A_4$ is a non-commutative subgroup of $S_4$.

*Example 3.3.5:* Let (S, +, •) be a bisemigroup. $S = Z_{12} \cup Z^+$, $Z_{12}$ is a semigroup under product modulo 12 and $Z^+$ is semigroup under multiplication; S is not a S-bisemigroup. But S is a commutative bisemigroup. All bisemigroups are not S-bisemigroup.

*Example 3.3.6:* Let (S, +, •) be a bisemigroup. $S = Z_{12} \cup Q^+$, $Z_{12}$ is a semigroup under multiplication modulo 12 and $Q^+$ is a semigroup under usual multiplication. S is a S-bisemigroup.

**DEFINITION 3.3.5:** *Let (S, +, •) be a S-bisemigroup. A proper subset A of S, be a sub-bigroup of S which cannot be contained in any proper sub-bisemigroup of S. Then we* *call A the largest sub-bigroup of S. Let S be a S-bisemigroup. If A be a proper subset of S which is a S-sub-bisemigroup of S and contains the largest bigroup of S then we say A to be the Smarandache hyper bisubsemigroup (S-hyper bisubsemigroup).*

**DEFINITION 3.3.6:** *Let $(S, +, \bullet)$ be a bisemigroup. An element $x \in S \setminus \{1\}$ is said to have a Smarandache inverse (S-inverse) y in G if $xy = 1$ and for $a, b \in S \setminus \{x, y, 1\}$ we have $xa = y$ (or $ax = y$), $yb = x$ (or $by = x$) with $ab = 1$.*

*Let $(S, +, \bullet)$ be a bisemigroup, x be the S-inverse of y and (x, y) the S-inverse pair with the related pair (a, b). If the pair (a, b) happens to be a S-inverse pair not necessarily with (x, y) as the related pair then we say (a, b) is the Smarandache co inverse pair (S-co inverse pair).*

**Note**: *It is to be noted all elements in the bisemigroup in general need not have S-inverse or S-co inverse.*

**DEFINITION 3.3.7:** *Let $(S, +, \bullet)$ be a bisemigroup; if no element in S has a S-inverse then we call S a S-inverse free group.*

**DEFINITION 3.3.8:** *Let $x \in S$, $(S, +, \bullet)$ be a bisemigroup. Let x have a S-inverse y for the pair (x, y), let (a, b) is the related pair. We say (x, y) is a Smarandache self-inversed pair (S-self-inversed pair) if (a, b) has the S-inverse and the related pair is (x, y).*

**DEFINITION 3.3.9:** *Let $(S, +, \bullet)$ be a bisemigroup. Let $x \in S$, x is said to have a Smarandache conjugate (S-conjugate) y in G if*

    i.   *x is conjugate to y i.e. there exists $a \in G$ such that $x = aya^{-1}$.*
   ii.   *a is conjugate with x and a is conjugate with y.*

**DEFINITION 3.3.10:** *Let $S_1, \ldots, S_n$ be n-S-bisemigroups. $S = S_1 \times \ldots \times S_n = \{(s_1, s_2, \ldots, s_n) \mid s_i \in S_i, i = 1, 2, \ldots, n\}$ is called the Smarandache direct product (S-direct product) of the S-bisemigroups; $S_1, S_2, \ldots, S_n$ if S is a maximal bisemigroup and G is got from the $S_1, \ldots, S_n$ as $G = G_1 \times \ldots \times G_n$ where each $G_i$ is the maximal sub-bigroup of the S bisemigroup $S_i$ for $i = 1, 2, \ldots, n$.*

*Here by a maximal bisubgroup of S we mean a proper subset $M \subset S$ (S a bisemigroup); M is the maximal sub-bigroup of S that is N is a sub-bigroup such that $M \subset N$ then $M = N$ is the only possibility.*

*If S the bisemigroup has only maximal sub-bigroup we call S a Smarandache maximal bisemigroup (S-maximal bisemigroup).*

**DEFINITION 3.3.11:** *Let S be a S-bisemigroup; $A_1, A_2, \ldots, A_n$ be non-empty subsets of S. We say $A_1, A_2, \ldots, A_s$ is the Smarandache internal direct product (S-internal direct product) of S if $S = A_1, A_2, \ldots, A_n = \{(a_1, \ldots, a_n) \mid a_i \in A_i, i = 1, 2, \ldots, n\}$ and this accounts for all elements of S. If $S = B * A_1 * \ldots * A_n$ where B is S-bisemigroup and $A_1, \ldots, A_n$ are maximal bisubgroup of $S_i$ then we say S is a Smarandache strong internal direct product (S-strong internal direct product).*

We by these S-strong internal direct product, S-internal direct product etc., we get more number of S-bisemigroups.

We recall some of its basic applications to the bi-semi-automaton and biautomaton.

**DEFINITION 3.3.12:** *A Smarandache bisemi-automaton (S-bisemi-automaton) is a triple $X = (Z, A, \delta)$ if $X = X_1 \cup X_2$ where $X_1 = (Z, A_1, \delta)$ and $X_2 = (Z, A_2, \delta)$ are two distinct bisemi-automatons (i.e. $X_1$ is not a sub-bisemi-automaton of $X_2$ or vice versa).*

*Thus every S-bisemi-automaton is a bisemi-automaton but every bisemi-automaton in general is not a S-bisemi-automaton.*

The reader is advised to obtain interesting properties about these newly constructed structures.

**DEFINITION 3.3.13:** *A Smarandache biautomaton (S-biautomaton) is a quintuple $Y = (Z, A, B, \delta, \lambda)$ where $(Z, A, \delta)$ is a S-bisemi-automaton and $Y = Y_1 \cup Y_2$ where both $Y_1$ and $Y_2$ are distinct bi-automatons of Y.*

*Thus we see every biautomaton is a S-biautomaton and in general every biautomaton is not a S-biautomaton.*

*We define new Smarandache bisemi-automaton and Smarandache biautomaton.*

**DEFINITION 3.3.14:** *Let $X_s = (Z, A, \delta)$ be a S-bisemi-automaton; if the free semigroup generated by A is a S-free semigroup then we call $\overline{X}_s = (Z, \overline{A}, \overline{\delta})$ a new S-bisemi-automaton.*

Note all new S-bisemi-automaton are S-bisemi-automatons.

**DEFINITION 3.3.15:** *Let $Y_s = (Z, A, B, \delta, \lambda)$ be a S-biautomaton. We call $\overline{Y}_s = (Z, \overline{A}, \overline{B}, \overline{\delta}, \overline{\lambda})$ to be a Smarandache new bi-automaton (S-new-bi-automaton) if $\overline{A}$ and $\overline{B}$ are S-free semigroups generated by A and B respectively.*

*Thus we see all new S-bi-automaton are new biautomaton and not conversely.*

In the last section we have defined the concept of bisets. Now we define the new notion called Smarandache bisets.

**DEFINITION 3.3.16:** *Let A be a biset i.e. $A = A_1 \cup A_2$ we call A Smarandache biset (S-biset) if in each of $A_1$ and $A_2$ we can find some proper subsets $P_1$ and $P_2$ which can be made into semigroups.*

*Example 3.3.7:* Let $A = \{2, 3, 9, -1, 1, 0\}$ be a set. Clearly $A = \{2, 3, 9, 0\} \cup \{-1, 1, 0, 2\}$ is a biset. Now A can be made into a S-biset for take $P_1 = \{0, 9, 3\}$ is a semigroup under multiplication modulo 27

| × | 0 | 3 | 9 |
|---|---|---|---|
| 0 | 0 | 0 | 0 |
| 3 | 0 | 9 | 0 |
| 9 | 0 | 0 | 0 |

and $P_2 = \{0, 1, 2\}$ is a semigroup under multiplication modulo 3 for

| × | 0 | 1 | 2 |
|---|---|---|---|
| 0 | 0 | 0 | 0 |
| 1 | 0 | 1 | 2 |
| 2 | 0 | 2 | 1 |

Thus A is a S-biset. It is important to note that all bisets in general are not S-bisets. For take $A = \{0, 1\}$. Clearly $A = \{0\} \cup \{1\}$ is biset but A can never be a S-biset.

Now we proceed on to define Smarandache biquasi groups.

**DEFINITION 3.3.17:** *Let $(T, +, \bullet)$ be a non-empty set with two binary operations $+$ and '$\bullet$'. We call T a Smarandache biquasi group (S-biquasi group) if $(T, +)$ is a S-semigroup and $(T, \bullet)$ is just a semigroup.*

Thus we see every biquasi group is a S-biquasi group; but a S-biquasi group in general is not a biquasi group.

**THEOREM 3.3.2:** *Every S-biquasi group is a bisemigroup.*

*Proof:* Straightforward by the very definitions.

**THEOREM 3.3.3:** *A bisemigroup in general need not be a S-biquasi group.*

*Proof:* By an example. Consider the bisemigroup $S = S_1 \cup S_2 = \{0, a \mid a^2 = 0\} \cup \{x, y, z \mid xy = z, xz = 0, zx = 0, z^2 = x^2 = y^2 = 0, zy = yz = 0\}$. Clearly S is a bisemigroup but is not a S-biquasi group but is not a S-biquasi group. Hence the claim.

**DEFINITION 3.3.18:** *Let $X = X_1 \cup X_2$ be a biquasi group. A proper subset $P = P_1 \cup P_2 \subset X$ is said to be a Smarandache sub-biquasi group (S-sub-biquasi group) if P itself is a S-biquasi group.*

**THEOREM 3.3.4:** *Let $X = X_1 \cup X_2$ be a biquasi group, then if $P \subset X$ is a S-subquasi group then X is a S-biquasi group.*

*Proof:* Direct; left for the reader to prove.

**DEFINITION 3.3.19:** *Let $X = X_1 \cup X_2$ be a biquasi group. A proper subset $P \subset X$ is said to be a S-normal sub-biquasi group if*

  i.  *P is a S-sub-biquasi group.*
  ii. *$P = P_1 \cup P_2$ and both $P_1$ and $P_2$ are ideals of $X_1$ and $X_2$ respectively.*

*We say the order of X to be of finite order if the number of elements in X is finite; other wise X is said to be of infinite order; we denote the order of X by o(X) or |X|.*

**DEFINITION 3.3.20:** *Let (X, +, •) be a S-biquasi group of finite order. If the order of every S sub-biquasi group divides the order of X; then we call the S biquasi group to be a Smarandache Lagrange biquasi group (S-Lagrange biquasi group). If there exists atleast one S-sub-biquasi group whose order divides the order of X then we call X a Smarandache weakly Langranges biquasi group (S-weakly Lagrange biquasi group).*

**DEFINITION 3.3.21:** *Let (X, +, •) with X = $X_1 \cup X_2$ be a S biquasi group of finite order. We call an element x ∈ X a Smarandache Cauchy element (S-Cauchy element) if there exists an integer n such that if x ∈ X, where $X_1$ is a S-semigroup and $x^n = e$ and n / |$A_1$| where $A_1 \subset X_1$ is a proper subset, which is a group under the operations of $X_1$.*

***Example 3.3.8:*** Let X = S(3) ∪ {0, a, b | $a^2 = b^2 = 0$ and ab = ba = 0}. Clearly X is a S-biquasi group. |X| = 27 + 3 = 30. Now

$$x = \begin{pmatrix} 1 & 2 & 3 \\ 2 & 3 & 1 \end{pmatrix} \in S(3)$$

with

$$x^3 = \begin{pmatrix} 1 & 2 & 3 \\ 1 & 2 & 3 \end{pmatrix}.$$

Clearly 3 / |$S_3$| where |$S_3$| = 6. Thus all elements in X are S-Cauchy elements.

**DEFINITION 3.3.22:** *Let (X, +, •) be a biquasi group of finite order. If p is a prime such that $p^\alpha$ | o(X) but $p^{\alpha+1} \nmid$ o(X) but X has a S-sub-biquasi group A of order $p^\alpha$ then, we call X a Smarandache p-Sylow sub-biquasi group (S-p-Sylow sub-biquasi group). If for every prime p such that $p^\alpha$ | o(X) and $p^{\alpha+1} \nmid$ o(X) we have an associated S-sub-biquasi group B of order $p^\alpha$ then we call X a Smarandache strong p-Sylow bi quasi group (S-strong p-Sylow biquasi group).*

**PROBLEMS:**

1. Find the order of the smallest S- bisemigroup.
2. Is (X, +, •) with X = S(3) ∪ $Z_{18}$ where $Z_{18}$ is a semigroup under '×' modulo 18, a S-bisemigroup?
3. Can a S-bisemigroup of order 11 exist?
4. Find all S-cyclic sub-bisemigroups of X, where X = $X_1 \cup X_2$ with $X_1$ = S(4) and $X_2 = Z_{27}$, semigroup under '×' modulo 27.
5. Give an example of a S-bisemi-automaton.
6. Illustrate by an example all biautomaton need not in general be a S-biautomaton.
7. Does the S bisemigroup given by X = $Z_{10} \cup$ S(5) have S-cyclic bisemigroup?

**Chapter 4**

# BILOOPS AND SMARANDACHE BILOOPS

This chapter has two sections. In section one we introduce the concept of biloops. The study of biloops started in the year 2001 by the author [110]. Here we bring in some of its properties and its basic applications. In section two for the first time define Smarandache biloops and give several interesting properties about them.

## 4.1 Biloops and its properties

In this section we give the definition of biloops and enumerate some of its properties.

**DEFINITION 4.1.1:** *Let $(L, +, \bullet)$ be a non empty set with two binary operations. L is said to be a biloop if L has two nonempty finite proper subsets $L_1$ and $L_2$ of L such that*

    i.   $L = L_1 \cup L_2$.
    ii.  $(L_1, +)$ is a loop.
    iii. $(L_2, \bullet)$ is a loop or a group.

*For a biloop we demand atleast one of $(L_1, +)$ or $(L_2, \bullet)$ to be a loop which is not a group.*

*Also $L_1$ and $L_2$ must be proper subsets of L i.e. $L_1 \not\subset L_2$ or $L_2 \not\subset L_1$.*

***Example 4.1.1:*** Let $(L, +, \bullet)$ be a biloop with $L = L_1 \cup L_2$ where $(L_1, +)$ is the loop given by the following table:

| +   | e   | $g_1$ | $g_2$ | $g_3$ | $g_4$ | $g_5$ |
|-----|-----|-------|-------|-------|-------|-------|
| e   | e   | $g_1$ | $g_2$ | $g_3$ | $g_4$ | $g_5$ |
| $g_1$ | $g_1$ | e   | $g_3$ | $g_5$ | $g_2$ | $g_4$ |
| $g_2$ | $g_2$ | $g_5$ | e   | $g_4$ | $g_1$ | $g_3$ |
| $g_3$ | $g_3$ | $g_4$ | $g_1$ | e   | $g_5$ | $g_2$ |
| $g_4$ | $g_4$ | $g_3$ | $g_5$ | $g_2$ | e   | $g_1$ |
| $g_5$ | $g_5$ | $g_2$ | $g_4$ | $g_1$ | $g_3$ | e   |

and $(L_2, \bullet)$ be the cyclic group of order 7.

**DEFINITION 4.1.2:** *Let $(L, +, \bullet)$ be a biloop. L is said to be a commutative biloop if $L = L_1 \cup L_2$ then $(L_1, +)$ and $(L_2, \bullet)$ are both commutative loops.*

***Example 4.1.2:*** Let $(L, +, \bullet)$ be a biloop. $L = L_1 \cup L_2$, here the loops $(L_1, +)$ and $(L_2, \bullet)$ are given by the following tables:

$L_1 = \{e, g_1, g_2, g_3, g_4, g_5\}$ and $L_2 = \{e, a_1, a_2, \ldots, a_7\}$.

Table for $L_1$

| + | e | $g_1$ | $g_2$ | $g_3$ | $g_4$ | $g_5$ |
|---|---|---|---|---|---|---|
| e | e | $g_1$ | $g_2$ | $g_3$ | $g_4$ | $g_5$ |
| $g_1$ | $g_1$ | e | $g_4$ | $g_2$ | $g_5$ | $g_3$ |
| $g_2$ | $g_2$ | $g_4$ | e | $g_5$ | $g_3$ | $g_1$ |
| $g_3$ | $g_3$ | $g_2$ | $g_5$ | e | $g_1$ | $g_4$ |
| $g_4$ | $g_4$ | $g_5$ | $g_3$ | $g_1$ | e | $g_2$ |
| $g_5$ | $g_5$ | $g_3$ | $g_1$ | $g_4$ | $g_2$ | e |

Table for $L_2$

| * | e | $a_1$ | $a_2$ | $a_3$ | $a_4$ | $a_5$ | $a_6$ | $a_7$ |
|---|---|---|---|---|---|---|---|---|
| e | e | $a_1$ | $a_2$ | $a_3$ | $a_4$ | $a_5$ | $a_6$ | $a_7$ |
| $a_1$ | $a_1$ | e | $a_5$ | $a_2$ | $a_6$ | $a_3$ | $a_7$ | $a_4$ |
| $a_2$ | $a_2$ | $a_5$ | e | $a_6$ | $a_3$ | $a_7$ | $a_4$ | $a_1$ |
| $a_3$ | $a_3$ | $a_2$ | $a_6$ | e | $a_7$ | $a_4$ | $a_1$ | $a_5$ |
| $a_4$ | $a_4$ | $a_6$ | $a_3$ | $a_7$ | e | $a_1$ | $a_5$ | $a_2$ |
| $a_5$ | $a_5$ | $a_3$ | $a_7$ | $a_4$ | $a_1$ | e | $a_2$ | $a_6$ |
| $a_6$ | $a_6$ | $a_7$ | $a_4$ | $a_1$ | $a_5$ | $a_2$ | e | $a_3$ |
| $a_7$ | $a_7$ | $a_4$ | $a_1$ | $a_5$ | $a_2$ | $a_6$ | $a_3$ | e |

Clearly $L = L_1 \cup L_2$. is a commutative biloop.

**DEFINITION 4.1.3:** *Let L be a biloop. $L = L_1 \cup L_2$. The number of distinct elements in L will be called as the order of the biloop denoted by o(L) or |L|.*

**DEFINITION 4.1.4:** *A biloop $L = L_1 \cup L_2$ is said to be a Moufang biloop if both the loops $(L_1, +)$ and $(L_2, \bullet)$ are Moufang loops i.e. all elements of $L_1$ and $L_2$ satisfy the one of the following identities:*

  i.   *(xy) (zx) = (x (yz))x.*
  ii.  *((xy) z) y = x(y(zy)).*
  iii. *x (y (xz)) = ((xy) x)z.*

**DEFINITION 4.1.5:** *A biloop $L = L_1 \cup L_2$. is called a Bruck biloop if both the loops $(L_1,+)$ and $(L_2, \bullet)$ are Bruck loops i.e. both of them satisfy the identities (x(yz))z = x(y(xz)) and $(xy)^{-1} = x^{-1}y^{-1}$ for all x, y $\in L_1$ and x, y $\in L_2$. Similarly we can define a biloop to be a Bol biloop, WIP biloop, right (or left) or alternative biloops.*

**DEFINITION 4.1.6:** *Let $(L, +, \bullet)$ be a biloop. A non-empty subset P of L is said to be a sub-biloop of L if P is itself a biloop under the operations of L. i.e. $P = P_1 \cup P_2$ and $(P_1,+)$ and $(P_2, \bullet)$ are loops.*

**THEOREM 4.1.1:** *Let L be a biloop. A non empty subset P of L is a sub-biloop of L if and only if*

  i.  *$P_1 = P \cap L_1$ is a subloop of $L_1$ where $L = L_1 \cup L_2$.*
  ii. *$P_2 = P \cap L_2$ is a subloop of $L_2$.*

*Proof:* Left as an exercise for the reader to prove.

**DEFINITION 4.1.7:** *Let $(L, +, \bullet)$ be a biloop where $L = L_1 \cup L_2$. A proper subset $P$ of $L$ is a normal sub-biloop of $L$ if*

  i.  *$P$ is a sub-biloop of $L$.*
  ii. *$xP_1 = P_1x$ for all $x, \in L_1$ (where $P_1 = P \cap L_1$ is a subloop of $L_1$).*
  iii. *$xP_2 = P_2x$ for all $x, \in L_2$ (where $P_2 = P \cap L_2$ is a subloop of $L_2$).*
  iv. *$y(xP_1) = (yx)P_1$ for all $x \in L_1$ and $y(xP_2) = (yx)P_2$ for all $x \in L_2$.*

*We call a biloop $(L, +, \bullet)$ to be a simple biloop if $L$ has no proper normal sub-biloops.*

**DEFINITION 4.1.8:** *Let $(L, +, \bullet)$ be a biloop; we say $C(L)$ is the Moufang center of this biloop where $C(L) = C(L_1) \cup C(L_2)$ where $C(L_i) = \{x \in L \mid xy = yx \text{ for all } y \in L_i\}$; $i = 1, 2$ where $L_1$ and $L_2$ are loops such that $L = L_1 \cup L_2$.*

**DEFINITION 4.1.9:** *Let $L$ and $L'$ be two biloops. We say a map $\theta$ from $L$ to $L'$ is biloop homomorphism if $\theta = \theta_1 \cup \theta_2$, '$\cup$' is just a symbol and $\theta_1$ is a loop homomorphism from $L_1$ to $L'_1$ and $\theta_2$ is a loop homomorphism from $L_2$ and $L'_2$ where $L = L_1 \cup L_2$ and $L' = L'_1 \cup L'_2$. i.e. $\theta_1(xy) = \theta_1(x)\theta_1(y)$ for all $x, y \in L_1$ and $\theta_2(xy) = \theta_2(x)\theta_2(y)$ for all $x, y \in L_2$.*

**THEOREM (LAGRANGE THEOREM FOR BILOOPS):** *Let $L$ be a finite biloop. If $H$ is a proper sub-biloop of $L$ then $o(H)$ in general is not a divisor of $o(L)$.*

*Proof:* By an example consider the biloop given in example 4.1.1. $|L| = 13$ but $L$ has sub-biloops $H$ where $H = \{e, a_i, 1\} = H_1 \cup H_2$ with $H_2 = \{1\}$ and $H_1 = \{e, a_1\}$ and $o(H) = 3$ but $3 \nmid 13$. Hence the claim.

**THEOREM (CAUCHY THEOREM FOR BILOOPS):** *Let $L$ be a biloop of finite order and $p \mid o(L)$ where $p$ is a prime number. Then in general the biloop $L$ does not contain an element $a \neq e \in L$ where $a^p = e$, $e$ is the identity in $L_1$ or $L_2$ according as $a \in L_1$ or $L_2$.*

*Proof:* Consider the biloop $L = L_1 \cup L_2$. $L_1$ is given by the following table and $L_2 = \{g \mid g^9 = 1\}$ the cyclic group of order 9.

Table for $L_1$

$L_1 = \{e, a_1, a_2, a_3, a_4, a_5\}$.

| *   | e   | $a_1$ | $a_2$ | $a_3$ | $a_4$ | $a_5$ |
|-----|-----|-----|-----|-----|-----|-----|
| e   | e   | $a_1$ | $a_2$ | $a_3$ | $a_4$ | $a_5$ |
| $a_1$ | $a_1$ | e   | $a_5$ | $a_4$ | $a_3$ | $a_2$ |
| $a_2$ | $a_2$ | $a_3$ | e   | $a_1$ | $a_5$ | $a_4$ |
| $a_3$ | $a_3$ | $a_5$ | $a_4$ | e   | $a_2$ | $a_1$ |
| $a_4$ | $a_4$ | $a_2$ | $a_1$ | $a_5$ | e   | $a_3$ |
| $a_5$ | $a_5$ | $a_4$ | $a_3$ | $a_2$ | $a_1$ | e   |

Clearly |L| = 15. Now p = 5 is prime and p/15 but L has no element of order 5. Hence the claim.

So now we define two new types of biloops.

**DEFINITION 4.1.10:** *Let $L = L_1 \cup L_2$ be a finite biloop. If the order of every sub-biloop of the biloop divides |L| then we call L a Lagrange biloop.*

**DEFINITION 4.1.11:** *Let $L = L_1 \cup L_2$ be a biloop of finite order. If for every prime p, p/|L| we have elements $a \neq e \in L$ such that $a^p = e$ then we call L a Cauchy biloop.*

*Example 4.1.3:* Let $L = L_1 \cup L_2$ be a biloop. $L_1$ is a loop of order 6 given by the following table and $L_2 = \langle g \,/\, g^2 = 1 \rangle$ i.e. $L_2$ is a cyclic group of order 2.

Table for $L_1$

| • | e | $x_1$ | $x_2$ | $x_3$ | $x_4$ | $x_5$ |
|---|---|---|---|---|---|---|
| e | e | $x_1$ | $x_2$ | $x_3$ | $x_4$ | $x_5$ |
| $x_1$ | $x_1$ | e | $x_3$ | $x_5$ | $x_2$ | $x_4$ |
| $x_2$ | $x_2$ | $x_5$ | e | $x_4$ | $x_1$ | $x_3$ |
| $x_3$ | $x_3$ | $x_4$ | $x_1$ | e | $x_5$ | $x_2$ |
| $x_4$ | $x_4$ | $x_3$ | $x_5$ | $x_2$ | e | $x_1$ |
| $x_5$ | $x_5$ | $x_2$ | $x_4$ | $x_1$ | $x_3$ | e |

o(L) = 8. 2 is the only prime that divides 8 and every element in L is of order 2. Hence the claim.

Thus the class of Cauchy biloops is non-trivial. The next natural question is about Sylow theorems.

**THEOREM (SYLOW THEOREM FOR BILOOPS):** *Let $L = L_1 \cup L_2$ be a biloop of finite order. If p is a prime such that $p^\alpha / o(L)$ but $p^{\alpha+1} \nmid o(L)$ then the biloop in general does not have a sub-biloop of order $p^\alpha$ (i.e. L does not have a p-Sylow sub-biloop).*

*Proof*: This is proved by the following example. Let $L = L_1 \cup L_2$ where $L_1 = \{g \,/\, g^7 = 1\}$ and $L_2 = \{e, a_1, a_2, a_3, \ldots, a_7\}$, |L| = 15. $L_1$ is a cyclic group of order 7 and $L_2$ is given by the following table:

| • | e | $a_1$ | $a_2$ | $a_3$ | $a_4$ | $a_5$ | $a_6$ | $a_7$ |
|---|---|---|---|---|---|---|---|---|
| e | e | $a_1$ | $a_2$ | $a_3$ | $a_4$ | $a_5$ | $a_6$ | $a_7$ |
| $a_1$ | $a_1$ | e | $a_5$ | $a_2$ | $a_6$ | $a_3$ | $a_7$ | $a_4$ |
| $a_2$ | $a_2$ | $a_5$ | e | $a_6$ | $a_3$ | $a_7$ | $a_4$ | $a_1$ |
| $a_3$ | $a_3$ | $a_2$ | $a_6$ | e | $a_7$ | $a_4$ | $a_1$ | $a_5$ |
| $a_4$ | $a_4$ | $a_6$ | $a_3$ | $a_7$ | e | $a_1$ | $a_5$ | $a_2$ |
| $a_5$ | $a_5$ | $a_3$ | $a_7$ | $a_4$ | $a_1$ | e | $a_2$ | $a_6$ |
| $a_6$ | $a_6$ | $a_7$ | $a_4$ | $a_1$ | $a_5$ | $a_2$ | e | $a_3$ |
| $a_7$ | $a_7$ | $a_4$ | $a_1$ | $a_5$ | $a_2$ | $a_6$ | $a_3$ | e |

5/15 but L has no 5-Sylow biloops. Also 3/15 and L has 3-Sylow biloops given by H = {e, a$_3$} ∪ {1}. Now we define p-Sylow biloops.

**DEFINITION 4.1.12:** *Let $L = L_1 \cup L_2$ be a biloop of finite order. If p is a prime such that $p^\alpha/|L|$ but $p^{\alpha+1} \nmid |L|$ and if L has a sub-biloop of order $p^\alpha$ but not of order $p^{\alpha+1}$. We call the sub-biloop of order $p^\alpha$ as p-Sylow sub-biloops.*

*If L is a biloop of finite order; L has for every prime p, p/|L| with $p^\alpha/|L|$, $p^{\alpha+1} \nmid |L|$ and if L has p-Sylow subloops then we call L a Sylow biloop.*

Thus we see biloops do not in general obey all classical theorems for groups.

Using the new class of loops $L_n(m) \in L_n$; $n > 3$, n odd and $(m, n) = 1$ and $(m - 1, n) = 1$, we get a nice class of biloops. Say $L_k(m) \cup L_t(p) \in L_k \cup L_t$, k and t, two distinct odd integers greater than 3. These new class of biloops will also be only of even order.

Unlike in loops we can using these loops from $L_n$ get both odd and prime order biloops. For example $L = L_{11}(7) \cup \{g / g^7 = 1\}$ gives a biloop of order $(12 + 7) = 19$ a prime order biloop.

$L = L_{13}(3) \cup \{g / g^7 = 1\}$, we get an odd order biloop of order $14 + 7 = 21$ as we know every loop in the class of loops $L_n$ is of order $(n + 1)$ as n is odd order of all the loops in $L_n$ are even. Thus we can get biloops of any desired order even, odd or prime using the appropriate loops from the class of loops $L_n$ and any group of finite order.
We define nuclei, commutator, associator etc. of a biloop and finally the concept of representation of biloops and isotopes of biloops.

**DEFINITION 4.1.13:** *Let $L = L_1 \cup L_2$ be a biloop. If x and y $\in L_1 \cup L_2$ are elements of $L_1$ or $L_2$ the commutator (x, y) is defined as $xy = (yx)(x, y)$ [It is very important to note that 'or' is used in the sense that x, y $\in L_1$ 'or' x, y $\in L_2$ only and not in any other manner].*

*The commutator sub-biloop of L denoted by BL' is the {subloop generated by commutators of $L_1$} ∪ {subloop generated by the commutator in $L_2$} i.e. $L'_1 \cup L'_2$. Clearly BL' is a sub-biloop of L.*

**DEFINITION 4.1.14:** *Let $L = L_1 \cup L_2$ be a biloop. If x, y and z are elements of a loop $L_1$ or $L_2$ an associator (x, y, z) is defined by $(xy)z = (x(yz))(x, y, z)$ i.e. the triple is in $L_1$ or in $L_2$, the 'or' used as in the case of commutator.*

*The associator sub-biloop of $L = L_1 \cup L_2$ denoted by B(AL) is the {subloop of $L_1$ generated by all its associators} ∪ {subloop of $L_2$ generated by all of its associators} $= A(L_1) \cup A(L_2)$ i.e. B(AL) is a sub-biloop of L.*

**DEFINITION 4.1.15:** *Let $L = L_1 \cup L_2$ be a biloop. The left nucleus of the biloop $BN_\lambda = N_{\lambda_1} \cup N_{\lambda_2}$ where $N_{\lambda_1} = \{a \in L_1 /(a, x, y) = e$ for all $x, y \in L_1\}$ is a subloop of $L_1$. $N_{\lambda_2} = \{a \in L_2 / (a, x, y) = e$ for all $x, y \in L_2\}$ is a subloop of $L_2$. So $BN_\lambda = N_{\lambda_1} \cup$* $N_{\lambda_2}$ is a sub-biloop of L. The middle nucleus of the biloop $BN_\mu = N_{\mu_1} \cup N_{\mu_2}$ where $N_{\mu_1} = \{a \in L_1 / (x, a, y) = e$ for all $x, y \in L_1\}$ is a subloop of $L_1$. $N_{\mu_2} = \{a \in L_2 / (x, a, y) = e$ for all $x, y \in L_2\}$ is a subloop of $L_2$. Thus $BN_\mu = N_{\mu_1} \cup N_{\mu_2}$ is a sub-biloop of L.

The right nucleus of a biloop $L = L_1 \cup L_2$ is defined as $BN_\rho = N_{\rho_1} \cup N_{\rho_2}$ where $N_{\rho_1} = \{a \in L_1 / (x, y, a) = e$ for all $x, y \in L_1\}$ is a subloop of $L_1$. $N_{\rho_2} = \{a \in L_2 / (x, y, a) = e$ for all $x, y \in L_2\}$ is a subloop of $L_2$. Thus $BN_\rho = N_{\rho_1} \cup N_{\rho_2}$ is a sub-biloop of $L = L_1 \cup L_2$. The nucleus of the biloop L is $BN(L) = B(N(L_1 \cup L_2)) = BN_\lambda \cap BN_\mu \cap BN_\rho = (N_{\lambda_1} \cap N_{\mu_1} \cap N_{\rho_1}) \cup (N_{\lambda_2} \cap N_{\mu_2} \cap N_{\rho_2}) = N(L_1) \cup N(L_2)$. Clearly $BN(L)$ is a sub-biloop of $L = L_1 \cup L_2$.

**DEFINITION 4.1.16:** *Let $L = L_1 \cup L_2$ be a biloop. The Moufang center $BC(L) = BC(L_1 \cup L_2)$ is {the set of elements of the loop $L_1$ which commute with every element of $L_1$} $\cup$ {the set of elements of the loop $L_2$ which commute with every element of $L_2$} $= C(L_1) \cup C(L_2)$, $BC(L) = BC(L_1 \cup L_2) = C(L_1) \cup C(L_2)$. The center $BZ(L)$ of the biloop L is $B(Z(L)) = (C(L_1) \cap N(L_1)) \cup (C(L_2) \cap N(L_2))$.*

Now we proceed on to define the representation of a finite biloop $L = L_1 \cup L_2$.

**DEFINITION 4.1.17:** *Let $(L, \bullet, +)$ be a biloop i.e. $L = L_1 \cup L_2$ and $(L_1, \bullet)$ and $(L_2, +)$ are finite loops. For $\alpha_1 \in L_1$ define $R_{\alpha_1}$ as a permutation of the loop $(L_1, \bullet)$. $R_{\alpha_1} : x \to x \bullet \alpha_1$. For $\alpha_2 \in L_2$ define $R_{\alpha_2}$ as a permutation of the loop $(L_2, +)$ as follows: $R_{\alpha_2} : x \to x + \alpha_2$. We will call the set $(R_{\alpha_1} \cup R_{\alpha_2} | \alpha_1 \in L_1, \alpha_2 \in L_2\}$ the right regular representation of $(L, +, \bullet)$ i.e. $(L_1, \bullet)$ and $(L_2, +)$ or briefly the representation of the biloop L.*

**DEFINITION 4.1.18:** *Let $(L, +, \bullet)$ be a biloop, i.e. $L = L_1 \cup L_2$. For any pre-determined pairs $a_1, b_1 \in L_1$ and $a_2, b_2 \in L_2$ in the biloop L, a principal isotope $(L, \oplus', \odot')$ of the biloop $(L, +, \bullet)$ is defined by $x_1 \oplus' y_1 = X_1 + Y_1$ where $X_1 + a_1 = x_1$ and $b_1 + Y_1 = y_1$ and $x_2 \odot' y_2 = X_2 + Y_2$ where $X_2 \bullet a_2 = x_2$ and $b_2 \bullet Y_2 = y_2$. $(L, +, \bullet)$ is said to be a G-biloop if it is isomorphic to all of its principal isotopes. We illustrate this by the following example.*

**Example 4.1.4:** Let $L = L_5(2) \cup L_5(4)$. Clearly L is a biloop. The tables of $(L_5(2), \bullet)$ and its principal isotope $(L_5(2), *)$ is given below:

| • | e | 1 | 2 | 3 | 4 | 5 |
|---|---|---|---|---|---|---|
| e | e | 1 | 2 | 3 | 4 | 5 |
| 1 | 1 | e | 3 | 5 | 2 | 4 |
| 2 | 2 | 5 | e | 4 | 1 | 3 |
| 3 | 3 | 4 | 1 | e | 5 | 2 |
| 4 | 4 | 3 | 5 | 2 | e | 1 |
| 5 | 5 | 2 | 4 | 1 | 3 | e |

Composition table of $(L_5(2), *)$

| * | e | 1 | 2 | 3 | 4 | 5 |
|---|---|---|---|---|---|---|
| e | 3 | 2 | 5 | e | 1 | 4 |
| 1 | 5 | 3 | 4 | 1 | e | 2 |
| 2 | 4 | e | 3 | 2 | 5 | 1 |
| 3 | e | 1 | 2 | 3 | 4 | 5 |
| 4 | 2 | 5 | 1 | 4 | 3 | e |
| 5 | 1 | 4 | e | 5 | 2 | 3 |

Composition table of principal isotope of $(L_5(2), •)$ is $(L_5(2), *)$. To obtain $(L_5(2), *)$ we have taken $a_1 = e$ and $b_1 = 3$. Now for the principal isotope we see the identity is 3. $x * x = 3$ for all $x \in (L_5(2), *)$.

Table for $(L_5(4), •)$ is

| • | e | 1 | 2 | 3 | 4 | 5 |
|---|---|---|---|---|---|---|
| e | e | 1 | 2 | 3 | 4 | 5 |
| 1 | 1 | e | 5 | 4 | 3 | 2 |
| 2 | 2 | 3 | e | 1 | 5 | 4 |
| 3 | 3 | 5 | 4 | e | 2 | 1 |
| 4 | 4 | 2 | 1 | 5 | e | 3 |
| 5 | 5 | 4 | 3 | 2 | 1 | e |

Composition table for $(L_5(4), *)$ is

| * | e | 1 | 2 | 3 | 4 | 5 |
|---|---|---|---|---|---|---|
| e | 1 | e | 5 | 4 | 3 | 2 |
| 1 | e | 1 | 2 | 3 | 4 | 5 |
| 2 | 4 | 2 | 1 | 5 | e | 3 |
| 3 | 2 | 3 | e | 1 | 5 | 4 |
| 4 | 5 | 4 | 3 | 2 | 1 | e |
| 5 | 3 | 5 | 4 | e | 2 | 1 |

Here we have taken $a_2 = 1$ and $b_2 = e$ the identity is 1. Thus $(L, *, *')$ where $L = (L_5(2), *) \cup (L_5(4), *')$ is the principal isotopes of the biloop. Clearly the principal isotope is also a biloop.

Several results enjoyed by the loops in this regard can be easily obtained for the case of biloops with suitable modifications.

Now we proceed on to define and find different representations of the biloops of order 2n where $L = L_1 \cup L_2$ where each $L_i$ is an even order loop should we paste or have common paths between the two $K_{2n}$ graphs. This question is an interesting one and is left for further analysis to the reader. In case of loops we have proved refer [109].

## 4.2 Smarandache biloops and its properties

In this section we introduce the notion of Smarandache biloops. Study of biloops is very recent and the concept of Smarandache biloops is introduced for the first time in this text. We enumerate several of its interesting properties enjoyed by them.

**DEFINITION 4.2.1:** *Let $(L, +, \bullet)$ be a biloop we call L a Smarandache biloop (S-biloop) if L has a proper subset P which is a bigroup.*

Now by this new structure, we study in a non-assocaitive structure an associative sub-bistructure.

*Example 4.2.1:* Let $(L, +, \bullet)$ be a biloop given by the following tables where $L = L_1 \cup L_2$ with $L_1 = \{g \,/\, g^9 = 1\}$ the cyclic group of degree 9 and $L_2$ is the loop given by the following table:

| $\bullet$ | e | $a_1$ | $a_2$ | $a_3$ | $a_4$ | $a_5$ |
|---|---|---|---|---|---|---|
| e | e | $a_1$ | $a_2$ | $a_3$ | $a_4$ | $a_5$ |
| $a_1$ | $a_1$ | e | $a_3$ | $a_5$ | $a_2$ | $a_4$ |
| $a_2$ | $a_2$ | $a_5$ | e | $a_4$ | $a_1$ | $a_3$ |
| $a_3$ | $a_3$ | $a_4$ | $a_1$ | e | $a_5$ | $a_2$ |
| $a_4$ | $a_4$ | $a_3$ | $a_5$ | $a_2$ | e | $a_1$ |
| $a_5$ | $a_5$ | $a_2$ | $a_4$ | $a_1$ | $a_3$ | e |

Take $P = \{1, g^3, g^6\} \cup \{e, a_5\}$. Clearly P is a bigroup. Hence L is a S-biloop.

Now we see that all biloops in general need not be S-biloops. Let $(L, +, \bullet)$ be a S-biloop, the order of L is denoted by $o(L)$ or $|L|$ is the number of elements in L; if L has finite number of elements we call L a S-biloop of finite order if L has infinite number of elements we call L a S-biloop of infinite order.

**DEFINITION 4.2.2:** *Let $(L, +, \bullet)$ be a biloop. If $P \subset L$ is a proper subset of L which is a sub-biloop of L, is said to be a S-sub-biloop of L if P itself is a S-biloop under the operations of L.*

*Example 4.2.2:* Let $L = L_1 \cup L_2$ where $L_1 = L_5(3) \times L_9(8)$ and $L_2 = S_6$. Clearly L is a S-biloop. Take $P = P_1 \cup P_2$ where $P_1 = \{e\} \times L_9(8)$ and $P_2 = A_6$. Clearly P is a sub-biloop of L which is a S-sub-biloop of L.

**THEOREM 4.2.1:** *Let $(L, +, \bullet)$ be a biloop which has a S-sub-biloop then L is a S-biloop.*

*Proof*: Straightforward by the very definitions; hence left for the reader to prove.

**DEFINITION 4.2.3:** *Let $(L, +, \bullet)$ be a biloop. L is said to be a Smarandache commutative biloop (S-commutative biloop) if every proper subset P of L which are bigroups, is a commutative bigroup of L. If $(L, +, \bullet)$ has only one proper subset P which is a bigroup, is commutative then we call L a Smarandache weakly commutative biloop (S-weakly commutative biloop).*

**DEFINITION 4.2.4:** *Let $(L, +, \bullet)$ be a biloop. We say $(L, +, \bullet)$ is a Smarandache cyclic biloop (S-cyclic biloop) if every proper subset which is a bigroup is a cyclic bigroup. If $(L, +, \bullet)$ has at least one proper subset which is a bigroup to be a cyclic bigroup then we call L a Smarandache weakly cyclic biloop (S-weakly cyclic biloop).*

**DEFINITION 4.2.5:** *Let $(L, +, \bullet)$ be a biloop. A proper S-sub-biloop, $P = P_1 \cup P_2$ of L is said to be a Smarandache normal biloop (S-normal biloop) if*

i. $xP = Px$.
ii. $(Px)y = P(xy)$.
iii. $y(xP) = (yx)P$

*for all $x, y \in L$.*

*If the biloop $(L, +, \bullet)$ has no proper S-normal sub-biloop then we call the biloop to be Smarandache simple (S-simple). Smarandache simple biloop and simple biloop may or may not in general have proper relations.*

**DEFINITION 4.2.6:** *Let $(L, +, \bullet)$ be a biloop and $P = P_1 \cup P_2$ be a S-sub-biloop of L; $(P \subset L)$ for the pair of elements $x, y$ in P, the commutator $(x, y)$ is defined by $xy = (yx)(x, y)$. The Smarandache commutator subloop (S-commutator subloop) of L denoted by $S(L^s)$ is the S-sub-biloop generated by all its commutators i.e. $\langle\{x \in P / x = (y, z)$ for some $y, z \in L\}\rangle$. Clearly if this set does not generate a S-sub-biloop we then say the S-commutator subloop generated relative to the S-sub-biloop P of L is empty. Thus it is important and interesting to note that in case of S-commutators in biloops we see there can be more than one S-commutator S-sub-biloop which is a marked difference between the commutator in biloops and S-commutator in biloops.*

**DEFINITION 4.2.7:** *Let $(L, +, \bullet)$ be a biloop. Let $P \subset L$ be a S-sub-biloop of L. If x, y and z are elements of the S-sub-biloop P of L; an associator $(x, y, z)$ is defined by $(xy)z = (x(yz))(x, y, z)$. The associator S-sub-biloop of the S-sub-biloop P of L denoted by $A(L^s)$ is the S-sub-biloop generated by all the associators, that is, $\langle\{x \in P / x = (a, b, c)$ for some $a, b, c \in P\}\rangle$. If $A(L^s)$ happens to be a S-sub-biloop then only we call $A(L^s)$ the Smarandache associator (S-associator) of L related to the S-sub-biloop P. Here also as in case of S-commutators we see there can be several S-associators unlike a unique associator in case of biloops.*

**DEFINITION 4.2.8:** *Let $(L, +, \bullet)$ be a biloop with $L = L_1 \cup L_2$. If for $a, b \in L$ with $ab = ba$ in L if we have $(ax)b = (bx)a$ (or $b(xa)$) for all $x \in L_1$ if $a, b \in L_1$ ($x \in L_2$ if $a, b \in L_2$, then we say the pair is pseudo commutative. If every commutative pair is pseudo commutative then we call L a pseudo commutative biloop.*

*If we have for every pair of commuting elements $a, b$ in a proper subset P of the biloop L (which is a S-sub-biloop) $a, b$ is a pseudo commutative pair for every $x \in P$, then we call the pair a Smarandache pseudo commutative pair (S-pseudo commutative pair).*

*(Thus we see if a pair is a pseudo commutative pair then it is obviously a S-pseudo commutative pair. But we see on the contrary that all S-pseudo commutative pairs in general need not be a pseudo commutative pair).*

*If in a S-sub-biloop P of L every pair in P is a S-pseudo commutative pair then we call the biloop L is a S-pseudo commutative biloop.*

**DEFINITION 4.2.9:** *Let (L, +, •) be a biloop. The pseudo commutator of L denoted by $P(L) = \langle\{p \in L / a(xb) = p[bx)a]\}\rangle$ we define the Smarandache pseudo commutator (S-pseudo commutator) of $P \subset L$ (P a S-sub-biloop of L) to be the S-sub-biloop generated by $\langle\{ p \in P / a(xb) = p([bx]a), a, b \in P\}\rangle$ denoted by $P_S(L))$. If $(P_S(L))$ is not a S-sub-biloop then we say $P_S(L) = \phi$. (Thus for a given biloop (L, +, •) we may have more than one $P_S(L)$).*

Now we proceed on to define pseudo associative biloop and Smarandache pseudo associative biloop.

**DEFINITION 4.2.10:** *Let (L, +, •) be a biloop. An associative triple a, b, c $\in P \subset L$ where P is a S-sub-biloop of L is said to be Smarandache pseudo associative (S-pseudo associative) if (ab)(xc) = (ax)(bc) for all x $\in$ P. If (ab)(xc) = (ax)(bc) for some x $\in$ P we say the triple is Smarandache pseudo associative (S-pseudo associative) relative to those x in P. If in particular every associative triple in P is S-pseudo associative then we say the biloop is a Smarandache pseudo associative biloop (S-pseudo associative biloop). Thus for a biloop to be a S-pseudo associative biloop it is sufficient that one of its S-sub-biloops are a S-pseudo associative biloop.*

*The Smarandache pseudo associator of a biloop (L, +, •) denoted by $S(AL_p) = \langle\{t \in P / (ab)(tc) = (at)(bc)$ where $a(bc) = (ab)c$ for $a, b, c \in P\}\rangle$; where P is a S-sub-biloop of L.*

**DEFINITION 4.2.11:** *A biloop (L, +, •) is said to be a Smarandache Moufang biloop (S-Moufang biloop) if there exists S-sub-biloop P of L which is a Moufang biloop.*

*A biloop L is said to be a Smarandache Bruck biloop (S-Bruck biloop) if L has a proper S-sub-biloop P where P is a Bruck biloop. Similarly a biloop L is said to be a Smarandache Bol biloop (S-Bol biloop) if L has a proper S-sub-biloop P where P is a Bol biloop. We call a biloop L to be a Smarandache Jordan biloop (S-Jordan biloop) if L has a S-sub-biloop P where P is a Jordan biloop.*

In a similar way we define Smarandache right (left) alternative biloop.

**DEFINITION 4.2.12:** *Let (L, +, •) and (L', +, •) be two biloops. We say a map $\phi$ from L to L' is a Smarandache biloop homomorphism (S-biloop homomorphism) if $\phi$ is a bigroup homomorphism from P to P' where P and P' are proper subsets of L and L' respectively such that they are bigroups of L and L'. Thus for us we need not even have $\phi$ to be defined on the whole of L. It is sufficient if $\phi$ is defined on a proper subset P of L which is a bigroup.*

Now we proceed on to define a Smarandache Moufang center.

**DEFINITION 4.2.13:** *Let L be a biloop. Let P be a S-sub-biloop of L. The Smarandache Moufang center (S-Moufang center) of L is defined to be the Moufang center of the S-sub-biloop P = $P_1 \cup P_2$. Thus for a given biloop L we can have several Moufang centres or no Moufang centre if the biloop L has no S-sub-biloops.*

**DEFINITION 4.2.14:** *Let L be a biloop. P be a S-sub-biloop of L. The Smarandache center (S-center) SZ(L) is the center of the biloop P $\subset$ L. Thus even in case of S-centers for a biloop we may have several S-centers depending on the number of S-sub-biloops.*

Now we define the Smarandache middle, left and right nucleus of a biloop.

**DEFINITION 4.2.15:** *Let L be a biloop. P be a S-sub-biloop of L. To define Smarandache nucleus related to P we consider the Smarandache left nucleus (S-left nucleus) of L to be the left nucleus of the S-sub-biloop P denoted by $N_\lambda^P$. Similarly the Smarandache right nucleus (S-right nucleus) of L to the right nucleus of the S-sub-biloop P denoted by $N_P^P$ and the Smarandache middle nucleus (S-middle nucleus) of the S-sub-biloop P.*

*Thus the Smarandache nucleus (S-nucleus) is SN(L) = $N_\mu^P \cap N_\lambda^P \cap N_P^P$.*

Now we proceed on to define Smarandache Lagrange biloop, Smarandache Cauchy element in a biloop and Smarandache p-Sylow sub-biloop.

**DEFINITION 4.2.16:** *Let (L, +, •) be a biloop of finite order. A proper subset where S is a S-sub-biloop of L is said to be a Smarandache Lagrange sub-biloop (S-Lagrange sub-biloop) of L if o(S)/o(L). If every S-sub-biloop of the biloop L happens to be a S-Lagrange sub-biloop then we call L a Smarandache Lagrange biloop (S-Lagrange biloop). If no S-sub-biloop of L is a S-Lagrange sub-biloop then we call L a Smarandache non-Lagrange biloop (S-non-Lagrange biloop).*

*Example 4.2.3:* Let L = $L_9(2) \cup L_5(3)$. Clearly L is a S-biloop. But L has no S-sub-biloops say P = $L_9(2) \cup \{e, a\}$. Clearly |L| = 16 whereas |P| = 12, 12 $\neq$ 16. Thus P is a S-sub-biloop, which is not a S-Lagrange sub-biloop.

**DEFINITION 4.2.17:** *Let (L, +, •) be a biloop of finite order. An element x $\in$ P, where P is a S-sub-biloop of L is said to be a Smarandache Cauchy element (S-Cauchy element) of the biloop L (S-Cauchy element of the biloop L) if $x^n$ = e, n > 1 and n / |P|. If every S-sub-biloop of L has a S-Cauchy element then we call the biloop L a Smarandache Cauchy biloop (S-Cauchy biloop).*

**DEFINITION 4.2.18:** *Let (L, +, •) be a biloop of finite order. If p is a prime such that $p^\alpha$/ |L| but $p^{\alpha+1} \neq$ |L| and if L has a S-sub-biloop P of order $p^\alpha$ then we call P the Smarandache p-Sylow sub-biloop (S-p-Sylow sub-biloop). If for every prime p, $p^\alpha$/|L|, $p^{\alpha+1} \neq$ |L| we have S-p-Sylow sub-biloop then we call L a Smarandache p-Sylow biloop (S-p-Sylow biloop).*

**DEFINITION 4.2.19:** *Let L be a biloop we call L a Smarandache pseudo biloop (S-pseudo biloop) if L a proper subset P where P is a bigroupoid.*

All properties about S-pseudo biloops can be defined as in case of S-biloops.

**DEFINITION 4.2.20:** *Let (L, +, •) and (L', +, •) be two biloops. A map $\phi: L \to L'$ is called a Smarandache biloop homomorphism (S-biloop homomorphism) if there exists a bigroup homomorphism from P to P' where P and P' are proper subsets of L and L' respectively which are groups.*

**DEFINITION 4.2.21:** *The notion of Smarandache biloop isomorphism (S-biloop isomorphism), automorphism (S-biloop automorphism) are defined as in case of other algebraic structures.*

**PROBLEMS:**

1. Give an example of a biloop, which is not a S-biloop.
2. Find the order of the smallest S-biloop.
3. Is the biloop $L = L_5(2) \cup S(7)$ a S-biloop? Justify your claim.
4. Find the commutator and S-commutator biloop of L given in problem 3.
5. Give an example of a biloop, which has a non-trivial S-associator sub-biloop.
6. Find the Moufang centre of the biloop (L, +, •) where $L = S_9 \cup L_{13}(5)$.
7. Does the biloop (L, +, •) where $L = S_3 \cup L_{15}(8)$ have S-Moufang centre?
8. Find how many S-Moufang centres does the biloop (L, +, •) where $L = G \cup L_{21}(17)$, $G = (g = \langle g \,/\, g^{24} = 1 \rangle$ have?
9. Give an example of a S-Cauchy biloop.
10. Does there exist a S-Lagrange biloop?
11. Is the biloop $L = L_9(5) \cup L_{11}(7)$ a S-Lagrange biloop?
12. Give an example of a biloop, which has no S-Cauchy element.
13. Do we have a biloop, which is a S-p-Sylow biloop?
14. Give an example of a S-pseudo biloop.
15. Let (L, +, •) be a biloop given by $L = L_{15}(8) \cup L_{13}(6)$. Find S-nucleus of L.

**Chapter 5**

# BIGROUPOIDS AND SMARANDACHE BIGROUPOIDS

This chapter has four sections. In the first section we define the notion of bigroupoids. Bigroupoids for the first time are defined and analyzed here. In section two we define Smarandache bigroupoids. Smarandache bigroupoids are very rare and they are defined. These structures find their applications in the study of automaton and semi-automaton which is discussed in section three. In section four we define the notion of direct product of S-automaton.

## 5.1 Bigroupoids and its Properties

In this section we define bigroupoids. Study of groupoids and S-groupoids started in 2002 [114]; as books solely on groupoids is very meager. Here we define bigroupoids introduce and analyze several of its interesting properties. Bigroupoids are also non-associative structure like biloops.

**DEFINITION 5.1.1:** *Let $(G, +, \bullet)$ be a non-empty set. We call G a bigroupoid if $G = G_1 \cup G_2$ and satisfies the following:*

  i. *$(G_1, +)$ is a groupoid (i.e. the operation + is non-associative).*
  ii. *$(G_2, \bullet)$ is a semigroup.*

*Example 5.1.1:* Let $(G, +, \bullet)$ be a bigroupoid, where $G = G_1 \cup G_2$ with $G_2 = S(3)$ and $G_1$ is a groupoid given by the following table:

| + | $x_1$ | $x_2$ | $x_3$ |
|---|---|---|---|
| $x_1$ | $x_1$ | $x_3$ | $x_2$ |
| $x_2$ | $x_2$ | $x_1$ | $x_3$ |
| $x_3$ | $x_3$ | $x_2$ | $x_1$ |

Clearly G is a bigroupoid as $(G_1, +)$ is a groupoid and $(G_2, \bullet)$ is a semigroup.

*Example 5.1.2:* Let $(G, +, \bullet)$ be a bigroupoid where $G = G_1 \cup G_2$ with $G_1 = S(4)$ and $G_2 = \{Z^+$ under the operation '$\bullet$', for any $x, y \in Z^+$ we define $x \bullet y = 7x + 5y\}$. Clearly G is a bigroupoid.

**DEFINITION 5.1.2:** *Let $(G, +, \bullet)$ be a bigroupoid. We say the order of the bigroupoid G is finite if the number of elements in G is finite. If the number of elements in G is infinite we say the groupoid G is an infinite groupoid. We denote the order of the groupoid G by $o(G)$ or $|G|$.*

The bigroupoid in example 5.1.1 is a finite bigroupoid whereas the bigroupoid given in example 5.1.2 is an infinite bigroupoid. A natural question that may arise in anyone's mind is that will the union of a groupoid and a semigroup in general give a bigroupoid. The answer to this question is no. For take $G = G_1 \cup G_2$ where $G_1 = (Z^+$ where the operation on $Z^+$ is for a, b $\in Z^+$ define a • b = pa + qb p, q $\in Z^+$ p ≠ q). Clearly $G_1$ is a groupoid, as the operation '•' is non-associative. $G_2 = \{Q^+$ semigroup under usually addition}. Thus $G_2$ is a semigroup.

Clearly $G = G_1 \cup G_2$ is not a bigroupoid as the set $G_1 \subset G_2$.

**DEFINITION 5.1.3:** *Let (G, +, •) be a bigroupoid with $G = G_1 \cup G_2$. A proper subset P of G said to be a sub-bigroupoid of G if $P = P_1 \cup P_2$ and $P_1$ is a subgroupoid of $(G_1, +)$ and $P_2$ is a subsemigroup of $(G_2, •)$.*

*Example 5.1.3:* Let $G = G_1 \cup G_2$; (G, +, •) a bigroupoid with $(G_1, +)$ a groupoid and $(G_2, •)$ a semigroup. Take $G_1 = Z_8(2, 6)$ and $G_2 = S(5)$. $G_1$ is the groupoid from the new class of groupoids, $G_2$ is the symmetric semigroup. The groupoid $Z_8$ (2, 6) is given by the following table:

| * | 0 | 1 | 2 | 3 | 4 | 5 | 6 | 7 |
|---|---|---|---|---|---|---|---|---|
| 0 | 0 | 6 | 4 | 2 | 0 | 6 | 4 | 2 |
| 1 | 2 | 0 | 6 | 4 | 2 | 0 | 6 | 4 |
| 2 | 4 | 2 | 0 | 6 | 4 | 2 | 0 | 6 |
| 3 | 6 | 4 | 2 | 0 | 6 | 4 | 2 | 0 |
| 4 | 0 | 6 | 4 | 2 | 0 | 6 | 4 | 2 |
| 5 | 2 | 0 | 6 | 4 | 2 | 0 | 6 | 4 |
| 6 | 4 | 2 | 0 | 6 | 4 | 2 | 0 | 6 |
| 7 | 6 | 4 | 2 | 0 | 6 | 4 | 2 | 0 |

Let $P = P_1 \cup P_2$ where $P_2 = S(3)$ and take $P_1 = \{0, 2, 4, 6\}$ to be the subgroupoid of $Z_8(2, 6)$. Thus P is a sub-bigroupoid.

**DEFINITION 5.1.4:** *Let (G, +, •) be a bigroupoid. $G = G_1 \cup G_2$. If $(G_1, +)$ and $(G_2, •)$ are idempotent groupoid and idempotent semigroup respectively, then (G, +, •) is an idempotent bigroupoid.*

**DEFINITION 5.1.5:** *Let (G, +, •) where $G = G_1 \cup G_2$ be a bigroupoid. A non-empty proper subset P of G is said to be a left bi-ideal of the bigroupoid, if $P = P_1 \cup P_2$ and $(P_1, +)$ is a left ideal of $G_1$ and $(P_2, •)$ is a left ideal of $G_2$.*

*On similar lines one can define right bi-ideal. $P = P_1 \cup P_2$, (P, +, •) is called as a bi-ideal of the bigroupoid (G, +, •) if P is simultaneously both a right bi-ideal and a left bi-ideal of G.*

**DEFINITION 5.1.6:** *Let $G = (G_1 \cup G_2, +, •)$ be a bigroupoid. A sub-bigroupoid (P = $P_1 \cup P_2$, +, •) of G is said to be a normal sub-bigroupoid of G if*

   i.   $aP_1 = P_1a$ or $aP_2 = P_2 a$ (according as a $\in P_1$ or $P_2$).

ii.   $P_1(xy) = (P_1x) y$ or $P_2 (xy) = (P_2x) y$ (according as $x, y \in P_1$ or $x, y \in P_2$).
iii.  $y(xP_1) = (xy) P_1$ or $y (xP_2) = (xy) P_2$ (according as $x, y \in P_1$ or $x, y \in P_2$).

for all $x, y, a \in G$.

*We call a bigroupoid simple if it has no non-trivial normal sub-bigroupoids. We define when bigroupoids are normal.*

**DEFINITION 5.1.7:** *Let $(G, +, \bullet)$ be a bigroupoid. $G = G_1 \cup G_2$ with $(G_1, +)$ is a groupoid and $(G_2, \bullet)$ is a semigroup. We say G is a normal bigroupoid if*

i.   $xG = Gx$ i.e. ($x_1 G_1 = G_1 x_1$ if $x_1 \in G$, and $x_2 G_2 = G_2 x_2$ if $x_2 \in G_2$).
ii.  $G(xy) = (Gx) y$ (i.e. $G_1 (xy) = (G_1 x) y$ if $x, y \in G_1$ and $(G_2 x) y = G_2(xy)$ if $x, y \in G_2$).
iii. $y (xG) = (yx) G$ (i.e. $(yx) G_1 = y (x G_1)$ if $x, y \in G_1$ and $(yx) G_2 = y (x G_2)$ if $x, y \in G_2$).

Several examples can be constructed with some innovative ideas.

**DEFINITION 5.1.8:** *Let $(G, +, \bullet)$ be a bigroupoid with $G = G_1 \cup G_2$ where $(G, +)$ is a groupoid and $(G_2, \bullet)$ is a semigroup. Let H and K be two proper sub-bigroupoids of $(G, +, \bullet)$ (where $H = H_1 \cup H_2$, $(H_1, +)$ is a subgroupoid of $(G_1, +)$ and $(H_2, \bullet)$ is a subsemigroup of $(G_2, \bullet)$. Similarly for $K = K_1 \cup K_2$ with $K \cap H = \phi$. We say H is bi conjugate with K if there exist $x_1 \in H_1$ and $x_2 \in H_2$ with $x_1 K_1 = H_1$ or $(K_1 x_1 = H_1)$ and $x_2 K_2 = H_2$ (or $K_2 x_2 = H_2$) or in the mutually exclusive sense.*

**DEFINITION 5.1.9:** *Let G be a bigroupoid. We say G is a inner commutative bigroupoid if every sub-bigroupoid of G is commutative.*

**DEFINITION 5.1.10:** *Let G be a bigroupoid. An element $x \in G$ is said to be a zero divisor if there exist $y \in G$ such that $xy = 0$ (we put xy as it can be the operation in $G_1$ or it can be the operation in $G_2$).*

**DEFINITION 5.1.11:** *Let $(G, +, \bullet)$ be a bigroupoid. The centre or bicentre of $(G, +, \bullet)$ denoted by $BC(G) = \{x \in G_1 \mid xa = ax$ for all $a \in G_1\} \cup \{ y \in G_2 \mid yb = by$ for all $b \in G_2\}$ i.e. why we call it as bicentre.*

Here it is pertinent to mention that by default of notation we put '+' and '•' it can be any binary operation not the usual addition '+' and the usual multiplication '•'.

**DEFINITION 5.1.12:** *Let $(G, +, \bullet)$ be a bigroupoid. $G = G_1 \cup G_2$, we say $a, b \in G$ is a bi conjugate pair if $a = b \bullet x_1$ (if $a, b \in G_1$ for some $x_1 \in G_1$ or $a = x_1 \bullet b$) and $b = a \bullet y_1$ (or $y_1 \bullet a$ for some $y_1 \in G_1$) (if $a, b \in G_2$ we choose x and y from $G_2$). On similar lines we define right bi-conjugate only for $a, b \in G_1$ (or $G_2$) as $a \bullet x = b$ and $b \bullet y = a$ for $x \in G_1$ ($x \in G_2$).*

We define left bi conjugate also in an analogous way.

Now we define a notion called biquasi loops analogous to biquasi groups where groups are replaced by loops and semigroups by groupoids.

**DEFINITION 5.1.13:** *Let $(Y, +, \bullet)$ be a set with $Y = Y_1 \cup Y_2$ where $(Y_1, +)$ is a loop and $(Y_2, \bullet)$ is a groupoid then we call $(Y, +, \bullet)$ as a biquasi loop.*

It is worthwhile to note all biquasi groups are biquasi loops as every group is trivially a loop and every semigroup is a groupoid.

*Example 5.1.4:* Let $(Y, +, \bullet)$ be a biquasi loop where $Y = Y_1 \cup Y_2$ with $Y_1 = L_7(3)$ the loop and $Y_2 = Z^+ \cup \{0\}$ under '~' where '~' denotes difference. Clearly Y is a biquasi loop.

**DEFINITION 5.1.14:** *Let $(Y, +, \bullet)$ be a biquasi loop. A proper subset $P \subset Y$ is said to be a sub-biquasi loop in P itself under the operations '+' and '$\bullet$' is a biquasi loop.*

**THEOREM 5.1.1:** *Every biquasi group is a biquasi loop and not conversely.*

*Proof*: Left as an exercise for the reader to prove.

*Example 5.1.5:* Let $(Y, +, \bullet)$ be a biquasi loop where $Y = Y_1 \cup Y_2$ where $(Y_1, +)$ is a loop and $(Y_2, \bullet)$ is a groupoid given by the following tables:

Table for $(Y_1, +)$

| + | e | a | b | c | d |
|---|---|---|---|---|---|
| e | e | a | b | c | d |
| a | a | e | c | d | b |
| b | b | d | a | e | c |
| c | c | b | d | a | e |
| d | d | c | e | b | a |

Table for $(Y_2, \bullet)$

| $\bullet$ | 0 | 1 | 2 | 3 |
|---|---|---|---|---|
| 0 | 0 | 2 | 0 | 2 |
| 1 | 3 | 1 | 3 | 1 |
| 2 | 2 | 0 | 2 | 0 |
| 3 | 1 | 3 | 1 | 3 |

Y is a biquasi loop of order 9.

**DEFINITION 5.1.15:** *Let $(X, +, \bullet)$ be a biquasi loop, $X = X_1 \cup X_2$; $(X_1, +)$ is a loop and $(X_2, \bullet)$ is a groupoid. $P = P_1 \cup P_2$ a sub-biquasi loop of X. We say P is a normal sub-biquasi loop of X if $P_1$ is a normal subloop of $(X_1, +)$ and $P_2$ is an ideal of the groupoid $(X_2, \bullet)$.*

**DEFINITION 5.1.16:** *Let (Y, +, •) be a biquasi loop of finite order. If the order of every proper sub-biquasi loop P of X divides the order of X, then we say X is a Lagrange biquasi loop.*

*If the order of at least one of the proper sub-biquasi loop P of X divides the order of X then we say X is a weakly Lagrange biquasi loop. If sub-biquasi loop exists and none of its order divides the order of X then X is said to be non-Lagrange biquasi loop.*

***Example 5.1.6:*** *Consider the biquasi loop (Y, +, •) with Y = $Y_1 \cup Y_2$ where $Y_1$ = $L_7(3)$ is a loop of order 8.*

$Y_2$ is the groupoid given by the following table:

| + | $x_1$ | $x_2$ | $x_3$ |
|---|---|---|---|
| $x_1$ | $x_1$ | $x_3$ | $x_2$ |
| $x_2$ | $x_2$ | $x_1$ | $x_3$ |
| $x_3$ | $x_3$ | $x_2$ | $x_1$ |

Clearly |Y| = 11. Y has no non-trivial sub-biquasi loops but it has sub-biquasi group. P = $P_1 \cup P_2$ where $P_1$ = {e, $a_6$} is a group and $P_2$ = {$x_1$} is a semigroup.

**DEFINITION 5.1.17:** *Let (Y, +, •) be a biquasi loop, of finite order. If p is a prime such that $p^\alpha$ / |Y| and $p^{\alpha+1} \neq$ |Y| we have a sub-biquasi loop P of order $p^\alpha$ then we say P is a p-Sylow sub-biquasi loop. If for every prime p / |Y| we have a p-Sylow sub-biquasi loop, then we call Y a p-Sylow bi-quasiloop.*

Now we proceed on to define the Cauchy elements in biquasi loops.

**DEFINITION 5.1.18:** *Let (Y, +, •) be a biquasi loop of finite order, an element x ∈ Y such that $x^t$ = e (e identity element of Y ) t >1 with t / |Y| is called the Cauchy element of Y. We see in case of groups every element is a Cauchy element but in case of biquasi loops only a few of the elements may be Cauchy elements; if the order of the biquasi loop is a prime, p none of the elements in it would be a Cauchy element even though we will have x ∈ Y with $x^t$ = e; t > 1, t ≠ p. as t ∤ p.*

We do not know whether there exists other biquasi loops in which no element is a Cauchy element. Several other interesting study can be carried out as in case of any other bistructure as a matter of routine.

**DEFINITION 5.1.19:** *Let (P, +, •) be a non-empty set . If (P, +) is a loop and (P, •) is a semigroup then we call P a biquasi semigroup.*

Like biquasi loops, biquasi semigroups are also non-associative structures. Only biquasi structure, which is associative, is biquasi groups.

**THEOREM 5.1.2:** *All biquasi semigroups are biquasi loops.*

*Proof:* Follows directly by the definitions.

The notions of sub-biquasi semigroups, ideals in biquasi semigroups, Lagrange biquasi semigroups, Cauchy elements in biquasi semigroups etc can be defined in an analogous way.

**DEFINITION 5.1.20:** *Let $(T, +, \bullet)$ be a non-empty set, T is a biquasi groupoid if $(T, +)$ is a group and $(T, \bullet)$ is a groupoid.*

*Example 5.1.7:* Let $(T, +, \bullet)$ be given by $T = T_1 \cup T_2$ where $T_1 = S_3$ and $T_2$ is a groupoid given by the following table:

| *     | $a_0$ | $a_1$ | $a_2$ | $a_3$ | $a_4$ | $a_5$ | $a_6$ |
|-------|-------|-------|-------|-------|-------|-------|-------|
| $a_0$ | $a_0$ | $a_4$ | $a_1$ | $a_5$ | $a_2$ | $a_6$ | $a_3$ |
| $a_1$ | $a_3$ | $a_0$ | $a_4$ | $a_1$ | $a_5$ | $a_2$ | $a_6$ |
| $a_2$ | $a_6$ | $a_3$ | $a_0$ | $a_4$ | $a_1$ | $a_5$ | $a_2$ |
| $a_3$ | $a_2$ | $a_6$ | $a_3$ | $a_0$ | $a_4$ | $a_1$ | $a_5$ |
| $a_4$ | $a_5$ | $a_2$ | $a_6$ | $a_3$ | $a_0$ | $a_4$ | $a_1$ |
| $a_5$ | $a_1$ | $a_5$ | $a_2$ | $a_6$ | $a_3$ | $a_0$ | $a_4$ |
| $a_6$ | $a_4$ | $a_1$ | $a_5$ | $a_2$ | $a_6$ | $a_3$ | $a_6$ |

Clearly T is a biquasi groupoid. The order of T is 13. T has elements of finite order, for

$$x = \begin{pmatrix} 1 & 2 & 3 \\ 1 & 3 & 2 \end{pmatrix} \in T$$

is such that

$$x^2 = \begin{pmatrix} 1 & 2 & 3 \\ 1 & 2 & 3 \end{pmatrix} = e;$$

but $|T| = 13$, so a prime order, hence T cannot have Cauchy elements. This T cannot have Lagrange sub-biquasi groupoids. T is a non-Lagrange bi quasi groupoid and T has no p-Sylow sub-biquasi groupoids.

**THEOREM 5.1.3:** *All biquasi groups are biquasi groupoids and not conversely.*

*Proof:* Direct; hence left for the reader to prove.

The above example is a biquasi groupoid, which is not a biquasi group. We can have a biquasi groupoid to have a substructure, which is a biquasi group and so on.

Such study will be carried out in the next section when we go on for the study of Smarandache bigroupoids.

**PROBLEMS:**

1. Give an example of a bigroupoid of order 7.

2. Show by an example a bigroupoid of order 13 can have proper sub-biquasi groupoids.

3. Find the Cauchy elements of the bigroupoid $G = G_1 \cup G_2$ given by the following tables:

   Table for $G_1$

   | *     | $a_0$ | $a_1$ | $a_2$ | $a_3$ | $a_4$ |
   |-------|-------|-------|-------|-------|-------|
   | $a_0$ | $a_0$ | $a_1$ | $a_2$ | $a_3$ | $a_4$ |
   | $a_1$ | $a_1$ | $a_1$ | $a_3$ | $a_1$ | $a_3$ |
   | $a_2$ | $a_2$ | $a_4$ | $a_2$ | $a_4$ | $a_2$ |
   | $a_3$ | $a_3$ | $a_3$ | $a_1$ | $a_3$ | $a_1$ |
   | $a_4$ | $a_4$ | $a_2$ | $a_4$ | $a_2$ | $a_4$ |

   $G_2 = S(3)$ the symmetric semigroup of order 27. Clearly $|G| = 32$. Does $G = G_1 \cup G_2$ have Lagrange sub-bigroupoids? Find sub-bigroupoids of G.

4. Define homomorphism of bigroupoids. Find a homomorphism from X to X' where $X = G_1 \cup S(2)$ and $X' = G_2 \cup S(3)$ where $G_1$ and $G_2$ are groupoids given by the following tables:

   Table of $G_1$

   | +     | $x_0$ | $x_1$ | $x_2$ | $x_3$ |
   |-------|-------|-------|-------|-------|
   | $x_0$ | $x_0$ | $x_1$ | $x_2$ | $x_3$ |
   | $x_1$ | $x_3$ | $x_0$ | $x_1$ | $x_2$ |
   | $x_2$ | $x_2$ | $x_3$ | $x_0$ | $x_1$ |
   | $x_3$ | $x_1$ | $x_2$ | $x_3$ | $x_0$ |

   Table of $G_2$

   | * | 0 | 1 | 2 | 3 | 4 | 5 | 6 | 7 |
   |---|---|---|---|---|---|---|---|---|
   | 0 | 0 | 6 | 4 | 2 | 0 | 6 | 4 | 2 |
   | 1 | 2 | 0 | 6 | 4 | 2 | 0 | 6 | 4 |
   | 2 | 4 | 2 | 0 | 6 | 4 | 2 | 0 | 6 |
   | 3 | 6 | 4 | 2 | 0 | 6 | 4 | 2 | 0 |
   | 4 | 0 | 6 | 4 | 2 | 0 | 6 | 4 | 2 |
   | 5 | 2 | 0 | 6 | 4 | 2 | 0 | 6 | 4 |
   | 6 | 4 | 2 | 0 | 6 | 4 | 2 | 0 | 6 |
   | 7 | 6 | 4 | 2 | 0 | 6 | 4 | 2 | 0 |

5. Give an example of a bigroupoid of order 12 which has no Cauchy elements.

6. Give an example of a biquasi loop.

7. Find a biquasi semigroup of order 8.

8. Give an example of a biquasi groupoid of order 16.

9. Is $X = L_5(3) \cup S(4)$ a
   a. biquasi groupoid?
   b. biquasi loop? Justify your claim.

10. Find for the bigroupoid $G = G_1 \cup G_2$ where $G_1$ is given by the following table

| •   | $a_0$ | $a_1$ | $a_2$ | $a_3$ | $a_4$ |
|-----|-------|-------|-------|-------|-------|
| $a_0$ | $a_0$ | $a_4$ | $a_3$ | $a_2$ | $a_1$ |
| $a_1$ | $a_2$ | $a_1$ | $a_0$ | $a_4$ | $a_3$ |
| $a_2$ | $a_4$ | $a_3$ | $a_2$ | $a_1$ | $a_0$ |
| $a_3$ | $a_1$ | $a_0$ | $a_4$ | $a_3$ | $a_2$ |
| $a_4$ | $a_3$ | $a_2$ | $a_1$ | $a_0$ | $a_4$ |

and $G_2 = \{a, b, 0, c, d \mid a^2 = b^2 = c^2 = d^2 = 0; ab = ba = c; ac = ca = 0; ad = da = 0; cd = dc = 0\}$. Does this bigroupoid have

a. sub-bigroupoid?
b. ideals?

## 5.2 Smarandache bigroupoids and its properties

In this section we for the first time introduce the notion Smarandache bigroupoids and also give several of its properties. Many interesting results about Smarandache bigroupoids are given and it finds its application in semi-automaton and automaton which will be dealt in the next section of this chapter.

**DEFINITION 5.2.1:** *Let $(G, +, \bullet)$ be a non-empty set with $G = G_1 \cup G_2$, we call G a Smarandache bigroupoid (S-bigroupoid) if*

   i. *$G_1$ and $G_2$ are distinct proper subsets of G such that $G = G_1 \cup G_2$ ($G_1 \not\subset G_2$ or $G_2 \not\subset G_1$).*
   ii. *$(G_1, +)$ is a S-groupoid.*
   iii. *$(G_2, \bullet)$ is a S-semigroup.*

*Now it is worthy to mention that '+' does not mean always usual '+' but only a binary operation which is non-associative in this case.*

**Example 5.2.1:** Let $(G, +, \bullet)$ denote a S-bigroupoid where $G = G_1 \cup G_2$ with the S-groupoid $G_1$ given by the following table:

| • | 0 | 1 | 2 | 3 | 4 | 5 |
|---|---|---|---|---|---|---|
| 0 | 0 | 3 | 0 | 3 | 0 | 3 |
| 1 | 1 | 4 | 1 | 4 | 1 | 4 |
| 2 | 2 | 5 | 2 | 5 | 2 | 5 |
| 3 | 3 | 0 | 3 | 0 | 3 | 0 |
| 4 | 4 | 1 | 4 | 1 | 4 | 1 |
| 5 | 5 | 2 | 5 | 2 | 5 | 2 |

and $G_2 = S(3)$ the S-semigroup. Clearly G is a S-bigroupoid.

**DEFINITION 5.2.2:** *Let (G, +, •) with $G = G_1 \cup G_2$ be a S-bigroupoid. We say G is of finite order if the number of elements in G is finite, if the number of elements in G is infinite we call G an infinite bigroupoid. We denote the order of the groupoid G by o(G) or |G|.*

**DEFINITION 5.2.3:** *Let (G, +, •) be a bigroupoid. A non-empty proper subset (P, +, •) is said to be a Smarandache sub-bigroupoid (S-sub-bigroupoid) if P itself is a S-bigroupoid.*

One may think why we have not taken a S-bigroupoid but only a bigroupoid.

The answer to this question is given by the following theorem.

**THEOREM 5.2.1:** *Let (G, +, •) be a bigroupoid. If G has a proper subset (P, +, •), which is a S-sub-bigroupoid, then G is a S-bigroupoid.*

*Proof*: Straightforward by the definitions.

**DEFINITION 5.2.4:** *Let (G, +, •) be a bigroupoid. We say G is a Smarandache commutative bigroupoid (S-commutative bigroupoid) if every S-sub-bigroupoid of G is commutative. If (G, +, •) has atleast one S-sub-bigroupoid which is commutative then we call G a Smarandache weakly commutative bigroupoid (S-weakly commutative bigroupoid).*

The following result is left for the reader to prove.

**THEOREM 5.2.2:** *Every S-commutative bigroupoid is S-weakly commutative bigroupoid.*

**DEFINITION 5.2.5:** *Let (G, +, •) be a bigroupoid of finite order. We call a S-sub-bigroupoid P of G to be a Smarandache Lagrange sub-bigroupoid (S- Lagrange sub-bigroupoid) of G if o(P) / o(G).*

*If every S-sub-bigroupoid P of G is a S-Lagrange sub-bigroupoid then we call the bigroupoid G to be a Smarandache Lagrange bigroupoid (S-Lagrange bigroupoid). If no S-sub-bigroupoid P of G is a S-Lagrange sub-bigroupoid then we call G a Smarandache non-Lagrange bigroupoid (S- non-Lagrange bigroupoid).*

It is worthwhile to mention here that G may be a S- non-Lagrange bigroupoid still G may have a sub-bigroupoid whose order divides the order of G i.e. G may not have only S-sub-bigroupoids.

**DEFINITION 5.2.6:** *Let $(G, +, \bullet)$ be a bigroupoid of finite order. If p is a prime such that $p^\alpha/o(G)$ but $p^{\alpha+1} \not\mid o(G)$ further if G has a S-sub-bigroupoid P of order $p^\alpha$ then we call P a Smarandache p-Sylow sub-bigroupoid (S- p-Sylow sub-bigroupoid).*

*If for every prime p, dividing the o(G) we have corresponding associative S-p-Sylow sub-bigroupoid then we call the bigroupoid G to be a Smarandache p-Sylow bigroupoid (S-p-Sylow bigroupoid).*

**DEFINITION 5.2.7:** *Let $(G, +, \bullet)$ be a S-bigroupoid of finite order. An element $x \in G$ is said to be a Smarandache Cauchy element (S- Cauchy element) if $x^n = e$ and $n/o(G)$.*

*If the S-bigroupoid has atleast a Cauchy element then we call G a Smarandache Cauchy bigroupoid (S-Cauchy bigroupoid).*

**DEFINITION 5.2.8:** *Let $(G, +, \bullet)$ be a S-bigroupoid. We call G an idempotent bigroupoid if every element in G is an idempotent.*

*We call the S-bigroupoid to be a Smarandache idempotent bigroupoid (S- idempotent bigroupoid) if every element in every S-sub-bigroupoid in G is a S-idempotent.*

**DEFINITION 5.2.9:** *Let $(G, +, \bullet)$ a bigroupoid, $G = G_1 \cup G_2$. A non-empty proper subset P of G is said to be a Smarandache left bi-ideal (S-left bi-ideal) of the bigroupoid G if*

    i.    *$(P, +, \bullet)$ is a S- sub-bigroupoid.*
    ii.   *Let $P = P_1 \cup P_2$; $(P_1, +)$ is a left ideal of $G_1$ and $(P_2, \bullet)$ is a right ideal of $G_2$.*

*On similar lines we define S-right bi-ideal. If P is simultaneously both the S-left bi-ideal and S-right bi-ideal then we say P is a Smarandache bi-ideal (S-bi-ideal) of G.*

A Smarandache normal sub-bigroupoid; is defined as follows:

**DEFINITION 5.2.10:** *Let $G = G_1 \cup G_2$ be a bigroupoid. A S-sub-bigroupoid $(P = P_1 \cup P_2, +, \bullet)$ is said to be a Smarandache normal sub-bigroupoid (S-normal sub-bigroupoid) of G if*

    i.    *$aP_1 = P_1 a$ or $aP_2 = P_2 a$ (according as $a \in P_1$ or $P_2$).*
    ii.   *$P_1(xy) = (P_1 x)y$ or $P_2(xy) = (P_2 x) y$ (according as $x, y \in P_1$ or $x, y \in P_2$)*

*for all $x, y, a \in G$.*

*We call the bigroupoid Smarandache simple (S-simple) if it has no nontrivial S-normal sub-bi groupoids.*

Next we define when are S-bigroupoids normal.

**DEFINITION 5.2.11:** *Let $(G, +, \bullet)$ be a finite bigroupoid. We say $G$ to be a Smarandache normal bigroupoid (S-normal bigroupoid) if the largest S-sub-bigroupoid of $G$ is a normal bigroupoid.*

We say the largest S-sub-bigroupoid to mean the number of elements in it is the maximum. So only the term 'largest' is used. That is we take all S-sub-bigroupoids of G and consider the S-sub-bigroupoid, P of G such that o(P) is the greatest. If we have several such P's and if at least one of the P's is a normal bigroupoid then we accept G to be a S-normal bigroupoid.

**DEFINITION 5.2.12:** *Let $(G, +, \bullet)$ be a bigroupoid with $G = G_1 \cup G_2$. Let H and K be any two S-sub-bigroupoids of G; i.e. $H = H_1 \cup H_2$ and $K = K_1 \cup K_2$. We say H is Smarandache biconjugate (S-biconjugate) with K if there exists $x_1 \in H_1$ and $x_2 \in H_2$ with $x_1 K_1 = H_1$ (or $K_1 x_1$) and $x_2 K_2 = H_2$ (or $K_2 x_2$) 'or' in the mutually exclusive sense.*

**DEFINITION 5.2.13:** *Let $(G, +, \bullet)$ be a bigroupoid. We say G is Smarandache inner commutative bigroupoid (S-inner commutative bigroupoid) if every S-sub-bigroupoid of G is commutative.*

The concept of S-zero divisor, zero divisor, S-units, units, S-idempotents and idempotents in a S-bigroupoid is the same as that of a bigroupoid.

We do not differentiate the structure of the elements in a bigroupoid and a S-bigroupoid.

**DEFINITION 5.2.14:** *Let $(G, +, \bullet)$ be a bigroupoid. We say G is a Smarandache Moufang bigroupoid (S-Moufang bigroupoid) if G has a proper subset P where P is a S-sub-bigroupoid of G and all elements of P satisfy the Moufang identify i.e. P is a Moufang bigroupoid.*

*We call a bigroupoid $(G, +, \bullet)$ to be a Smarandache Bol bigroupoid (S-Bol bigroupoid) if G has a S-sub-bigroupoid P, where P satisfies the Bol identity.*

*Similarly if the bigroupoid G has a S-sub-bigroupoid, which is Bruck bigroupoid, then we call G a Smarandache Bruck bigroupoid (S-Bruck bigroupoid).*

*The notions of Smarandache right (left) alternative (S-right(left) alternative) or alternative bigroupoids is defined in a similar way. Thus for a bigroupoid to satisfy the Smarandache identities we demand a proper S-sub-bigroupoid of G must satisfy the identities.*

**DEFINITION 5.2.15:** *Let $(G, +, \bullet)$ be a bigroupoid. P be a S-sub-bigroupoid of G. We say G is a Smarandache P-bigroupoid (S-P-bigroupoid) if $(x \bullet y) \bullet x = x \bullet (y \bullet x)$ for all $x, y \in P$. If every S-sub-bigroupoid of G satisfies the identity $(x \bullet y) \bullet x = x \bullet (y \bullet x)$ we call G a Smarandache strong P-bigroupoid (S-strong P-bigroupoid).*

Several interesting properties about S-bigroupoids can be defined and studied. This task is left for the reader as an exercise.

Now we proceed onto define 3 types of Smarandache quasi bistructures. viz. Smarandache quasi biloops, Smarandache quasi groupoids and Smarandache quasi semigroups; all the three bi-structures are non-associative.

**DEFINITION 5.2.16:** *Let $(G, +, \bullet)$ be a biquasi loop. We say G is a Smarandache biquasi loop (S-biquasi loop) if it has a proper subset $(P, +, \bullet)$. $P \subset G$ such that P is a biquasi semigroup.*

*Example 5.2.2:* Let $(G, +, \bullet)$ be a biquasi loop. $G = G_1 \cup G_2$; $G_1 = L_5(2)$ and $G_2$ be a groupoid given by the following table:

| *   | $a_0$ | $a_1$ | $a_2$ | $a_3$ | $a_4$ | $a_5$ |
|-----|-------|-------|-------|-------|-------|-------|
| $a_0$ | $a_0$ | $a_4$ | $a_2$ | $a_0$ | $a_4$ | $a_2$ |
| $a_1$ | $a_2$ | $a_0$ | $a_4$ | $a_2$ | $a_0$ | $a_4$ |
| $a_2$ | $a_4$ | $a_2$ | $a_0$ | $a_4$ | $a_2$ | $a_0$ |
| $a_3$ | $a_0$ | $a_4$ | $a_2$ | $a_0$ | $a_4$ | $a_2$ |
| $a_4$ | $a_2$ | $a_0$ | $a_4$ | $a_2$ | $a_0$ | $a_4$ |
| $a_5$ | $a_4$ | $a_2$ | $a_0$ | $a_4$ | $a_2$ | $a_0$ |

Take $P = L_5(2) \cup \{0\}$ is a biquasi semigroup. So G is a S-biquasi loop.

*Example 5.2.3:* Let $(G, +, \bullet)$ be a biquasi loop where $G = G_1 \cup G_2$ with $G_1 = L_{11}(3)$ and $G_2$ given by the following table:

| * | 0 | 1 | 2 | 3 | 4 | 5 | 6 | 7 |
|---|---|---|---|---|---|---|---|---|
| 0 | 0 | 6 | 4 | 2 | 0 | 6 | 4 | 2 |
| 1 | 2 | 0 | 6 | 4 | 2 | 0 | 6 | 4 |
| 2 | 4 | 2 | 0 | 6 | 4 | 2 | 0 | 6 |
| 3 | 6 | 4 | 2 | 0 | 6 | 4 | 2 | 0 |
| 4 | 0 | 6 | 4 | 2 | 0 | 6 | 4 | 2 |
| 5 | 2 | 0 | 6 | 4 | 2 | 0 | 6 | 4 |
| 6 | 4 | 2 | 0 | 6 | 4 | 2 | 0 | 6 |
| 7 | 6 | 4 | 2 | 0 | 6 | 4 | 2 | 0 |

P be a proper set with $P = P_1 \cup P_2$ where $P = L_{11}(3) \cup \{0, 4\}$. Clearly P is a biquasi semigroup. Hence $(G, +, \bullet)$ is a S-biquasi groupoid.

**DEFINITION 5.2.17:** *Let $(G, +, \bullet)$ be a biquasi semigroup. We say $(G, +, \bullet)$ is Smarandache biquasi semigroup (S-biquasi semigroup) if G has a proper subset P such that P is biquasi group under the operations of G.*

*Example 5.2.4:* Let $(G, +, \bullet)$ where $G = G_1 \cup G_2$ be a biquasi semigroup with $G_1 = L_5(3)$ and $G_2 = S(3)$. Take $P = P_1 \cup P_2$ where $P_1 = \{e, 3\}$ and $P_2 = S(3)$. Clearly $(P, +, \bullet)$ is a biquasi group. Hence G is a biquasi semigroup.

**DEFINITION 5.2.18:** *Let $(G, +, \bullet)$ be a biquasi groupoid. If P is a proper subset of G and $(P, +, \bullet)$ is such that $(P, +)$ is a group and $(P, \bullet)$ is a S-groupoid then we call G a Smarandache biquasi groupoid (S-biquasi groupoid).*

*Example 5.2.5:* Let $(G, +, \bullet)$ be a biquasi groupoid $G = G_1 \cup G_2$ where $G_1 = S_7$ and $G_2$ is a groupoid given by the following table:

| * | 0 | 1 | 2 | 3 | 4 | 5 |
|---|---|---|---|---|---|---|
| 0 | 0 | 0 | 0 | 0 | 0 | 0 |
| 1 | 3 | 3 | 3 | 3 | 3 | 3 |
| 2 | 0 | 0 | 0 | 0 | 0 | 0 |
| 3 | 3 | 3 | 3 | 3 | 3 | 3 |
| 4 | 0 | 0 | 0 | 0 | 0 | 0 |
| 5 | 3 | 3 | 3 | 3 | 3 | 3 |

Clearly $P = A_7 \cup \{0, 2, 3\}$ is such that $A_7$ is a group and $\{0, 2, 3\}$ is a S-groupoid. Thus G is a S-biquasi groupoid.

For these three Smarandache biquasi structure we can define Smarandache bi-ideals. Smarandache Moufang, Bol, Bruck, alternative biquasi structures, Smarandache normal sub-biquasi structure, Smarandache homomorphism and several other properties as in case of any Smarandache algebraic bistructure. All these research is left for the reader.

**PROBLEMS:**

1. Give an example of a S-biquasi groupoid and find

    i. S-sub-biquasi groupoid
    ii. S-normal sub-biquasi groupoid.

2. For the S-biquasi loop $G = L_{15}(8) \cup \{Z_9 \,(3, 6)$ i.e. for a, b $\in Z_9$ define $a \bullet b = 3a + 6b \pmod 9\}$ find whether G satisfies any one of the well known identities so that the biquasi loop is a S-(Moufang / Bol / Bruck / alternative / P) biquasi loop.

3. Let $(G, +, \bullet)$ be such that $G = G_1 \cup G_2$ where $G_1 = L_9(8)$ is a loop, $G_2 = S(7)$, the semigroup.

    i. Is G a S-biquasi loop?
    ii. Is G a S-biquasi groupoid?
    iii. Is G a S-biquasi semigroup?

4. Give an example of a S-biquasi groupoid of order 18.

5. Give an example of a S-biquasi semigroup of order 12?

6. Can there exist a S-biquasi loop of order 11? Justify your claim.

7. Give an example of a S-biquasi loop of order 13, which has proper S-sub-biquasi loops.

8. For the S-biquasi semigroups (G, +, •) where G = $L_7(3) \cup S(5)$ and G' = $L_{11}(3) \cup S(7)$ find a S-biquasi semigroup homomorphism.

9. Find the number of S-sub-biquasi semigroups in (G, +, •) given by G = $G_1 \cup G_2$ = $L_5(2) \cup S(3)$.

10. Can the S-sub-biquasi semigroup given in problem 9 have S-normal sub-biquasi semigroup?

11. Is (G, +, •) given by G = $G_1 \cup G_2$ where $G_1$ = S(7) a semigroup and $G_2$ = $Z_6(4, 2)$ a groupoid a S-quasi bi-semigroup? S-quasi bi-groupoid? S-quasi bi-loop?

12. Can a S-biquasi biloop be a S-quasi group? Justify your claim.

## 5.3 Applications of bigroupoids and S-bigroupoids

In this section we introduce the applications of bigroupoids as S-groupoids are used in the application of semi-automaton and automaton. We find the application of bigroupoids to semi-automaton and automaton. The notion of Smarandache semi-automaton and Smarandache automaton was defined only in the year 2002 [114].

**DEFINITION 5.3.1:** $Y_s = (Z, \overline{A}_s, \overline{\delta}_s)$ *is said to be a Smarandache semi-automaton if* $\overline{A}_s = \langle A \rangle$ *is the free groupoid generated by A with Λ the empty element adjoined with it and* $\overline{\delta}_s$ *is the function from* $Z \times \overline{A}_s \to Z$. *Thus the S-semi-automaton contains Y = (Z, $\overline{A}$, $\overline{\delta}$) as a new semi-automaton which is a proper sub-structure of $Y_s$.*

*Or equivalently, we define a S-semi-automaton as one, which has a new semi-automaton as a sub-structure.*

The advantages of the S-semi-automaton are if for the triple Y = (Z, A, δ) is a semi-automaton with Z, the set of states, A the input alphabet and δ : Z × A → Z is the next state function.

When we generate the S-free groupoid by A and adjoin with it the empty alphabet Λ then we are sure that $\overline{A}$ has all free semigroups. Thus, each free semigroup will give a new semi-automaton. Thus by choosing a suitable A we can get several new semi-automaton using a single S-semi-automaton.

We now give some examples of S-semi-automaton using finite groupoids.

When examples of semi-automaton are given usually the books use either the set of modulo integers $Z_n$ under addition or multiplication we build using groupods in $Z_n$.

**DEFINITION 5.3.2:** $\overline{Y}_s' = (Z_1, \overline{A}_s, \overline{\delta}'_s)$ is called the Smarandache subsemi-automaton (S-subsemi-automaton) of $\overline{Y}_s = (Z_2, \overline{A}_s, \overline{\delta}'_s)$ denoted by $\overline{Y}_s' \leq \overline{Y}_s$ if $Z_1 \subset Z_2$ and $\overline{\delta}'_s$ is the restriction of $\overline{\delta}_s$ on $Z_1 \times \overline{A}_s$ and $\overline{Y}_s'$ has a proper subset $\overline{H} \subset \overline{Y}_s'$ such that $\overline{H}$ is a new semi-automaton.

*Example 5.3.1:* Let $Z = Z_4(2, 1)$ and $A = Z_6(2, 1)$. The S-semi-automaton $(Z, A, \delta)$ where $\delta : Z \times A \to Z$ is given by $\delta(z, a) = z \bullet a \pmod 4$.

We get the following table:

| δ | 0 | 1 | 2 | 3 |
|---|---|---|---|---|
| 0 | 0 | 1 | 2 | 3 |
| 1 | 2 | 3 | 0 | 1 |
| 2 | 0 | 1 | 2 | 3 |
| 3 | 2 | 3 | 0 | 1 |

We get the following graph for this S-semi-automaton.

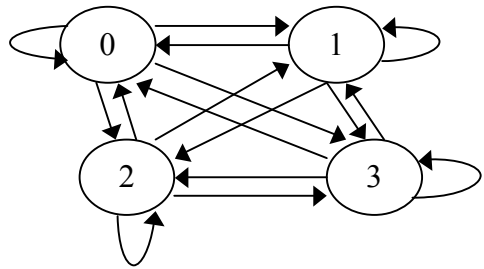

**Figure 5.3.1**

This has a nice S-subsemi-automaton given by the table. In the example 5.3.1 we see this S-semi-automaton has a S-subsemi-automaton given by the table:

| δ | 0 | 2 |
|---|---|---|
| 0 | 0 | 2 |
| 2 | 0 | 2 |

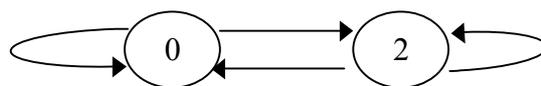

**Figure 5.3.2**

Thus this has a S-subsemi-automaton $Z_1$ given by $Z_1 = \{0, 2\}$ states.

**DEFINITION 5.3.3:** $\overline{K}_s = (Z, \overline{A}_s, \overline{B}_s, \overline{\delta}_s, \overline{\lambda}_s)$ is defined to be a Smarandache automaton (S-automaton) if $\overline{K} = (Z, \overline{A}_s, \overline{B}_s, \overline{\delta}_s, \overline{\lambda}_s)$ is the new automaton and $\overline{A}_s$ and $\overline{B}_s$, the S-free groupoids so that $\overline{K} = (Z, \overline{A}_s, \overline{B}_s, \overline{\delta}_s, \overline{\lambda}_s)$ is the new automaton got from K and $\overline{K}$ is strictly contained in $\overline{K}_s$.

Thus S-automaton enables us to adjoin some more elements which is present in A and freely generated by A, as a free groupoid; that will be the case when the compositions may not be taken as associative.

Secondly, by using S-automaton we can couple several automaton as

$$\begin{aligned} Z &= Z_1 \cup Z_2 \cup ... \cup Z_n \\ A &= A_1 \cup A_2 \cup ... \cup A_n \\ B &= B_1 \cup B_2 \cup ... \cup B_n \\ \lambda &= \lambda_1 \cup \lambda_2 \cup ... \cup \lambda_n \\ \delta &= \delta_1 \cup \delta_2 \cup ... \cup \delta_n. \end{aligned}$$

where the union of $\lambda_i \cup \lambda_j$ and $\delta_i \cup \delta_j$ denote only extension maps as '$\cup$' has no meaning in the composition of maps, where $K_i = (Z_i, A_i, B_i, \delta_i, \lambda_i)$ for i = 1, 2, 3, ..., n and $\overline{K} = \overline{K}_1 \cup \overline{K}_2 \cup ... \cup \overline{K}_n$. Now $\overline{K}_s = (\overline{Z}_s, \overline{A}_s, \overline{B}_s, \overline{\lambda}_s, \overline{\delta}_s)$ is the S-automaton.

A machine equipped with this S-automaton can use any new automaton as per need. We give some examples of S-automaton using S-groupoids.

*Example 5.3.2:* Let $Z = Z_4$ (3, 2), $A = B = Z_5$ (2, 3). $K = (Z, A, B, \delta, \lambda)$ is a S-automaton defined by the following tables where $\delta (z, a) = z * a \pmod 4$ and $\lambda (z, a) = z * a \pmod 5$.

| λ | 0 | 1 | 2 | 3 | 4 |
|---|---|---|---|---|---|
| 0 | 0 | 3 | 1 | 4 | 2 |
| 1 | 2 | 0 | 3 | 1 | 4 |
| 2 | 4 | 2 | 0 | 3 | 1 |
| 3 | 1 | 4 | 2 | 0 | 3 |

| δ | 0 | 1 | 2 | 3 | 4 |
|---|---|---|---|---|---|
| 0 | 0 | 2 | 0 | 2 | 0 |
| 1 | 3 | 1 | 3 | 1 | 3 |
| 2 | 2 | 0 | 2 | 0 | 2 |
| 3 | 1 | 3 | 1 | 3 | 1 |

We obtain the following graph:

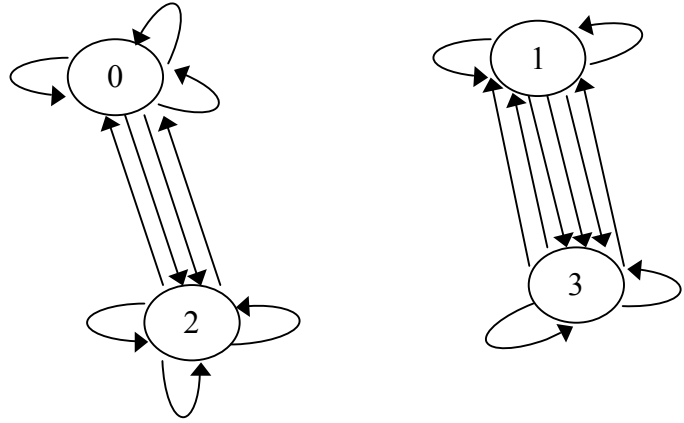

**Figure 5.3.3**

Thus we see this automaton has 2 S-sub automatons given by the states {0, 2} and {1, 3}.

**DEFINITION 5.3.4:** $\overline{K}_s' = (Z_1, \overline{A}_s, \overline{B}_s, \overline{\delta}_s, \overline{\lambda}_s)$ is called Smarandache sub-automaton (S-sub-automaton) of $\overline{K}_s = (Z_2, \overline{A}_s, \overline{B}_s, \overline{\delta}_s, \overline{\lambda}_s)$ denoted by $\overline{K}_s' \leq \overline{K}_s$ if $Z_1 \subseteq Z_2$ and $\overline{\delta}_s'$ and $\overline{\lambda}_s'$ are the restriction of $\overline{\delta}_s$ and $\overline{\lambda}_s$ respectively on $Z_1 \times \overline{A}_s$ and has a proper subset $\overline{H} \subset \overline{K}_s'$ such that $\overline{H}$ is a new automaton.

**DEFINITION 5.3.5:** Let $\overline{K}_1$ and $\overline{K}_2$ be any two S-automaton where $\overline{K}_1 = (Z_1, \overline{A}_s, \overline{B}_s, \overline{\delta}_s, \overline{\lambda}_s)$ and $\overline{K}_2 = (Z_2, \overline{A}_s, \overline{B}_s, \overline{\delta}_s, \overline{\lambda}_s)$. A map $\phi : \overline{K}_1$ to $\overline{K}_2$ is a Smarandache automaton homomorphism (S-automaton homomorphism) if $\phi$ restricted from $K_1 = (Z_1, A_1, B_1, \delta_1, \lambda_1)$ and $K_2 = (Z_2, A_2, B_2, \delta_2, \lambda_2)$ denoted by $\phi_r$ is a automaton homomorphism from $K_1$ to $K_2$. $\phi$ is called a monomorphism (epimorphism or isomorphism) if there is an isomorphism $\phi_r$ from $K_1$ to $K_2$.

**DEFINITION 5.3.6:** Let $\overline{K}_1$ and $\overline{K}_2$ be two S-automatons, where $\overline{K}_1 = (Z_1, \overline{A}_s, \overline{B}_s, \overline{\delta}_s, \overline{\lambda}_s)$ and $\overline{K}_2 = (Z_2, \overline{A}_s, \overline{B}_s, \overline{\delta}_s, \overline{\lambda}_s)$.

*The Smarandache automaton direct product (S-automaton direct product) of $\overline{K}_1$ and $\overline{K}_2$ denoted by $\overline{K}_1 \times \overline{K}_2$ is defined as the direct product of the automaton $K_1 = (Z_1, A_1, B_1, \delta_1, \lambda_1)$ and $K_2 = (Z_2, A_2, B_2, \delta_2, \lambda_2)$ where $K_1 \times K_2 = (Z_1 \times Z_2, A_1 \times A_2, B_1 \times B_2, \delta, \lambda)$ with $\delta((z_1, z_2), (a_1, a_2)) = (\delta_1(z_1, a_1), \delta_2(z_2, a_2))$, $\lambda((z_1, z_2), (a_1, a_2)) = (\lambda_1(z_1, a_2), \lambda_2(z_2, a_2))$ for all $(z_1, z_2) \in Z_1 \times Z_2$ and $(a_1, a_2) \in A_1 \times A_2$.*

**Remark:** Here in $\overline{K}_1 \times \overline{K}_2$ we do not take the free groupoid to be generated by $A_1 \times A_2$ but only free groupoid generated by $\overline{A}_1 \times \overline{A}_2$

Thus the S-automaton direct product exists wherever a automaton direct product exists.

We have made this in order to make the Smarandache parallel composition and Smarandache series composition of automaton extendable in a simple way.

**DEFINITION 5.3.7:** *A Smarandache groupoid $G_1$ divides a Smarandache groupoid $G_2$ if the corresponding semigroups $S_1$ and $S_2$ of $G_1$ and $G_2$ respectively divides, that is, if $S_1$ is a homomorphic image of a sub-semigroup of $S_2$.*

*In symbols $G_1 \mid G_2$. The relation divides is denoted by '|'.*

**DEFINITION 5.3.8:** *Let $\overline{K}_1 = (Z_1, \overline{A}_s, \overline{B}_s, \overline{\delta}_s, \overline{\lambda}_s)$ and $\overline{K}_2 = (Z_2, \overline{A}_s, \overline{B}_s, \overline{\delta}_s, \overline{\lambda}_s)$ be two S-automaton. We say the S-automaton $\overline{K}_1$ divides the S-automaton $\overline{K}_2$ if in the automatons $K_1 = (Z_1, A, B, \delta_1, \lambda_1)$ and $K_2 = (Z_2, A, B, \delta_2, \lambda_2)$, if $K_1$ is the homomorphic image of a sub-automaton of $K_2$.*

*Notationally $K_1 \mid K_2$.*

**DEFINITION 5.3.9:** *Two S-automatons $\overline{K}_1$ and $\overline{K}_2$ are said to be equivalent if they divide each other. In symbols $\overline{K}_1 \sim \overline{K}_2$.*

For more about automaton, semi-automaton, S-automaton and S-semi automaton please refer [114].

## 5.4 Direct Product of S-automaton

We proceed on to define direct product of S-automaton and study about them. We can extend the direct product of semi-automaton to more than two S-automatons. Using the definition of direct product of two automaton $K_1$ and $K_2$ with an additional assumption we define Smarandache series composition of automaton.

**DEFINITION 5.4.1:** *Let $K_1$ and $K_2$ be any two S-automatons where $\overline{K}_1 = (Z_1, \overline{A}_s, \overline{B}_s, \overline{\delta}_s, \overline{\lambda}_s)$ and $\overline{K}_2 = (Z_2, \overline{A}_s, \overline{B}_s, \overline{\delta}_s, \overline{\lambda}_s)$ with an additional assumption $A_2 = B_1$.*

*The Smarandache automaton composition series (S-automaton composition series) denoted by $\overline{K}_1 \Vdash \overline{K}_2$ of $\overline{K}_1$ and $\overline{K}_2$ is defined as the series composition of the automaton $K_1 = (Z_1, A_1, B_1, \delta_1, \lambda_1)$ and $K_2 = (Z_2, A_2, B_2, \delta_2, \lambda_2)$ with $\overline{K}_1 \Vdash \overline{K}_2 = (Z_1 \times Z_2, A_1, B_2, \delta, \lambda)$ where $\delta((z_1, z_2), a_1) = (\delta_1(z_1, a_1), \delta_2(z_2, \lambda_1(z_1, a_1))$ and $\lambda((z_1, z_2), a_1) = (\lambda_2(z_2, \lambda_1(z_1, a_1))$ $((z_1, z_2) \in Z_1 \times Z_2, a_1 \in A_1)$.*

*This automaton operates as follows: An input $a_1 \in A_1$ operates on $z_1$ and gives a state transition into $z_1' = \delta_1(z_1, a_1)$ and an output $b_1 = \lambda_1(z_1, a_2) \in B_1 = A_2$. This output $b_1$ operates on $Z_2$ transforms a $z_2 \in Z_2$ into $z_2' = \delta_2(a_2, b_1)$ and produces the output $\lambda_2(z_2, b_1)$.*

*Then $\overline{K}_1 \Vdash \overline{K}_2$ is in the next state $(z_1', z_2')$ which is clear from the following circuit:*

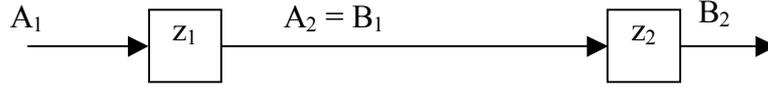

Now a natural question would be do we have a direct product, which corresponds to parallel composition, of the 2 Smarandache automatons $\overline{K}_1$ and $\overline{K}_2$.

Clearly the Smarandache direct product of automatons $\overline{K}_1 \times \overline{K}_2$ since $Z_1$, and $Z_2$ can be interpreted as two parallel blocks. $A_1$ operates on $Z_i$ with output $B_i$ ($i \in \{1, 2\}$), $A_1 \times A_2$ operates on $Z_1 \times Z_2$, the outputs are in $B_1 \times B_2$.

For more about S-automaton and S-semi-automaton refer [33, 114]. Now we proceed on to define the notion of semi-automaton and automaton using bigroupoids.

We have already defined the concept of bisemi-automaton and biautomaton in chapter 3 using bisemigroups. Now we define using bigroupoids the concept of bisemi-automaton and biautomaton.

**DEFINITION 5.4.2:** $\overline{Y} = (Z, \overline{A}_s, \overline{\delta}_s)$ *is said to be a Smarandache bisemi-automaton (S-bisemi-automaton) if* $\overline{A}_s = \langle A \rangle$ *is the free groupoid generated by A with $\Lambda$ the empty element adjoined with it and $\overline{\delta}_s$ is the function from $Z \times \overline{A}_s \to Z$. Thus the S-bisemi-automaton contains $Y = (Z, \overline{A}, \overline{\delta})$ as a new bi semi-automaton which is a proper substructure of $\overline{Y}$.*

*Or equivalently, we define a S-semi-automaton as one, which has a new bisemi-automaton as its substructure.*

**DEFINITION 5.4.3:** $X_t = (Z_t, \overline{A}_s, \overline{\delta}_s)$ *is called the Smarandache subsemi-automaton (S-subsemi-automaton) of $X = (Z_p, \overline{A}_s, \overline{\delta}_s)$ denoted by $X_t \leq X$ if $Z_t \subset Z_p$ and $\overline{\delta}'_s$ is the restriction of $\overline{\delta}_s$ on $Z_t \times A_s$ and $X_t$ has a proper subset $\overline{H} \subset X_t$ such that H is a new semi-automaton.*

**DEFINITION 5.4.4:** $\overline{X} = (Z, \overline{A}_s, \overline{B}_s, \overline{\delta}_s, \overline{\lambda}_s)$ *is said to be Smarandache biautomaton (S-biautomaton) if $\overline{A}_s = \langle A \rangle$ is the free bigroupoid generated by A with $\Lambda$ the empty element adjoined with it and $\overline{\delta}_s$ is the function from $Z \times \overline{A}_s \to Z$ and $\overline{B}_s$. Thus the S-biautomaton contains $X = (Z, \overline{A}, \overline{B}, \overline{\delta}, \overline{\lambda})$ as the new biautomaton which is a proper substructure of $\overline{X}$.*

*The notion of Smarandache sub-bi-automaton (S-sub-bi-automaton) and Smarandache sub-bi-semi-automaton (S-sub-bisemi-automaton) are defined in a similar way as that of S-subautomaton and S-subsemi-automaton.*

Several important results and new definitions can be introduced and studied. The S-bisemi-automaton and S-biautomaton can certainly do multifold of problems in multidimensional ways. Thus these finite machines can do several work ingrained in a single piece.

Thus in the opinion of the author these S-biautomaton and S-bisemi-automaton when used practically in the construction of finite machines will yield better results with lesser economy.

**PROBLEMS:**

1. Give an example of S-bisemi-automaton.
2. Illustrate by graph a biautomaton using the set of states to be just 7.
3. Prove every S-bisemi-automaton is a new bi-semi-automaton.
4. Can we define bicascades?
5. Is it possible to define Smarandache bicascades? Justify your claim.
6. Describe a biautomaton with only 2 states.
7. Can two finite automaton be combined to get a bi-automaton? Justify your claim.

**Chapter 6**

# BIRINGS AND SMARANDACHE BIRINGS

Since there is no book on the concepts of bialgebraic structures like bigroups, biloops, birings … , we felt it very essential to define and describe these systems and illustrate them with examples. To the best of the authors knowledge the author has not seen any paper on birings. Only bigroups have been introduced in the year 1994 [34]. Further the notions like bialgebras are introduced mainly not totally in the algebraic flavour. They are built for Lie algebras and Jordan algebras as a moniod object in the category of co-moniod. But as our book is on bialgebraic structures and above all Smarandache bialgebraic structures and as there is no solid Smarandache category theory in literature we restrict ourselves only to the case of pure algebraic treatment. The chapter has three sections. In the first section we deal with the introduction and study of birings and in the second section we define non-associative bi-ring and in section three we define Smarandache birings and analyse them.

## 6.1 Birings and its properties

In this section we define the concept of birings and introduce the notions like bi-ideals, biring homomorphism etc. The study of birings is very new and is also very interesting and innovative. All the possible ways of analyzing birings is carried out in this section with illustrative examples.

**DEFINITION 6.1.1:** *A nonempty set $(R, +, \bullet)$ with two binary operations '+' and '$\bullet$' is said to be a biring if $R = R_1 \cup R_2$ where $R_1$ and $R_2$ are proper subsets of $R$ and*

      i. *$(R_1, +, \bullet)$ is a ring.*
      ii. *$(R_2, +, \bullet)$ is a ring.*

*Example 6.1.1*: Let $R = \{0, 5, 2, 4, 6, 8\}$ be a set under '+' and '$\bullet$' modulo 10. Clearly $R = R_1 \cup R_2$ is a biring; as $R_1 = \{0, 5\}$ and $R_2 = \{0, 2, 4, 6, 8\}$ are rings under '+' and '$\bullet$' modulo 10. Hence $R = R_1 \cup R_2$ is a biring.

**DEFINITION 6.1.2:** *A biring R is said to be finite if it has a finite number of elements if R has infinite number of elements R is said to be an infinite biring.*

*We denote the order of the biring R by either o(R) or |R|.*

The biring given in example 6.1.1 is a biring of finite order and $o(R) = 6$.

**DEFINITION 6.1.3:** *A biring $R = R_1 \cup R_2$ is said to be commutative biring if both $R_1$ and $R_2$ are commutative rings otherwise i.e. if one of $R_1$ or $R_2$ is a non-commutative ring then we say R is a non-commutative biring. We say the biring R has a monounit if a unit exists which is common to both $R_1$ and $R_2$ we call the common unit as the* monounit. If $R_1$ and $R_2$ are rings which has separate units than we say the biring R has a unit.

**Example 6.1.2:** Let R be the collection of all $2 \times 2$ upper triangular and lower triangular matrices with entries from Z; the ring of integers.

$$R = \left\{ \begin{pmatrix} a & 0 \\ b & c \end{pmatrix}, \begin{pmatrix} x & y \\ 0 & z \end{pmatrix} \middle| a, b, c, x, y, z \in Z \right\}.$$

$R = R_1 \cup R_2$ where,

$$R_1 = \left\{ \begin{pmatrix} a & 0 \\ b & c \end{pmatrix} \middle| a, b, c \in Z \right\} \text{ and } R_2 = \left\{ \begin{pmatrix} x & y \\ 0 & z \end{pmatrix} \middle| x, y, z \in Z \right\}$$

$R_1$ and $R_2$ are rings under usual matrix addition and multiplication;

$$\begin{pmatrix} 1 & 0 \\ 0 & 1 \end{pmatrix} = I_{2 \times 2}$$

is the mono-unit of R.

**DEFINITION 6.1.4:** *Let $R = R_1 \cup R_2$ be a biring. A non-empty subset S of R is said to be a sub-biring of R if $S = S_1 \cup S_2$ and S itself is a biring and $S_1 = S \cap R_1$ and $S_2 = R_2 \cap S$.*

**Example 6.1.3:** Let $R = R_1 \cup R_2$ be a biring given in example 6.1.2.
Take $S = S_1 \cup S_2$ where

$$S_1 = \left\{ \begin{pmatrix} 0 & 0 \\ a & 0 \end{pmatrix} \middle| a \in Z \right\} \text{ and}$$

$$S_2 = \left\{ \begin{pmatrix} 0 & a \\ 0 & 0 \end{pmatrix} \middle| a \in Z \right\}.$$

Clearly S is a sub-biring of R.

**THEOREM 6.1.1:** *Let R be a biring where $R = R_1 \cup R_2$. A non-empty subset $S = S_1 \cup S_2$ of R is a sub-biring of R if and only if $R_1 \cap S = S_1$ and $R_2 \cap S = S_2$ are subrings of $R_1$ and $R_2$ respectively.*

*Proof:* Let $R = R_1 \cup R_2$ be a biring. S a proper subset of R. If S is a sub-biring of R then S is a biring, so $S = S_1 \cup S_2$, where $S_1$ and $S_2$ are subrings of $R_1$ and $R_2$ respectively and $S_1 = S \cap R_1$ and $S_2 = S \cap R_2$. Clearly, the converse is true by the very definitions of birings and sub-birings.

**DEFINITION 6.1.5:** *Let $R = R_1 \cup R_2$ be a biring. A non-empty subset I of R is said to be a right bi-ideal of R if $I = I_1 \cup I_2$ where $I_1$ is a right ideal of $R_1$ and $I_2$ is a right* ideal of $R_2$. I is said to be a left bi-ideal of R if $I = I_1 \cup I_2$ where $I_1$ and $I_2$ are left ideals of $R_1$ and $R_2$ respectively. If $I = I_1 \cup I_2$ is such that both $I_1$ and $I_2$ are ideals of $R_1$ and $R_2$ respectively then we say I is a bi-ideal of R.

Now it may happen when $I = I_1 \cup I_2$; $I_1$ may be a right ideal of $R_1$ and $I_2$ may be a left ideal of $R_2$ then how to define ideal structures. For this case we give the following definition.

**DEFINITION 6.1.6:** Let $R = R_1 \cup R_2$ be a biring. We say the set $I = I_1 \cup I_2$ is a mixed bi-ideal of R if $I_1$ is a right (left) ideal of $R_1$ and $I_2$ is a left (right) ideal of $R_2$.

Thus we see only in case of birings we can have the concept of mixed bi-ideals i.e. an ideal simultaneously having a section to be a left ideal and another section to be right ideal.

We give an example of right bi-ideal of a biring.

**Example 6.1.4:** Let $R = R_1 \cup R_2$ be a biring, where

$$R_1 = \left\{ \begin{pmatrix} a & 0 \\ c & d \end{pmatrix} \middle| a, c, d \in Z \right\} \text{ and}$$

$$R_2 = \left\{ \begin{pmatrix} x & y \\ 0 & z \end{pmatrix} \middle| x, y, z \in Z \right\}$$

be rings under usual addition and multiplication respectively. $I = I_1 \cup I_2$ where

$$I_1 = \left\{ \begin{pmatrix} a & 0 \\ b & 0 \end{pmatrix} \middle| a, b \in Z \right\} \text{ and}$$

$$I_2 = \left\{ \begin{pmatrix} 0 & a \\ 0 & 0 \end{pmatrix} \middle| a \in Z \right\}$$

are left ideals of $R_1$ and $R_2$ respectively. Thus I is a left bi-ideal of R.

Now take $J = J_1 \cup J_2$ where

$$J_1 = \left\{ \begin{pmatrix} 0 & 0 \\ a & b \end{pmatrix} \middle| a, b \in z \right\}$$

is a right ideal of $R_1$ and

$$J_2 = \left\{ \begin{pmatrix} 0 & b \\ 0 & 0 \end{pmatrix} \middle| b \in Z \right\} \text{ is}$$

a right ideal of $R_2$ so J is a right bi-ideal of R.

Take $K = K_1 \cup K_2$ where

$$K_1 = \left\{ \begin{pmatrix} a & 0 \\ b & 0 \end{pmatrix} \middle| a, b \in Z \right\} \text{ and}$$

$$K_2 = \left\{ \begin{pmatrix} 0 & b \\ 0 & 0 \end{pmatrix} \middle| b \in Z \right\} \text{ are}$$

left and right ideals of $R_1$ and $R_2$ respectively; so K is a mixed bi-ideal of R.

**DEFINITION 6.1.7:** *Let $R = R_1 \cup R_2$ be a biring. A bi-ideal $I = I_1 \cup I_2$ is called a maximal bi-ideal of R if $I_1$ is a maximal ideal of $R_1$ and $I_2$ is a maximal ideal of $R_2$. Similarly we can define the concept of minimal bi-ideal; as $J = J_1 \cup J_2$ is a minimal bi-ideal, if $J_1$ is a minimal ideal of $R_1$ and $J_2$ is a minimal ideal of $R_2$. It may happen in a bi-ideal I of a ring $R = R_1 \cup R_2$ that one of $I_1$ or $I_2$ may maximal or minimal then what do we call the structure $I = I_1 \cup I_2$. We call $I = I_1 \cup I_2$, a bi-ideal in which only $I_1$ or $I_2$ is maximal ideal as quasi maximal bi-ideal. Similarly we can define quasi minimal bi-ideal.*

**DEFINITION 6.1.8:** *Let $R = R_1 \cup R_2$ and $S = S_1 \cup S_2$ be two birings. We say a map $\phi$ from R to S is a biring homomorphism if $\phi = \phi_1 \cup \phi_2$ where $\phi_1 = \phi / R_1$ from $R_1$ to $S_1$ is a ring homomorphism and $\phi_2 = \phi / R_2$ is a map from $R_2$ to $S_2$ is a ring homomorphism.*

*We for notational convenience denote by $\phi = \phi_1 \cup \phi_2$ through this '$\cup$' is not the set theoretic union. We define for the homomorphism $\phi: R \to S$, where $R = R_1 \cup R_2$ and $S = S_1 \cup S_2$ are brings the kernel of the homomorphism $\phi$ as ker $\phi$ = ker $\phi_1 \cup$ ker $\phi_2$; here ker $\phi_1 = \{a_1 \in R_1 \mid \phi_1(a_1) = 0\}$ and ker $\phi_2 = \{a_2 \in R_2 \mid \phi_2(a_2) = 0\}$ i.e. ker $\phi = \{a_1 \in R_1, a_2 \in R_2 \mid \phi_1(a_1) = 0 \text{ and } \phi_2(a_2) = 0\}$.*

**THEOREM 6.1.2:** *Let $R = R_1 \cup R_2$ and $S = S_1 \cup S_2$ be two birings with $\phi$ a biring homomorphism from R to S then ker $\phi$ = ker $\phi_1 \cup$ ker $\phi_2$ is a bi-ideal of the biring R.*

*Proof:* Straightforward and hence left for the reader to prove.

**Remark:** The notion of biring isomorphism and automorphism are extendable as in case of ring homomorphism.

**DEFINITION 6.1.9:** *A biring $R = R_1 \cup R_2$ is said to be a bifield if $(R_1, +, \bullet)$ and $(R_2, +, \bullet)$ are fields. If R is finite i.e. $|R| < \infty$ we call R a finite field. If the characteristic of both $R_1$ and $R_2$ are finite then we say $R = R_1 \cup R_2$ is a bifield of finite characteristic.*

*If in R = $R_1 \cup R_2$ one of $R_1$ or $R_2$ is a field of characteristic 0 and one of $R_1$ or $R_2$ is of some finite characteristic we do not associate any characteristic with it.*

*If both $R_1$ or $R_2$ in R = $R_1 \cup R_2$ is of zero characteristic then we say R is a field of characteristic zero.*

Thus unlike in fields we see in case of bifields we can have characteristic p, p a prime or zero characteristic or no characteristic associated with it.

***Example 6.1.5:*** Let R = $R_1 \cup R_2$ where $R_1 = Z_{11}$ and $R_2 = Q$ then the bifield R has no characteristic associated with it.

***Example 6.1.6:*** Let R = $R_1 \cup R_2$, where $R_1 = Q(\sqrt{2})$ and $R_2 = Q(\sqrt{3})$ then R is a bifield of characteristic zero.

***Example 6.1.7:*** Let R = $Z_2 \cup Z_3$ where $Z_2$ and $Z_3$ are prime fields of characteristic 2 and 3 respectively then R is a bifield of finite characteristic.

**DEFINITION 6.1.10:** *Let F = $F_1 \cup F_2$ be a bifield we say a proper subset S of F to be a sub-bifield if S = $S_1 \cup S_2$ and $S_1$ is a subfield of $F_1$ and $S_2$ is a subfield of $F_2$. If the bifield has no proper sub-bifield then we call F a prime bifield.*

***Example 6.1.8:*** Let F = $F_1 \cup F_2$ where $F_1 = Z_7$ and $F_2 = Q$. Clearly F is a prime bifield as F has no proper sub-bifields. Now we proceed on to define polynomial birings.

**DEFINITION 6.1.11:** *Let R = $R_1 \cup R_2$ where R is either a bifield or a biring which is commutative with unit. The polynominal biring R [x] = $R_1$ [x] $\cup$ $R_2$ [x]. Here $R_1$ [x] and $R_2$ [x] are polynomial rings over the rings or fields in the same indeterminate x. If is easily verified that polynomial biring is a biring.*

**DEFINITION 6.1.12:** *Let D = $D_1 \cup D_2$ where $D_1$ and $D_2$ are integral domains then we call D a bidomain.*

***Example 6.1.9:*** Let D = Z $\cup$ Q[x]. Clearly D is a bidomain.

**THEOREM 6.1.3:** *Let R = $R_1 \cup R_2$ where $R_1$ and $R_2$ are fields; then the polynomial ring R[x] is an bidomain.*

*Proof*: The result is straightforward hence left for the reader as an exercise.

Several theorems about polynomial rings can also be obtained for polynomial birings. The concept of degree of a polynomial in birings is defined only as in the case of polynomial rings.

Let p(x) ∈ R [x] where, R [x] = $R_1$ [x] $\cup$ $R_2$ [x], p(x) ∈ $R_1$ [x] or $R_2$ [x] or it may be in R [x]. Let I be the ideal generated by p(x) then I in general is not a bi-ideal of R [x] for it may or may not be represented as the union of two ideals, so we define bi-ideals in polynomial birings as follows:

**DEFINITION 6.1.13:** *Let $R[x] = R_1[x] \cup R_2[x]$ be a polynomial biring. Let p(x) be a polynomial in R[x]. We say p(x) generates a bi-ideal I, if I can be written as $I_1 \cup I_2$, where $I_1 = I \cap R_1[x]$ and $I_2 = I \cap R_2[x]$ are proper ideals generated by p(x) in $R_1[x]$ and $R_2[x]$ respectively. Clearly $I_1 \neq (0)$ and $I_2 \neq (0)$. $I_1 \subseteq R_1[x]$ and $I_2 \subseteq R_2[x]$, for other wise there does not exist any bi-ideal generated by p(x).*

*Example 6.1.10:* Let $R[x] = R_1[x] \cup R_2[x]$ be a polynomial biring over the bifield $R = R_1 \cup R_2$ consider $R_1$ a ring of characteristic two given by $R_1 = \{0, 1, a, b \mid a^2 = 1, b^2 = 0, ab = ba = b\}$ and $R_2 = Z_3 = \{0, 1, 2\}$.

We have $R[x] = R_1[x] \cup R_2[x]$. Let $p(x) = x^2 + 1 \in R[x]$. The ideal I generated by p(x) is $I = \langle p(x) \rangle = \{(x^2 + 1), 0, a(x^2 + 1), b(x^2 + 1), p(x)(x^2 + 1)\} \cup \{0, x^2 + 1, 2(x^2 + 1), p(x)(x^2 + 1)\}$.

Thus we see I is a bi-ideal of the polynomial biring.

It is important to note that unlike the polynomial rings in which every polynomial will generate an ideal we see in case of polynomial birings a polynomial may or may not generate a polynomial bi-ideal.

**THEOREM 6.1.4:** *Let $R = R_1 \cup R_2$ be a commutative biring with unit or a bifield. R[x] be the polynomial biring. Then every polynomial p(x) in R[x] need not generate a polynomial bi-ideal.*

*Proof*: By an example $R = Z_2 \cup Z_3$; $R[x] = Z'_2[x] \cup Z'_3[x]$. Here $Z'_2[x]$ denotes the set of all polynomials in x of degree less than or equal to 7 and $Z'_3[x]$ = set of all polynomials of degree less than or equal to 5.

Let $x^2 + 1 \in R[x]$; the bi-ideal generated by $x^2 + 1$ is given by $I = \{0, (x^2 + 1), p(x)(x^2 + 1)\} \cup \{0, x^2 + 1, 2(x^2 + 1), p(x)(x^2 + 1)\} = I_1 \cup I_2$.

Since we have $I_1 \subset I_2$ so I is not a polynomial bi-ideal. Hence the claim.

**THEOREM 6.1.5:** *Let $R[x] = R_1[x] \cup R_2[x]$ be a polynomial biring. A polynomial $p(x) \in R[x]$ generates a bi-ideal if and only if I can be written as $I_1 \cup I_2$, where $I_1 = I \cap R_1[x]$ and $I_2 = I \cap R_2[x]$ are proper ideals generated by p(x) in $R_1[x]$ and $R_2[x]$ respectively.*

*Proof*: Direct; hence left as an exercise to the reader.

**DEFINITION 6.1.14:** *Let $R[x] = R_1[x] \cup R_2[x]$ be a polynomial biring. A polynomial $p(x) \in R[x]$ is said to be reducible in R[x] if and only if it is reducible in both $R_1[x]$ and $R_2[x]$. It is said to be irreducible if and only if it is irreducible in both $R_1[x]$ and $R_2[x]$. Unlike in case of polynomial rings, we have in case of polynomial birings we have for any $p(x) \in R[x] = R_1[x] \cup R_2[x]$ if p(x) is reducible in one of $R_1[x]$ or $R_2[x]$ and not in the other, we cannot of conclude anything about reducibility or irreducibility which is yet another marked difference between the polynomial rings and the polynomial birings.*

The following example gives us a polynomial which is reducible nor irreducible.

*Example 6.1.11:* Consider the polynomial biring $R[x] = Z_2[x] \cup Z_3[x]$ and the polynomial $x^2 +1 \in R[x]$. $\langle x^2 + 1 \rangle$ is reducible in $Z_2[x]$ but is irreducible on in $Z_3[x]$. Thus the polynomial $p(x) = x^2 +1$ is neither reducible nor irreducible in $R[x]$.

**DEFINITION 6.1.15:** *Let $R = R_1 \cup R_2$ be a biring we say R is a simple biring if R has no nontrivial bi-ideals; i.e. the only bi-ideals of R are (0) and R.*

**DEFINITION 6.1.16:** *Let $R = R_1 \cup R_2$ be a biring. We say an element $x \in R$, $(x \neq 0)$ is a zero divisor in R if we have a $y \in R$ $(y \neq 0)$ with $xy = 0$.*

*A non-commutative biring, which has no zero divisors will be called as a bidivision ring.*

*Example 6.1.12:* Let $R = Z_3 \cup Q'$ where $Q'$ is the ring of quaternions and $Z_3 = \{0, 1, 2\}$ prime field of characteristic 3. Clearly R is a bidivision ring which is not a bifield.

**DEFINITION 6.1.17:** *Let $R = R_1 \cup R_2$ be a biring an element $x \in R$ is an idempotent if $x^2 = x$ and we say an element $a \in R$ is nil potent if $a^n = 0$ for some $n \geq 2$.*

Similarly we can define semiidempotent in bi-ideals as follows:

**DEFINITION 6.1.18:** *Let $R = R_1 \cup R_2$ be a biring. We say an element $\alpha \in R$ ($\alpha \neq 0$ or $\alpha \neq 1$) is said to be a semiidempotent of R if $\alpha \notin$ {the ideal generated by $\alpha^2 - \alpha$} or the ideal generated by $\alpha^2 - \alpha$ is whole of R.*

One of the very interesting concepts which one can define using birings is group birings, semigroup birings, loop birings and groupoid birings. Here on wards by a biring R we mean only a biring, which is commutative, and one which contains the identity element 1 or a bifield.

**DEFINITION 6.1.19:** *$R = R_1 \cup R_2$ be a commutative biring with 1 or a bifield. Let G be a group. The group biring $RG = R_1G \cup R_2G$, where $R_1G$ and $R_2G$ are group rings of the group G over the rings $R_1$ and $R_2$ respectively. Thus RG is called the group biring.*

Almost all the questions studied in case of group rings can also be studied in case of group birings.

**THEOREM 6.1.6:** *Let $R = R_1 \cup R_2$ be a biring of characteristic 0 and G be a group having elements of finite order. The group biring RG has nontrivial divisors of zero.*

*Proof*: Let $1 \neq g \in G$ with $g^n = 1$ where $n > 1$. Let RG be the group biring. Take $\alpha = (1 - g)$ and $\beta = 1 + g + g^2 + \ldots + g^{n-1}$. Clearly both $\alpha, \beta \in RG$; we see $\alpha\beta = 0$. Thus the group biring RG has nontrivial divisors of zero.

Similarly the study of idempotents, units, semiidempotents, nilpotents can be carried out as a matter of routine in case of group birings.

It is easily verified that all group birings are also birings.

***Example 6.1.13:*** Let $R = Z_3 \cup Z_4$ be a biring and $G = \{g \mid g^2 = 1\}$ be the cyclic group of order 2. $RG = Z_3G \cup Z_4 G$ is the bigroup ring.

This has nontrivial units, idempotents and zero divisors. A natural question would be suppose we have a biring R which is the union of two grouprings then what will its structure be called. Suppose the group rings uses the same group G then we can call the biring as group biring; if they use two distinct groups then we only call them as birings.

***Example 6.1.14:*** Let $R = QS_3 \cup Z_2S_4$ be the group ring. We only say R is biring and not as a group biring.

Now what is the advantage of birings or bigroups? The chief advantage of it is as follows:

Suppose $S_3 = \{1, p_1, p_2, p_3, p_4, p_5\}$ be the symmetric group of degree 3. Consider the subgroups of $S_3$ say $H_1 = \{1, p_2\}$ and $H_2 = \{1, p_4, p_5\}$. We have the set $H = H_1 \cup H_2$ this is only a set for if we have not introduced the notion of bigroups.

So $H = H_1 \cup H_2$ is a bigroup i.e. union of two subgroups is not a group but a bigroup. Thus we have the nice theorem, which in general holds for all algebraic structures be it groups, loops, semigroups, groupoids, rings or semirings.

**THEOREM 6.1.7:** *Let G be any group $H_1$ and $H_2$ be any two distinct subgroups of G i.e. $H_1 \not\subset H_2$ and $H_2 \not\subset H_1$ then $H = H_1 \cup H_2$ is a bigroup.*

*Proof:* Let $H = H_1 \cup H_2$. Clearly as $H_1$ and $H_2$ are groups under the same operation so we have H to be a bigroup.

**Remark:** It is very important to note that if R is a ring or a semiring or a loop or a semigroup or a groupoid having two distinct subrings or subsemirings or subloops or subsemigroups or subgroupoids such that they are not related then we have their union to be birings or bisemirings or biloops or bisemigroups or bigroupoids respectively. By the term not related we mean if $R_1$ and $R_2$ are two substructure we have $R_1 \not\subset R_2$ or $R_2 \not\subset R_1$ but we can have $R_1 \cap R_2 \neq \phi$.

***Example 6.1.15:*** Let $G = \langle g \mid g^{12}=1 \rangle$ be the cyclic group of order 12. Take $H_1 = \{g^2, g^4, g^6, g^8, g^{10}, 1\}$ and $H_2 = \{1, g^3, g^6, g^9\}$. We see $H = H_1 \cup H_2$ is not a group only just a set but H is bigroup, we see $H_1 \not\subset H_2$ and $H_2 \not\subset H_1$ but $H_1 \cap H_2 = \{1, g^6\}$.

Thus the concept of bigroups enables us to recognize the structure of union of two subgroups as a bigroup. One of the unsolved question would be, suppose G is a group and $R = R_1 \cup R_2$ be a biring when is the bigroup ring a bidomain or a bidivision ring. The answer to this question still remains open. But we can give partial answers.

**THEOREM 6.1.8:** *Let $R = R_1 \cup R_2$ be a biring and G any group. The group biring RG is a bidomain if and only if both the group rings $R_1G$ and $R_2G$ are integral domains or division rings i.e. they do not have zero divisors.*

*Proof*: Straightforward by using the results in group rings.

**Note:** It is pertinent to mention here that the group biring will be a bidomain when in the biring $R = R_1 \cup R_2$, $R_1$ and $R_2$ are both fields of characteristic zero. Further the group G taken as torsion free abelian group; then certainly the group biring is a bidomain.

**DEFINITION 6.1.20:** *A biring $R = R_1 \cup R_2$ is semiprime.*

  i.   *If each of the $R_i$'s are semiprime i = 1, 2 and*
  ii.  *If R has prime ideals common in $R_1$ and $R_2$, then intersection is zero.*

*The condition (ii) becomes superfluous if R has no ideals common in both $R_1$ and $R_2$.*

Characterization of semiprime birings is not an easy task.

**THEOREM 6.1.9:** *Let $R = R_1 \cup R_2$ be a biring. If R has no ideals common in $R_1$ and $R_2$ and if $R_1$ and $R_2$ are semiprime then R is a semiprime biring.*

*Proof*: Follows from simple arguments, hence left for the reader to solve.

So we do not conclude the biring $R = R_1 \cup R_2$ to be semiprime even if $R_1$ and $R_2$ are semiprime rings.

A natural question would be if we take the equivalent formulation of a semiprime ring viz. "A ring R is semiprime if and only if R contains no non-zero ideal with square zero".

So we define semiprime birings as follows which can be taken as an equivalent formation of the definition.

**DEFINITION 6.1.21:** *Let $R = R_1 \cup R_2$ be a biring we say R is a semiprime biring if and only if R contain no non-zero bi-ideal with square zero.*

Using this definition we give a nice result about group rings, which are semiprime.

**THEOREM 6.1.10:** *Let $K = K_1 \cup K_2$ be the bifield of characteristic 0. Then KG the group biring is semiprime where G is any group.*

*Proof*: Follows from the fact that if K is a field of characteristic zero and G any group then the group ring is semiprime. So $KG = K_1G \cup K_2G$, is the biring in which both $K_1G$ and $K_2G$ has no non-zero ideal with square zero. Hence the claim.

**THEOREM 6.1.11:** *Let G be a group having an abelian subgroup A with $[G : A] = n < \infty$. Let $K = K_1 \cup K_2$ be a bifield. Then $KG = K_1 G \cup K_2 G$ satisfies the standard polynomial identity of degree 2n only when $K_1 G \cap K_2 G = G$.*

*Proof*: Follows from the results of [41, 42].

**THEOREM 6.1.12:** *Let $K = K_1 \cup K_2$ be a bifield of characteristic 0. KG be the group biring which is semiprime, satisfying a polynomial identity of degree n. Then G has an abelian subgroup A with $[G: A] \leq n! \, u(n)$.*

*Proof*: Follows from [41, 42].

**Note:** $u(n) = (n!)^{n(n!)^2}$ is the notation used in the above theorem. Several interesting results about semiprime group birings and those group birings which satisfy the polynomial identity can be studied and analyzed. There is an ocean of results, which can be got in this direction. Now we proceed on to define semisimple birings.

**DEFINITION 6.1.22:** *Let $R = R_1 \cup R_2$ be a biring, a non-empty set M is said to be a R-bimodule (or a bimodule over the biring R) if M is an abelian bigroup, under addition '+' say with $M = M_1 \cup M_2$ such that for every $r_1 \in R_1$ and $m_1 \in M_1$ there exists an element $r_1 m_1$ in $M_1$ and for every $r_2 \in R_2$ and $m_2 \in M_2$ there exists an element $r_2 m_2$ in $M_2$ subject to*

1. $r_1 (a_1+b_1) = r_1 a_1 + r_1 b_1$ ; $r_1 \in R_1$ and $a_1, b_1 \in M_1$.
2. $r_2 (a_2+b_2) = r_2 a_2 + r_2 b_2$ ; $r_2 \in R_2$ and $a_2, b_2 \in M_2$.
3. $r_1 (s_1 a_1) = (r_1 s_1) a_1$ for $r_1, s_1 \in R_1$ and $a_1 \in M_1$.
4. $r_2 (s_2 a_2) = (r_2 s_2) a_2$ for $r_2, s_2 \in R_2$ and $a_2 \in M_2$.
5. $(r_1 + s_1)a_1 = r_1 a_1 + s_1 a_1$, $r_1, s_1 \in R_1$ and $a_1 \in M_1$.
6. $(r_2 + s_2)a_2 = r_2 a_2 + s_2 a_2$ where $r_2, s_2 \in R_2$ and $a_2 \in M_2$.

*In short if M is a bigroup ($M = M_1 \cup M_2$) under '+' then $M_1$ is a $R_1$-module and $M_2$ is a $R_2$-module then we say M is a R-bimodule over the biring $R = R_1 \cup R_2$.*

*If the biring R has unit 1 and if $1.m_1 = m_1$ and $1.m_2 = m_2$ for every $m_1 \in M_1$ and $m_2 \in M_2$ then we call M a unital R-bimodule. Thus the concept of R-bimodule forces both the structures to be bistructures i.e. we demand the group should be a bigroup and also the ring must be a biring then only we can speak of a R-bimodule.*

*Let M be a R-bimodule an additive sub-bigroup, A of M i.e. $A = A_1 \cup A_2$ where $A_1$ is a subgroup of $M_1$ and $A_2$ is a subgroup of $A_2$ is called the sub-bimodule of the bimodule M if when ever $r_1 \in R_1$ and $r_2 \in R_2$ and $a_1 \in A_1$ and $a_2 \in A_2$ we have $r_1 a_1 \in A_1$ and $r_2 a_2 \in A_2$.*

*A bimodule M is cyclic if there is an element $m_1 \in M_1$ and $m_2 \in M_2$ such that for every $m \in M_1$ is of the form $m = r_1 m_1$ where $r_1 \in R_1$ and for every $m' \in M_2$ is of the form $m' = r_2 m_2$ where $r_2 \in R_2$. Thus cyclic bimodules is nothing but bicyclic groups, that is the bigroup $G = G_1 \cup G_2$ is a bicyclic group if both $G_1$ and $G_2$ are cyclic groups.*

Several other analogous results about bimodules can be defined in case of birings and studied. It is left as a piece of research for the reader.

**DEFINITION 6.1.23:** *Let $R = R_1 \cup R_2$ be a biring. We say R is a biregular biring if for any $\alpha \in R$ there exists $\beta \in R$ with $\alpha \beta \alpha = \alpha$. If both $R_1$ and $R_2$ are regular rings then the biring R is biregular.*

In view of this definition we have the following theorem:

**THEOREM 6.1.13:** *Let $R = R_1 \cup R_2$ be a biring, $R_1 \cap R_2 = \phi$ or $\{0\}$ or $\{1\}$ or $\{0, 1\}$ if it has common unit. The biring R is biregular if and only if every finitely generated right ideal is generated by an idempotent.*

*Proof*: Exactly as in case of rings.

Several other interesting results about biregular birings and biregular group birings can be obtained. All these type of research is left for the reader.

Now leaving all the analogous ideas which can be developed in case of birings to the reader we proceed on to define the concept of semigroup birings in the following.

**DEFINITION 6.1.24:** *Let R be a biring; $R = R_1 \cup R_2$. S be any semigroup. We call RS the semigroup biring if $RS = R_1S \cup R_2S$ where $R_1S$ and $R_2S$ are semigroup rings of the semigroup S over $R_1$ and $R_2$ respectively.*

**Note:** If R is a biring such that $R = R_1 \cup R_2$ where $R_1$ and $R_2$ are any two semigroup rings; that is $R_1 = R'_1 S_1$ and $R_2 = R'_2 S_2$ two different semigroup rings over two rings $R'_1$ and $R'_2$ respectively; then we do not call R a semigroup biring it will be known as just a biring.

*Example 6.1.16:* Let $R = Z_2, S_1 \cup QS_2$ where $S_1 = $ {the semigroup under multiplication of the set of modulo integers $Z_6$} and $S_2 = $ {S(4), the set of all maps of the set (1, 2, 3, 4) to itself is the symmetric semigroup}; then R is not a semigroup biring but only a biring.

*Example 6.1.17:* Let $R = Z_3 \cup Z$. Let S(6) be the symmetric semigroup. Then $RS(6) = Z_3S(6) \cup ZS(6)$ is the semigroup biring of the semigroup S(6) over R. Clearly RS(6) is also a biring. Further $Z_3S(6)$ and $ZS(6)$ are semigroup rings hence rings.

All properties defined and studied for group birings wherever relevant can be defined and extended in case of semigroup birings. The task of studying, characterizing and obtaining very many interesting results about semigroup birings are left for the reader.

Now we proceed on to define yet another new type of birings using bigroups.

**DEFINITION 6.1.25:** *Let $G = G_1 \cup G_2$ be a bigroup. R any commutative ring with unit or a field. The bigroup ring RG consisting of $RG = RG_1 \cup RG_2$ where $RG_1$ and $RG_2$ are group rings of the groups $G_1$ and $G_2$ over the ring R.*

**THEOREM 6.1.14:** *Let $G = G_1 \cup G_2$ be a bigroup and R any ring. The bigroup ring RG is a biring.*

*Proof*: Follows from the very definition as $RG = RG_1 \cup RG_2$ so both $RG_1$ and $RG_2$ are rings so RG is a biring.

*Example 6.1.18*: Let Z be the ring of integers and $G = \{1, g^2, g^4, g^6, g^8, g^{10}, g^3, g^6, g^9 \mid g^{12} = 1\}$. Clearly G is a bigroup. The bigroup ring ZG is a biring $ZG = ZG_1 \cup ZG_2$ where $G_1 = \{1, g^2, g^4, g^6, g^8, g^{10}\}$ and $G_2 = \{1, g^3, g^6, g^9\}$. It is easily verified ZG is a biring.

The study of bigroup rings is very new. No one has, till date introduced or studied this concept so we may not be able to give any literature in this direction but leave it for the reader to develop notions about bigroup rings using concepts from group rings and birings.

**THEOREM 6.1.15:** *Let K be any field; G any finite bigroup. Then the bigroup ring KG has nontrivial divisors of zero.*

*Proof*: Straightforward; hence left for the reader to prove.

**THEOREM 6.1.16:** *Let G be a group having proper subgroups H and K such that $H \not\subset K$ or $K \not\subset H$; and R be a ring. Then the group ring RG has a proper subset, which is a biring.*

*Proof*: Given RG is a group ring of the group G over the ring R. It is further given G has two proper subgroups H and K such that $H \not\subset K$ or $K \not\subset H$. Consider RH and RK. Clearly $RH \subset RG$ and $RK \subset RG$ and RH and RK are subrings of RG. Now let $RH \cup RK = S$, S is a proper subset of RG and S is a biring. Hence the claim.

**THEOREM 6.1.17:** *Let K be any field of characteristic 0 and $G = G_1 \cup G_2$ be a bigroup. The bigroup ring KG is semiprime.*

*Proof*: We know $KG = KG_1 \cup KG_2$. By theorems of [41, 42], $KG_1$ and $KG_2$ are semiprime, so KG the bigroup ring is a semiprime biring.

**THEOREM 6.1.18:** *Let K be any field and G a finite bigroup i.e. $G = G_1 \cup G_2$. Then the bigroup ring KG is a union of two Artinian ring.*

*Proof*: Follows from the fact if G is finite then KG is artinian. Now $G = G_1 \cup G_2$ since G is finite so are the groups $G_1$ and $G_2$. Thus the bigroup ring $KG = KG_1 \cup KG_2$ is artinian as $KG_1$ and $KG_2$ are artinian. The theorem is true.

**THEOREM 6.1.19:** *Let $G = G_1 \cup G_2$ be a bigroup such that $G_1$ has an abelian subgroup $A_1$ such that $[G_1 A_1] = n < \infty$ and $G_2$ has an abelian subgroup $A_2$ such that $[G_2, A_2] = n < \infty$. Then KG the bigroup ring, where K is a field satisfies the standard polynomial identity.*

*Proof*: Follows from the results of [41, 42].

**Note:** We demand both $[G_1; A_1] = n$ and $[G_2; A_2] = n$ with $n < \infty$. Recall $\Delta = \Delta(G) = \{x \in G \mid [G : C_G(x)] < \infty\}$ where G is any group and $C_G(x)$ is the centralizer of x in G. Now if G is a bigroup say $G = G_1 \cup G_2$ then we denote $\Delta(G) = \Delta(G_1) \cup \Delta(G_2)$, clearly $\Delta(G)$ is a bigroup as $\Delta(G_1)$ and $\Delta(G_2)$ are normal subgroups of $G_1$ and $G_2$ respectively. In fact $\Delta(G)$ is the sub-bigroup of G or to be more precise $\Delta(G)$ is a normal sub-bigroup of G.

**THEOREM 6.1.20:** *The following are equivalent; here $G = G_1 \cup G_2$ is a bigroup and K any field*

   i.   *K[G] is prime.*
   ii.  *$\Delta(G) = \Delta(G_1) \cup \Delta(G_2)$ is torsion free abelian bigroup.*
   iii. *$G = G_1 \cup G_2$ has no nonidentity finite normal sub-bigroup.*

*Proof*: Let $G = G_1 \cup G_2$ be a bigroup, K a field and KG the bigroup ring. Clearly $KG = KG_1 \cup KG_2$, we know if $KG_1$ and $KG_2$ are prime then so is KG. Let $\Delta(G_1)$ and $\Delta(G_2)$ be normal subgroups so that $\Delta(G)$ is the normal sub-bigroup. Then by the results of [16] the three statements are equivalent.

Several results in this direction can be stated and proved as in the case of group rings. So study of bigroup rings itself can be treated as a separate research.

Now we proceed on to define the concept of bisemigroup rings.

**DEFINITION 6.1.26:** *Let K be a field or a commutative ring with 1 and S be a bisemigroup i.e. $S = S_1 \cup S_2$. Then $KS = KS_1 \cup KS_2$ is the bisemigroup ring which is the set theoretic union of the semigroup rings $KS_1$ and $KS_2$.*

**THEOREM 6.1.21:** *Let $KS = KS_1 \cup KS_2$ be the bisemigroup ring of the bisemigroup $S = S_1 \cup S_2$ over the field K or the ring K. Then KS is a biring.*

*Proof*: Direct by the very definition; hence left for the reader to prove.

It is important to mention here that, as semigroups are the generalized structure of groups so are the bisemigroups, the generalized structure of bigroups. So we request the reader to develop the new concept of bisemigroup rings analogous to bigroup rings. We only propose some interesting research problems for the reader in chapter ten.

Now we proceed on to define still a newer concept called bigroup birings. So far no algebraist have defined or discussed about this new notion. The author feels that this new concept would certainly find a lot of application in several fields like finite automata, control theory and in linear logic.

**DEFINITION 6.1.27:** *Let $G = G_1 \cup G_2$ be a bigroup and $R = R_1 \cup R_2$ be a biring. The bigroup biring $RG = (R_1 G_1 \cup R_1 G_2 \cup R_2 G_1 \cup R_2 G_2)$ consists of union of grouprings. Thus a bigroup biring has four grouprings.*

What is its structure? What shall we call bigroup biring? To this we define a new algebraic structure as follows:

**DEFINITION 6.1.28:** *Let $(R, +, \bullet)$ be a ring, we call R a quad ring if $R = R_1 \cup R_2 \cup R_3 \cup R_4$, where each $R_i$ is a ring under the same operation '+' and '$\bullet$' of R (i = 1, 2, 3, 4).*

**DEFINITION 6.1.29:** *Let $(G, \bullet)$ be a group we call G a quad group if $G = G_1 \cup G_2 \cup G_3 \cup G_4$ where each $G_i$ is a group under the operation of G for i = 1, 2, 3, 4.*

In view of these definitions we have the following theorem:

**THEOREM 6.1.22:** *Let $G = G_1 \cup G_2$ be a bigroup and $R = R_1 \cup R_2$ be a biring. The bigroup biring RG is a quad ring.*

*Proof*: Obvious by the very definition.

Now we can give an equivalent formulation of quad groups and quad rings using bigroups and birings.

**DEFINITION 6.1.30:** *Let $G = G_1 \cup G_2$ where $(G, \bullet)$ is set with a closed binary operation. If each of $G_1$ and $G_2$ are bigroups then we call G a quad group.*

**DEFINITION 6.1.31:** *Let $(R, +, \bullet)$ where $R = R_1 \cup R_2$; we call R a quad ring if $R_1$ and $R_2$ are birings under the same operations '+' and '$\bullet$'.*

Now we proceed on to define quad ring using two birings.

**DEFINITION 6.1.32:** *Let $(R, +, \bullet)$ be a non-empty set with two binary operations. If $R = R_1 \cup R_2$ where $R_1$ and $R_2$ are two proper subset of R such that $R_1$ and $R_2$ under the operations of R are birings. i.e. $R_1$ is a biring and $R_2$ is a biring then we call R a quad ring.*

One can study all the properties and notions studied for birings and rings in case of quad rings. Similarly study of bigroup birings boils down to the study of quad rings.

*Example 6.1.19*: Let $G = \{1, g^{105}, g^{70}, g^{140}, g^{42}, g^{84}, g^{126}, g^{168}, g^{30}, g^{60}, g^{90}, g^{120}, g^{150}, g^{180}, / g^{210} = 1\}$ be a set under multiplication. $G = G_1 \cup G_2 \cup G_3 \cup G_4$ where

$$\begin{aligned} G_1 &= \{1, g^{105}\}, \\ G_2 &= \{1, g^{70}, g^{140}\}, \\ G_3 &= \{1, g^{42}, g^{84}, g^{126}, g^{168}\} \text{ and} \\ G_4 &= \{1, g^{30}, g^{60}, g^{90}, g^{120}, g^{150}, g^{180}\} \end{aligned}$$

each $G_i$ are groups under the operations of using following the condition $g^{210} = 1$. Thus G is a quad group.

It is to be noted that G is just a set not even closed under multiplication but G happen to be a nice quad group.

*Example 6.1.20*: Let R = {0, 105, 70, 140, 30, 60, 90, 120, 150, 180, 42, 84, 126, 168} the integers under '+' and '•' modulo 210. Clearly R is a quad ring as R = $R_1 \cup R_2 \cup R_3 \cup R_4$ where $R_1$ = {0, 105}, $R_2$ = {0, 70, 140}, $R_3$ = {0, 30, 60, 90, 120, 150, 180} and $R_4$ = {0, 42, 84, 126, 168} where each $R_i$ is a ring under addition and multiplication modulo 210. But clearly the set R is not even closed under product or sum.

Thus the study of bigroup birings is nothing but the study of quad rings. Hence the reader is expected to find innovative properties about these new algebraic structure. Also the notions of quad subrings, quad ideals, quad module etc can be defined analogous to birings or rings. Further we just define the notions of quad domain, quad division ring and quad fields.

**DEFINITION 6.1.33:** *Let (R, +, •) be a non-empty set endowed with two binary operations '+' and ' • '. Let R = $R_1 \cup R_2 \cup R_3 \cup R_4$. We say R is a quad domain if each ($R_i$, +, •) is an integral domain for i = 1, 2, 3, 4.*

*Example 6.1.21*: The quad ring given in example 6.1.19 is a quad domain. It can be easily verified.

**DEFINITION 6.1.34:** *Let (R, +, •) be a non-empty set. We say R = $R_1 \cup R_2 \cup R_3 \cup R_4$ is a quad division ring if each $R_i$ is a division ring, i = 1, 2, 3, 4.*

**DEFINITION 6.1.35:** *Let (R, +, •) be a non-empty set. If R = $R_1 \cup R_2 \cup R_3 \cup R_4$ where each $R_i$ is a field for i =1, 2, 3, 4. Then we say R is a quad field.*

Now as in case of birings and rings we can define the notion of polynomial quad rings as; R[x] = $R_1[x] \cup R_2[x] \cup R_3[x] \cup R_4[x]$ will be a polynomial quad ring only if each $R_i$ is a commutative ring with unit (we just call a quad ring R to be commutative if R = $R_1 \cup R_2 \cup R_3 \cup R_4$ then each $R_i$ is commutative for i = 1, 2, 3, 4 we say R is quad ring with unit if R has a common unit i.e. 1 ∈ R such that 1 • r = r • 1 = r for all r ∈ R and 1 ∈ $R_i$ for i = 1, 2, 3, 4.

Thus quad ideals in the polynomial quad ring can be defined and the notions of reducibility or irreducibility of the polynomials can be defined. Here it has become important for us to define reducibility in two ways.

**DEFINITION 6.1.36:** *Let R = $R_1 \cup R_2 \cup R_3 \cup R_4$ where each ($R_i$, +, •) is a ring; that is R be a quad ring we assume R to be a commutative quad ring with 1. Let R[x] = $R_1[x] \cup R_2[x] \cup R_3[x] \cup R_4[x]$ be the polynomial quad ring. Let p(x) ∈ R(x) we say p(x)is reducible if p(x) = (x-$\alpha_1$) …(x - $\alpha_t$) where $\alpha_1, \alpha_2, \alpha_3, …, \alpha_t$ ∈ $R_i$ for some i, i fixed i.e. all the roots belong to $R_i$ for some fixed i. We call p(x) ∈ R[x] as quasi reducible; if P(x) = (x-$\alpha_1$) …(x-$\alpha_t$) and $\alpha_1, …, \alpha_t$ ∈ R where all the $\alpha_i$'s may not be in a single ring $R_1$ or $R_2$ or $R_3$ or $R_4$. So in case of polynomial quad rings we have two types of reducibility, unlike other rings: viz reducible and quasi reducible. We say a polynomial p(x) ∈ R[x] is irreducible if all the $\alpha_1, \alpha_2, …, \alpha_t$ does not belong to one of the $R_i$'s .We say p(x) is not even quasi irreducible if the roots $\alpha_1, \alpha_2, …, \alpha_t$ does not belong to R.*

The interested reader can develop a complete notion of polynomial quad rings and polynomial quad birings. The last notion about these polynomial quad rings is that, if p(x) is a polynomial, the ideal generated by p(x) in general need not be a quad ideal.

Now we proceed on to define the concept of bisemigroup birings.

**DEFINITION 6.1.37:** *Let $(R, +, \bullet)$ be a biring with unit and $R = R_1 \cup R_2$ is assumed to be commutative. Let $(S, \bullet)$ be a bisemigroup i.e. $S = S_1 \cup S_2$ where each $S_i$ is a semigroup. The bisemigroup biring $RS = R_1S_1 \cup R_1S_2 \cup R_2S_1 \cup R_2S_2$ is union of the four semigroup rings $R_i S_j$, $i = 1, 2$ and $j = 1, 2$.*

Thus we have all bisemigroup birings are quad rings.

Almost all results pertaining to bigroup birings can also be extended to the case of bisemigroup birings with some simple modifications.

**PROBLEMS:**

1. Give an example of a biring of order 10.
2. Does there exist a biring of order 3? Justify your claim.
3. Find ideals in the biring $R = Z \cup \{0, 2, 4, 6, 8, 10, 12, 14$ addition under modulo 16$\}$.
4. Is $R = Z \cup 3Z$ a biring? Justify your answer.
5. Is $R = Q \cup Z$ a biring?
6. Show $R = \{A_{n \times n}\} \cup \{B^{n \times n}\}$ where $A_{n \times n}$ denotes the set of all $n \times n$ lower triangular matrices and $B^{n \times n}$ set of all upper triangular $n \times n$ matrices under usual matrix addition and matrix multiplication is a biring.
   Find
   a. Bi-ideals of R.
   b. Sub-birings of R.
   c. Zero divisors in R.
   d. Idempotents in R.
   e. Nilpotents in R.
7. Give an example of a group biring which is prime.
8. Does there exist a group biring, which is not semiprime? (Substantiate your answer).
9. Find an example of isomorphic group birings.
10. Give an example of a semigroup biring, which is not prime.
11. Give an example of a semigroup biring, which is semiprime?
12. Show by an example that there exists group birings which satisfies standard polynomial identity.
13. Find nontrivial ideals in the semigroup biring RS (5) where $R = Z_8 \cup Z_{12}$ and S(5) is the symmetric semigroup using 5 elements.
14. Give an example of a quad ring.
15. What is the smallest order of the quad ring?
16. Give an example of a bigroup biring which is not a quad domain.
17. Give an example of a bisemigroup biring, which is not a quad domain.

18. Give an example of a non-commutative quad ring.
19. Do all finite quad rings commutative? Justify your claim.

## 6.2 Non-associative birings

In this section we define the notion of non-associative birings. Study of birings is itself very new so the concept is non-associative birings is absent in literature. Here we define loop birings, groupoid birings and finally the concept of quad birings which are non-associative.

It is important to mention here that by no means we say that these birings are Lie bialgebras or Jordan bialgebras, we deal in this book purely in an algebraic vein so intermingling of category or topological concepts are totally absent in this book.

**DEFINITION 6.2.1:** *Let $R = R_1 \cup R_2$ with two binary operation '+' and '•'. R is said to be a non-associative biring if $(R_1, +, •)$ or $(R_2, +, •)$ is a non-associative ring 'or' not in the mutually exclusive sense for both $(R_1, +, •)$ and $(R_2, +, •)$ can be also non-associative rings; so what we demand is atleast one of $R_1$ or $R_2$ must be a non-associative ring for the biring to be a non-associative biring.*

*Example 6.2.1:* Let $R = (ZL \cup Q)$ under usual addition and product, where ZL is the loop ring of the loop giving by the following table:

| • | e | $g_1$ | $g_2$ | $g_3$ | $g_4$ | $g_5$ |
|---|---|---|---|---|---|---|
| e | e | $g_1$ | $g_2$ | $g_3$ | $g_4$ | $g_5$ |
| $g_1$ | $g_1$ | e | $g_3$ | $g_5$ | $g_2$ | $g_4$ |
| $g_2$ | $g_2$ | $g_5$ | e | $g_4$ | $g_1$ | $g_3$ |
| $g_3$ | $g_3$ | $g_4$ | $g_1$ | e | $g_5$ | $g_2$ |
| $g_4$ | $g_4$ | $g_3$ | $g_5$ | $g_2$ | e | $g_1$ |
| $g_5$ | $g_5$ | $g_2$ | $g_4$ | $g_1$ | $g_3$ | e |

Clearly R is a non-associative biring.

*Example 6.2.2:* Let R = {0, 2, 4, 6, 8, 10, addition and multiplication modulo 12} ∪ KG, (where K = 3, 6, 9, 0, addition and multiplication modulo 12) and G a groupoid given by the following table:

| • | $a_1$ | $a_2$ | $a_3$ | $a_4$ |
|---|---|---|---|---|
| $a_1$ | $a_1$ | $a_4$ | $a_3$ | $a_2$ |
| $a_2$ | $a_3$ | $a_2$ | $a_1$ | $a_4$ |
| $a_3$ | $a_1$ | $a_4$ | $a_3$ | $a_2$ |
| $a_4$ | $a_3$ | $a_2$ | $a_1$ | $a_4$ |

KG is the groupoid ring, which is clearly a non-associative ring.

Thus R = {0, 2, 4, 6, 8, 10} ∪ KG is a non-associative biring.

**DEFINITION 6.2.2:** *Let R be a non-associative biring. We say R is a commutative non-associative biring if $R = R_1 \cup R_2$ were both $R_1$ and $R_2$ are non-associative and commutative ring. We call R a non-associative biring with unit if $1 \in R$ such that $1 \bullet r = r \bullet 1 = r$ for all $r \in R$. We say R is finite if R has finite number of elements in them other wise infinite.*

*Example 6.2.3*: Let $R = ZG \cup ZL$ be the loop rings where both G and L are commutative loops given by the following tables:

The table for the loop $(G, *)$ is given by

| *    | e    | $a_1$ | $a_2$ | $a_3$ | $a_4$ | $a_5$ |
|------|------|-------|-------|-------|-------|-------|
| e    | e    | $a_1$ | $a_2$ | $a_3$ | $a_4$ | $a_5$ |
| $a_1$ | $a_1$ | e    | $a_4$ | $a_2$ | $a_5$ | $a_3$ |
| $a_2$ | $a_2$ | $a_4$ | e    | $a_5$ | $a_3$ | $a_1$ |
| $a_3$ | $a_3$ | $a_2$ | $a_5$ | e    | $a_1$ | $a_4$ |
| $a_4$ | $a_4$ | $a_5$ | $a_3$ | $a_1$ | e    | $a_2$ |
| $a_5$ | $a_5$ | $a_3$ | $a_1$ | $a_4$ | $a_2$ | e    |

The loop $(L, \bullet)$ is given by the following table:

| $\bullet$ | e | $g_1$ | $g_2$ | $g_3$ | $g_4$ | $g_5$ | $g_6$ | $g_7$ |
|---|---|---|---|---|---|---|---|---|
| e | e | $g_1$ | $g_2$ | $g_3$ | $g_4$ | $g_5$ | $g_6$ | $g_7$ |
| $g_1$ | $g_1$ | e | $g_5$ | $g_2$ | $g_6$ | $g_3$ | $g_7$ | $g_4$ |
| $g_2$ | $g_2$ | $g_5$ | e | $g_6$ | $g_3$ | $g_7$ | $g_4$ | $g_1$ |
| $g_3$ | $g_3$ | $g_2$ | $g_6$ | e | $g_7$ | $g_4$ | $g_1$ | $g_5$ |
| $g_4$ | $g_4$ | $g_6$ | $g_3$ | $g_7$ | e | $g_1$ | $g_5$ | $g_2$ |
| $g_5$ | $g_5$ | $g_3$ | $g_7$ | $g_4$ | $g_1$ | e | $g_2$ | $g_6$ |
| $g_6$ | $g_6$ | $g_7$ | $g_4$ | $g_1$ | $g_5$ | $g_2$ | e | $g_3$ |
| $g_7$ | $g_7$ | $g_4$ | $g_1$ | $g_5$ | $g_2$ | $g_6$ | $g_3$ | e |

Clearly the non-associative biring R is a commutative ring with unit.

The biring $R = ZG \cup ZL$ is a ring of infinite order. The concept of bi-ideal and sub-biring are defined as in case of birings.

Now we define some special properties about these non-associative rings.

**DEFINITION 6.2.3:** *Let $R = R_1 \cup R_2$ is a non-associative biring. We say proper subset S of R is said to be a sub-biring of R if $S = S_1 \cup S_2$ where $S_1$ is a subring of $R_1$ and $S_2$ is a subring of $R_2$.*

**DEFINITION 6.2.4:** *Let $R = R_1 \cup R_2$ be a non-associative biring. We say R is a subcommutative biring if R has atleast a proper sub-biring, which is commutative.*

The notion of subcommutativity has meaning only when R is not a commutative bring. We define the notion of strongly subcommutative biring in what follows:

**DEFINITION 6.2.5:** *Let $R = R_1 \cup R_2$ be a non-associative biring. If every sub-biring of R is commutative then we call R a strongly subcommutative biring.*

Here also the notion of strongly subcommutative has meaning only when R is a non-commutative biring.

**DEFINITION 6.2.6:** *Let $R = R_1 \cup R_2$ be a biring. We say R is a subassociative biring if R has a sub-biring $S = S_1 \cup S_2$ such that S is an associative biring. If every sub-biring of $R = R_1 \cup R_2$ of a non-associative biring is associative then we call R a strongly subassociative biring.*

*Example 6.2.4:* Let $R = Z_5 G \cup Z_5 L$ where G and L are loops given by the following tables and R is a non-associative biring. The table for the loop G is given below:

| *   | e   | $g_1$ | $g_2$ | $g_3$ | $g_4$ | $g_5$ |
|-----|-----|-------|-------|-------|-------|-------|
| e   | e   | $g_1$ | $g_2$ | $g_3$ | $g_4$ | $g_5$ |
| $g_1$ | $g_1$ | e   | $g_5$ | $g_4$ | $g_3$ | $g_2$ |
| $g_2$ | $g_2$ | $g_3$ | e   | $g_1$ | $g_5$ | $g_4$ |
| $g_3$ | $g_3$ | $g_5$ | $g_4$ | e   | $g_2$ | $g_1$ |
| $g_4$ | $g_4$ | $g_2$ | $g_1$ | $g_5$ | e   | $g_3$ |
| $g_5$ | $g_5$ | $g_4$ | $g_3$ | $g_2$ | $g_1$ | e   |

The table for the loop L

| *   | e   | $g_1$ | $g_2$ | $g_3$ | $g_4$ | $g_5$ | $g_6$ | $g_7$ |
|-----|-----|-------|-------|-------|-------|-------|-------|-------|
| e   | e   | $g_1$ | $g_2$ | $g_3$ | $g_4$ | $g_5$ | $g_6$ | $g_7$ |
| $g_1$ | $g_1$ | e   | $g_4$ | $g_7$ | $g_3$ | $g_6$ | $g_2$ | $g_5$ |
| $g_2$ | $g_2$ | $g_6$ | e   | $g_5$ | $g_1$ | $g_4$ | $g_7$ | $g_3$ |
| $g_3$ | $g_3$ | $g_4$ | $g_7$ | e   | $g_6$ | $g_2$ | $g_5$ | $g_1$ |
| $g_4$ | $g_4$ | $g_2$ | $g_5$ | $g_1$ | e   | $g_7$ | $g_3$ | $g_6$ |
| $g_5$ | $g_5$ | $g_7$ | $g_3$ | $g_6$ | $g_2$ | e   | $g_1$ | $g_4$ |
| $g_6$ | $g_6$ | $g_5$ | $g_1$ | $g_4$ | $g_7$ | $g_3$ | e   | $g_2$ |
| $g_7$ | $g_7$ | $g_3$ | $g_6$ | $g_2$ | $g_5$ | $g_1$ | $g_4$ | e   |

Clearly R is a biring, which is both non-associative and non-commutative and contain 1. $R = Z_5 G \cup Z_5 L$ is a strongly subcommutative and strongly subassociative biring; for the order of the loops, though is even have no subloops other than the subgroups of order 2. Further $Z_5$ has no subrings as $Z_5$ is the prime field of characteristic 5. Thus $R = Z_5G \cup Z_5L$ has all its sub-birings to be associative ones and all of its sub-birings are also commutative; clearly R is not a associative or a commutative biring.

We can have examples of just subcommutative and subassociative birings also. This task is left as an exercise for the reader.

**THEOREM 6.2.1:** *Every strongly subcommutative (or subassociative) biring is subcommutative (or subassociative) biring.*

*Proof*: The proof is straightforward hence left for the reader to prove.

**DEFINITION 6.2.7:** *Let $R = R_1 \cup R_2$ be a non-associative biring we call R a Moufang biring if any one of the following identities is satisfied by all elements of $R_1$ and by all elements of $R_2$.*

    i.    *(xy) (zx) = (x(yz))x.*
    ii.    *((xy)z) y = x(y (zy)).*
    iii.    *x(y (xz)) = ((xy) x)z;*

*for all x, y, z ∈ $R_1$ and for all x, y, z ∈ $R_2$.*

**DEFINITION 6.2.8:** *Let $R = R_1 \cup R_2$, we call the biring a Bruck biring if (x (yx)) z = x(y (xz) ) and $(xy)^{-1} = x^{-1}y^{-1}$ for all x, y, z in $R_1$ and for all x, y, z in $R_2$.*

It is to be noted that the identity must be satisfied by both the non-associative rings $R_1$ and $R_2$.

**DEFINITION 6.2.9:** *Let $R = R_1 \cup R_2$ be a non-associative biring. R is said to be a Bol biring if ((xy)z) y = x ((yz) y) for all x, y, z ∈ $R_1$ and for all x, y, z ∈ $R_2$.*

*Similarly we say the biring $R = R_1 \cup R_2$ is an alternative biring if (xy) y = x (yy) and x (xy) = (xx) y are true for all x, y ∈ $R_1$ and x, y ∈ $R_2$. We call a biring to be a weak inverse property biring, WIP biring if (xy) z = e imply x (yz) = e for all x, y, z ∈ $R_1$ and x, y, z ∈ $R_2$.*

**DEFINITION 6.2.10:** *If x, y are elements of a non-associative biring $R = R_1 \cup R_2$ (x, y ∈ $R_1$ or x, y ∈ $R_2$) the commutator (x, y) is defined by xy = (yx) (x, y). The commutator sub-biring is the biring generated under '+' and '•' by 〈 {x ∈ R = $R_1 \cup R_2$ / x = (y, z) for some y, z ∈ R} 〉.*

**DEFINITION 6.2.11:** *If x, y and z are elements of a non-associative biring $R = R_1 \cup R_2$; an associator (x, y, z) is defined by, (xy) z = (x(yz)) (x, y, z).*

The associator sub-biring of the biring R denoted by A(R) is the sub-biring generated by all of its associators that is by 〈 {x ∈ R = $R_1 \cup R_2$ | x = (a, b, c) for some a, b, c ∈ R} 〉.

The natural study would be, is the commutator subring a sub-biring? Is the associator subring a sub-biring? We see it has become impossible to find non-associative rings except those non-associative rings built using loops and groupoids i.e. loop rings and groupoid rings. Apart from this we are not able to construct non-associative rings, barring the well known class of non-associative rings viz Lie algebras and Jordan algebras, but both these algebras makes use of only the tools of vector spaces over fields so they are not like rings built using a single set.

All the more even these non-associative rings cannot be got naturally using Z or Q or R or $Z_n$. These non-associative rings, which we investigate, are only using loops over rings i.e. loop rings analogous to group rings and groupoid over rings analogous to semigroup rings. Thus we are not in a position to give nontrivial examples of them.

**DEFINITION 6.2.12:** *Let $R = R_1 \cup R_2$ be a non-associative biring. Let $x \in R$, we say $x$ is said to be right quasi regular if there exists a $y \in R$ such that $x \circ y = 0$; $x \in R$ is said to be left quasi regular if there exists $y' \in R$ such that $y' \circ x = 0$. An element is quasi regular if it is both left and right quasi regular. $y$ will be known as the right quasi inverse of $x$ and $y'$ as the left quasi inverse of $x$.*

A right ideal or left ideal in R is said to be right quasi regular if each of its elements is right quasi regular. One cannot always say a left quasi regular ideal or a right quasi regular ideal of a biring in general need not be a bi-ideal of $R = R_1 \cup R_2$.

**THEOREM 6.2.2:** *Let $R = R_1 \cup R_2$ be a biring such that $R_1 \cap R_2 = \{0, 1\}$. If both $R_1$ and $R_2$ has quasi regular ideals then we say R has a quasi regular bi-ideals.*

*Proof*: Straightforward by the very definitions.

**DEFINITION 6.2.13:** *Let $R = R_1 \cup R_2$ be a biring the Jacobson radical of R is defined as $J(R) = J(R_1) \cup J(R_2) = \{a_1 \in R_1 / a_1 R_1$ is right quasi regular ideal of $R_1\} \cup \{a_2 \in R_2 / a_2 R_2$ is right quasi regular ideal of $R_2\}$. If $J(R) = 0$, then we say the biring is bisemi-simple; i.e. $J(R_1) = 0$ and $J(R_2) = 0$.*

Now we proceed on to define the notion of loop birings and groupoid birings where birings are taken to be associative commutative biring with unit or bifields.

**DEFINITION 6.2.14:** *Let $R = R_1 \cup R_2$ be a biring which is commutative, associative with 1 or a bifield. L any loop. The loop biring $RL = R_1 L \cup R_2 L$ is defined as the union of the loop rings $R_1 L$ and $R_2 L$.*

**THEOREM 6.2.3:** *The loop biring of the loop L over the biring $R = R_1 \cup R_2$ is a non-associative biring.*

**Proof:** Straightforward by the very definitions.

**Note:** Let $L_1$ and $L_2$ be any two loops $R_1$ and $R_2$ be any two rings, $R_1 L_1$ and $R_2 L_2$ be loop rings. If $R = R_1 L_1 \cup R_2 L_2$ then we do not call R a loop biring even if $R_1 \cup R_2$ happens to be a biring

*Example 6.2.5*: Let L be a loop given by the following table:

| * | e | a | b | c | d |
|---|---|---|---|---|---|
| e | e | a | b | c | d |
| a | a | e | c | d | b |
| b | b | d | a | e | c |
| c | c | b | d | a | e |
| d | d | c | e | b | a |

$R = R_1 \cup R_2$ be a commutative biring with 1 or a bifield then $RL = R_1 L \cup R_2 L$ is a loop biring.

Study of properties related to loop birings are very many and interesting.

**DEFINITION 6.2.15:** *Let $RL = R_1L \cup R_2L$ be a loop biring of the loop L over the biring $R = R_1 \cup R_2$. The augmentation bi-ideal of RL denoted by $W(RL) = W(R_1L) \cup W(R_2L)$ where*

$$W(R_iL) = \left\{ \alpha = \sum \alpha_j m_j \in R_iL \mid \sum \alpha_j = 0 \right\} \text{ for } i = 1, 2.$$

Clearly the augmentation ideal of a loop ring is a bi-ideal of RL. The study of the concepts of zero divisors, units, idempotents etc in case of loop birings are similar to those in case of birings.

**THEOREM 6.2.4:** *Let $R = Q(\sqrt{2}) \cup Q(\sqrt{3})$ be a bifield where $Q(\sqrt{2})$ and $Q(\sqrt{3})$ are fields. Let L be a loop having atleast one element $h_i$ such that $h_i^2 = 1$ ($h_i \neq 1$), then $x = m(1 + h_i) \in RL$ is quasi regular if and only if x is not an idempotent.*

*Proof*: Suppose $x = m(1 + h_i) \in RL$ is quasi regular. ($m \neq 0$ or $m \neq ½$), x is an idempotent only when $m = ½$ so x is not an idempotent in RL. Conversely if $x = m(1 + h_i) \in RL$, $m \neq 0$ is not idempotent that is $m \neq ½$. Take $y = m / 2m - 1 (1 + h_i) \in RL$ then

$$x + y - xy = m(1 + h_i) + \frac{m}{2m-1}(1 + h_i)$$

$$- m(1 + h_i)\frac{m}{2m-1}(1 + h_i) = 0$$

Hence x is right quasi regular with y as its right quasi inverse.

**THEOREM 6.2.5:** *Let L be a finite loop and $R = Q(\sqrt{2}) \cup Q(\sqrt{3})$ be a bifield of characteristic zero. Then $J(RL) \subset W(RL)$.*

*Proof*: We know $J(RL) = J(R_1L) \cup J(R_2L)$ where $R_1 = Q(\sqrt{2})$ and $R_2 = Q(\sqrt{3})$ similarly $W(RL) = W(R_1L) \cup W(R_2L)$. To prove $J(RL) \subseteq W(RL)$ it is sufficient to show that if an element $x = \Sigma a_i h_i \notin W(RL)$ then ($x \notin W(R_1L)$ and $x \notin W(R_2L)$), ($x \notin J(RL)$) (i.e. $x \notin J(R_1L)$ and $x \notin J(R_2L)$).

Since ($x \notin W(RL)$ we have $\Sigma a_i \neq 0$. But it is straightforward from the fact if $x \in RL$ is in J(RL) then $\Sigma a_i \neq 0$. So we see $x \notin J(RL)$ Hence the claim.

**DEFINITION 6.2.16:** *Let $R = R_1 \cup R_2$ be a non-associative biring such that both the rings $R_1$ and $R_2$ are rings without divisors of zero then we call R a non-associative bidivision rings. If in particular both $R_1$ and $R_2$ happen to be commutative without divisors of zero we call R a non-associative bidomain.*

Study of loop birings can be carried out as in the case of group birings. Several interesting results in this direction can be obtained.

Now we give an example of a non-associative biring using the biring $R = R_1 \cup R_2$ where R is an associative biring.

*Example 6.2.6*: Let Z be the ring of integers, $L_1$ and $L_2$ be any two distinct loops. $ZL_1$ and $ZL_2$ be the loop rings of the loops $L_1$ and $L_2$ over Z. Then $R = ZL_1 \cup ZL_2$ is a non-associative biring. Thus this method helps us to construct more and more non-associative birings.

*Example 6.2.7*: Let $Z_{12}$ be the ring of integers modulo 12. $G = \langle g \mid g^{21} = 1 \rangle$ and L be a loop given by the following table:

| • | e | $g_1$ | $g_2$ | $g_3$ | $g_4$ | $g_5$ |
|---|---|---|---|---|---|---|
| e | e | $g_1$ | $g_2$ | $g_3$ | $g_4$ | $g_5$ |
| $g_1$ | $g_1$ | e | $g_3$ | $g_5$ | $g_2$ | $g_4$ |
| $g_2$ | $g_2$ | $g_5$ | e | $g_4$ | $g_1$ | $g_3$ |
| $g_3$ | $g_3$ | $g_4$ | $g_1$ | e | $g_5$ | $g_2$ |
| $g_4$ | $g_4$ | $g_3$ | $g_5$ | $g_2$ | e | $g_1$ |
| $g_5$ | $g_5$ | $g_2$ | $g_4$ | $g_1$ | $g_3$ | e |

$R = Z_{12} G \cup Z_{12} L$ is the non-associative biring where $Z_{12}G$ is the group ring of G over $Z_{12}$ and $Z_{12}L$ is a loop ring of the loop L over $Z_{12}$.

*Example 6.2.8*: Let $Z_2 = \{0, 1\}$ be the prime field of characteristic 2 and $G = S_3$ and L be the loop given in the above example. Then $R = Z_2L \cup Z_2S_3$ is a non-associative biring which has nontrivial units, zero divisors and idempotents.

*Example 6.2.9*: Let Z be the ring of integers. S(5) be the symmetric semigroup on 5 elements and L be any loop, ZS(5) be the semigroup ring and ZL the loop ring. Then $R = ZL \cup ZS(5)$ is a biring which is non-associative. Thus we have an infinite number of birings, which are non-associative.

Now we proceed on to define the concept of groupoid biring.

**DEFINITION 6.2.17**: *Let $R = R_1 \cup R_2$ be any biring, commutative with unit and G any groupoid. The groupoid biring $RG = R_1G \cup R_2G$ is the biring by taking the union of the groupoid ring $R_1G$ with the groupoid ring $R_2G$.*

Clearly the groupoid ring RG is a non-associative biring. It is important to note that if $G_1$ and $G_2$ are any two groupoids and $R = R_1 \cup R_2$ be a biring then $K = R_1G_1 \cup R_2G_2$ is not a groupoid biring; but certainly K is a non-associative biring. Also if $R_1$ and $R_2$ are just any two distinct rings and $G_1$ and $G_2$ any arbitrary groupoids with $R_1G_1$ and $R_2G_2$ the groupoid rings then $R = R_1G_1 \cup R_2G_2$ is a non-associative biring and not a groupoid biring.

*Example 6.2.10*: Let $R = R_1 \cup R_2$ be any biring say $R_1 = Q(\sqrt{5})$ and $R_2 = Z$ (the ring of integers). G be any groupoid given by the following table:

| • | $a_0$ | $a_1$ | $a_2$ |
|---|---|---|---|
| $a_0$ | $a_0$ | $a_2$ | $a_1$ |
| $a_1$ | $a_1$ | $a_0$ | $a_2$ |
| $a_2$ | $a_2$ | $a_1$ | $a_0$ |

Clearly $R = R_1G_1 \cup R_2G_2$ is a groupoid biring which is a non-associative biring.

***Example 6.2.11:*** Let G and $G_1$ be two groupoids given by the following tables:

Table for G

| • | $a_0$ | $a_1$ | $a_2$ | $a_3$ | $a_4$ | $a_5$ |
|---|---|---|---|---|---|---|
| $a_0$ | $a_0$ | $a_2$ | $a_4$ | $a_0$ | $a_2$ | $a_4$ |
| $a_1$ | $a_2$ | $a_4$ | $a_0$ | $a_2$ | $a_4$ | $a_0$ |
| $a_2$ | $a_4$ | $a_0$ | $a_2$ | $a_4$ | $a_0$ | $a_2$ |
| $a_3$ | $a_0$ | $a_2$ | $a_4$ | $a_0$ | $a_2$ | $a_4$ |
| $a_4$ | $a_2$ | $a_4$ | $a_0$ | $a_2$ | $a_4$ | $a_0$ |
| $a_5$ | $a_4$ | $a_0$ | $a_2$ | $a_4$ | $a_0$ | $a_2$ |

Table for $G_1$

| * | a | b | c | d |
|---|---|---|---|---|
| a | a | c | a | c |
| b | a | c | a | c |
| c | a | c | a | c |
| d | a | c | a | c |

Let Z be the ring of integers and $Z_{11}$ be the prime field of characteristic 11, then ZG and $Z_{11}G_1$ be the groupoid ring. Now $R = ZG \cup Z_{11}G_1$ is a biring but not a groupoid biring.

***Example 6.2.12:*** Let $R = R_1G_1 \cup R_2G_2$ with groupoid rings $R_1G_1$ and $R_2G_2$ where $G_1$ and $G_2$ are groupoids given by the following tables:

Table for $G_1$

| • | $a_0$ | $a_1$ | $a_2$ | $a_3$ | $a_4$ | $a_5$ |
|---|---|---|---|---|---|---|
| $a_0$ | $a_0$ | $a_0$ | $a_0$ | $a_0$ | $a_0$ | $a_0$ |
| $a_1$ | $a_3$ | $a_3$ | $a_3$ | $a_3$ | $a_3$ | $a_3$ |
| $a_2$ | $a_0$ | $a_0$ | $a_0$ | $a_0$ | $a_0$ | $a_0$ |
| $a_3$ | $a_3$ | $a_3$ | $a_3$ | $a_3$ | $a_3$ | $a_3$ |
| $a_4$ | $a_0$ | $a_0$ | $a_0$ | $a_0$ | $a_0$ | $a_0$ |
| $a_5$ | $a_3$ | $a_3$ | $a_3$ | $a_3$ | $a_3$ | $a_3$ |

Table for $G_2$

| • | $g_0$ | $g_1$ | $g_2$ | $g_3$ | $g_4$ | $g_5$ | $g_6$ | $g_7$ |
|---|---|---|---|---|---|---|---|---|
| $g_0$ | $g_0$ | $g_6$ | $g_4$ | $g_2$ | $g_0$ | $g_6$ | $g_4$ | $g_2$ |
| $g_1$ | $g_2$ | $g_0$ | $g_6$ | $g_4$ | $g_2$ | $g_0$ | $g_6$ | $g_4$ |
| $g_2$ | $g_4$ | $g_2$ | $g_0$ | $g_6$ | $g_4$ | $g_2$ | $g_0$ | $g_6$ |
| $g_3$ | $g_6$ | $g_4$ | $g_2$ | $g_0$ | $g_6$ | $g_4$ | $g_2$ | $g_0$ |
| $g_4$ | $g_0$ | $g_6$ | $g_4$ | $g_2$ | $g_0$ | $g_6$ | $g_4$ | $g_2$ |
| $g_5$ | $g_2$ | $g_0$ | $g_6$ | $g_4$ | $g_2$ | $g_0$ | $g_6$ | $g_4$ |
| $g_6$ | $g_4$ | $g_2$ | $g_0$ | $g_6$ | $g_4$ | $g_2$ | $g_0$ | $g_6$ |
| $g_7$ | $g_6$ | $g_4$ | $g_2$ | $g_0$ | $g_6$ | $g_4$ | $g_2$ | $g_0$ |

Clearly R is a non-associative biring.

**Example 6.2.13:** Let $R = Z_2G \cup Z_2S_6$; clearly R is a non-associative biring where G is some groupoid and $Z_2G$ is a groupoid ring and $Z_2S_6$ is a group ring. The biring is a non-commutative biring having zero divisors, units, idempotents and nilpotents.

All results pertaining to loop biring and semigroup biring can also be derived as in case of groupoid birings with suitable modifications. Thus the work of doing research about groupoid birings is left for the reader to develop.

The conditions when a Moufang groupoid is used what is the structure of the groupoid biring.

Let $R = R_1 \cup R_2$ be a biring, which is commutative with 1 or a bifield. G be a Moufang groupoid and RG be the groupoid biring. Then $RG = R_1G \cup R_2G$ is a biring which is a Moufang biring.

Similar results can be used to obtain Bol, Bruck, alternative birings. Not only these still interesting notions which has not been excavated about groupoid birings is the associator sub-biring and the commutator sub-biring. Now we have given methods of constructing classes of non-associative birings using groupoids and loops. Thus if we take any groupoids $G_1$, $G_2$ or loops $L_1$, $L_2$ using them construct the groupoid rings $RG_1$ and $RG_2$ over the same ring R then $K = R G_1 \cup RG_2$ is a non-associative biring.

Similarly if $K_1 = RL_1 \cup RL_2$ then also $K_1$ will be a non-associative biring. Not only these even if some $P = RL_1 \cup RS$ where $RL_1$ is the loop ring of the loop $L_1$ over R and RS is the semigroup ring then also P is a non-associative biring. Thus we have given methods to construct non-associative birings.

Now we proceed on to define biloop birings, bigroupoid birings and consequent of this we define quad rings, which are non-associative.

**DEFINITION 6.2.18:** *Let (R, +, •) be a set. We call R a non-associative quad ring if R = $R_1 \cup R_2 \cup R_3 \cup R_4$ where $R_i$'s are distinct and each $R_i$ i.e. ($R_i$, +, •) is a ring of which atleast one of the $R_i$ is a non-associative ring. It can happen that all of the $R_i$'s are non-associative rings.*

*When we say the sets $R_1$, $R_2$, $R_3$, $R_4$ are distinct we basically assume that $R_i \not\subset R_j$ for $i \neq j$; i, j = 1, 2, 3, 4.*

Quad rings, which are non-associative, can also be defined in another way.

**DEFINITION 6.2.19:** *Let (R, +, •) be a set. We say (R, +, •) is a non-associative quad ring if $R = R_1 \cup R_2$ where $R_1$ and $R_2$ are birings under the operations of R and atleast one of the birings is a non-associative biring.*

*Example 6.2.14:* Let R be a commutative ring with 1 say the set of integers. Take a loop L, a groupoid G and the group $S_8$ and the semigroup S(4) where the loop L and the groupoid G are given by the following tables:

Table for L

| • | e | $g_1$ | $g_2$ | $g_3$ | $g_4$ | $g_5$ |
|---|---|---|---|---|---|---|
| e | e | $g_1$ | $g_2$ | $g_3$ | $g_4$ | $g_5$ |
| $g_1$ | $g_1$ | e | $g_3$ | $g_5$ | $g_2$ | $g_4$ |
| $g_2$ | $g_2$ | $g_5$ | e | $g_4$ | $g_1$ | $g_3$ |
| $g_3$ | $g_3$ | $g_4$ | $g_1$ | e | $g_5$ | $g_2$ |
| $g_4$ | $g_4$ | $g_3$ | $g_5$ | $g_2$ | e | $g_1$ |
| $g_5$ | $g_5$ | $g_2$ | $g_4$ | $g_1$ | $g_3$ | e |

Table for G

| • | $a_1$ | $a_2$ | $a_3$ |
|---|---|---|---|
| $a_1$ | $a_1$ | $a_3$ | $a_2$ |
| $a_2$ | $a_2$ | $a_1$ | $a_3$ |
| $a_3$ | $a_3$ | $a_2$ | $a_1$ |

Now take R = ZG ∪ ZL ∪ $ZS_8$ ∪ ZS(4) clearly R is a quad ring which is non-associative. We see that in general we can, by this method get several quad rings, which are non-associative. So we have the nontrivial existence of quad rings. Now we get yet another new classes of quad rings which are non-associative using biloops and bigroupoids.

**DEFINITION 6.2.20:** *Let $R = R_1 \cup R_2$ be a biring and L be a biloop. $L = L_1 \cup L_2$. The biloop biring $RL = R_1L_1 \cup R_1L_2 \cup R_2L_1 \cup R_2L_2$ is a quad ring which is non-associative.*

**THEOREM 6.2.6:** *Every biloop biring is a non-associative quad ring.*

*Proof*: Given $L = L_1 \cup L_2$ is a biloop and $R = R_1 \cup R_2$ is a biring. Suppose RL = $R_1L_1$ ∪ $R_1L_2$ ∪ $R_2L_1$ ∪ $R_2L_2$ be the biloop biring.

Now RL is a union of four rings certainly some of them are non-associative rings. Thus every biloop biring is a non-associative quad ring.

*Example 6.2.15:* Let $L = L_1 \cup L_2$ be a biloop where $L_1 = G = \langle g \mid g^7 = 1 \rangle$ be a cyclic group of order 7 and $L_2$ be the loop given by the following table:

| • | e | $g_1$ | $g_2$ | $g_3$ | $g_4$ | $g_5$ |
|---|---|---|---|---|---|---|
| e | e | $g_1$ | $g_2$ | $g_3$ | $g_4$ | $g_5$ |
| $g_1$ | $g_1$ | e | $g_4$ | $g_2$ | $g_5$ | $g_3$ |
| $g_2$ | $g_2$ | $g_4$ | e | $g_5$ | $g_3$ | $g_1$ |
| $g_3$ | $g_3$ | $g_2$ | $g_5$ | e | $g_1$ | $g_4$ |
| $g_4$ | $g_4$ | $g_5$ | $g_3$ | $g_1$ | e | $g_2$ |
| $g_5$ | $g_5$ | $g_3$ | $g_1$ | $g_4$ | $g_2$ | e |

Let R = Q(√2) ∪ Q(√3) be the bifield. The biloop biring RL = Q(√2) G ∪ Q(√2)L$_2$ ∪ Q(√3)G ∪ Q(√3)L$_2$ is a non-associative quad ring.

**DEFINITION 6.2.21:** *Let R = R$_1$ ∪ R$_2$ ∪ R$_3$ ∪ R$_4$ be a non-associative quad ring; we say R is a commutative quad ring if each of the rings R$_i$ are commutative. If 1 ∈ R such that 1• r = r • 1 = r for all r ∈ R we say R is a quad ring with unit. If the quad ring R has no zero divisors i.e. for no r ∈ R we have a s ≠ 0 in R such that r • s = 0 then we say R is a non-associative quad division ring. If R is a non-associative quad division ring, which is commutative, then we call the quad ring as a non-associative quad field.*

*All notions like semiprime, prime, semisimple can also defined as in case of biloop birings. For instance we define the Jacobson radical of the biloop biring RL = R$_1$L$_1$ ∪ R$_1$L$_2$ ∪ R$_2$L$_1$ ∪ R$_2$L$_2$ where L = L$_1$ ∪ L$_2$ and R = R$_1$ ∪ R$_2$ is defined as J(RL) = J(R$_1$L$_1$) ∪ J (R$_1$L$_2$) ∪ J (R$_2$L$_2$) ∪ J (R$_2$L$_2$).*

*Similarly we define the augmentation ideal W(RL) = W(R$_1$L$_1$) ∪ W(R$_1$L$_2$) ∪ W(R$_2$L$_1$) ∪ W(R$_2$L$_2$).*

Now we will define quad ideal I of quad ring R as I = I$_1$ ∪ I$_2$ ∪ I$_3$ ∪ I$_4$ where each I$_i$ is an ideal of R$_i$, i = 1, 2, 3, 4. The concepts of maximal quad ideal, minimal, principal etc can be defined in an analogous way. It is interesting to see that the Jacobson radical of the quad ring R is a quad ideal of R. Similarly the augmentation ideal of the quad ring is a quad ideal.

A marked difference between a quad ring and a non-associative quad ring is that in case of non-associative quad ring we will have several special types of quad rings which satisfies the identities like Moufang, Bol, Bruck, alternative.

So in case a quad ring R satisfies the Bol identity we call the quad ring R a Bol quad ring. Similarly we define quad rings which are alternative quad rings.

Study of these quad rings is really interesting and innovative. The concept of simple and semisimple quad rings is defined in an analogous way as in case of birings. Now we define bigroupoid birings.

**DEFINITION 6.2.22:** *Let G be a bigroupoid i.e. G = G$_1$ ∪ G$_2$ where one of G$_1$ or G$_2$ is a groupoid and R = R$_1$ ∪ R$_2$ is the associative commutative biring with unit. The bigroupoid biring RG = R$_1$G$_1$ ∪ R$_1$G$_2$ ∪ R$_2$G$_1$ ∪ R$_2$G$_2$ is defined as the union of groupoid rings R$_i$G$_j$, i = 1, 2 and j = 1, 2.*

**THEOREM 6.2.7:** *Let RG be a bigroupoid biring; G = G$_1$ ∪ G$_2$ is the bigroupoid and R = R$_1$ ∪ R$_2$ is a biring. RG is a non-associative quad ring.*

*Proof*: Follows from the fact that the bigroupoid biring. RG = R$_1$G$_1$ ∪ R$_1$G$_2$ ∪ R$_2$G$_1$ ∪ R$_2$G$_2$ which is the union of 4 rings of which at least two are non-associative. So RG is a quad ring.

Thus using the new class of groupoids and biring we can have several quad rings, which are non-associative. Study of bigroupoid biring is very new. Only in this book we have defined the concept of bigroupoid biring, which will form a new class of non-associative quad rings. The main features, which will be interesting, are

  i. When are two bigroupoid birings isomorphic?
 ii. When is a bigroupoid biring semiprime?
iii. When is a bigroupoid biring prime?

To this end we define the concept of the Jacobson radical and the augmentation ideal of the non-associative quad rings.

**DEFINITION 6.2.23:** *Let RG be a bigroupoid biring where $R = R_1 \cup R_2$ is an associative biring and $G = G_1 \cup G_2$ is a bigroupoid. The Augmentation ideal of RG denoted by $W(RG) = W(R_1G_1) \cup W(R_1G_2) \cup W(R_2G_1) \cup W(R_2G_2)$.*

Clearly the augmentation ideal of a bigroupoid biring is a quad ideal of RG.

**DEFINITION 6.2.24:** *Let RG be a bigroupoid biring of the bigroupoid $G = G_1 \cup G_2$ over the biring $R = R_1 \cup R_2$. The Jacobson radical of RG denoted by $J(RG) = J(R_1G_1) \cup J(R_1G_2) \cup J(R_2G_1) \cup J(R_2G_2)$.*

Clearly the Jacobson radical of RG is an ideal this ideal is also a quad ideal of RG. The analysis of bigroupoid birings to be a Moufang quad ring, Bol quad ring, Bruck quad ring and right alternative quad ring will yield a lot of results. Can we say if G is a Bol bigroupoid then the bigroupoid biring RG where R is a biring of characteristic zero be a Bol quad ring?

Such study can be extended to cases of Bruck quad ring, Moufang quad ring etc.

**PROBLEMS:**

1. What is the order of the smallest non-associative biring?
2. Give an example of a non-associative biring, which is commutative with 1.
3. Does there exist a non-associative biring of order 11? Justify your claim.
4. Can there be a non-associative biring of prime order p, p any prime? (p < 19).
5. Does there exist a non-associative biring of order less than 16?
6. Give an example of a non-associative bifield of finite order (i.e. having only finite number of elements in them).
7. Give an example of a quad ring of finite order.
8. Let $R = QL \cup QS_3$ be a biring where $S_3$ is the symmetric group of degree 3, Q the field of rationals and L a loop given by the following table:

| • | e | a | b | c | d |
|---|---|---|---|---|---|
| e | e | a | b | c | d |
| a | a | e | c | d | b |
| b | b | d | a | e | c |
| c | c | b | d | a | e |
| d | d | c | e | b | a |

Find
   i. J(R).
   ii. W(R).
   iii. Ideals.
   iv. sub-birings.

9. Give an example of a quad ring of order 120.
10. In the quad ring $R = Z_2 G \cup Z_2 S_3 \cup Z_2 L \cup Z_2 G_1$, where $Z_2 = \{0, 1\}$ the prime field of characteristic two and $S_3$ the symmetric group of degree 3, G and $G_1$ be two groupoids given by the tables. L is a loop such that $L \in L_7$.

| *     | e     | $a_0$ | $a_1$ | $a_2$ |
|-------|-------|-------|-------|-------|
| e     | e     | $a_0$ | $a_1$ | $a_2$ |
| $a_0$ | $a_0$ | e     | $a_2$ | $a_1$ |
| $a_1$ | $a_1$ | $a_2$ | e     | $a_0$ |
| $a_2$ | $a_2$ | $a_1$ | $a_0$ | e     |

| *     | $a_0$ | $a_1$ | $a_2$ | $a_3$ | $a_4$ | $a_5$ |
|-------|-------|-------|-------|-------|-------|-------|
| $a_0$ | $a_0$ | $a_2$ | $a_4$ | $a_0$ | $a_2$ | $a_4$ |
| $a_1$ | $a_0$ | $a_2$ | $a_4$ | $a_0$ | $a_2$ | $a_4$ |
| $a_2$ | $a_0$ | $a_2$ | $a_4$ | $a_0$ | $a_2$ | $a_4$ |
| $a_3$ | $a_0$ | $a_2$ | $a_4$ | $a_0$ | $a_2$ | $a_4$ |
| $a_4$ | $a_0$ | $a_2$ | $a_4$ | $a_0$ | $a_2$ | $a_4$ |
| $a_5$ | $a_0$ | $a_2$ | $a_4$ | $a_0$ | $a_2$ | $a_4$ |

Find
   i. J(R).
   ii. W(R).
   iii. Ideals of R.
   iv. Right ideals of R.
   v. Quad subrings, which are not quad ideals of R.
   vi. Does R satisfy any of the standard identities like Moufang, Bol, Bruck or alternative?
   vii. Does R have units, zero divisors and idempotents?

## 6.3 Smarandache birings and its properties

In this section we define the notion of Smarandache birings (S-birings) which are both associative as well as non-associative. Study of this type is being carried out for the first time in this book. Also in this section we define Smarandache bifields and study some of its basic properties. Several important properties about these new structures are introduced in this section.

**DEFINITION 6.3.1:** *A Smarandache biring (S-biring) (R, +, •) is a non-empty set with two binary operations '+' and '•' such that $R = R_1 \cup R_2$ where $R_1$ and $R_2$ are proper subsets of R and*

    i. $(R_1, +, \bullet)$ is a S-ring.
    ii. $(R_2, +, \bullet)$ is a S-ring.

*If only one of $R_1$ or $R_2$ is a S-ring then we call $(R, +, \bullet)$ a Smarandache weak biring (S-weak biring).*

**THEOREM 6.3.1:** *All S-birings are S-weak birings and not conversely.*

*Proof*: Follows directly by the very definition.

*Example 6.3.1:* Take $R = Q[x] \cup Z_4$ clearly R is a S-weak biring and not a S-biring.

*Example 6.3.2:* Let $R = Q[x] \cup Z_{12}$ is a S-biring. For both $\{0, 3, 8\} \subset Z_{12}$ and $Q \subset Q[x]$ are fields.

**DEFINITION 6.3.2:** *Let $(R, +, \bullet)$ be a biring. R is said to be a Smarandache conventional biring (S-conventional biring) if and only if R has a proper subset P that is a bifield.*

In view of this definition we can say a S-biring can have a characteristic associated with it. Also we see both definitions may not be equivalent.

For in this case a S-conventional biring is a S-biring but all S-birings are not S-conventional biring. We illustrate this by the following example.

*Example 6.3.3:* Let $(R, +, \bullet)$ be a non-empty set, where $R = R_1 \cup R_2$ with $R_1 = Q \times Q$ and $R_2 = Q[x]$. Clearly both $R_1$ and $R_2$ are S-rings but $Q[x]$ has Q to be the subfield and $Q \times Q$ has $Q \times \{0\}$ to be the field. So $P = Q \cup Q \times \{0\}$ is not a bifield as they are not proper subsets of P. Thus all S-conventional birings are S-birings but every S-biring in general need not be a S-conventional biring.

**DEFINITION 6.3.3:** *Let $(R, +, \bullet)$ be a biring. A proper subset P of R is said to be a Smarandache sub-biring (S-sub-biring) if $(P, +, \bullet)$ is itself a S-biring. Similarly we say for a biring $(R, +, \bullet)$ a proper subset P of R is said to be a Smarandache conventional sub-biring (S-conventional sub-biring) if $(P, +, \bullet)$ is itself a S-conventional biring under the operations of R.*

Thus we have the following very important observations.

**THEOREM 6.3.2:** *If a biring R has a S-sub-biring then R is itself a S-biring.*

*Proof*: Direct by the very definition.

On similar lines the reader can prove the following theorem.

**THEOREM 6.3.3:** *If a biring R has a S-conventional sub-biring then R itself is a S-conventional biring.*

**DEFINITION 6.3.4:** *Let R be a biring. We say R is a Smarandache commutative biring (S-commutative biring) if every S-sub-biring is commutative. If at least one of the S-sub-biring of R is commutative then we say R is a Smarandache weakly commutative biring (S-weakly commutative biring).*

*We define similarly in case of S-conventional birings to be S-commutative and S-weakly commutative conventional birings.*

**THEOREM 6.3.4:** *Let (R, +, •) be a biring. If R is a S-commutative biring then R is a S-weakly commutative biring.*

*Proof*: Left for the reader as an exercise.

Now we define of level II S-biring.

**DEFINITION 6.3.5:** *Let (R, +, •) be a non-empty set. We say (R, +, •) is a Smarandache biring II (S-biring II) if R contains a proper subset P such that P is a bidivision ring.*

Thus we see all S-conventional birings are S-biring II.

Hence we say a S-biring II, R is commutative if it has S-conventional biring. We say the S-biring is finite if R has finite number of elements otherwise R has infinite number of elements and the S-biring is of infinite order. In a biring we define S-units, S-zero divisors and S-idempotents as in case of rings for these notions, the biring need not be a S-biring or a S-conventional biring.

**DEFINITION 6.3.6:** *Let (R, +, •) be a commutative biring. If R has no S-zero divisors then we call R a Smarandache integral bidomain (S-integral bidomain).*

*If R is a non-commutative biring and R has no S-zero divisors then we call R a Smarandache division biring (S-division biring).*

**DEFINITION 6.3.7:** *Let (R, +, •) be a biring. The Smarandache bi-ideal (S-bi-ideal) P of R is defined as an ideal of the biring R where P is a S-sub-biring.*

**DEFINITION 6.3.8:** *Let (R, +, •) be a S-conventional biring. B ⊂ R is a bifield of R. A non-empty subset P of R is said to be Smarandache pseudo right bi-ideal (S-pseudo right bi-ideal) of R related to P if*

      i.   (P, +) a an abelian bigroup
     ii.   For b ∈ B and p ∈ P we have p • b ∈ P.

*On similar lines we define Smarandache pseudo left bi-ideal (S-pseudo left bi-ideal). A non-empty subset P of R is said to be a Smarandache pseudo bi-ideal (S-pseudo bi-ideal) if P is both a S-pseudo right bi-ideal and S-pseudo left bi-ideal.*

The concept of S-maximal and S-minimal ideals in birings are defined in an analogous way as in case of rings similarly one can define S-minimal pseudo bi-ideal,

S-maximal pseudo bi-ideals, S-cyclic pseudo bi-ideals, S-prime pseudo bi-ideals in a biring.

**DEFINITION 6.3.9:** *Let (R, +, •) be a biring . If R is a S-biring then so is the polynomial biring R[x]; as R ⊂ R[x], R is a S-biring.*

**DEFINITION 6.3.10**: *If a S-biring (R, +, •) has no S-bi-ideals then we call R a Smarandache simple biring (S-simple biring). If R has no S-pseudo bi-ideals then we call, R Smarandache pseudo simple biring (S-pseudo simple biring).*

*The Smarandache quotient biring (S-quotient biring) is defined as R/I = {a + I/a ∈ R} where R is a biring and I is an S-bi-ideal of R. The Smarandache pseudo quotient biring (S-pseudo quotient biring) is defined as R/J = {a + J/a ∈ R} where J is a S-pseudo bi-ideal of the biring R.*

The concept of Smarandache R-module for rings was introduced in [44, 49] we define Smarandache R-module for birings.

**DEFINITION 6.3.11:** *Let (R, +, •) be a biring. A non-empty subset M is said to be a Smarandache bimodule (S-bimodule) if M is an S-abelian semigroup with 0 under the operation '+' i.e. $M = M_1 \cup M_2$ and all the conditions enumerated in the definition 6.1.22 is satisfied.*

*The only change we make in the definition of S-bimodule is that we replace the abelian bigroups M by a S-abelian semigroup with 0 i.e. S-commutative monoid.*

Several interesting results in this direction can be obtained.

**DEFINITION 6.3.12:** *Let (R, +, •) be a biring. We say R is a Smarandache semiprime biring (S-semiprime biring) if R contains no non-zero S-bi-ideal with square zero.*

Now we proceed on to define Smarandache group birings, Smarandache semigroup birings and Smarandache quad birings.

**DEFINITION 6.3.13:** *Let R be a biring and S a semigroup. RS the semigroup biring is a Smarandache semigroup biring (S-semigroup biring) if and only if S is a S-semigroup.*

**DEFINITION 6.3.14:** *Let (R, +, •) be a biring and G any group. The group biring RG is a Smarandache group biring (S-group biring) if and only if R is a S-biring.*

It may so happen that a group biring may be a S-biring but still we do not call it a S-group biring. Similarly the semigroup biring may be a S-biring but we do not call it a S-semigroup biring unless S the semigroup is a S-semigroup.

Now we define still new concepts using a S-bisemigroup and a ring, S-bigroup and a ring.

**DEFINITION 6.3.15:** *Let R be a ring. $S = S_1 \cup S_2$ be a S-bisemigroup. The bisemigroup ring RS is defined as $RS_1 \cup RS_2$ we call this bisemigroup ring as a*

*Smarandache bisemigroup ring (S-bisemigroup ring). Clearly S-bisemigroup ring is a biring.*

*Similarly when we take a $G = G_1 \cup G_2$ a S-bigroup and define a bigroup ring we call it as the Smarandache bigroup ring (S-bigroup ring). Thus $RG = RG_1 \cup RG_2$ is a biring.*

The next interesting study is when are the concepts bigroup biring and bisemigroup biring S-birings.

**DEFINITION 6.3.16:** *Let R be a biring i.e. $R = R_1 \cup R_2$ and S be a bisemigroup i.e. $S = S_1 \cup S_2$ then $RS = R_1S_1 \cup R_1S_2 \cup R_2S_1 \cup R_2S_2$ is the union of four semigroup rings. We say RS is Smarandache quad ring (S-quad ring) if each of the four semigroup rings are S-semigroups ring.*

*Similarly when we take G to be a bigroup say $G = G_1 \cup G_2$ and $R = R_1 \cup R_2$ a biring then the bigroup biring RG is the union of four group rings $R_1 G_1 \cup R_1 G_2 \cup R_2 G_1 \cup R_2 G_2$.*

*We say RG is a Smarandache quad ring II if each of the four group rings are S-group rings. So we say a S-quad ring Q is the union of two S-birings i.e. $Q = Q_1 \cup Q_2$ where $Q_1$ and $Q_2$ are S-birings.*

Several inter-relations can be defined and studied in case of S-quad rings, S-quad rings I and II. The next study will be based on non-associative birings, which are Smarandache na-birings. We define Smarandache na-birings in the following.

**DEFINITION 6.3.17:** *Let $(R, +, \bullet)$ with $R = R_1 \cup R_2$ be a na-biring, R is said to be a Smarandache na-biring (S-na-biring) if R has a proper subset P such that $P = P_1 \cup P_2$ and P is an associative biring under the operations of R.*

*Example 6.3.4:* Let $(R, +, \bullet)$ be a na-biring. $R = Z_{10}G \cup Z$ where $Z_{10}G$ is a groupoid ring of the groupoid G; G is given by the following table:

| * | 0 | 1 | 2 | 3 | 4 |
|---|---|---|---|---|---|
| 0 | 0 | 2 | 4 | 1 | 3 |
| 1 | 2 | 4 | 1 | 3 | 0 |
| 2 | 4 | 1 | 3 | 0 | 2 |
| 3 | 1 | 3 | 0 | 2 | 4 |
| 4 | 3 | 0 | 2 | 4 | 1 |

Clearly $(R, +, \bullet)$ is a S-na-biring as $P = Z_{10} \cup Z$ is an associative biring.

**DEFINITION 6.3.18:** *Let $(R, +, \bullet)$ be a na-biring. A proper subset $P \subset R$ is said to be a Smarandache na sub-biring (S-na-sub-biring) if P is itself a na-sub-biring and P is a S-na-biring.*

**THEOREM 6.3.5:** *If (R, +, •) has a na-sub-biring which is a S-na sub-biring then we say R is a S-na-biring.*

*Proof*: Follows directly from the definitions, hence left for the reader to prove.

**DEFINITION 6.3.19:** *Let (R, +, •) be a na-biring. We say R is a Smarandache commutative na-biring (S-commutative na-biring) if every S- na sub-biring of R is commutative. If at least one S-na sub-biring is commutative then we call R a Smarandache weakly commutative na-biring (S-weakly commutative na-biring).*

*Example 6.3.5:* Let (R, +, •) be a na-biring where $R = R_1 \cup R_2$ with $R_1 = Z_7L$ and $R_2 = Z_3$ (the ring of integers modulo 3). Here $R_1$ is the loop ring of the loop L over $Z_7$.

The loop L is given by the following table:

| * | e | $g_1$ | $g_2$ | $g_3$ | $g_4$ | $g_5$ | $g_6$ | $g_7$ |
|---|---|---|---|---|---|---|---|---|
| e | e | $g_1$ | $g_2$ | $g_3$ | $g_4$ | $g_5$ | $g_6$ | $g_7$ |
| $g_1$ | $g_1$ | e | $g_4$ | $g_7$ | $g_3$ | $g_6$ | $g_2$ | $g_5$ |
| $g_2$ | $g_2$ | $g_6$ | e | $g_5$ | $g_1$ | $g_4$ | $g_7$ | $g_3$ |
| $g_3$ | $g_3$ | $g_4$ | $g_7$ | e | $g_6$ | $g_2$ | $g_5$ | $g_1$ |
| $g_4$ | $g_4$ | $g_2$ | $g_5$ | $g_1$ | e | $g_7$ | $g_3$ | $g_6$ |
| $g_5$ | $g_5$ | $g_7$ | $g_3$ | $g_6$ | $g_2$ | e | $g_1$ | $g_4$ |
| $g_6$ | $g_6$ | $g_5$ | $g_1$ | $g_4$ | $g_7$ | $g_3$ | e | $g_2$ |
| $g_7$ | $g_7$ | $g_3$ | $g_6$ | $g_2$ | $g_5$ | $g_1$ | $g_4$ | e |

Several of the S-sub-bigroup is commutative except $P = Z_7L \cup \{0\}$.

So to overcome this problem we will state a sub-biring of the form $P = R_1 \cup \{0\}$ or $\{0\} \cup R_2$ where $R = R_1 \cup R_2$ will be declared as trivial. Thus the na-biring given in this example is a S-commutative na-biring.

**THEOREM 6.3.6:** *Let $R = R_1 \cup R_2$ be a na-biring. If R is a strongly subcommutative biring then R is a S-commutative biring.*

*Proof:* Straightforward, hence left for the reader to prove.

Now we proceed on to define Bruck, Bol, Moufang birings.

**DEFINITION 6.3.20:** *Let (R, +, •) be a na-biring, we call R a Smarandache Moufang biring (S-Moufang biring) if R has a proper subset P; such that P is a S-sub-na-biring of P and every triple of R satisfies the Moufang identity*

$$(xy)(zx) = (x(yz))x$$

for all $x, y, z \in P$.

Similarly $R = R_1 \cup R_2$ is a Smarandache Bruck biring (S-Bruck biring) if R has proper subset P such that

    i.    P is a S-na sub-biring of R.
    ii.    $(x(yx))z = x(y(xz))$ and $(xy)^{-1} = x^{-1}y^{-1}$

is true for all $x, y, z \in P$. We call $(R, +, \bullet)$ a Smarandache Bol biring (S-Bol biring) if $((xy\,z)\,y = x((yx)\,y)$ for all $x, y, z \in P$, $P \subset R$ is a S-na sub-biring of R.

On similar lines we define Smarandache left (right) alternative birings and Smarandache WIP-biring.

**DEFINITION 6.3.21:** *Let $R = (R_1 \cup R_2, +, \bullet)$ be a na-biring. If every S-na sub-biring P of R satisfies Moufang, Bol, Bruck, alternative (right / left) or WIP then we call R a Smarandache strong Moufang, Bol, Bruck, alternative (right/ left) or WIP biring respectively.*

**THEOREM 6.3.7:** *Every S-strong Moufang, Bol, Bruck, ... biring is a S-Moufang, Bol, Bruck, ... biring.*

*Proof*: Obvious by the very definitions.

**DEFINITION 6.3.22:** *Let $(R, +, \bullet)$ be a na-biring, $R = R_1 \cup R_2$. Let $P \subset R$ be a subset of R where $P = P_1 \cup P_2$ is a S-sub-biring of R. The commutator for x, y, in P is defined as $xy = (yx)(x, y)$ [where $x, y \in P_1$ or $x, y \in P_2$]. The commutator sub-biring of $P = \langle\{x \in P\,/\,x = (y, z)$ for some $y, z \in P\}\rangle$ is defined as the Smarandache commutator (S-commutator) of R related to P.*

It is important to note we may have more than one S-commutator for a biring. In some cases it may so happen that R may contain as many S-commutator sub-biring as S-sub-birings of R.

In a similar way Smarandache associator biring (S-associator biring) can be defined.

**DEFINITION 6.3.23:** *Let $(R = R_1 \cup R_2, +, \bullet)$ be a na-biring. Let $x \in S$ we say x is said to be right quasi regular if there exists a $y \in S$ such $x \circ y = 0$ where S is a S-sub-biring of R. Let $x \in S$, x is said to be left quasi regular if there exists $y' \in S$ such that $y' \circ x = 0$. An element is quasi regular if it is both right and left quasi regular.*

*A S-right ideal or S-left ideal in S where S is a S-sub-biring of R is said to be a S-right quasi regular ideal or S-left quasi regular ideal if each of its elements is S-right quasi regular or S-left quasi regular respectively.*

*Thus the Jacobson radical of $S \subset R$ as $J(S) = J(S_1) \cup J(S_2) = \{a_1 \in S_1\,/\,a_1 S_1$ is right quasi regular ideal of $S_1\} \cup \{a_2 \in S_2\,/\,a_1 S_2$ is right quasi ideal of $S_2\}$. We call $J(S)$ the Smarandache right quasi regular bi-ideal (S-right quasi regular bi-ideal) related to S.*

*It is important to note that a biring can have several S-right quasi-regular bi-ideals. If each J(S) = 0 for every S a S-sub-biring of R = R$_1$ ∪ R$_2$ then we say R is Smarandache bisimple (S-bisimple).*

We define Smarandache strongly Jacobson radical.

**DEFINITION 6.3.24:** *Let R = R$_1$ ∪ R$_2$ be a na-biring. We have defined S-quasi regular elements. If SJ (R) = S[J(R$_1$)] ∪ S(J(R$_2$)) = {a$_1$ ∈ R / a$_1$ R$_1$ is S-right quasi regular)} ∪ { a$_2$ ∈ R / a$_2$ R$_2$ is S-right quasi regular} then SJ ( R) is called the Smarandache strong Jacobson radical (S-strong Jacobson radical).*

*It is worthwhile to note that S-strong Jacobson radical is unique and if J(R) = 0 we say R = R$_1$ ∪ R$_2$ is Smarandache strongly bisimple (S-strongly bisimple).*

Several interesting relations inter relating them can be had which is left as an exercise to the reader. The study of na-birings arising from groupoid birings and loop birings will be yielding several nice results. As in case of quad birings we can define na-quadrings.

**DEFINITION 6.3.25:** *Let R = R$_1$ ∪ R$_2$ be a biring G = G$_1$ ∪ G$_2$ be a S-bigroupoid. The Smarandache bigroupoid biring (S-bigroupoid biring) is defined as RG = R$_1$G$_1$ ∪ R$_1$G$_2$ ∪ R$_2$G$_1$ ∪ R$_2$G$_2$, a quad ring where each of R$_1$G$_1$, R$_1$G$_2$, R$_2$G$_1$ and R$_2$G$_2$ are groupoid rings. Clearly this quad ring is a na-quad ring.*

On similar lines we define biloop birings, which are na-quad rings. Clearly the birings, which we use, are only associative birings.

**DEFINITION 6.3.26:** *Let R = R$_1$ ∪ R$_2$ be a biring; L = L$_1$ ∪ L$_2$ be a S-biloop. The Smarandache biloop biring (S-biloop biring) RL = R$_1$L$_1$ ∪ R$_1$L$_2$ ∪ R$_2$L$_1$ ∪ R$_2$L$_2$ are loop rings. Clearly RL is a na-quad ring.*

Several properties about these na-quad rings can be defined and analyzed by any researcher.

Now we define when are these na-quad rings, Smarandache quad rings.

**DEFINITION 6.3.27:** *Let R be a na-quad ring i.e. R = R$_1$ ∪ R$_2$ ∪ R$_3$ ∪ R$_4$ where at least some of the R$_i$ are na-ring. We call R a Smarandache weakly na-quad ring (S-weakly na-quad ring) if at least one of the R$_i$'s is a S-ring.*

*We call R is a Smarandache strongly na-quad ring (S-strongly na-quad ring) if each of the R$_i$'s is a S-ring. R is said to be a Smarandache bistrongly na-quad ring (S-bistrongly na-quad ring) if each of the birings is a S-biring. (This has relevance as the quad ring is the union of two birings). R is said to be just a Smarandache biweakly na-quad ring (S-biweakly na-quad ring) if one of the birings is a S-biring.*

All inter-related results can be defined and deduced by any interested reader.

**DEFINITION 6.3.28:** *Let RL be a biloop biring we say RL is a Smarandache strong biloop biring (S-strong biloop biring) if the biloop is a S-biloop and the biring is a S-biring.*

*Similarly we define RG the bigroupoid biring to be a Smarandache strong bigroupoid biring (S-strong bigroupoid biring) if the bigroupoid is a S-bigroupoid and the biring is a S-biring.*

*Example 6.3.6:* Let $L = L_1 \cup L_2$ be a biloop where $L_1 = L_7(3)$ and $L_2 = S_3$. Clearly L is a S-biloop. Take $R = Z_7G \cup R$, R is the field of reals $Z_7G$ the group ring of the group G over the ring $Z_7$. RL is a S-quad ring.

The Smarandache subquad ring is defined as follows:

**DEFINITION 6.3.29:** *Let R be a na-quad ring. We say a proper subset $P \subset R$ is said to be a Smarandache na subquad ring (S-na-subquad ring) if P itself is a na-quad ring.*

*We define the na-quad ring R to be a Smarandache Moufang quad ring if R has a proper S-subquad ring P such that every element of P satisfies the Moufang identity.*

*The quad ring R is said to be Smarandache strong Moufang quad ring (S-strong Moufang quad ring) if every S-subquad ring P of R is a S-Moufang quad ring.*

*On similar lines we can define Smarandache Bol quad ring (S-Bol quad ring), Smarandache Bruck quad ring (S-Bruck quad ring), Smarandache alternative (right / left) quad rings as well as their Smarandache strong concept.*

Thus we see these structures can be easily built using groupoid quad rings and loop quad rings.

Finally we leave the following theorem for the reader to prove.

**THEOREM 6.3.8:** *If R is a quad ring which is S-Moufang quad ring (S-Bol / S-Bruck / S-alternative / S-WIP), then R is a S-quad ring.*

The notions of Smarandache commutators, Smarandache associators, Smarandache center, Smarandache Moufang center, Smarandache nuclei, etc. can be defined in case of S-quad ring.

The task of such study and introduction of these concepts is left as an exercise for the reader to prove.

**PROBLEMS:**

1. Find a S-sub-biring of the bi-ring, $R = Z_7G \cup Z_3S_3$ where $G = \langle g \mid g^2 = 1 \rangle$.

2. What is the smallest order of the S-biring?

3. Is $R = Z_3S_3 \cup Z_9S(7)$ a S-biring? Justify your claim.

4. Can R = $Z_3L_{11}(4) \cup Z_7G$ (where G = $Z_6(3, 2)$ is the groupoid) be a S-na-biring?

5. Does R given in problem 4 have proper S-sub-biring?

6. What is the order of the smallest probable

   i. na-quad ring?
   ii. S-na-quad ring?
   iii. Quad ring?
   iv. S-quad ring?

7. Define ideals of a quad ring.

8. Find any S-ideals in the quad ring RL (where R = $Z_{12} \cup Q$, L = $L_7(2) \cup S_7$) i.e. RL = $Z_{12} L_7(2) \cup Z_{12}S_7 \cup QS_7 \cup QL_7(2)$.

9. Can RL given in problem 8 have S-subna-quad rings which are not S-ideals?

10. Find S-Moufang centre of RL given in problem 8.

11. Is R = $Z_2L_{11}(5) \cup Z_3S_9 \cup Z_7S(5) \cup Z_6(Z_{10}(2, 8))$ a S-quad ring?

12. Find S-ideals, S-subquad ring of R, R given in problem 11.

# Chapter 7

# BISEMIRINGS, S-BISEMIRINGS AND BISEMIVECTOR SPACES AND S-BISEMIVECTOR SPACES

This chapter has four sections. In section one we just define bisemirings and give some of its basic properties. Section two is devoted to the study and introduction of non-associative bisemirings. In section three, which is the vital part of this chapter, the notions of Smarandache bisemirings both associative and non-associative are introduced. In the final chapter we define both the concept of bisemivector spaces and Smarandache bisemivector spaces and give some interesting results about them. In the opinion of the author all concepts studied and introduced in this chapter are new and hence the study and research on them will yield several new notions.

## 7.1 Bisemirings and its properties

In this section we introduce the concept of bisemirings which is carried out in a systematic way. To the best of authors knowledge till date this concept was unknown and unraveled by researchers. In this section we study only associative bisemirings and give methods and means to construct associative bisemirings, polynomial bisemirings and matrix bisemirings.

**DEFINITION 7.1.1:** *Let $(S, +, \bullet)$ be a non-empty set with two binary operations '+' and '$\bullet$'. We call $(S, +, \bullet)$ a bisemiring if*

  i.  $S = S_1 \cup S_2$ where both $S_1$ and $S_2$ are distinct subsets of S, $S_1 \not\subset S_2$ and $S_2 \not\subset S_1$.
  ii. $(S_1, +, \bullet)$ is a semiring.
  iii. $(S_2, +, \bullet)$ is a semiring.

*Example 7.1.1:* Let $S = Z^o \cup C_5$ where $(Z^o = Z^+ \cup \{0\}, +, \bullet)$ is a semiring and $C_5$ which is a chain lattice with five elements. S is a bisemiring.

**DEFINITION 7.1.2:** *Let $(S, +, \bullet)$ be a bisemiring. If the number of elements in S is finite we call S a finite bisemiring. If the number of elements in S is infinite we say S is an infinite bisemiring.*

Thus we see if one of $S_1$ or $S_2$ is an infinite semiring we get the bisemiring S to be an infinite bisemiring. The example 7.1.1 gives a bisemiring of infinite order.

**DEFINITION 7.1.3:** *Let $(S, +, \bullet)$ be a bisemiring. We say S is a commutative bisemiring if both the semirings $S_1$ and $S_2$ are commutative, otherwise we say the bisemiring S is a non-commutative bisemiring.*

**DEFINITION 7.1.4:** *Let $(S, +, \bullet)$ be a bisemiring. We say S is a strict bisemiring if both $S_1$ and $S_2$ are strict bisemirings where $S = S_1 \cup S_2$.*

**DEFINITION 7.1.5:** *Let $(S, +, \bullet)$ be a bisemiring; with $S = S_1 \cup S_2$; we call an element $0 \neq x \in S$ to be a zero divisor if there exist a $y \neq 0 \in S$ such that $x \bullet y = 0$. We say a bisemiring S has a unit 1 in S if $a \bullet 1 = 1 \bullet a = a$ for all $a \in S$.*

**DEFINITION 7.1.6:** *Let $(S, +, \bullet)$ be a bisemiring, we say S is a bisemifield if $S_1$ is a semifield and $S_2$ is a semifield where $S = S_1 \cup S_2$.*

*If both $S_1$ and $S_2$ are non-commutative semirings with no zero divisors then we call $S = S_1 \cup S_2$, to be a bisemidivision ring. It is to be noted that even if one $S_1$ or $S_2$ is non-commutative with no zero divisors but other is a semifield still S is a bisemidivision ring.*

*Example 7.1.2:* Let $(S, +, \bullet)$ be a non-empty set $S = S_1 \cup S_2$ where $S_1 = C_7$ the chain lattice of order 7 and $S_2 =$

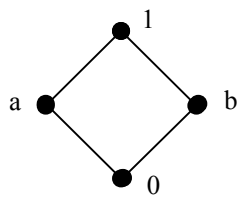

**Figure 7.1.1**

a distributive lattice of order 4. S is a finite bisemiring which is not a bisemifield or a bisemidivision ring but S has units and also S has zero divisors.

*Example 7.1.3:* Let $S = Q^o \cup C_2$ where $Q^o = Q^+ \cup \{0\}$ is a semifield and $C_2$ the chain lattice of order two. S is a bisemifield of infinite order.

Now we proceed on to give examples of non-commutative bisemirings. We will also follow this notation throughout this book.

*Example 7.1.4:* Let $Z^+_{m \times m} = \{(a_{ij}) \mid a_{ij} \in Z^+\} \cup (0)_{m \times m}$ i.e., the set of all $m \times m$ matrices with entries from the set of positive integers together with the zero $(m \times m)$ matrix. Take $S = Z^+_{m \times m} \cup C_3$, clearly S is a non-commutative bisemiring of infinite order.

*Notations:*

$Z^o_{m \times m} = \{(a_{ij}) \mid a_{ij} \in Z^o = Z^+ \cup \{0\}\}$; the set of $m \times m$ matrices.
$Z^+_{m \times m} = \{(a_{ij}) \mid a_{ij} \in Z^+\} \cup (0)_{m \times n}$
$Z'^+_{m \times m} = \{(a_{ij}) \mid a_{ij} \in Z^+\} \cup (0)_{m \times n} \cup I_{m \times m}$.

$Z^o_{m \times m}$ is a semiring non-commutative having zero divisors and unit, $Z^+_{m \times m}$ has no unit but it is a non--commutative semiring with no zero divisors. $Z'^+_{m \times m}$ has unit but is non-commutative with no zero divisors. We see $Z^+_{m \times m} \subset Z'^+_{m \times m} \subset Z^o_{m \times m}$. We can replace $Z^+$ by $Q^+$ or $R^+$ or $Z^o$ by $Q^o$ or $R^o$ and all the results hold good; thus we have several classes of non-commutative infinite semirings using these matrices.

*Example 7.1.5:* $S = Q^+_{m \times m} \cup C_{11}$ is a non-commutative bisemiring which is a bidivision semiring.

**DEFINITION 7.1.7:** *Let $(S = S_1 \cup S_2, +, \bullet)$ be a bisemiring, which is commutative and has unit element. The polynomial bisemiring is denoted by $S[x] = S_1[x] \cup S_2[x]$ where $S_1[x]$ and $S_2[x]$ are polynomial semirings.*

This notion will play a nice role when a polynomial in $S[x]$ may be reducible in $S_1[x]$ and irreducible in $S_2[x]$. Several interesting results in this direction can be derived which is left as an exercise for the reader to solve. Do we have yet any other classes of bisemirings? the answer is yes.

**DEFINITION 7.1.8:** *Let $S = S_1 \cup S_2$ be a bisemiring. Let G by any group. $S_1$ and $S_2$ commutative semirings with unit. $S_1G$ and $S_2G$ be the group semirings of the group G over $S_1$ and $S_2$ respectively.*

*We define the group bisemiring $SG = S_1G \cup S_2G$ which is clearly a bisemiring.*

**Remark:** If G happens to be a commutative group then the group bisemiring is a commutative bisemiring. If on the other hand G happens to be a non-commutative group then we get the group bisemiring to be a non-commutative bisemiring. Thus we get an infinite class of both non-commutative bisemirings of various orders both finite and infinite.

Now we proceed on to define semigroup bisemirings, which still happen to be a new class of bisemirings. This class happens to be a more generalized class then the class of group bisemirings as every group is trivially a semigroup.

**DEFINITION 7.1.9:** *Let $(S, +, \bullet)$ be a bisemiring, which is commutative with unit, P, any semigroup. We define the semigroup bisemiring SP as $SP = S_1P \cup S_2P$ were $S_1P$ and $S_2P$ are semigroup semirings of the semigroup P over the semirings $S_1$ and $S_2$ respectively.*

Depending on the nature of the semigroups the semigroup bisemiring will be finite or infinite, commutative or non-commutative.

Any innovative reader can study several nice results about these structures. As our aim is the study of bialgebraic structures and more on the study of Smarandache bialgebraic structures we do not indulge in this study, as they will deviate the path on some other research direction. But the reader is advised to purse in this direction, as it would help him to obtain some new theme for research. Having given methods to generate bisemirings we proceed on to define some basic properties of bisemirings.

**DEFINITION 7.1.10:** *Let $(S, +, \bullet)$ be a bisemiring. We call a proper subset P of S ($P \subset S$) to be a sub-bisemiring of S, if $(P, +, \bullet)$ itself a bisemiring under the operations of '+' and '$\bullet$'.*

*Example 7.1.6:* Let $S = (Q^o \cup C_7, +, \bullet)$ be a bisemiring. $P = (Z^o \cup \{0, 1\}, +, \bullet)$ is a sub-bisemiring of S.

**DEFINITION 7.1.11:** *Let $(S, +, \bullet)$ be a bisemiring. $(P, +, \bullet)$ be a proper sub-bisemiring of S. We call P an ideal of the bisemiring if for all $s \in S$ and $p \in P$, $s \bullet p$ and $p \bullet s \in P$.*

The concept of right ideal of a bisemiring and left ideal of a bisemiring can also be defined in an analogous way as in case of rings.

**DEFINITION 7.1.12:** *Let $(S, +, \bullet)$ be a bisemiring $S = S_1 \cup S_2$; if both $S_1$ and $S_2$ are semirings of characteristic 0 then we say S is a bisemiring of characteristic 0. If one of $S_i$ alone $(i = 1, 2)$ is of characteristic zero and other $S_i$ undefined or of finite characteristic we say that the semiring has no characteristic associated with it.*

*Example 7.1.7:* Let $S = Z^o \cup C_7$ be a bisemiring. S has no characteristic associated with it.

*Example 7.1.8:* Let $S = Z^o[x] \cup Q^o$ be the bisemiring. The characteristic of the bisemiring is 0.

*Example 7.1.9:* Let $S = C_8 \cup D$ be a bisemiring where $C_8$ is the chain lattice of order 8 and D is the distributive lattice given by the following figure:

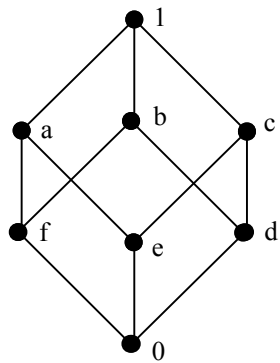

**Figure 7.1.2**

This bisemiring has no characteristic associated with it.

**DEFINITION 7.1.13:** *Let $(S, +, \bullet)$ and $(S', +, \bullet)$ be two bisemirings where $S = S_1 \cup S_2$ and $S' = S'_1 \cup S'_2$. We say a map $\phi : S \to S'$ is a bisemiring homomorphism if $\phi_1: S_1 \to S'_1$ and $\phi_2: S_2 \to S'_2$ is a semiring homomorphism or $\phi_1: S_1 \to S'_2$ and $\phi_2: S_2 \to S'_1$ is a semiring homomorphism where we denote $\phi$ just by default of notation as $\phi = \phi_1 \cup \phi_2$.*

Now we had defined the notion of bisemifield, we define here the concept of sub-bisemifield.

**DEFINITION 7.1.14:** *Let $(S, +, \bullet)$ be a bisemifield. If a proper subset $P \subset S$ is such that $(P, +, \bullet)$ is a bisemifield under the operations of S then we call P a sub-bisemifield. We also call S the extension bisemifield of the bisemifield P.*

***Example 7.1.10:*** Let $S = Z^o[x] \cup C_{20}$; clearly S is a bisemifield. Take $P = Z^o \cup \{0, 1\}$ a proper subset of S, we see P is a sub-bisemifield of S.

***Example 7.1.11:*** Let $S = Z^o \cup C_2$ be a bisemifield. We see S has no sub-bisemifield.

**DEFINITION 7.1.15:** *Let $(S, +, \bullet)$ be a bisemifield. If S has no sub-bisemifield then we call S a prime bisemifield.*

The bisemifield given in example 7.1.11 is a prime bisemifield.

***Example 7.1.12:*** Let $S = Z^o[x] \cup C_3$ be a bisemifield. Clearly S is not a prime bisemifield. For $P = Z^o \cup \{0, 1\}$ is a sub-bisemifield of S.

**PROBLEMS:**

1. Let $(S, +, \bullet)$ be a bisemiring where $S = S_1 \cup S_2$ with $S_1 = C_{10}$ and $S_2 = D$; here $S_1$ is a chain lattice given by the following diagram: Is S a bisemifeild? Justify your claim.

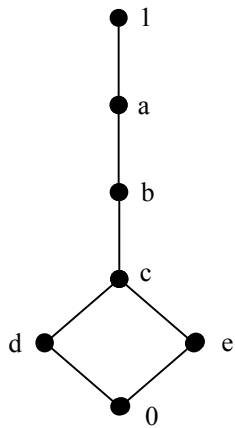

2. Let $S = Z^o[x] \cup Z^+_{3 \times 3}$ be the bisemiring.

   a. Find at least two proper sub-bisemirings.
   b. Does S have ideal?
   c. Find atleast two right ideals of S.
   d. Is S a bisemifield?

3. Let $S = C_2[x] \cup Z^o_{2 \times 2}$ be a bisemiring. Does S have nontrivial idempotents? units?

4. Give an example of a bisemiring of order 7.
5. Can we have bisemirings of any order? Justify your claim.
6. What is the order of the smallest non-prime bisemifield?
7. What is the order of the smallest prime bisemifield?
8. Does their exist a bisemifield of order 8? Justify your claim.

9. Let $(S = S_1 \cup S_2, +, \bullet)$ be a bisemiring where $S_1 = C_2S_3$ (the group semiring of the group $S_3$ over $C_2$) and $S_2 = Z^oS(3)$ (The semigroup semiring of the semigroup $S(3)$ over $Z^o$).

   a. Is S a commutative bisemiring?
   b. Does S have units?
   c. Find idempotents in S?
   d. Find right ideals in S.

10. Is the bisemiring given in problem (9) a group bisemiring? Justify!

11. Is $(S = S_1 \cup S_2, +, \bullet)$ where $S_1 = Z^oG$ and $S_2 = C_3G$ (with $G = \langle g \mid g^7 = 1 \rangle$) a group bisemiring? Is S commutative?

12. Is the bisemiring given in problem 11 a bisemifield?

## 7.2 Non-associative bisemirings and its properties

In this section we introduce the concept of non-associative bisemirings. At the first instant we state that natural construction of non-associative bisemirings is still an open problem. We construct non-associative bisemiring only by using loops or groupoids and associative bisemirings. Throughout this section we denote a non-associative bisemiring by na-bisemiring.

**DEFINITION 7.2.1:** *Let $(S, +, \bullet)$ be a bisemiring. We say S is a non-associative bisemiring if $S = S_1 \cup S_2$ where $S_1$ and $S_2$ are semirings where at least one of $S_1$ or $S_2$ is a non-associative semiring.*

**DEFINITION 7.2.2:** *Let $(S, +, \bullet)$ be an associative bisemiring. Let L be a loop, the loop bisemiring $SL = S_1L \cup S_2L$ (where $S = S_1 \cup S_2$; $S_1$ and $S_2$ are commutative semirings with unit) and $S_1L$ and $S_2L$ are the loop semirings of the loop L over the semirings $S_1$ and $S_2$ respectively. Clearly the loop bisemiring is a non-associative bisemiring.*

*Example 7.2.1:* Let $(S, +, \bullet)$ be a bisemiring. $S = S_1 \cup S_2$ where $S_1 = Z^oL_7(3)$ and $S_2 = C_3S_3$. Clearly the bisemiring is a non-associative bisemiring of infinite order. Clearly this bisemiring is not a loop bisemiring or a group bisemiring.

**DEFINITION 7.2.3:** *Let $(S, +, \bullet)$ be a commutative bisemiring with 1 and G any groupoid. The groupoid bisemiring $SG = S_1G \cup S_2G$ is a non-associative bisemiring. ($S_1G$ and $S_2G$ are the groupoid semirings of the groupoid G over the semirings $S_1$ and $S_2$ respectively).*

It is pertinent to mention here that if we take $S = S_1G \cup S_2H$ where G and H are two distinct groupoids with $S_1G$ and $S_2H$ the groupoid semirings. S is not a groupoid bisemiring but only a na-bisemiring.

The notions like sub-bisemiring, ideals (right or left or two-sided), zero-divisors, idempotents, units can be defined in an analogous way with least modifications, as associativity does not play any major role in the study of these concepts.

**DEFINITION 7.2.4:** *Let $(S, +, \bullet)$ be a non-associative bisemiring. We say S is a strict bisemiring if $a + b = 0$ then $a = 0$ and $b = 0$. We call $(S, +, \bullet)$ a na-bisemiring, to be a commutative bisemiring if both $S_1$ and $S_2$ are commutative semirings. $(S, +, \bullet)$ is said to be a semiring with unit if there exists $1 \in S$ such that $1 \bullet s = s \bullet 1 = s$ for all $s \in S$. An element $0 \neq x \in S$ is said to be a zero divisor if there exists $0 \neq y \in S$ such that $x \bullet y = 0$. A na-bisemiring which has no zero-divisors but which is commutative with unit is called a non-associative bisemifield. If the operation in S is non-commutative we call S a na-bidivision semiring.*

Now we proceed on to define notions like special identities satisfied by bisemirings, which are non-associative.

**DEFINITION 7.2.5:** *Let $(S, +, \bullet)$ be a na-bisemiring. We say S is a Moufang bisemiring if all elements of S satisfy the Moufang identity i.e. $(xy)(zx) = (x(yz))x$ for all $x, y, z \in S$. A na-bisemiring S is said to be a Bruck bisemiring if $(x(yz))z = x(y(xz))$ and $(xy)^{-1} = x^{-1}y^{-1}$ for all $x, y, z \in S$. A na-bisemiring S is called a Bol bisemiring if $((xy)z)y = x((yz)y)$ for all $x, y, z \in S$. We call a na-bisemiring S to be right alternative if $(xy)y = x(yy)$ for all $x, y \in S$; left alternative if $(xx)y = x(xy)$ and alternative if it is both right and left alternative.*

**DEFINITION 7.2.6:** *Let $(S, +, \bullet)$ be a na-bisemiring. $x \in S$ is said to be right quasi regular (r.q.r) if there exists a $y \in R$ such that $x \circ y = 0$ and $x$ is said to be left quasi regular (l.q.r) if there exists a $y' \in R$ such that $y' \circ x = 0$.*

*An element is quasi regular (q.r) if it is both right and left quasi-regular.*

*y is known as the right quasi inverse (r.q.i) of x and y' as the left quasi inverse (l.q.i) of x. A right ideal or left ideal in R is said to be right quasi regular (l.q.r or q.r respectively) if each of its elements is right quasi regular (l.q.r or q.r respectively).*

**DEFINITION 7.2.7:** *Let S be a na-bisemiring. An element $x \in S$ is said to be right regular (left regular) if there exists a $y \in S$ ($y' \in S$) such that $x(yx) = x$ (($xy'$)$x = x$).*

**DEFINITION 7.2.8:** *Let S be a na-bisemiring. The Jacobson radical $J(S)$ of a bisemiring S is defined as follows: $J(S) = \{x \in S \, / \, xS$ is right quasi-regular ideal of $S\}$. A bisemiring is said to be semisimple if $J(S) = \{0\}$ where $J(S)$ is the Jacobson radical of S.*

**DEFINITION 7.2.9:** *Let S be a na-bisemiring. We say S is prime if for any two ideals A, B in S, $AB = 0$ implies $A = 0$ or $B = 0$.*

***Example 7.2.2:*** Let $S = S_1G \cup S_2G$ be a bisemiring where $S_1$ and $S_2$ are semirings given by the following figures:

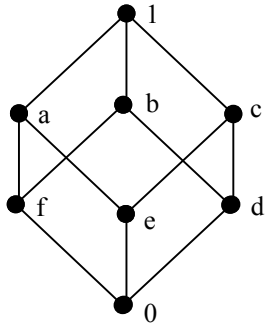

$S_1$

**Figure 7.2.1**

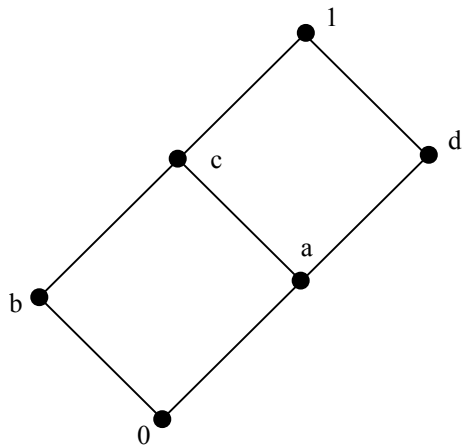

$S_2$

**Figure 7.2.2**

and G is the groupoid given by the following table:

| · | e | $a_0$ | $a_1$ | $a_2$ |
|---|---|---|---|---|
| e | e | $a_0$ | $a_1$ | $a_2$ |
| $a_0$ | $a_0$ | e | $a_2$ | $a_1$ |
| $a_1$ | $a_1$ | $a_2$ | e | $a_0$ |
| $a_2$ | $a_2$ | $a_1$ | $a_0$ | e |

$SG = S_1G \cup S_2G$ is a na-bisemiring having non-trivial zero divisors and idempotents.

**PROBLEMS:**

1. Find the order of the smallest na-bisemiring.
2. Is $S = Z^\circ L_5(3) \cup C_7$ a na-bisemiring?
3. Is S given in problem 2 a na-loop bisemiring? Justify your claim.

4. Is $S = C_7G \cup Q^oG$ where G is a groupoid given by the following table:

   | • | $a_o$ | $a_1$ | $a_2$ | $a_3$ |
   |---|---|---|---|---|
   | $a_0$ | $a_0$ | $a_2$ | $a_0$ | $a_2$ |
   | $a_1$ | $a_2$ | $a_0$ | $a_2$ | $a_0$ |
   | $a_2$ | $a_0$ | $a_2$ | $a_0$ | $a_2$ |
   | $a_3$ | $a_2$ | $a_0$ | $a_2$ | $a_0$ |

   a groupoid bisemiring? Justify your claim.

5. Let $S = S_1 \cup S_2$ be a bisemiring where $S_2$ and $S_2$ are given by the following diagrams:

   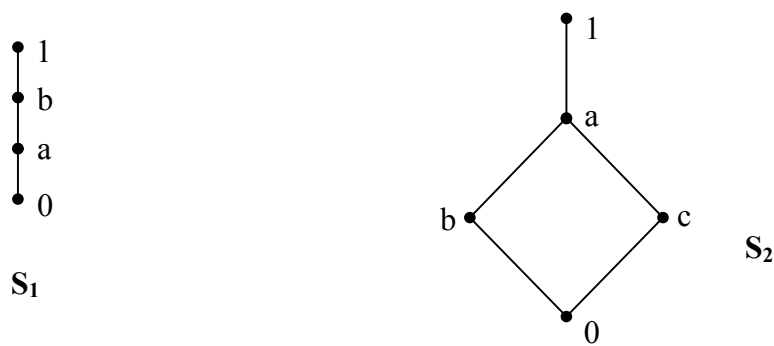

   $SL_7(3) = S_1L_7(3) \cup S_2L_7(3)$ is a loop bisemiring which is a na-bisemiring.

   i. Does $SL_7(3)$ have zero divisors?
   ii. Find idempotents if any in $SL_7(3)$.
   iii. Is $SL_7(3)$ a commutative na-bisemiring?
   iv. Find the order of $SL_7(3)$.
   v. Does $SL_7(3)$ have ideals?
   vi. Can $SL_7(3)$ have sub-bisemirings, which are not ideals?
   vii. Does it satisfy any of the standard identities like Bol, Bruck, and Moufang?

6. Let $S = Z^0L_{11}(2) \cup C_3L_5(2)$ be the na-bisemiring. Find ideals in S.

   a. What is the characteristic of S?
   b. Is S a strict na-bisemiring.

7. Give an example of a finite na-bisemiring of order 14 having no zero divisors, no units but having idempotents and non-trivial sub-bisemirings.

## 7.3 S-bisemirings and its properties

In this section we introduce the notion of Smarandache bisemirings and several of its properties. Studies of S-semirings have been introduced in 2002 [122, 123]. Here we give several types of results about S-bisemirings, which are associative as well as non-associative.

**DEFINITION 7.3.1:** *Let (S, +, •) be a bisemiring. We call (S, +, •) a Smarandache bisemiring (S-bisemiring) if S has a proper subset P such that P under the operations of S is a bisemifield.*

*Example 7.3.1:* Let $(S, +, •)$ be a bisemiring, where $S = S_1 \cup S_2$ with $S_1 = Q^o[x]$ and $S_2 = C_3[x]$. Then S is a S-bisemiring as $P = P_1 \cup P_2$ is a bisemifield where $P_1 = Q^o$ and $P_2 = C_3$.

All bisemirings are not in general S-bisemirings.

*Example 7.3.2:* Let $(S, +, •)$ be a bisemiring. $S = S_1 \cup S_2$ where $S_1 = Z^o$ and $S_2 = C_2$. Clearly the bisemiring is not a S-bisemiring.

If a S-bisemiring has only a finite number of elements we say the S-bisemiring is finite otherwise infinite.

**DEFINITION 7.3.2:** *Let S be a bisemiring. A non-empty proper subset A of S is said to be a Smarandache sub-bisemiring (S-sub-bisemiring) if A is a S-bisemiring i.e. A has a proper subset P such that P is a bisemifield under the operations of S.*

**THEOREM 7.3.1:** *If (S, +, •) is a bisemiring having a S-sub-bisemiring then S is a S-bisemiring.*

*Proof:* Follows from the very definitions.

It is left for the reader to prove that if $(S, +, •)$ is S-bisemiring then every sub-bisemiring of S need not be a S-sub-bisemiring of S.

*Example 7.3.3:* Let $(S, +, •)$ be a S-bisemiring. $S = Q^o \cup C_7$. Take $P = \{0, 2, 4, 6, \ldots\} \cup \{0, 1\} \subset S$; P is a sub-bisemiring of S which is clearly not a S-sub-bisemiring of S.

**DEFINITION 7.3.3:** *Let S be a S-bisemiring. We say S is a Smarandache commutative bisemiring (S-commutative bisemiring) if S has a S-sub-bisemiring, which is commutative. If S has no commutative S-sub-bisemiring then we call S to be a Smarandache non-commutative bisemirings (S-non-commutative bisemirings). If every S-sub-bisemiring of S is commutative then we call S a Smarandache strongly commutative bisemirings (S-strongly commutative bisemirings).*

**DEFINITION 7.3.4:** *Let (S, +, •) be a bisemiring. A non-empty subset P of S is said to be a Smarandache right (left) bi-ideal (S-right (left) bi-ideal) of S if the following conditions are satisfied.*

  i.  *P is a S-sub-bisemiring.*
  ii. *For every $p \in P$ and $A \subset P$ where A is a bisemifield of P we have for all $a \in A$ and $p \in P$, $a • p$ ($p • a$) is in A. If P is simultaneously both a S-right bi-ideal and S-left bi-ideal then we say P is a Smarandache bi-ideal (S-bi-ideal) of S.*

*Example 7.3.4:* Let

$$M_{2\times 2} = \left\{ \begin{pmatrix} a & b \\ c & d \end{pmatrix} \middle/ a,b,c,d \in C_2 = [0,1] \right\}$$

= set of all $2 \times 2$ matrices with entries from the chain lattice $C_2$. Take $S = M_{2\times 2} \times C_6$, $P = A \cup P_1$ where

$$A = \left\{ \begin{pmatrix} 1 & 0 \\ 0 & 0 \end{pmatrix}, \begin{pmatrix} 1 & 0 \\ 0 & 1 \end{pmatrix}, \begin{pmatrix} 0 & 0 \\ 0 & 0 \end{pmatrix} \right\} \text{ and } P_1 = \{0, a_1, 1\}$$

is a S-sub-bisemiring of S, which is not an S-bi-ideal of S.

**DEFINITION 7.3.5:** *Let $(S, +, \bullet)$ be a bisemiring. A proper subset A of S is said to be Smarandache pseudo sub-bisemiring (S-pseudo sub-bisemiring) if the following conditions are true:*

*If there exists a subset of P of S such that $A \subset P$, where P is a S-sub-bisemiring, i.e. P has a subset B such that B is a bisemifield under the operations of S or P itself is a bisemifield under the operations of S.*

**DEFINITION 7.3.6:** *Let S be a bisemiring. A non-empty subset P of S is said to be a Smarandache pseudo right (left) bi-ideal (S-pseudo right (left) bi-ideal) of the bisemiring S if the following conditions are true:*

  i.  *P is a S-pseudo sub-bisemiring i.e. $P \subset A$, A is a bisemifield in S.*
  ii. *For every $p \in P$ and every $a \in A$, $a \bullet p \in P$ ($p \bullet a \in P$). If P is simultaneously both a S-pseudo right bi-ideal and S-pseudo left bi-ideal then we say P is a Smarandache pseudo bi-ideal (S-pseudo bi-ideal).*

**DEFINITION 7.3.7:** *Let S be a bisemiring. A nonempty subset P of S is said to be a Smarandache dual bi-ideal (S-dual bi-ideal) of S if the following conditions hold good:*

  i.  *P is a S-sub-bisemiring.*
  ii. *For every $p \in P$ and $a \in A \setminus \{0\}$; $a + p$ is in A, where $A \subset P$.*

**DEFINITION 7.3.8:** *Let S be a bisemiring. A nonempty subset P of S is said to be a Smarandache pseudo dual bi-ideal (S-pseudo dual bi-ideal) of S if the following conditions are true:*

  i.  *P is a S-pseudo sub-bisemiring i.e. $P \subset A$, A is a bisemifield in S or A contains a bisemifield.*
  ii. *For every $p \in P$ and $a \in A$, $p + a \in P$. Clearly P is simultaneously left and right S-pseudo dual bi-ideal of S as S is additively commutative.*

**DEFINITION 7.3.9:** *Let S be a S-bisemiring. S is said to be a Smarandache bisemidivision ring (S-bisemidivision ring) if the proper subset $A \subset S$ is such that*

> i. *A is a S-sub-bisemiring.*
> ii. *A contains a subset P such that P is a bisemidivision ring, that is P has no zero divisors and P is a non-commutative bisemiring.*

The concept of S-zero divisors, S-units and S-idempotents in bisemirings is defined analogous to rings and for this we do not even demand the bisemiring to be a S-bisemiring. Several interesting results can be got in this direction.

**DEFINITION 7.3.10:** *Let S be a bisemiring we say S is a Smarandache right chain bisemiring (S-right chain bisemiring) if the S-right bi-ideals of S are totally ordered by inclusion.*

*Similarly we define Smarandache left chain bisemirings (S-left chain bisemirings). If all the S-bi-ideals of the bisemiring are ordered by inclusion we say S is a Smarandache chain bisemiring (S-chain bisemiring).*

**DEFINITION 7.3.11:** *Let S be a bisemiring. If $S_1 \subset S_2 \subset ...$ is a monotonic ascending chain of S-bi-ideals $S_i$ and there exists a positive integer r such that $S_r = S_n$ for all $r \geq n$, then we say the bisemiring S satisfies the Smarandache ascending chain conditions (S-acc) for S-bi-ideals in the bisemiring S.*

*We say S satisfies Smarandache descending chain conditions (S-dcc) on S-bi-ideals $S_i$, if every strictly decreasing sequence of S-ideals $N_1 \supset N_2 \supset ...$ in S is of finite length. The Smarandache min conditions (S-min conditions) for S-bi-ideals holds in S if given any set P of S-bi-ideals in S, there is a bi-ideal of P that does not properly contain any other bi-ideal in the set P. Similarly one can define Smarandache maximum condition (S-maximum condition) for S-bi-ideals in case of bisemirings.*

**DEFINITION 7.3.12:** *Let S be a bisemiring. S is said to be a Smarandache idempotent bisemiring (S-idempotent bisemiring) if a proper subset P of S that is a sub-bisemiring of S satisfies the following conditions:*

> i. *P is a S-sub-bisemiring.*
> ii. *P is an idempotent bisemiring.*

**DEFINITION 7.3.13:** *Let S be a bisemiring. S is said to be a Smarandache e-bisemiring (S-e-bisemiring) if S contains a proper subset A satisfying the following conditions:*

> i. *A is a sub-bisemiring.*
> ii. *A is a S-sub-bisemiring.*
> iii. *A is a e-bisemiring.*

**DEFINITION 7.3.14:** *Let S be any bisemiring. G be a S-semigroup. Consider the semigroup bisemiring SG, we call SG the Smarandache group bisemiring (S-group bisemiring) i.e. $SG = S_1G \cup S_2G$ where each of $S_1G$ and $S_2G$ are S-group semirings.*

**DEFINITION 7.3.15:** *Let G be a group, G is a Smarandache antigroup (S-antigroup), that is if G contains a proper subset, which is a semigroup. The group bisemiring FG where F is any bisemiring and G a S-antigroup is called the Smarandache semigroup bisemiring (S-semigroup bisemiring).*

**DEFINITION 7.3.16:** *A bisemiring S is said to be Smarandache bisemiring of level II (S-bisemiring of level II) if S contains a proper subset P that is a bifield.*

*Example 7.3.5:* A bisemiring S is said to be S-bisemiring II if $S = S_1 \cup S_2$ where $S_1$ and $S_2$ are S-semiring of level II.

**DEFINITION 7.3.17:** *Let R be a biring. R is said to be a Smarandache anti-bisemiring (S-anti-bisemiring) if R contains a subset S such that S is just a bisemiring.*

Several interesting results can be derived from these concepts, which is left for the reader.

**DEFINITION 7.3.18:** *Let S be a bisemifield. $S = S_1 \cup S_2$ is said to be a Smarandache bisemifield (S-bisemifield) if a proper subset $P = P_1 \cup P_2$ of S is a S-bisemialgebra with respect to the same induced operations and an external operator (i.e. $P_1$ is a S-semialgebra and $P_2$ is a k-semialgebra, $P = P_1 \cup P_2$ is a S-bisemialgebra).*

We define Smarandache bisemifields of level II.

**DEFINITION 7.3.19:** *Let S be a bisemifield. S is said to be a Smarandache bisemifield of level II (S-bisemifield of level II) if S has a proper subset P where P is a bifield.*

**DEFINITION 7.3.20:** *Let S be a bisemifield. A proper subset P of S is said to be a Smarandache sub-bisemifield I (II) (S-sub-bisemifield I (II)) if P is a S-bisemifield of level I (or level II).*

**DEFINITION 7.3.21:** *Let S be a bifield or a biring. S is said to be a Smarandache anti-bisemifield (S-anti-bisemifield) if S has a proper subset, A which is a bisemifield.*

Several innovative results can be developed in this direction.

**DEFINITION 7.3.22:** *Let S be a biring or a bifield. A proper subset P in S is said to be a Smarandache anti sub-bisemifield (S-anti sub-bisemifield) of S if P is itself a S-anti bisemifield.*

**DEFINITION 7.3.23:** *Let S be a bifield or a biring, which is a S-anti bisemifield. If we can find a subset P in the sub-bisemifield T in S such that*

   i.     *P is a bisemiring.*
   ii.    *for all $p \in P$ and $t \in T$, $pt \in P$, then P is called the Smarandache anti bi-ideal (S-anti bi-ideal) of the S-anti-bisemifield.*

**Note**: We cannot have the concept of right or left ideals as the sub-bisemifield is commutative.

Now we proceed on to define and develop non-associative S-bisemirings.

The concept of non-associative semirings is very rare and they cannot be defined on $Q^o$ or $Z^o$ or L or $C_n$. We need at least two algebraic structures to define non-associative semirings. Thus all the more difficult is the concept of non-associative bisemirings. This book gives the definition of non-associative bisemirings and now we study the notion of non-associative Smarandache bisemirings.

**DEFINITION 7.3.24:** *Let (S, +, •) be a na-bisemiring. We say S is a Smarandache na bisemiring (S-na bisemiring) if S has a proper subset P, (P $\subset$ S) such that P is an associative bisemiring.*

*Example 7.3.6:* Let (S, +, •) be a na-bisemiring, where S = $S_1 \cup S_2$ with $S_1 = Z^o L_5(3)$ and $S_2 = C_8$. Take P = $Z^o \cup \{0, 1\}$. Clearly P is an associative bisemiring of S. Hence S is a S-na-bisemiring.

**DEFINITION 7.3.25:** *Let (S, +, •) be a na- bisemiring. A proper subset P of S is said to be a Smarandache na sub-bisemiring (S-na sub-bisemiring) if P itself under the operations of S is a S-na-bisemiring.*

**DEFINITION 7.3.26:** *Let (S, +, •) be a na bisemiring. We say S is a Smarandache commutative na bisemiring (S-commutative na bisemiring) if S has a S-sub-bisemiring P that is commutative. If every S-sub-bisemiring P of S is commutative we call S a Smarandache strongly commutative bisemiring (S-strongly commutative bisemiring) even if S is commutative and S has no S-sub-bisemiring then also S is not Smarandache commutative (S-commutative).*

**THEOREM 7.3.2:** *If (S, +, •) is a S-strongly commutative bisemiring then S is a S-commutative bisemiring and not conversely.*

*Proof:* Follows directly by the definitions.

The concept of S-ideal, S-zero divisors, S-idempotents, S-units and S-anti zero divisors can be defined as in case of associative bisemirings. The only notions, which is not present in any associative bisemiring, is the concept of special identities.

**DEFINITION 7.3.27:** *Let (S, +, •) be a na bisemiring. We call S a Smarandache Moufang bisemiring (S-Moufang bisemiring) if S has a proper subset P where P is a S-sub-bisemiring of S and P satisfies the Moufang identities i.e. P is a Moufang bisemiring.*

**Note**: When we say P is a S-sub-bisemiring we assume P is a S-sub-bisemiring we assume P is not the associative bisemiring but P only contains a proper subset L such that only L is an associative bisemiring of P.

On similar lines we define Smarandache Bol bisemiring, Smarandache Bruck bisemiring, Smarandache right (left) alternative bisemiring, Smarandache alternative bisemiring and so on.

We also define a stronger concept called Smarandache strong Moufang bisemiring, Smarandache strong Bol bisemiring and so on.

**DEFINITION 7.3.28:** *Let (S, +, •) be a non-associative bisemiring.*

*If every S-sub-bisemiring P of S satisfies the Moufang identity i.e. every P is a Moufang bisemiring then we call S a Smarandache strong Moufang bisemiring (S-strong Moufang bisemiring).*

**THEOREM 7.3.3:** *Let (S, +, •) be a S-strong Moufang bisemiring, then S is a S-Moufang bisemiring and not conversely.*

*Proof:* Follows from the very definitions directly. To prove the converse the reader is requested to construct a counter example.

The notions of Smarandache P-bisemiring and Smarandache Jordan bisemiring can also be defined in a similar way. We have not discussed how the set of S-bi-ideals of the bisemiring look like, i.e. do they form a distributive lattice or a modular lattice and so on.

Further the study of the set of bi-ideals of a bisemiring is an unseen problem in case of bistructures.

**PROBLEMS:**

1. Is (S, +, •) a S-bisemiring where S = $Z^oG \cup C_2S_3$ where G = $\langle g / g^2 = 1 \rangle$?

2. For the S-bisemiring (S, +, •) where S = $C_2S(2) \cup DG$ where D is the lattice -

    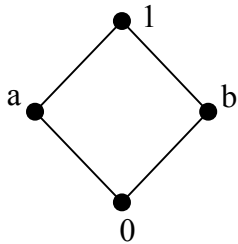

    and G = $\langle g / g^3 = 1 \rangle$. Find S-subsemiring.

3. Does S in problem 2 have

    i. S-units?
    ii. S-zero divisors?
    iii. S-ideals?

4. Let (S, +, •) be a na-bisemiring given by S= $S_1 \cup S_2$ where $S_1$ = $C_2L_5(2)$ and $S_2 = C_3S(3)$.

    i. Is S a S-na bisemiring?
    ii. Is S a S-Moufang bisemiring?
    iii. Find all S-ideals of S.

    iv. How many S-sub-bisemiring does S contain?

5. Give an example of na-bisemiring, which satisfies the S-Bol identity.
6. Show by an example that all bisemirings are not in general S-bisemiring.
7. Give an example of S-bisemiring II.
8. Give an example of a finite S-bisemifield.
9. Can a S-bisemiring of order 3 exist? Justify your claim.
10. What is order of the smallest S-bisemiring, which is non-associative?
11. Can a S-bisemiring of order 11 exist?
12. Find an example of a bisemiring, which has S-pseudo dual bi-ideal.

## 7.4 Bisemivector spaces and S-bisemivector spaces

In this section we introduce the concept of bisemivector spaces and their Smarandache analogue. The study of semivector space is very meager, so the notion of bisemivector spaces is practically not defined or studied till date. Also as each chapter of this book can be written as a separate book we mainly give the basic notions and a few of its fundamental properties. The development of the subject is left for the reader. As it may not be practically possible to dwell into each and every one of its properties, as this book mainly tries to introduce all bialgebraic structures and their Smarandache analogue.

**DEFINITION 7.4.1:** *A strict semigroup under '+' is one which is a semigroup in which $a + b = 0$ implies $a = 0$ and $b = 0$. A strict semigroup is commutative if $a + b = b + a$. We assume throughout this section by a strict semigroup it is a commutative semigroup with 0.*

**Example 7.4.1:** $Z^o$ is a strict semigroup under '+'.

**Example 7.4.2:** $Z^o[x]$ is a strict semigroup under '+'.

**Example 7.4.3:** $Z^o_{m \times n}$ under '+' is a strict semigroup.

**DEFINITION 7.4.2:** *Let $S = S_1 \cup S_2$ be a bisemigroup. We say S is a strict bisemigroup if both $S_1$ and $S_2$ are strict semigroups under '+' and contains zero which is commutative.*

**DEFINITION 7.4.3:** *Let S be a strict semigroup. S is a semivector space over the semifield F if for all $s \in S$ and $\alpha \in F$*

  i. $\alpha s \in S$.
  ii. $(\alpha \beta)s = \alpha(\beta s)$.
  iii. $\alpha(s_1 + s_2) = \alpha s_1 + \alpha s_2$;

*for all $\alpha, \beta \in F$ and $s, s_1, s_2 \in S$. Further $(\alpha + \beta)s = \alpha s + \beta s$ and $1.s = s$ where $1 \in F$ and $s \in S$.*

Now we are in a position to define bisemivector spaces.

**DEFINITION 7.4.4:** *Let $V = V_1 \cup V_2$ be a strict bisemigroup. S be a semifield. We say V is a bisemivector space over the semifield F if both $V_1$ and $V_2$ are semivector spaces over the semifield F.*

*Example 7.4.4:* $V = Z^o[x] \cup Z^o_{3\times 2}$ be a strict bisemigroup. $Z^o$ is a semifield. V is a bisemivector space over the semifield $Z^o$.

It is worthwhile to note that V is not a bisemivector space over the semifield $Q^o$ or $R^o$. Thus as in case of vector spaces we see even in case of bisemivector spaces, the definition depends on the related semifield over which we try to define. Now we proceed on to define several new concepts in these bisemivector spaces.

**DEFINITION 7.4.5:** *Let $V = V_1 \cup V_2$ be a bisemivector space over a semifield S. A set of vectors $\{v_1, ..., v_n\} \in V$ is said to be linearly dependent if $A = \{v_1, ..., v_k\}$, $B = \{v_k, ..., v_n\}$ are subsets of $V_1$ and $V_2$ respectively and if there exists a non-trivial relation among them i.e. some $v_j$ is expressible as a linear combination from the vectors in A and $v_i$ is expressible as a linear combination of vectors from B. A set of vectors, which is not linearly dependent, is called linearly independent.*

Several results in this direction can be got as in case of semivector spaces.

**DEFINITION 7.4.6:** *Let $V = V_1 \cup V_2$ be a bisemivector space over the semifield F. A linearly independent set $P = A \cup B$ spanning the bisemivector space V is called the basis of V.*

It is interesting to note that unlike vector spaces, the basis in certain bisemivector spaces are unique as evidenced by the following example.

*Example 7.4.5:* Let $V = V_1 \cup V_2$, where $V_1 = Z^o \times Z^o$ and $V_2 = Z^o_8[x]$. Clearly V is a bisemivector space over $Z^o$ and the unique basis for this bisemivector space is $P = \{(1, 0), (0, 1)\} \cup \{1, x, x^2, \ldots, x^8\}$. $Z^o_8[x]$ denotes the collection of all polynomials of degree less than or equal to 8 with coefficient from $Z^o$. This semivector space has no other basis; P is its only basis.

**THEOREM 7.4.1:** *In a bisemivector space V the number of elements in a set, which spans V, need not be an upper bound for the number of vectors that can be linearly independent in V.*

*Proof*: By an example, Let $V = V_1 \cup V_2$, $V_1 = Z^o \times Z^o$ and $V_2 = Z^o_3[x] = \{$all polynomials of degree less than or equal to 3$\}$. V is a bisemivector space over $Z^o$. $P = \{(0, 1), (1, 0)\} \cup \{1, x, x^2, x^3\}$ is a basis which span V. But if we take $P_1 = \{(1, 1), (2, 1), (3, 0)\} \cup \{1, x, x^2, x^3\}$. Clearly it is easily verified, $P_1$ is a linearly independent set in V; but not a basis of V as it does not span V. This property is a very striking difference between a bivector space and a bisemivector space.

The following theorem is left for the reader to prove.

**THEOREM 7.4.2:** *In a n-dimensional bisemivector space we need not have in general, every set of (n + 1)- vectors to be linearly dependent.*

Now several important properties can be built using bisemivector spaces over semifields. Now we are going to build several nice results, which we are going to call as bipseudo semivector spaces.

**DEFINITION 7.4.7:** *Let V be a strict semigroup. Let $S = S_1 \cup S_2$ be a bisemifield. If V is a semivector space over both the semifields $S_1$ and $S_2$, then we call V a bipseudo semivector space.*

**Example 7.4.6:** $Q^o[x]$ is a bipseudo semivector space over the bisemifield $Z^o[x] \cup Q^o$.

**Example 7.4.7:** $R^o[x]$ is a bipseudo semivector space over the bisemifield $Z^o[x] \cup Q^o$.

**DEFINITION 7.4.8:** *Let $V = V_1 \cup V_2$ be a bisemivector space over the semifield S. A proper subset $P \subset V$ is said to be a sub-bisemivector space if P is itself a bisemivector space over S.*

**Example 7.4.8:** Take $V = R^o[x] \cup \{Q^o \times Q^o\}$ is a bisemivector space over $Z^o$. Clearly $P = R^o \cup \{Z^o \times Z^o\}$ is a sub-bisemivector space over $Z^o$.

Several such examples can be had. The concept of isomorphism does not exist for the same dimensional bisemivector spaces even over same semifields.

**DEFINITION 7.4.9:** *Let $V = V_1 \cup V_2$ and $W = W_1 \cup W_2$ be any two bisemivector spaces defined over the same semifield S. A map $T: V \to W$ is called the bisemivector space homomorphism or linear transformation if $T = T_1 \cup T_2$ where $T_1: V_1 \to W_1$ and $T_2: V_2 \to W_2$ are linear transformations of semivector spaces $V_1$ to $V_2$ and $W_1$ to $W_2$.*

Several results in this direction can be made as the field is at the dormant state.

Now the sub-bipseudo semivector space is defined and analogously its transformations are introduced.

**DEFINITION 7.4.10:** *Let V be a bipseudo semivector space over the bisemifield $S = S_1 \cup S_2$. A proper subset P of V is said to be a sub-bipseudo semivector space if P is itself a bipseudo semivector space over S.*

**DEFINITION 7.4.11:** *Let V and W be two bipseudo semivector spaces over the same bisemifield $S = S_1 \cup S_2$. A map $T: V \to W$ is the bipseudo semivector space homomorphism (or a transformation $T: V \to W$) if T is a semivector transformation from V onto W as semivector space over $S_1$ and T a transformation as semivector spaces over V onto W as semivector spaces over $S_2$.*

The notions of Smarandache bisemivector spaces are defined in the following. Now the notion of Smarandache bisemivector spaces will be introduced. To do so first we define the concept of strict S-semigroup.

**DEFINITION 7.4.12:** *Let (S, +) be a semigroup. We say S is a Smarandache strict semigroup (S-strict semigroup) if S has a proper subset P such that (P, +) is a strict semigroup, so P is commutative and has unit.*

**DEFINITION 7.4.13:** *Let $V = V_1 \cup V_2$, we say the bisemigroup is a Smarandache strict bisemigroup (S-strict bisemigroup) if V has a proper subset P such that $P \subset V$, $P = P_1 \cup P_2$ where both $P_1$ and $P_2$ are strict semigroups, i.e. P is a strict bisemigroup.*

*Example 7.4.9:* Let $V = Z[x] \cup Z_{2 \times 3}$ i.e. the polynomials over Z under '+' and the set of all $2 \times 3$ matrices with entries from the integers under '+'. Clearly V is not a strict bisemigroup but V is a S-strict bisemigroup as take $P \subset V$ with $P = P_1 \cup P_2$ where $P_1 = Z^o[x]$ and $P_2 = Z^0_{2\times 3}$. Then P is a strict semigroup so V is a S-strict semigroup.

**DEFINITION 7.4.14:** *Let $V = V_1 \cup V_2$ be a bisemigroup. V is said to be a Smarandache bisemivector space (S-bisemivector space) over the semifield F if*

  i.  *both $V_1$ and $V_2$ are such that for each $v_i \in V_i$ and $\alpha \in F$, $\alpha v_i \in V_i$ for i = 1, 2.*
  ii. *both $V_1$ and $V_2$ has proper subspaces say $P_1$ and $P_2$ respectively such that $P_1 \cup P_2 = P$ is a strict bisemivector space over F;*

*or equivalently we can say that in $V = V_1 \cup V_2$ we have $P = P_1 \cup P_2$ where P is a strict bisemigroup and P is a bisemivector space over F. The concept of Smarandache sub-bisemivector space (S-sub-bisemivector space) can be introduced in a similar way.*

**DEFINITION 7.4.15:** *Let V be a S-strict semigroup. S be a S-bisemifield. If $P \subset V$ where P is a strict semigroup and is a bipseudo semivector space over S then we say P is a Smarandache bipseudo semivector space (S-bipseudo semivector space).*

*A Smarandache sub-bipseudo semivector space (S-sub-bipseudo semivector space) can be defined in a similar way.*

*Example 7.4.10:* Let R[x] be the semigroup under '+'. Consider the S-bisemifield $Z^o[x] \cup Q^o$. Take $P = R^o[x]$, P is a S-strict semigroup and P is a bipseudo semivector space over $Z^o[x] \cup Q^o$. Hence R[x] is a S-bipseudo semivector space.

*Example 7.4.11:* Let $Q_7[x]$ be the set of all polynomials of degree less than or equal to 7. $S = Z^o_7[x] \cup Q^o$ is a bisemifield. Clearly $Q_7[x]$ is a S-bipseudo semivector space over S. The basis for $Q_7[x]$ is $\{1, x, x^2, x^3, x^4, x^5, x^6, x^7\}$ over S.

**PROBLEMS:**

1. Give an example of a bisemivector space with a unique basis.
2. Illustrate by an example a bisemivector space having more than one basis.
3. Is the semigroup Z[x] under '+' a S-strict semigroup?
4. Give an example of a S-strict semigroup with finite number of elements in it.

5. Give an example of a bipseudo semivector space, which is never a S-bipseudo semivector space.
6. Is $V = Z \cup Q[x]$ over the bisemifeild $S = R^o \cup Q^o[x]$ be a S-bipseudo semivector space?
7. Is $V = R^o[x]$ be a S-bipseudo semivector space over $S = Q^o[x] \cup R^o$? Justify your claim.

# Chapter 8

# BINEAR-RINGS AND SMARANDACHE BINEAR-RINGS

The study of bistructures is very recent that too the study of Smarandache bistructures is being introduced here for the first time. We in this chapter introduce the concept of binear-rings and Smarandache binear-rings. This chapter is organized into three sections; in section one we introduce the notion of binear-rings. Smarandache binear-rings are introduced in section two. Generalizations like biquasi near-rings, biquasi rings, biquasi semirings and their Smarandache analogue are introduced in section three. Several results and problems are proposed in this chapter.

## 8.1 Binear-rings and its properties

The study of S-near-rings is very recent (2002) [117, 118]. So the notion of binear-rings does not find its place in any literature. Here we define binear-rings, biseminear-rings and introduce several new and innovative concepts about them.

**DEFINITION 8.1.1:** *Let $(N, +, \bullet)$ be a non-empty set. We call N a binear-ring if $N = N_1 \cup N_2$ where $N_1$ and $N_2$ are proper subsets of N i.e. $N_1 \not\subset N_2$ or $N_2 \not\subset N_1$ satisfying the following conditions:*

*At least one of $(N_i, +, \bullet)$ is a right near-ring i.e. for preciseness we say*

    i.    *$(N_1, +, \bullet)$ is a near-ring*
    ii.   *$(N_2, +, \bullet)$ is a ring.*

*We say that even if both $(N_i, +, \bullet)$ are right near-rings still we call $(N, +, \bullet)$ to be a binear-ring. By default of notation by binear-ring we mean only right binear-ring unless explicitly stated.*

**Example 8.1.1:** *Let $(N, +, \bullet)$ is a non-empty set. $N = Z \cup Z_{12}$ where Z is a near-ring and $Z_{12}$ is a ring. Thus $(N, +, \bullet)$ is a binear-ring.*

**Example 8.1.2:** *Let $(N, +, \bullet)$ be a non-empty set. $N = Z_{12} \cup Q$ where $(Z_{12}, +, \bullet)$ is a near-ring and Q a ring. $(N, +, \bullet)$ is a binear-ring.*

**DEFINITION 8.1.2:** *Let $(N, +, \bullet)$ be a binear-ring. $|N|$ or $o(N)$ denotes the order of the binear-ring, that is number of elements in N.*

**DEFINITION 8.1.3:** *Let $(N, +, \bullet)$ be a binear-ring. We call $(N, +, \bullet)$ as abelian if $(N, +)$ is abelian i.e. if $N = N_1 \cup N_2$ then $(N_1, +)$ and $(N_2, +)$ are both abelian. If both $(N_1, \bullet)$ and $(N_2, \bullet)$ are commutative then we call N a commutative binear-ring. If $N = N_d$, i.e. $N_1 = (N_1)_d$ and $N_2 = (N_2)_d$ then we say N is a distributive binear-ring. If all non-zero elements of N are left (right) cancellable we say that N fulfils the left (right)* cancellation law. N is a bi-integral domain if both $N_1$ and $N_1$ has no zero-divisors. If $N\setminus\{0\} = N_1\setminus\{0\}$ and $N_2\setminus\{0\}$ are both groups then we call N a binear-field.

**DEFINITION 8.1.4:** *Let $(P, +)$ be a bigroup (i.e. $P = P_1 \cup P_2$) with 0, and let N be a binear-ring. $\mu : N \times P \to P$ is called the N-bigroup if for all $p_i \in P_i$ and for all $n, n_i \in N_i$ we have $(n + n_i)p = np + n_ip$ and $(nn_i)p = n(n_ip)$ for $i = 1, 2$. $N^P = N_1^P \cup N_2^P$ stands for N-bigroups.*

**DEFINITION 8.1.5:** *A sub-bigroup M of a binear-ring N with $M.M \subset M$ is called a sub-binear-ring of N. A bisubgroup S of $N^P$ with $NS \subset S$ is a N-sub-bigroup of P.*

**DEFINITION 8.1.6:** *Let $(N, +, \bullet)$ and $(N', +, \bullet)$ be two binear-rings. ($N = N_1 \cup N_2$ and $N' = N_1' \cup N_2'$ ). P and P' be two N –bigroups.*

    i.    *$h : N \to N'$ is called a binear-ring homomorphism if for all $m_1, n_1 \in N_1$ and $m_2, n_2 \in N_2$*

$$h(m_1 + n_1) = h(m_1) = h(n_1)$$
$$h(m_2 + n_2) = h(m_2) + h(n_2)$$

        *and*

$$h(m_1n_1) = h(m_1)h(n_1)$$
$$h(m_2n_2) = h(m_2)h(n_2).$$

    ii.    *$h : P \to P'$ ($P = P_1 \cup P_2$ and $P' = P_1' \cup P_2'$ is called N-bigroup homomorphism if for all $p_1, q_1 \in P_1$; $p_2, q_2 \in P_2$ and for all $n_1 \in N_1$ and $n_2 \in N_2$ we have*

$$h(p_1 + q_1) = h(p_1) + h(q_1)$$
$$h(n_1p_1) = n_1h(p_1)$$

        *and*

$$h(p_2 + q_2) = h(p_2) + h(q_2)$$
$$h(n_2p_2) = n_2h(p_2).$$

**DEFINITION 8.1.7:** *Let $N = N_1 \cup N_2$ be a binear-ring and P a N-bigroup. A binormal subgroup (or normal bisubgroup) I of ($I = I_1 \cup I_2$) $(N, +)$ is called a bi-ideal of N if*

    i.    $I_1N_1 \subset I_1$ and $I_2N_2 \subset I_2$
    ii.    $n, n_1 \in N_1, i_1 \in I_1, n(n_1 + i_1) - nn_1 \in I_1$ and $n, n_2 \in N_2; i_2 \in I_2; n(n_2 + i_2) - nn_2 \in I_2$.

*Normal sub-bigroup T of $(N, +)$ with (i) is called right bi-ideal of N while normal sub-bigroup L of $(N, +)$ with (ii) are called left-bi-ideals.*

*A normal sub-bigroup S of P is called bi-ideal of $N^P$ if for all $s_i \in P_i$ ($i = 1, 2$) ($P = P_1 \cup P_2$) and $s \in S_i$ ($i = 1, 2$; and $S = S_1 \cup S_2$) for all $n_i \in N_i$ ($i = 1, 2; N = N_1 \cup N_2$). $n_i(s + s_i) - ns \in S_i; i = 1, 2$. Factor binear-ring N/I and factor N bigroup P/S are defined as in case of birings.*

**DEFINITION 8.1.8:** *A sub-binear-ring M of the binear-ring N is called bi-invariant if $MN_1 \subset M_1$ and $M_2N_2 \subset M_2$ (where $M = M_1 \cup M_2$ and $N = N_1 \cup N_2$) and $N_1M_1 \subset M_1$ and $N_2M_2 \subset M_2$.*

A minimal bi-ideal, minimal right bi-ideal and minimal left bi-ideal and dually the concept of maximal right bi-ideal, left bi-ideal and maximal bi-ideal are defined as in case of bi-rings.

**DEFINITION 8.1.9:** *Let N be a binear-ring ($N = N_1 \cup N_2$, +, •) and S a sub-bisemigroup of ($N = N_1 \cup N_2$, +) (where $S = S_1 \cup S_2$). A binear-ring $N_S$ is called a binear-ring of left (right) quotients with respect to S if*

i. $N_S = (N_1)_{S_1} \cup (N_2)_{S_2}$ *has identity*
ii. *N is embeddable in $N_S$ by a binear-ring homomorphism h*
iii. *For all $s_i \in S_i$ ($i = 1, 2; S = S_1 \cup S_2$). $h(s_i)$ is invertible in $((N_i)_{S_i}, •)$, $i = 1, 2$.*
iv. *For all $q_i \in N_{S_i}$ there exists $s_i \in S_i$ and there exists $n_i \in N_i$ such that $q_i = h(n_i)h(s_i)^{-1}$; ($q_i = h(s_i)^{-1}h(n_i)$), $i = 1, 2$.*

**DEFINITION 8.1.10:** *The binear-ring N ($N = N_1 \cup N_2$) is said to fulfill the left (right) ore condition with respect to a given sub-bisemigroup $S_i$ of ($N_i$, •) if for all $(s, n) \in S_i \times N_i$ there exist $n • s_1 = s • n_1$ ($s_1 • n = n_1 • s$); $i = 1, 2$.*

Let V denote the collection of all binear-rings and X be any non-empty subset.

**DEFINITION 8.1.11:** *A binear-ring $F_X \in V$ is called a free binear-ring in V over a binear-ring N if there exists $f: X \to F_X$ (where X is any non-empty set) for all $N \in V$ and for all $g: X \to N$ there exists a homomorphism $h \in Hom (F_X, N)$; $h \circ f = g$*

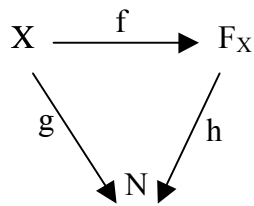

*Let V = {set of all binear-rings}, we simply speak about free binear-ring on X. A binear-ring is called free if it is free over some set X.*

**DEFINITION 8.1.12:** *A finite sequence $N = N_0 \supset N_1 \supset N_2 \supset ... \supset N_t = \{0\}$ of sub-binear-rings $N_i$ of the binear-ring N is called the binormal sequence of N if and only if for all $i \in \{1, 2, ..., t\}$, $N_i$ is a bi-ideal of $N_{i-1}$.*

*In the special case when all the $N_i$'s are bi-ideals of the binear-ring N, we call the binormal sequence an bi-invariant sequence.*

**DEFINITION 8.1.13:** *Let P be a bi-ideal of a binear-ring N. P is called the prime bi-ideal if for all bi-ideals I and J of N; IJ ⊂ N implies I ⊂ P or J ⊂ P. The binear-ring N is called a prime binear-ring if {0} is a prime ideal.*

**DEFINITION 8.1.14:** *Let S be a bi-ideal of a binear-ring N. S is semiprime if and only if for all bi-ideals I on N, $I^2 \subset S$ implies I ⊂ S. N is called a semiprime binear-ring if {0} is a semiprime bi-ideal.*

**DEFINITION 8.1.15:** *Let L be a left bi-ideal in a binear-ring N. L is called modular if and only if there exists e ∈ $N_1 \cup N_2$ and for all n ∈ $N_1 \cup N_2$ and for all n ∈ $N_1 \cup N_2$, n – ne ∈ L. In this case we also say that L is modular by e and that e is a right identity modulo L, since for all n ∈ N, ne = n (mod L).*

*Notation*: For z ∈ N, denote the bi-ideal generated by the set {n – nz / n ∈ N} by $L_Z$. $L_Z$ = N if z = 0. z ∈ N is called quasi regular if z ∈ $L_Z$, S ⊂ N is called quasi-regular if and only if for all s ∈ S, s is quasi regular.

Let N be a binear-ring. An idempotent e ∈ N = $N_1 \cup N_2$ is called central if it is in the center of ($N_1$, •) or ($N_2$, •) i.e. for all n ∈ $N_1$ (or n ∈ $N_2$) we have ne = ne.

**DEFINITION 8.1.16:** *A binear-ring is said to be biregular if there exists some set E of central idempotents with*

  i.  *For all e ∈ $N_1$ (e ∈ $N_2$), $N_1$e ($N_2$e) is an ideal of $N_1$ (or $N_2$). i.e. $N_1$e ∪ {0} or {0} ∪ $N_2$e or if e ∈ $N_1 \cap N_2$ then $N_1$e ∪ $N_2$e is a bi-ideal of N.*

**DEFINITION 8.1.17:** *Let I be a bi-ideal of the biring N. The intersection of all prime bi-ideals P, which contain I is called the prime bi-radical of I.*

**DEFINITION 8.1.18:** *Let (N, +, •) be a binear-ring and x an indetermeniate. N[x] is called the binear-polynomial ring if N[x] = $N_1$[x] ∪ $N_2$[x] where $N_1$[x] and $N_2$[x] are near-polynomial rings.*

Several differences exists between binear-polynomial rings and near-polynomial rings.

**DEFINITION 8.1.19:** *Let (N, +, •) be a binear-ring. We define the binear-matrix ring as follows: N = $N_1 \cup N_2$ where $N_1$ and $N_2$ are near-rings. The n × n binear-ring matrix ring with entries from the binear-ring N denoted by N(n × n) = {($a_{ij}^1$)/ $a_{ij}^1$ ∈ $N_1$} ∪ {($a_{ij}^2$)/ $a_{ij}^2$ ∈ $N_2$} where ($a_{ij}$) denotes the n × n matrix with entries from $N_1$ and $N_2$.*

**DEFINITION 8.1.20:** *Let N be a binear-ring; N = $N_1 \cup N_2$ where $N_1$ and $N_2$ are near-rings. The binear-ring N is said to fulfill the insertion factors property (IFP) provided for all a, b, n ∈ $N_1$ (or a, b, n ∈ $N_2$) we have ab = 0 implies anb = 0. The biring has strong IFP property if every homomorphic image of N has the IFP. The binear-ring N has the strong IFP if and only if for all I in N and for all a, b, n ∈ $N_1$ (a, b, n ∈ $N_2$), a, b ∈ $I_1$ implies anb ∈ $I_1$ (ab ∈ $I_2$ implies anb ∈ $I_2$) where I = $I_1 \cup I_2$.*

**DEFINITION 8.1.21:** *Let p be a prime. A binear-ring N is called a p-binear-ring provided that for all $x \in N$, $x^p = x$ and $px = 0$.*

**DEFINITION 8.1.22:** *A biright ideal I of a binear-ring N is called right quasi bireflexive if whenever A and B are bi-ideals of the binear-ring N with $AB \subset I$ then $b(b' + a) - bb' \in I$ for all $a \in A$ and $b, b' \in B$. A binear-ring N is strongly sub-bicommutative if it is right quasi bireflexive.*

**DEFINITION 8.1.23:** *Let N be a binear-ring S a subnormal sub-bigroup of (N, +). S is called a quasi bi-ideal of N if $SN \subset NS \subset S$ whereby NS we mean elements of the form $\{n(n' + s) - nn' / \text{for all } s \in S \text{ and for } n, n' \in N\} = NS$.*

**DEFINITION 8.1.24:** *A left binear-ring N is said to satisfy left self distributive identity if $abc = abac$ for all $a, b, c \in N_1$ or for all $a, b, c \in N_2$. A binear-ring N is said to be left permutable if $abc = bac$, right permutable if $abc = acb$, medial if $abcd = acbd$, right self distributive $abc = acbc$ for all $a, b, c, d \in N = N_1 \cup N_2$.*

**DEFINITION 8.1.25:** *Let $(S, +, \bullet)$ be a binear-ring; $S = S_1 \cup S_2$ and M a right S-module. Let $R = S \times M = (S_1 \times M) \cup (S_2 \times M)$, define $(\alpha, s) \odot (\beta, t) = (\alpha\beta, s\beta + t)$. Then $(R, +, \odot)$ is a left binear-ring, the abstract affine binear-ring inducted by S and M. If S is a right self-distributive binear-ring and $MS^2 = 0$, then R is a right self-distributive binear-ring.*

If $MS = 0$ and S is a medial, left permutable or left-self distributive then $(R, \odot)$ is a medial, left permutable or left-self distributive respectively.

**DEFINITION 8.1.26:** *A binear-ring R is defined to be equiprime if for all $0 \neq a \in R = R_1 \cup R_2$ and $x_1, y_1 \in R_1$ ($x_2, y_2 \in R_2$) if $a \in R_1$, $ar_1x_1 = ar_1y_1$ for all $r_1 \in R_1$ ($ar_2x_2 = ar_2y_2$ if $a \in R_2$ and $r_2 \in R_2$) implies $x_1 = y_1$ ($x_2 = y_2$). If P is a bi-ideal of R, P is called an equiprime bi-ideal of R, if R/P is an equiprime binear-ring. It can be proved if P is a bi-ideal of R if and only if for all $a \in R \backslash P$, $x, y \in P$, $arx - ary \in P$ for all $r \in R$ implies $x - y \in P$. Now we define the new notion of infra binear-ring.*

**DEFINITION 8.1.27:** *An infra binear-ring (INR) is a triple $(N, +, \bullet)(N = N_1 \cup N_2)$ where*

    i.   *(N, +) is a bigroup*
    ii.  *(N, •) is a bisemigroup.*
    iii. *$(x + y) \bullet z = x \bullet z - 0 \bullet z + y \bullet z$ for all $x, y, z \in N$.*

**DEFINITION 8.1.28:** *Let I be a bi-ideal of a binear-ring N. Suppose that I satisfies the following conditions*

  i.   *$a, x, y \in N$, $anx - any \in I$ for all $n \in N$ implies $x - y \in I$. ($I = I_1 \cup I_2$ and $x - y \in I_1$ or $x - y \in I_2$ according as $anx - any \in I_1$ or $I_2$)*
  ii.  *I is left bi-invariant*
  iii. *$0N \subset I$.*

*Then I is called an equiprime left bi-ideal of N. We call an idempotent e of a binear-ring to be central if it is in the center of (N, •) i.e. for all $n_1 \in N_1$ ($n_2 \in N_2$) we have $n_1e = ex$ or $n_2e = ex$ according as $e \in N_1$ or $N_2$ respectively.*

**DEFINITION 8.1.29:** *Let (N, +, •) be a binear-ring. N is said to be partially ordered by ≤ if*

i. *≤ makes the bigroup (N, +) i.e. ($N_1$, +) and ($N_2$, +) into a partially ordered group.*
ii. *for all $n_1, n_1' \in N_1$ ($n_2, n_2' \in N_2$), $n_1 \geq 0$ and $n_1' \geq 0$ implies $n_1 n_1' \geq 0$ ($n_2 n_2' \geq 0$)*

Several interesting results about these concepts can be defined and derived, this task is left for the reader.

The concept of group binear-ring, bigroup near-ring, bigroup binear-ring, semigroup binear-ring, bisemigroup near-ring and bisemigroup binear-ring will be defined some of these definitions lead to the concept of quad near-ring or bi-binear-ring.

**DEFINITION 8.1.30:** *Let G be any group. (R, +, •) be a binear-ring. The group binear-ring, RG is defined as $RG = R_1G \cup R_2G$ where $R_iG$ for i = 1, 2 are group near-rings. Clearly the group binear-ring is a binear-ring.*

*Example 8.1.3:* Let $G = S_3$ be the symmetric group of degree 3 and $R = Z_{10} \cup Z$ be the binear-ring. Group binear-ring $RG = RS_3 = Z_{10}S_3 \cup ZS_3$ is a binear-ring where $Z_{10}S_3$ is a group near-ring and $ZS_3$ is a group ring.

**DEFINITION 8.1.31:** *Let S be a semigroup. $N = N_1 \cup N_2$ be a binear-ring. The semigroup binear-ring NS is defined to be $NS = N_1S \cup N_2S$ where $N_1S$ and $N_2S$ are semigroup near-rings. Thus NS is a binear-ring.*

*Example 8.1.4:* Let $S = S(7)$ be a semigroup. $N = N_1 \cup N_2$ be a binear-ring say $N = Z_2 \cup Z_{14}$. $Z_2$ is a near-field and $Z_{14}$ is the ring; $NS = Z_2S \cup Z_{14}S$ is a binear-ring.

**DEFINITION 8.1.32:** *Let N be a near-ring and $G = G_1 \cup G_2$ be the bigroup. The bigroup near-ring $NG = NG_1 \cup NG_2$ is a binear-ring as both $NG_1$ and $NG_2$ are near-rings.*

*Example 8.1.5:* Let Z be a near-ring and $G = S_3 \cup G_1 = \langle g / g^7 = 1 \rangle$ be a bigroup. The bigroup near-ring $ZG = ZS_3 \cup ZG_1$ is a binear-ring.

**DEFINITION 8.1.33:** *Let N be a near-ring. $S = S_1 \cup S_2$ be a bisemigroup. The bisemigroup near-ring $NS = NS_1 \cup NS_2$ is the binear-ring which is the union of two near-rings viz. semigroup near-rings $NS_1$ and $NS_2$.*

*Example 8.1.6:* Let $Z_2$ be a near-ring. The bisemigroup $S = S(3) \cup Z_{25}$. The bisemigroup near-ring $Z_2S = Z_2S(3) \cup Z_2(Z_{25})$ is a binear-ring which is the union of near-rings.

**DEFINITION 8.1.34:** *Let $(N, +, \bullet)$ be a non-empty set, where $N = N_1 \cup N_2 \cup N_3 \cup N_4$ where each $N_i$ is a near-ring or a ring for each $i = 1, 2, 3, 4$. $(N, +, \bullet)$ is defined as the quad near-ring or bi-binear-ring.*

*Example 8.1.7:* Let $(N, +, \bullet)$ be a non-empty set, where $N = Z_2 \cup Z_{12} \cup Z_7 \cup Z$ where $Z_2$ and $Z$ are near-rings and $Z_{12}$ is the ring of integers modulo 12 and $Z_7$ is the prime field of characteristic 7. N is a quad near-ring.

**DEFINITION 8.1.35:** *Let $(N, +, \bullet)$ be a binear-ring where $N = N_1 \cup N_2$ and $G = G_1 \cup G_2$ be a bigroup. The bigroup binear-ring $NG = N_1G_1 \cup N_1G_2 \cup N_2G_1 \cup N_2G_2$ is a quad near-ring where each $N_iG_j$ is a group near-ring or group ring, $i = 1, 2$ and $j = 1, 2$. Similarly we can define bisemigroup binear-ring $NS = N_1S_1 \cup N_2S_1 \cup N_1S_2 \cup N_2S_2$ where S is a bisemigroup ($S = S_1 \cup S_2$) and $N = N_1 \cup N_2$ is a binear-ring. It is easily verified that the bisemigroup binear-ring is also a quad near-ring where $N_iS_j$ are semigroup near-rings or semigroup rings, $i = 1, 2$ and $j = 1, 2$. All results connected with binear-rings can be extended in case of quad near-rings.*

**DEFINITION 8.1.36:** *Let $Z_p$ be a near-ring on modulo integers p. Take $N = Z_p \cup Z_q$ the binear-ring where both are near-rings and $p \neq q$. G be any group, NG be the group binear-ring. Define the bimod (p, q) envelope of NG as $1 + U = (1 + U_1) \cup (U_2 + 1)$ where $U_1 = \{\Sigma\alpha_i g_i / \Sigma\alpha_i = 0, \alpha_i \in Z_p\}$, and $U_2 = \{\Sigma\beta_j g_j / \Sigma\beta_j = 0, \beta_j \in Z_q\}$, $G^* = 1 + U = (1 + U_1) \cup (1 + U_2)$ is called the bimod (p, q) envelope of the group binear-ring NG. Similarly we define for any semigroup S the bimod (p, q) envelope of NS denoted by $S^*$. The bimod (p, q) envelope may be a bisemigroup under multiplication or a bigroup under multiplication.*

Now we proceed onto define mod p envelope of a bigroup near-ring over the near-ring $Z_p$.

**DEFINITION 8.1.37:** *Let $Z_p$ be a near-ring of the integers modulo p and $G = G_1 \cup G_2$ be a bigroup of finite order. $Z_pG$ be the bigroup near-ring. The mod p envelope of $Z_pG$ is $1 + U = (1 + U_1) \cup (1 + U_2)$ where $U_1 = \{\Sigma\alpha_i g_i / \Sigma\alpha_i = 0, \alpha_i \in Z_p \text{ with } g_i \in G_1\}$ and $U_2 = \{\Sigma\beta_j h_j / \Sigma\beta_j = 0, \beta_j \in Z_p \text{ with } h_j \in G_2\}$.*

*Then $1 + U = G^*$ is the mod p envelop of $Z_pG$. Clearly $G^*$ is a bigroup or a bisemigroup under '$\bullet$'. Such study can be made even in case of bisemigroup near-ring $Z_pS$ where $S = S_1 \cup S_2$ is a bisemigroup. What is the structure of $S^*$? The study in this direction is interesting and remains unexplored.*

Now we proceed on to define non-associative binear-rings. It is like all other structures we don't have a method by which we can have non-assocative binear-rings.

**DEFINITION 8.1.38:** *Let $(N, +, \bullet)$ be a non-empty set we say N is a non-associative binear-ring if*

  i.  $N = N_1 \cup N_2$ *where $N_1$ and $N_2$ are proper subsets of N i.e. $N_1 \not\subset N_2$ or $N_2 \not\subset N_1$.*
  ii. *$(N_1, +, \bullet)$ is a non-associative near-ring.*
  iii. *$(N_2, +, \bullet)$ can be a near-ring or a ring.*

Throughout this book we denote a non-associative binear-ring by na-binear-ring.

***Example 8.1.8:*** Let $(N, +, \bullet)$ be a na-binear-ring where $N = N_1 \cup N_2$ with $N_1 = Z_8 L_5(3)$ and $N_2 = Z$. (i.e. $N_1$ is the loop near-ring of the loop $L_5(3)$ over the near-ring $Z_8$) and Z is the ring of integers. Clearly $(N, +, \bullet)$ is a na-binear-ring. Thus na-binear-rings can be constructed using either loop near-rings and groupoid near-rings as these classes of near-rings are always non-associative.

We can define all properties defined and studied in case of associative binear-rings in case of na-binear-ring. Thus this work is a matter of routine so we leave it for the reader as an exercise. Now we proceed on to define notions which are not common with associative binear-rings.

**DEFINITION 8.1.39:** *Let $(N, +, \bullet)$ be a binear-ring which is associative. Let L be a loop. The loop binear-ring $NL = N_1L \cup N_2L$ is a na-binear-ring. Here $N_iL$ are nothing but loop near-rings for i = 1, 2.*

**DEFINITION 8.1.40:** *Let $(N, +, \bullet)$ be an associative binear-ring i.e. $N = N_1 \cup N_2$. G be any groupoid; the groupoid binear-ring. $NG = N_1G \cup N_2G$ is a non-associative binear-ring where $N_iG$ are just the groupoid near-rings of the groupoid G over the near-ring $N_i$, i = 1, 2.*

**DEFINITION 8.1.41:** *Let N be a near-ring and $L = L_1 \cup L_2$ be a biloop. The biloop near-ring $NL = NL_1 \cup NL_2$ is a binear-ring which is non-assocative, where $NL_i$ are just loop near-rings.*

**DEFINITION 8.1.42:** *Let N be a near-ring and $G = G_1 \cup G_2$ be a bigroupoid. The bigroupoid near-ring $NG = NG_1 \cup NG_2$ is a non-associative binear-ring.*

Thus we see using loops, biloops, groupoids and bigroupoids we can get four distinct classes of non-associative binear-rings. The only special property that is not enjoyed by associative binear-rings are the notions of special identities like Moufang, Bol, Bruck, Alternative (right/ left) P-identity and so on. These are properties solely attached with na-binear-rings.

We define Moufang binear-ring.

**DEFINITION 8.1.43:** *Let $(N, +, \bullet)$ be a na-binear-ring. We call N a Moufang binear-ring if $(xy)(zx) = (x(yz))x$ for all x, y, z in N. That is we say a binear-ring $(N, +, \bullet)$ is called a Moufang binear-ring if each of the near-rings satisfy the Moufang identity $(xy)(zx) = (x(yz))x$ for all x, y, z $\in N_1$ and x, y, z in $N_2$ where $N = N_1 \cup N_2$.*

*On similar lines we define Bol binear-ring, Bruck binear-ring, alternative (left or right) binear-ring and so on. The notions of bi-ideals, sub-binear-rings etc. are defined as in case of associative binear-rings.*

Yet another type of binear-rings which are based on the non-associative near-rings introduced by [46] is studied and analysed in the following.

**DEFINITION 8.1.44:** *A right loop binear-ring is a system $(N, +, \bullet)$ where $N = N_1 \cup N_2$ satisfying the following conditions:*

 i. *$(N, +)$ is a biloop*
 ii. *$(N, \bullet)$ is a bisemigroup*
 iii. *The multiplication is right-distributive over addition '+' on both $N_1$ and $N_2$ i.e. for all $n, m, s \in N_i$ we have $(n + m)s = ns + ms$, for $i = 1, 2$.*

*Thus we can roughly define a loop binear-ring $(N, +, \bullet)$ as for $N = N_1 \cup N_2$ we have $(N_1, +, \bullet)$ is a right loop near-ring and $(N_2, +, \bullet)$ may be a right loop near-ring or can be non-right loop near-ring.*

**DEFINITION 8.1.45:** *Let $(N, +, \bullet)$ be a right loop binear-ring i.e. $N = N_1 \cup N_2$. An additive sub-biloop A of N ($A = A_1 \cup A_2$, $A_1$ a subloop of $N_1$ and $A_2$ a subloop of $N_2$) is a N-sub-biloop if $NA \subset A$ ($AN \subset A$) where $Na = \{n_1a_1 / n_1 \in N_1, a_1 \in A_1\} \cup \{n_2a_2 / n_2 \in N_2 \text{ and } a_2 \in A_2\}$.*

*A non-empty subset I of N where $I = I_1 \cup I_2$ is called the left-bi-ideal in N if*

 i. *$(I_1, +)$ and $(I_2, +)$ are normal subloop of $(N_1, +)$ and $(N_2, +)$ respectively.*
 ii. *$n(n_1 + i) + n_r n_1 \in I_i$, $i = 1, 2$ and $n, n_1 \in N_i$, $i = 1, 2$ where $n_r$ denotes the unique right inverse of n.*

*A non-empty subset I of N is called a bi-ideal in N if*

 i. *I is a left-bi-ideal*
 ii. *$IN \subset I$.*

**DEFINITION 8.1.46:** *A loop binear-ring $N$ ($N = N_1 \cup N_2$) is said to be left bipotent if $N_i a = N_i a^2$ for every $a \in N_i$, $i = 1, 2$. A loop binear-ring $(N, +, \bullet)$ where $N = N_1 \cup N_2$ is said to be a s-loop near-ring if $a \in N_i a$ for every a in $N_i$, $i = 1, 2$.*

*Let N be a s-loop binear-ring, then N is regular if and only if for each $a$ ($\neq 0$) in $N_i$ there is an idempotent $e_i \in N_i$ such that $N_i a = N_i e$; $i = 1, 2$.*

Now we proceed on to define quad near-rings which are non-associative.

**DEFINITION 8.1.47:** *Let $(N, +, \bullet)$ be an associative binear-ring and $L = L_1 \cup L_2$ be a biloop. The biloop binear-ring $NL = N_1L_1 \cup N_1L_2 \cup N_2L_1 \cup N_2L_2$ is a bi-binear-ring.*

**DEFINITION 8.1.48:** *Let $(N, +, \bullet)$ be a non-empty set. $N = N_1 \cup N_2 \cup N_3 \cup N_4$ where each $N_i$ is a near-ring and at least one of the $N_i$'s is a non-associative near-ring (We assume all the proper subsets $N_1, N_2, N_3, N_4$ are distinct and $N_i \not\subset N_j$ if $i \neq j$, $i, j \in \{1, 2, 3, 4\}$). We call N a non-associative quad near-ring. Thus the biloop binear-ring is a non-associative quad near-ring or bi-binear-ring.*

**DEFINITION 8.1.49:** *Let G be a bigroupoid and $(N = N_1 \cup N_2, +, \bullet)$ be an associative biring. The bigroupoid binear-ring $NG = N_1G_1 \cup N_1G_2 \cup N_2G_1 \cup N_2G_2$ where $G =$ $G_1 \cup G_2$ is the bigroupoid and each of the $N_1G_1$, $N_2G_1$, $N_2G_2$ and $N_1G_2$ are groupoid near-rings. Thus a bigroupoid binear-ring which is non-associative.

Thus we have got a new classes of quad near-rings viz. biloop binear-rings and bigroupoid binear-rings. All properties studied in case of non-associative binear-rings can be defined and studied in case of quad near-rings as they are nothing but non-associative bi-binear-rings.

We can define binear-fields and quad near-fields as in case of near-rings. Such study is also a matter of routine by replacing near-rings by near-fields. The only notion which we have not introduced is the concept of biseminear-rings but it can be associative as well as non-associative.

**DEFINITION 8.1.50:** *Let $(N, +, \bullet)$ is a non-empty set where $N = N_1 \cup N_2$. We say N is a biseminear-ring if one of $(N_1, +, \bullet)$ or $(N_2, +, \bullet)$ are seminear-rings. Thus for N to be biseminear-ring we need one of $(N_1, +, \bullet)$ or $(N_2, +, \bullet)$ or both to be a seminear-ring.*

We can have several such biseminear-rings. All notions studied in case of binear-rings can also be easily studied and generalized for biseminear-rings. The only means to generate several new classes of biseminear-rings are as follows:

**DEFINITION 8.1.51:** *Let $(S, +, \bullet)$ be a biseminear-ring with $S = S_1 \cup S_2$. G any group the group biseminear-ring $SG = S_1G \cup S_2G$ where $S_iG$ are group seminear-rings. Clearly every group biseminear-ring is a biseminear-ring.*

By using this definition we can have infinite class of both commutative or non-commutative and infinite or finite classes of biseminear-rings. If in the definition of group biseminear-rings if we replace the group by a semigroup we get yet another new class of biseminear-rings.

*Example 8.1.9:* Let $(S, +, \bullet)$ be a biseminear-ring with $S = S_1 \cup S_2$. For $S(9)$ the semigroup we obtain the semigroup biseminear-ring $S(S(9))$ which is a non-commutative biseminear-ring.

**DEFINITION 8.1.52:** *Let G be a bigroup and S a seminear-ring. The bigroup seminear-ring $SG = S_1G \cup S_2G$ is a union of two group seminear-rings; hence a bigroup seminear-ring is a biseminear-ring.*

On similar lines, if we take a seminear-ring S and T a bisemigroup i.e. $T = T_1 \cup T_2$ where $T_1$ and $T_2$ are semigroups then the bisemigroup seminear-ring $ST = ST_1 \cup ST_2$ is a biseminear-ring.

Thus we have so far given four methods of getting biseminear-rings. It is left for the reader to study biseminear-rings in a similar way as that of birings or binear-rings. Another new concept analogous to quad near-rings are quad seminear-rings. Here we define quad seminear-rings which are nothing but bi-bi-seminear-rings.

**DEFINITION 8.1.53:** *Let $(N, +, \bullet)$ be a non-empty set. We say $(N, +, \bullet)$ is a quad seminear-ring or bi-biseminear-ring if the following conditions hold good.*

i. $N = N_1 \cup N_2 \cup N_3 \cup N_4$ where $N_i$ are proper subsets of $N$ such that $N_i \not\subset N_j$ for any $i \neq j$, $i, j = 1, 2, 3, 4$.
ii. $(N_i, +, \bullet)$ is at least a seminear-ring for $i = 1, 2, 3, 4$ under the operations of $N$. Then we call $(N, +, \bullet)$ a quad seminear-ring or bi-biseminear-ring.

Now we shall show that this new class of quad seminear-rings is non-empty.

**DEFINITION 8.1.54:** *Let $(N, +, \bullet)$ be a biseminear-ring i.e. $N = N_1 \cup N_2$ and $G = G_1 \cup G_2$ be a bigroup. The bigroup biseminear-ring $NG = N_1G_1 \cup N_2G_2 \cup N_1G_2 \cup N_2G_1$ where each of $N_iG_j$ are group seminear-rings. $1 \leq i, j \leq 2$.*

**THEOREM 8.1.1:** *Every bigroup biseminear-ring is a quad seminear-ring or a bi-biseminear-ring.*

*Proof:* Left as an exercise for the reader to prove.

**DEFINITION 8.1.55:** *Let $(N, +, \bullet)$ be a biseminear-ring i.e. $N = N_1 \cup N_2$ and $S = S_1 \cup S_2$ be a bisemigroup. Let $NS$ be the bisemigroup biseminear-ring of the bisemigroup $S$ over the biseminear-ring $N$ i.e. $NS = N_1S_1 \cup N_1S_2 \cup N_2S_1 \cup N_2S_2$.*

**THEOREM 8.1.2:** *Every bisemigroup biseminear-ring is a bi-biseminear-ring or a quad seminear-ring.*

*Proof:* Straightforward. Hence left for the reader to prove.

***Example 8.1.10:*** Let $S = S(3) \cup Z_{10}$ be the bisemigroup. $N = Z^o \cup Z_7$ the biseminear-ring. The bisemigroup biseminear-ring $NS$ is a quad seminear-ring.

As it is a well-known fact that we do not have a direct method of constructing non-associative seminear-rings it is more difficult to get non-associative biseminear-rings. Now we give the definition of non-associative biseminear-rings.

**DEFINITION 8.1.56:** *Let $(N, +, \bullet)$ be a non-empty set. $N$ is said to be a non-associative biseminear-ring if $N = N_1 \cup N_2$ is a union of proper subsets of $N$ such that at least one of $(N_1, +, \bullet)$ or $(N_2, +, \bullet)$ are non-associative seminear-rings.*

We can build na-biseminear-rings in a very easy way using loops or groupoids and associative seminear-rings. We just define several new classes of na-biseminear-rings and leave it for the reader to develop all notions as in case of na-binear-rings.

**DEFINITION 8.1.57:** *Let $G = G_1 \cup G_2$ be a bigroupoid. $(N, +, \bullet)$ be a seminear-ring. The bigroupoid seminear-ring $NG = NG_1 \cup NG_2$ is a biseminear-ring where $NG_1$ and $NG_2$ are groupoid seminear-ring. Clearly both or at least one of $NG_1$ or $NG_2$ is a non-asssociative seminear-ring so $NG$ the groupoid seminear-ring is a biseminear-ring.*

**DEFINITION 8.1.58:** *Let $L = L_1 \cup L_2$ be a biloop and $N$ be any associative seminear-ring. $NL$ the biloop seminear-ring is a na-biseminear-ring. For $NL = NL_1 \cup NL_2$ and at least one of $NL_1$ or $NL_2$ is a non-associative seminear-ring.*

Now we define another class of biseminear-ring.

**DEFINITION 8.1.59:** *Let $(N, +, \bullet)$ be an associative biseminear-ring. L be any loop. The loop biseminear-ring $NL = N_1L \cup N_2L$ where $N = N_1 \cup N_2$ and $N_1L$ and $N_2L$ are loop seminear-rings. Hence NL is a na-biseminear-ring.*

**DEFINITION 8.1.60:** *Let G be any groupoid. $(N, +, \bullet)$ be an associative biseminear-ring. The groupoid biseminear-ring $NG = N_1G \cup N_2G$ where both $N_1G$ and $N_2G$ are na-seminear-rings so NG is a na-biseminear-ring.*

All properties studied for na-binear-rings can be studied in case of na-biseminear-rings. Now we proceed on to define na-quad seminear-rings or bi-biseminear-rings.

**DEFINITION 8.1.61:** *Let $(N, +, \bullet)$ be a non-empty set such that the following conditions are true:*

  i. $N = N_1 \cup N_2 \cup N_3 \cup N_4$ *where each of the $N_i$'s are proper subsets of N and $N_i \not\subset N_j$ if $i \neq j$, $1 \leq i, j \leq 4$*
  ii. *At least one of the $(N_i, +, \bullet)$ is a na-seminear-ring and others are just seminear-rings under the operations of N.*

*Then we call $(N, +, \bullet)$ a non-associative quad seminear-ring.*

We have got new class of quad seminear-rings built using bigroupoids and biloops.

**DEFINITION 8.1.62:** *Let $(N, +, \bullet)$ be a biseminear-ring (i.e. $N = N_1 \cup N_2$) and $G = G_1 \cup G_2$ be a bigroupoid. The bigroupoid biseminear-ring $NG = N_1G_1 \cup N_1G_2 \cup N_2G_1 \cup N_2G_2$ is the union of groupoid seminear-rings of the groupoids $G_1$ or $G_2$ over $N_1$ or $N_2$. Thus NG is a quad seminear-ring.*

**DEFINITION 8.1.63:** *Let $L = L_1 \cup L_2$ be a biloop and $(N, +, \bullet)$ be a biseminear-ring ($N = N_1 \cup N_2$). The biloop biseminear-ring NL of the biloop L over the biseminear-ring N is given by $NL = N_1L_1 \cup N_1L_2 \cup N_2L_1 \cup N_2L_2$.*

**THEOREM 8.1.3:** *Every biloop biseminear-ring is a na-quad seminear-ring.*

*Proof:* Left for the reader as an exercise. Thus we have introduced na-quad seminear-rings.

***Example 8.1.11:*** Let $L = L_7(3) \cup S_3$ be the biloop. $N = Q^o \cup Z_{11}$ be the biseminear-ring; NL is a na-quad seminear-ring.

**PROBLEMS:**

1. Let $N = Z_{20} \cup Z$ where $Z_{20}$ is a near-ring and Z is a ring. $(N, +, \bullet)$ be a binear-ring.
   i. Find bi-ideals of N.
   ii. Find sub-binear-rings of N which are not bi-ideals of N.

2. Let N be a binear-ring $\{I_k\}_{k \in K}$ be the collection of all bi-ideals of N. Prove

   i. the set of all finite sums of elements of the $I_k$'s is equivalent to
   ii. the set of all finite sums of elements of different $I_k$'s.the bi-ideal of N generated by $\bigcup_{k \in K} I_k$ .

3. Give an example of a free binear-ring.
4. Is the binear-ring (N, +, •) where N = Z ∪ R[x] where Z is the near-ring and R[x] the ring of polynomials over the reals R, a free binear-ring?
5. What is the smallest order of quad near-ring?
6. Can the smallest order of a na-quad near-ring be less than 8? Justify your claim.
7. Is NG where G = $L_{11}(3)$ ∪ $S_7$ and N = Z ∪ $Z_{11}$ be a quad near-ring? Justify your claim.
8. For NG given in problem 7. Find ideals, sub-biseminear-rings.
9. Can NG given in problem 7 have idempotents? zero divisors?
10. What is the order of the smallest quad seminear-ring?
11. Find the order of the smallest na-quad seminear-ring. Can it be less than 6?
12. Is NG where G = $L_{13}(3)$ ∪ $D_{27}$ and N = Z ∪ $Z_{12}$ (the biloop G and the binear-ring N) be a quad near-ring? Justify.
13. For NG given in problem 12 find quad ideals and N-sub-biloops.
14. Does there exist a quad seminear-ring of order 12?
15. Give an example of a quad seminear-ring of order 72.

## 8.2 S-binear-rings and its generalizations

In this section we introduce the notion of Smarandache binear-rings and Smarandache biseminear-ring. Several new concepts are introduced and studied. We also deal with Smarandache na binear-rings and Smarandache na biseminear-rings. In the study of Smarandache na binear-rings and na biseminear-rings the special identities are also introduced and studied.

**DEFINITION 8.2.1:** *Let (N, +, •) be a binear-ring. (N = $N_1$ ∪ $N_2$). We say N is a Smarandache binear-ring (S-binear-ring) if N contains a proper subset P such that P under the operations '+' and '•' is a binear field; i.e., (P = $P_1$ ∪ $P_2$, +, •) is a binear field.*

*Example 8.2.1:* Let (N, +, •) be a binear-ring; where $N_1$ = $Z_2$G and $N_2$ = $Z_3S_3$ with N = $N_1$ ∪ $N_2$. Clearly (N, +, •) is a S-binear-ring.

**DEFINITION 8.2.2:** *Let (N, +, •) be a binear-ring. A proper subset P of N is said to be a Smarandache subbinear-ring (S-subbinear-ring) if P itself is a S-binear-ring.*

**THEOREM 8.2.1:** *If (N, +, •) is a binear-ring and if P is a S-subbinear-ring then we have N to be a S-binear-ring.*

*Proof:* Follows from the very definitions.

**DEFINITION 8.2.3:** *Let P be a S-bisemigroup with 0 and let N be a S-binear-ring. A map $\mu: N \times Y \to Y$ where Y is a proper subset of P which is a bigroup under the operations of P, $(P, \mu)$ is called the Smarandache N-bigroup (S-N-bigroup) if for all $y \in Y$ and for all $n, n_l \in N$ we have $(n + n_l)y = ny + n_l y$ and $(nn_l)y = n(n_l y)$. $S(N^P)$ stands for the S-N-bigroups.*

**DEFINITION 8.2.4:** *A S-bisemigroup M of the binearring N is called a Smarandache quasi sub-binear-ring (S-quasi sub-binear-ring) of N if $X \subset M$ where X is a bisubgroup of M which is such that $XX \subset X$.*

**THEOREM 8.2.2:** *If N has a S-quasi sub-binear-ring then N has a sub-binear-ring.*

*Proof:* Follows directly by definitions.

**DEFINITION 8.2.5:** *A S-sub-bisemigroup Y of $S(N^P)$ with $NY \subset Y$ is said to be a Smarandache N-sub-bigroup (S-N-sub-bigroup) of P.*

**DEFINITION 8.2.6:** *Let N and N' be two S-binear-rings P and P' be S-N-sub-bigroups*

  i. $h : N \to N'$ *is called a Smarandache binear-ring homomorphism (S-binear-ring homomorphism) if for all $m, n \in M$ we have $h(m + n) = h(m) + h(n)$, $h(mn) = h(m)h(n)$ where $h(m), h(n) \in M'$ (M' is a proper subset of N' which is a binearfield). It is to be noted that h need not even be defined on whole of N.*

  ii. $h: P \to P'$ *is called the Smarandache N-sub-bigroup homomorphism (S-N-sub-bigroup homomorphism) if for all $p, q$ in S (S the proper subset of P which is a S-N-sub-bigroup of the S-bisemigroup P) and for all $m \in M \subset N$ (M a nearfield of N); $h(p + q) = h(p) + h(q)$ and $h(mp) = mh(p)$, $h(p), h(q)$ and $mh(p) \in S'$ (S' is a proper subset of P' which is S-N'-sub-bigroup of S-bisemigroup P').*

*Here also we do not demand h to be defined on whole of P.*

**DEFINITION 8.2.7:** *Let $(N, +, \bullet)$ be a S-binear-ring. A normal sub-bigroup I of $(N, +)$ is called a Smarandache bi-ideal (S-bi-ideal) of N related to X, where X is a binearfield contained in N if ($X = X_1 \cup X_2$, $X_1$ and $X_2$ are near fields).*

  i. $I_1 X_1 \subset I_1$; $I_2 X_2 \subset I_2$ where $I = I_1 \cup I_2$.
  ii. $\forall x_i, y_i \in X_i$ and for all $k_i \in I_i$, $x_i(y_i + k_i) - x_i y_i \in I_i$; $i = 1, 2$.

**DEFINITION 8.2.8:** *A proper subset S of P is called a S-bi-ideal of $S(N^P)$ related to Y if*

  i. *S is a S-normal sub-bigroup of the S-bisemigroup.*
  ii. *For all $s_1 \in S$ and $s \in Y$ and for all $m \in M$ (M the binear field of N) $n(s + s_1) - ns \in S$.*

*A S-binear-ring is S-bi-ideal if it has no S-bi-ideals. $S(N^P)$ is called Smarandache N-bisimple (S-N-bisimple) if it has no S-normal sub-bigroups expect 0 and P.*

**DEFINITION 8.2.9:** *A S-sub-binear-ring M of a binear-ring N is called Smarandache bi-invariant (S-bi-invariant) related to the binearfield X in N if $MX \subset M$ and $XM \subset M$ where X is a S-binear field of N. Thus in case of Smarandache bi-invariance (S-bi-invariance) it is only a relative concept as a S-sub-binear-ring M may not be invariant related to every binear field in the binear-ring N.*

**DEFINITION 8.2.10:** *The binear-ring N is said to fulfill the Smarandache left ore condition (S-left ore condition) with respect to a given S-sub-bisemigroup P of (N, •) if for $(s, n) \in S \times N$ there exists $n \bullet s_1 = s \bullet n_1$ ($s_1 \bullet n = n_1 \bullet s$).*

**DEFINITION 8.2.11:** *Let S(BV) denote the set of all Smarandache binear-rings. A S-binear-ring $F_X \in S(BV)$ is called a Smarandache free binear-ring (S-free binear-ring) in BV over X if there exists $f: X \to F_X$ and for all $N \in BV$ and for all $g: X \to N$ there exists a S-binear-ring homomorphism $h \in S(Hom\ F_X, N)$ such that $h \circ f = g$.*

**DEFINITION 8.2.12:** *A finite sequence $N = N_0 \supset N_1 \supset N_2 \supset ... \supset N_t = \{0\}$ of S-sub-binear-rings $N_i$ of N is called a Smarandache normal sequence (S-normal sequence) of N if and only if for all $i \in \{1, 2, ..., n\}$, $N_i$ is an S-bi-ideal of $N_{i-1}$. In the special case that all $N_i$ is an S-bi-ideal of N then we call the normal sequence a Smarandache bi-invariant sequence (S-bi-invariant sequence) and t is called the Smarandache length (S-length) of the sequence $N_{i-1}/N_i$ are called the Smarandache bi-factors (S-bi-factors) of the sequence as*

$$N_{i-1} = N_{i-1}^1 \cup N_{i-1}^2, N_i = N_i^1 \cup N_i^2.$$

*So*

$$N_{i-1}/N_i = \left(N_{i-1}^1/N_i^1\right) \cup \left(N_{i-1}^2/N_i^2\right).$$

**DEFINITION 8.2.13:** *Let P be a S-bi-ideal of the binear-ring N. P is called Smarandache prime bi-ideal (S-prime bi-ideal) if for all S-ideals I and J of N, $IJ \subset P$ implies $I \subset P$ or $J \subset P$.*

**DEFINITION 8.2.14:** *Let P be a S-left bi-ideal of a binear-ring N. P is called Smarandache bimodular (S-bimodular) if and only if there exists a S-idempotent $e \in N$ and for all $n \in N$; $n - ne \in P$.*

The notion of S-idempotents, S-units, S-zero divisors and S-nilpotents are defined in binear-rings as in case of near-rings as every binear-ring is roughly the union of two near-rings.

**DEFINITION 8.2.15:** *Let N be a binear-ring ($N = N_1 \cup N_2$). $z \in N$ is called Smarandache quasi regular (S-quasi regular) if $z \in S(L_Z)$ where $S(L_Z) = \{n_1 - n_1z \mid n_1 \in N_1\} \cup \{n_2 - n_2z \mid n_2 \in N_2\}$. An S-bi-ideal $P \subset N$ is called S-quasi regular if and only if for all $s \in P s$ is S-quasi regular. An S-bi-ideal P is a Smarandache bi-principal bi-ideal (S-bi-principal bi-ideal) if $P = P_1 \cup P_2$ and each of $P_1$ and $P_2$ are S-principal ideals of $N_1$ and $N_2$ where $N = N_1 \cup N_2$ is a binear-ring. A binear-ring N is*

*Smarandache biregular (S-biregular)* if there exists some set E of central S-idempotents with

    i.       *For all $e \in N$, Ne is a S-bi-ideal of N.*
    ii.     *For all $n \in N$ there exists $e \in E$; $Ne = (n)$.*
    iii.    *For all $e, f \in E$, $e + f = f + e$.*
    iv.    *For all $e, f \in E$, $ef$ and $e + f - ef \in E$.*

Now we define Smarandache binear-ring of level II.

**DEFINITION 8.2.16:** *Let N be a binear-ring. We say N is a Smarandache binear-ring of level II (S-binear-ring of level II) if N contains a proper subset P such that P is a biring under the operations of N.*

**Example 8.2.2:** *Let $N = T \times U$ where T is a binear-ring i.e., $T = T_1 \cup T_2$ ($T_1$, $T_2$ near-rings) and $U = U_1 \cup U_2$ where U is a biring where $U_1$ and $U_2$ are rings. Take $P = \{0\} \times U$; P is a S-binear-ring II.*

Now we proceed on to define Smarandache biseminear-rings of level I and level II.

**DEFINITION 8.2.17:** *Let $(N, +, \bullet)$ be a biseminear-ring. We say N is a Smarandache biseminear-ring (S-biseminear-ring) if N contains a proper subset P, such that P is a binear-ring under the operations of N. We call this S-biseminear-ring a S-biseminear-ring of level I.*

Now we define Smarandache biseminear-ring of level II.

**DEFINITION 8.2.18:** *Let $(N, +, \bullet)$ be a seminear-ring. We call N a Smarandache biseminear-ring of level II (S-biseminear-ring of level II) if N contains a proper subset P, such that P under the operations of N is a bisemiring.*

We can obtain classes of S-biseminear-ring II by defining a Smarandache mixed direct product of a bisemirings, i.e., $W = V \times Z$ where $V = V_1 \cup V_2$ is a biseminear-ring and $Z = Z_1 \cup Z_2$ is a bisemiring.

**Example 8.2.3:** *Let $(T, +, \bullet)$ be a non-empty set where T is the S-mixed direct product. $T = V \times Y$ where $V = Z_{10} \cup Z^o$ is a biseminear-ring and $Y = C_2 \cup Z^o$ bisemiring. Clearly T is a S-biseminear-ring II.*

**DEFINITION 8.2.19:** *Let N be a S-binear-ring. N is said to fulfill the Smarandache bi-insertions factors property (S-bi-insertions factors property) if for all $a, b \in N = N_1 \cup N_2$ we have $a \bullet b = 0$ implies $anb = 0$ for all $n \in P$, $P \subset N$ is a P-binear field. (if $a, b \in N_1$ then for all $n \in P_2$ we have $ab = 0$ implies $anb = 0$ and if $a, b \in N_2$ with $ab = 0$, then $\forall n \in P_2$ we have $anb = 0$).*

*We say the S-binear-ring satisfies the Smarandache strong bi-insertion factors property (S-strong bi-insertion factors property) if for all $a, b \in N_i$ ($N = N_1 \cup N_2$), $ab \in I_i$ ($I = I_1 \cup I_2$ a S-bi-ideal of N) implies $anb \in I_i$ where $n \in P_i$ ($P \subset N$, P a binear-ring with $P = P_1 \cup P_2$); $i = 1, 2.$*

**DEFINITION 8.2.20:** *Let p be a prime. A S-binear-ring N is called a Smarandache p-binear-ring (S-p-binear-ring) provided for all $x \in P_i$, $x^p = x$ and $px = 0$, where $P = P_1 \cup P_2$, $P \subset N$; P is a bisemifield.*

**DEFINITION 8.2.21:** *Let N be a binear-ring. A S-rightbi-ideal I of N is called Smarandache right biquasi reflexive (S-right biquasi reflexive) if whenever A and B are S-bi-ideals of N with $AB \subset I$, then $b(b' + a) - bb' \in I$ for all $a \in A$ and for all b, $b' \in B$. A binear-ring N is Smarandache strongly bi-subcommutative (S-strongly bi-commutative) if every S-right bi-ideal of it is S-right biquasi reflexive. Let N be a binear-ring P a S-subnormal bi-subgroup of (N, +). P is called Smarandache quasi bi-ideal (S-quasi bi-ideal) of N if $PN \subset P$ and $NP \subset P$ where by NP we mean elements of the from $\{n(n' + s) - nn'$ for all $s \in P$ and for all n, $n' \in N\} = NP$.*

**DEFINITION 8.2.22:** *A left binear-ring N is said to be a Smarandache left biself distributive (S-left-biself distributive) if the identity abc = abac is satisfied for every a, b, c in $A_i$, i = 1, 2 where $A = A_1 \cup A_2 \subset N$, A a S-sub-binear-ring N to be Smarandache bileft permutable (S-bileft permutable) if abc = bac for all a, b, c $\in A_i$, i = 1, 2 where $A = A_1 \cup A_2$ is a S-sub-binear-ring of N.*

*A binear-ring N is said to be Smarandache biright permutable (S-biright permutable) if abc = acd for all a, b, c $\in A_i$, i = 1, 2, where $A = A_1 \cup A_2 \subset N$ is a S-sub-binear-ring of N. The binear-ring N is said to be Smarandache bimedian (S-bimedian) if abcd = acbd for all a, b, c, d $\in A_i$, i = 1, 2; where $A = A_1 \cup A_2$ is a S-sub-binear-ring of N.*

*Similarly one can define Smarandache biright self distributive (S-biright self distributive) if abc = acbc for all a, b, c in $A_i$; i = 1, 2; $A = A_1 \cup A_2$ is a S-sub-binear-ring of N.*

**DEFINITION 8.2.23:** *Let (R, +, •) be a S-biring and M a S-right R-bimodule. Let $W = R \times M$ and define $(a, s) \odot (\beta, t) = (\alpha\beta, s\beta + t)$. Then $(W, +, \odot)$ is a S-left binear-ring the abstract Smarandache affine binear-ring (S-affine binear-ring) inducted by R and M. All other notions introduced for bi-ideals can also be defined and extended in case of bi-R-modules.*

**DEFINITION 8.2.24:** *Let R be a binear-ring. R is said to be Smarandache bi-equiprime (S-bi-equiprime) if for all $0 \neq a \in P_i$, i = 1, 2; where $P = P_1 \cup P_2 \subset R$ is a binear field in R and for x, y $\in R$, arx = ary for all r $\in R$ implies x = y. If B is a S-biideal of R, B is called a Smarandache equiprime bi-ideal (S-equiprime bi-ideal) if R/B is an equiprime binear-ring.*

Let us define the notion of Smarandache infra biseminear-ring.

**DEFINITION 8.2.25:** *Let (N, +, •) be a triple. N is said to be Smarandache infra binear-ring (S-infra binear-ring) where*

  i. *(N, +) is a S-bisemigroup.*
  ii. *(N, •) is a bisemigroup*
  iii. *$(x + y) \cdot z = x \cdot z - 0 \cdot z + y \cdot z$ for all x, y, z $\in N$.*

**DEFINITION 8.2.26:** *Let I be a S-left bi-ideal of N. Suppose I satisfies the following conditions:*

i. *a, x, y ∈ N, anx – any ∈ $I_k$; k = 1, 2; where I = $I_1 \cup I_2$ for all n ∈ N implies x – y ∈ $I_k$; k = 1, 2.*
ii. *I is left invariant.*
iii. *0N ⊂ I.*

*Then I is called a Smarandache equiprime bileft ideal (S-equiprime bileft ideal) of N. Now we define Smarandache group binear-ring (S-group binear-ring).*

**DEFINITION 8.2.27:** *Let G be any group ; N a binear-ring. Then the group binear-ring NG = $N_1G \cup N_2G$ where N = $N_1 \cup N_2$ is called the S-group binear-ring if and only if N is a S-binear-ring. It is pertinent to note that even if the group binear-ring NG is a S-binear-ring but N is not a S-binear-ring then the group binear-ring NG is not a S-group binear-ring.*

*Example 8.2.4:* Let G = $\langle g \mid g^3 = 1 \rangle$ be the cyclic group of order 3, N = $Z_2 \cup Z$ be a binear-ring. NG is not a S-group binear-ring.

*Example 8.2.5:* Let $Z_3S_3$ and ZG (G = $\langle g \mid g^{11} = 1 \rangle$) be group near-rings; N = $Z_3S_3 \cup$ ZG is a S-binear-ring. It has nothing to do with group binear-rings.

But these methods help us to define more and more classes of S-binear-rings.

**DEFINITION 8.2.28:** *Let T be a semigroup and N a binear-ring the semigroup binear-ring NT is a Smarandache semigroup binear-ring (S-semigroup binear-ring) if and only if T is a S-semigroup. If T is not a S-semigroup still NT happens to be a S-binear-ring then also we do not call NT a S-semigroup binear-ring.*

*Example 8.2.6:* Let $Z_2S(7)$ and $Z_7S$ (where S = {a, b, c, 0 | $a^2 = b^2 = c^2 = 0$, ab = ba = c, ac = ca = b, bc = cb = a} is a semigroup) be semigroup near-rings. Then N = $Z_2S(7) \cup Z_7S$ is a S-binear-ring.

This is certainly not a S-semigroup binear-ring, but only a S-binear-ring.

*Example 8.2.7:* Let $Z_2S(7)$ be a semigroup near-ring and $Z_3D_{27}$ be the group near-ring; then N = $Z_2S(7) \cup Z_3D_{27}$ is a binear-ring infact this is an example of a S-binear-ring.

It is neither a group binear-ring nor a semigroup binear-ring.

**DEFINITION 8.2.29:** *Let N = $N_1 \cup N_2$ be a binear-ring. If the S-bi-ideals of the binear-ring N fulfills bi D.C.C conditions (A.C.C. conditions) then we say N satisfies Smarandache D.C.C condition (S-D.C.C. condition) i.e., if {$N_{1i}$} are the S-bi-ideals of $N_1$ and {$N_{2j}$} are the S-bi-ideals of $N_2$ and if both {$N_{1i}$} and {$N_{2i}$} satisfies the D.C.C condition (A.C.C) then we say the S-bi-ideals of N satisfy Smarandache bi-D.C.C. (Smarandache bi-A.C.C. conditions).*

Such definitions can be made for right bi-ideals, left bi-ideals and ideals, S-bi-D.C.C.R ideals, S-bi-A.C.C.R ideals, S-bi-D.C.C.L ideals, S-bi-A.C.C.L ideals and so on. This notion is also applicable to N-bigroups.

**THEOREM 8.2.3:** *Let NS be a S-semigroup binear-ring; if N is a S-binear-ring, then NS has a S-group binear-ring.*

*Proof:* Follows from the very fact that if S is a S-semigroup then S has a subset P which is a group; since N is a S-binear-ring, NP is a S-group binear-ring.

**DEFINITION 8.2.30:** *Let N be a biseminear-ring i.e. $N = N_1 \cup N_2$ and S be a semigroup. The semigroup biseminear-ring NS is a S-semigroup biseminear-ring if and only if S is a S-semigroup i.e. $NS = N_1S \cup N_2S$ is also a biseminear-ring.*

**DEFINITION 8.2.31:** *Let N be a biseminear-ring and G any group. The group biseminear-ring NG is a Smarandache group biseminear-ring (S-group biseminear-ring) if and only if N is a S-biseminear-ring.*

**THEOREM 8.2.4:** *Let NS be a S-semigroup biseminear-ring, then NS has a sub-biseminear-ring which is a S-group biseminear-ring only if N is a S-biseminear-ring.*

*Proof:* Follows from the very definitions of these concepts, hence left for the reader as an exercise.

**DEFINITION 8.2.32:** *Let N and N' be any two S-biseminear-rings. We say a map $\phi$ from N to N' is a Smarandache biseminear-ring homomorphism (S-biseminear-ring homomorphism) from A to A' where $A \subset N$ and $A' \subset N'$ is a binear-ring if $\phi(x + y) = \phi(x)\phi(y)$, $\phi(xy) = \phi(x)\phi(y)$ where $x, y \in A_i$; $i = 1, 2$; and $\phi(x), \phi(y) \in A'_i$; $i = 1, 2$.*

*Thus for a Smarandache biseminear-ring homomorphism we do not require $\phi$ to be defined on whole of N, it is sufficient if it is defined on a proper near-ring A which is a subset of N. We say the biseminear-ring N is a Smarandache strict biseminear-ring if the subset $A \subset N$ where A is a strict bisemiring.*

*So we can define Smarandache biseminear-ring homomorphism II between two S-biseminear-rings N and N' as $\phi : N \to N'$ is a S-biseminear-ring homomorphism, $\phi : A \to A'$ where A and A' are bisemirings contained in N and N' respectively is a bisemiring homomorphism. Clearly as in case of S-homomorphism I. $\phi$ need not in general be defined on the whole of N.*

**DEFINITION 8.2.33:** *Let N be a Smarandache pseudo biseminear-ring (S-pseudo biseminear-ring) if N is a binear-ring and has proper subset A of N, which is a biseminear-ring under the operations of N.*

The notions of S-pseudo biseminear-rings homomorphism can be defined without any difficulty and the concept of S-pseudo sub-biseminear-rings and S-pseudo bi-ideals can be easily carried out.

**DEFINITION 8.2.34:** *Let N be a S-binear-ring an element $a \in P = P_1 \cup P_2$ (P a binear-field contained in N) is called Smarandache binormal element (S-binormal element) of N, if $aN_1 = N_1a$, $aN_2 = N_2a$, where $N = N_1 \cup N_2$; if $aN_1 = N_1a$ and $aN_2 = N_2a$ for every $a \in P = P_1 \cup P_2$ then N is called the Smarandache binormal near-ring (S-binormal near-ring). Let B(N) denote the set of all S-binormal bielements of N. N is called a S-binormal near-ring if and only if $B(N) = P$ where $P \subset N$ and P is a binear-field.*

The notions of Smarandache binormal sub-bisemigroup and Smarandache bi-invariant can be defined.

**DEFINITION 8.2.35:** *A commutative S-binear-ring N with identity is called a Smarandache Marot binear-ring (S-Marot binear-ring) if each regular S-bi-ideal of N is generated by regular elements where by regular S-bi-ideal we mean a bi-ideal with regular elements that is the non-zero divisors of the binear-ring. We say two S-bi-ideals I and K of the S-N-sub-bigroup G are Smarandache bi-complements (S-bi-complements) if*

  i. $K \cap I = \{0\}$ *i.e.* $K_1 \cap I_1 = \{0\}$ *and* $K_2 \cap I_2 = \{0\}$ *where* $K = K_1 \cup K_2$ *and* $I = I_1 \cup I_2$.
  ii. *If K' is a S-bi-ideal of G such that $K \subset K'$ then $K' \neq K$ implies $K' \cap I \neq \{0\}$.*

*Similar conditions hold for any I' a S-bi-ideal of G.*

**DEFINITION 8.2.36:** *Suppose H and K be any two Smarandache N-sub-bisemigroups (S-N-sub-bisemigroups) of G. K is said to be Smarandache bi-supplement (S-bi-supplement) for H if $H + K = G$ i.e. $H_1 + K_1 = G_1$; $H_2 + K_2 = G_2$ and $H = H_1 \cup H_2$, $K = K_1 \cup K_2$ and $G = G_1 \cup G_2$ and $H + K' \neq G$ for any proper S-bi-ideal K' of K. We say a binear-ring $N = N_1 \cup N_2$ is Smarandache bi-ideally strong (S-bi-ideally strong) if for every S-sub-binear-ring of N is a S-bi-ideal of N.*

**DEFINITION 8.2.37:** *Let N be a binear-ring. $\{I^k\} = \{I_1^k \cup I_2^k\}$ be the collection of all S-bi-ideals of N, N is said to be a Smarandache $I^*$-binear-ring (S-$I^*$-binear-ring) if for every pair of bi-ideals $I^s$, $I^t \in \{I^k\}$ and for every $x \in N \setminus \{I^s \cup I^t\}$ the S-bi-ideal generated by x and $I^s$ and x and $I^t$ are equal i.e. $\langle x \cup I^s \rangle = \langle x \cup I^t \rangle$; here $\langle x \cup I_1^s \rangle = \langle x \cup I_1^t \rangle$; if $x \in N_1$ and if $x \in N_2$ then $\langle x \cup I_2^s \rangle = \langle x \cup I_2^t \rangle$ where $I^s = I_1^s \cup I_2^s$ and $I^t = I_1^t \cup I_2^t$.*

Now we proceed on to define Smarandache quad near-ring and Smarandache quad seminear-ring.

**DEFINITION 8.2.38:** *Let $(N, +, \bullet)$ be a non-empty set. We say N is a Smarandache quad near-ring (S-quad near-ring) if the following conditions are true:*

  i. $N = N_1 \cup N_2 \cup N_3 \cup N_4$ *where $N_i$'s are proper subsets of N such that $N_i \not\subset N_j$ if $i \neq j$, $1 \leq i, j \leq 4$.*

ii. Each $(N_i, +, \bullet)$ is a near-ring, $i = 1, 2, 3, 4$ under the same operations of N.
iii. At least one of the $N_i$'s is a S-near-ring.

**DEFINITION 8.2.39:** *Let $(N, +, \bullet)$ be a binear-ring and G a bigroup. The bigroup binear-ring $NG = N_1G_1 \cup N_1G_2 \cup N_2G_1 \cup N_2G_2$ where each of the $N_iG_j$, $1 \le i, j \le 2$ are group near-rings. NG is a Smarandache bigroup binear-ring (S-bigroup binear-ring) if and only if G is a S-bigroup and N is a S-binear-ring.*

**THEOREM 8.2.5:** *Every bigroup binear-ring, which is a S-bigroup binear-ring, is a S-quad near-ring.*

*Proof:* Left for the reader as an exercise.

**DEFINITION 8.2.40:** *Let $(N, +, \bullet)$ be a binear-ring and S be a bisemigroup. The bisemigroup binear-ring NS is a Smarandache bisemigroup binear-ring (S-bisemigroup binear-ring) if and only if S is a S-bisemigroup and N is a S-binear-ring. That is if $S = S_1 \cup S_2$ and $N = N_1 \cup N_2$ then $NS = N_1S_1 \cup N_1S_2 \cup N_2S_1 \cup N_2S_2$.*

**THEOREM 8.2.6:** *Let NS be a S-bisemigroup biseminear-ring with N a S-seminear-ring. Then NS has a S-group biseminear-ring.*

*Proof:* Straightforward by the definitions hence left for the reader to prove.

**DEFINITION 8.2.41:** *Let $(N, +, \bullet)$ be a non-empty set. N is said to be a Smarandache quad seminear-ring (S-quad seminear-ring) if*

i. $N = N_1 \cup N_2 \cup N_3 \cup N_4$. $N_i$'s proper subsets of N with $N_i \not\subset N_j$ $(i \ne j)$.
ii. Each $(N_i, +, \bullet)$ is a seminear-ring.
iii. At least one of the $N_i$'s is a S-seminear-ring I (II).

**THEOREM 8.2.7:** *Every S-semigroup biseminear-ring is a S-quad seminear-ring.*

*Proof*: Left as an exercise for the reader to prove.

Now we proceed onto define non-associative (na) Smarandache binear-ring, na-Smarandache biseminear-ring, na-Smarandache quad near-ring and na-Smarandache quad seminear-ring.

**DEFINITION 8.2.42:** *Let $(N, +, \bullet)$ be a na-binear-ring. We say N is a Smarandache non-associative binear-ring (S-non-associative binear-ring) if*

i. $(N, +)$ is a S-bigroup.
ii. $(N, \bullet)$ is a S-bigroupoid.
iii. $(a + b) \bullet c = a \bullet c + b \bullet c$ for all $a, b, c \in N$.

**DEFINITION 8.2.43:** *Let $(N, +, \bullet)$ be a non-empty set, we say N is a Smarandache non-associative biseminear-ring (S-non-associative biseminear-ring) if*

i. $(N, +)$ is a S-bisemigroup

ii. $(N, \bullet)$ is a S-bigroupoid.
  iii. $(a + b) \bullet c = a \bullet c + b \bullet c$ for all $a, b, c \in N$.

**DEFINITION 8.2.44:** *Let $(N, +, \bullet)$ be a na-binear-ring if $N$ has a proper subset $P$, where $P$ is an associative binear-ring then we call $N$ a Smarandache non-associative binear-ring of level II (S-non-associative binear-ring of level II).*

**DEFINITION 8.2.45:** *Let $(N, +, \bullet)$ be a na-biseminear-ring, if $N$ has a proper subset $P$ such that $P$ is an associative biseminear-ring then we call $N$ a Smarandache non-associative biseminear-ring of level II (S-non-associative biseminear-ring of level II).*

**DEFINITION 8.2.46:** *Let $(N, +, \bullet)$ be a na-quad near-ring we call $N$ a Smarandache non-associative quad near-ring (S-non-associative quad near-ring) if $N$ has a proper subset $P$ such that $P$ is an associative quad near-ring.*

**DEFINITION 8.2.47:** *Let $(N, +, \bullet)$ be a na quad-seminear-ring, we call $N$ a Smarandache non-associative quad seminear-ring (S-non-associative quad seminear-ring) if $N$ has a proper subset $P$ such that $P$ is an associative quad seminear-ring.*

All properties studied and defined in case of Smarandache associative binear-rings and Smarandache associative biseminear-rings can be very easily studied and defined for Smarandache non-associative biseminear-rings and Smarandache non-associative binear-rings. The only concept non-existent in case of Smarandache associative binear-rings and S-associative biseminear-rings is the notions of identities like Moufang, Bol, Bruck, WIP, alternative (right/ left) etc. Here we just define for one of them and the reader can define and study in case of other identities.

Before we go for the study of these identities we proceed on to give methods by which we can get S-na-binear-rings and S-na-biseminear-rings. Natural examples of these are non-existent.

**DEFINITION 8.2.48:** *Let $L$ be a loop, $N$ a binear-ring, $NL$ be the loop binear-ring i.e. $NL = N_1L \cup N_2L$ where $N = N_1 \cup N_2$ with $N_iL$ loop near-rings, $i = 1, 2$. Clearly $NL$ is a S-na-near-ring.*

**DEFINITION 8.2.49:** *Let $G$ be a groupoid and $N$ a binear-ring. $NG$ be the groupoid binear-ring i.e. $NG = N_1G \cup N_2G$ where $N_iG$ are groupoid near-rings for $i = 1, 2$. $NG$ is a S-na binear-ring.*

**DEFINITION 8.2.50:** *Let $L$ be a biloop and $N$ a near-ring, $NL$ be the biloop near-ring i.e. $NL = NL_1 \cup NL_2$ where $L = L_1 \cup L_2$. Clearly $NL$ is a non-associative binear-ring.*

**THEOREM 8.2.8:** *If $NL$ is a bilooop near-ring, $NL$ in general is not a S-binear-ring.*

*Proof*: If $L = L_1 \cup L_2$ is a biloop such that none of the loops $L_1$ and $L_2$ are S-loops then $NL$ the biloop near-ring is not a S-binear-ring.

**THEOREM 8.2.9:** *Let $N$ be a near-ring and $L = L_1 \cup L_2$ be a S-biloop, then the biloop near-ring is a S-na-binear-ring.*

*Proof:* Follows from the fact that if $L = L_1 \cup L_2$ is a S-biloop then L contains a bigroup $G = G_1 \cup G_2$. Take $NG = NG_1 \cup NG_2$ where $NG \subset NL$ so NL is a S-na-binear-ring.

**DEFINITION 8.2.51:** *Let N be a near-ring and G be a bigroupoid. The bigroupoid near-ring is a binear-ring.*

**DEFINITION 8.2.52:** *Let N be a binear-ring and G be a S-bigroupoid. The bigroupoid near-ring is a S-na-binear-ring.*

**DEFINITION 8.2.53:** *Let L be a biloop and N be a binear-ring. The biloop binear-ring $NL = N_1L_1 \cup N_1L_2 \cup N_2L_1 \cup N_2L_2$.*

**DEFINITION 8.2.54:** *Let $(N, +, \bullet)$ be a non-empty set. We call N a Smarandache non-associative quad near-ring (S-non-associative quad near-ring) if the following conditions are true:*

  i. $N = N_1 \cup N_2 \cup N_3 \cup N_4$ where $N_i$'s are proper subsets of N and $N_i \not\subset N_j$ if $i \neq j$, $1 \leq i, j \leq 4$.
  ii. All $(N_i, +, \bullet)$ are near-rings under the operations of N, associative or non-associative.
  iii. At least one of the $(N_i, +, \bullet)$ is a S-non-associative near-ring.

**THEOREM 8.2.10:** *The biloop binear-ring NL is a S-na-quad near-ring where L is a biloop and N is an associative near-ring.*

*Proof:* Follows from the definitions, hence left as an exercise for the reader to prove.

**DEFINITION 8.2.55:** *Let $(N, +, \bullet)$ be an associative biseminear-ring. L any biloop. The biloop biseminear-ring $NL = N_1L_1 \cup N_1L_2 \cup N_2L_1 \cup N_2L_2$.*

**DEFINITION 8.2.56:** *Let $(N, +, \bullet)$ be a non-empty set. We call N a Smarandache na-quad seminear-ring (S-na-quad seminear-ring) if the following conditions are true:*

  i. $N = N_1 \cup N_2 \cup N_3 \cup N_4$ where each $N_i$ is a proper subset of N with $N_i \not\subset N_j$ if $i \neq j$.
  ii. $(N_i, +, \bullet)$ are seminear-rings associative or non-associative, $i = 1, 2, 3, 4$.
  iii. At least one of $(N_i, +, \bullet)$ is a S-na-seminear-ring.

**THEOREM 8.2.11:** *Let $(N, +, \bullet)$ be a biseminear-ring. L be a biloop. The biloop biseminear-ring NL is a S-na quad seminear-ring.*

*Proof:* Follows from the definitions.

Thus S-biloop binear-rings, S-biloop near-rings, S-biloop seminear-rings, S-bigroupoid near-rings, S-bigroupoid seminear-rings etc. leads to several classes of examples of S-na binear-ring, S-na biseminear-rings, S-na-quad near-rings and so on.

Now we proceed on to define special identities in S-na-binear-rings.

**DEFINITION 8.2.57:** *Let (N, +, •) be a na-binear-ring. We call N a Smarandache Moufang binear-ring (S-Moufang binear-ring) if N has a proper subset $P = P_1 \cup P_2$ where P is a S-sub-binearring of N and for all x, y, z $\in P_i$ we have (xy)(zx) = (x(yz))x, i = 1, 2.*

If every proper subset P of N which are S-sub-binear-rings of N satisfies the Moufang identity then we call the binear-ring a Smarandache strong Moufang binear-ring. On similar lines we can define Smarandache Bol binear-ring, Smarandache Bruck binear-ring etc. and their Smarandache strong Bol binear-ring, Smarandache strong Bruck binear-ring and so on.

The concept of S-zero divisors, S-idempotents and S-units are defined as in case of S-na-near-rings and S-na-binear-rings in a similar way. The concept of S-sub-binear-rings and S-bi-ideals are defined in an analogous way.

**PROBLEMS:**

1. Find the order of the smallest S-binear-ring.
2. Can a S-biseminear-ring of order ten exist?
3. Is $N = Z_5L_5(3) \cup Z_7(Z_9(3, 5))$ a S-na-binear-ring, where $Z_2$ and $Z_7$ are near-rings, $L_5(3)$, the loop of order 6 and $Z_9(3, 5)$ is a groupoid?
4. Find the sub-binear-ring of N given in problem 3.
5. Give an example of a Moufang binear-ring.
6. Give an example of a S-Bol biseminear-ring.
7. Can a S-alternative biseminear-ring of order 12 exist?
8. Using biloop and binear-ring can we find biloop binear-ring of order 16?
9. Using a bisemigroup S and a near-ring N, can NS be a biseminear-ring?
10. Using a bigroup G and a seminear-ring N, can NG be a biseminear-ring?
11. What is the order of the smallest S-quad near-ring?
12. Find the smallest order of the S-quad seminear-ring.

## 8.3 Generalizations, S-analogue and its applications

In this section we define some generalizations called biquasi structures and we give the Smarandache analogue. Further this section gives the applications of binear-ring and bisemi near-rings to automatons and the use of binear-rings to error correcting codes.

**DEFINITION 8.3.1:** *Let (N, +, •) be a non-empty set we say N is a biquasi ring if the following conditions are true.*

  i.   *$N = N_1 \cup N_2$ ($N_1$ and $N_2$ are proper subsets of N)*
  ii.  *$(N_1, +, •)$ is a ring.*
  iii. *$(N_2, +, •)$ is a near-ring.*

***Example 8.3.1:*** *Let $N = Z_{12} \cup Z$ where $Z_{12}$ is a ring and Z is a near-ring. Clearly N is a biquasi ring.*

**DEFINITION 8.3.2:** *Let $(P, +, \bullet)$ be a non-empty set with two binary operations '+' and '$\bullet$' P is said to be a biquasi semiring if*

  i.  $P = P_1 \cup P_2$ *are proper subsets.*
  ii. $(P_1, +, \bullet)$ *is a ring.*
  iii. $(P_2, +, \bullet)$ *is a semiring.*

*Example 8.3.2:* Let $N = Q \cup C_6$ where Q is a ring and $C_6$ is a chain lattice; N is a biquasi semiring.

**DEFINITION 8.3.3:** *Let $(T, +, \bullet)$ be a non-empty set with two binary operation. We call T a biquasi near-ring if*

  i.  $T = T_1 \cup T_2$ *where $T_1$ and $T_2$ are proper subsets of T.*
  ii. $(T_1, +, \bullet)$ *is a semiring.*
  iii. $(T_2, +, \bullet)$ *is a near-ring.*

*Example 8.3.3:* Let $N = C_7 \cup Z_2$ where $C_7$ is a chain lattice which is a semiring and $Z_2$ is a near-ring. Thus N is a biquasi near-ring.

**DEFINITION 8.3.4:** *Let $(N, +, \bullet)$ be a biquasi ring. A proper subset $P \subset N$, $P = P_1 \cup P_2$ is said to be a sub-biquasi ring if $(P, +, \bullet)$ is itself a biquasi ring.*

*Example 8.3.4:* Let $(N, +, \bullet)$ be a biquasi ring, where $N = Q \cup Z$. Clearly N has a proper subset $P = 2Z \cup 3Z$ to be sub-biquasi ring.

**DEFINITION 8.3.5:** *Let $(N, +, \bullet)$ be a biquasi semiring. A proper subset $P \subset N$ where $P = P_1 \cup P_2$ is a sub-bi-quasi semiring if P itself is a biquasi semiring.*

*Example 8.3.5:* Let $N = Z \cup L$ where L is a lattice given by the following diagram:

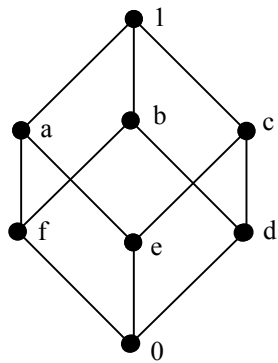

**Figure 8.3.1**

Take $P = 5Z \cup \{0, 1\}$. N is a biquasi semiring. P is a sub-biquasi semiring.

Similarly we define sub-biquasi near-ring for biquasi near-ring. All results related to biring, bisemiring and binear-ring can be extended to biquasi rings, biquasi semirings and biquasi near-rings.

Now we define non-associative biquasi rings as follows:

**DEFINITION 8.3.6:** *Let $(N, +, \bullet)$ be a biquasi ring we call N a non-associative biquasi ring if $(N_1, +, \bullet)$ which is a ring is a non-associative ring and $(N_2, +, \bullet)$ is a near-ring.*

**DEFINITION 8.3.7:** *Let $(N, +, \bullet)$ be a biquasi semiring we say N a non-associative biquasi semiring if $(N_2, +, \bullet)$ which is a semiring is a non-associative semiring.*

*Example 8.3.6:* Let $N = Q \cup Z^o L_5(3)$ be the biquasi semiring. It is easily verified that N is a non-associative biquasi semiring.

*Example 8.3.7:* Let $N = Z_{12} L_9 (8) \cup Z$ be a biquasi ring; N is a non-associative biquasi ring.

**DEFINITION 8.3.8:** *Let $(N, +, \bullet)$ be a biquasi near-ring. We say N is a non-associative biquasi near-ring where $N = Z^0 \cup ZL_{11}(3)$ where $Z^0$ is a semiring and $ZL_{11}(3)$ is the loop near-ring which is non-associative, hence N is na biquasi near-ring.*

The study of when these na biquasi structures satisfy the typical identities like Moufang, Bol, Bruck, W.I.P can be made as in case of any other non-associative structures. Substructures like biquasi ideals etc can also be defined.

Now we proceed for their Smarandache analogue.

**DEFINITION 8.3.9:** *Let $(N, +, \bullet)$ where $N = N_1 \cup N_2$ be a biquasi ring. We call N a Smarandache biquasi ring (S-biquasi ring) if $(N_1, +, \bullet)$ is a S-ring.*

**DEFINITION 8.3.10:** *Let $(N, +, \bullet)$ be a biquasi semiring. N is said to be a Smarandache biquasi semiring (S-biquasi semiring) if the semiring $(N_2, +, \bullet)$ is a S-semiring.*

**DEFINITION 8.3.11:** *Let $(T, +, \bullet)$ be a non-empty set which is a biquasi near-ring. We call T a Smarandache biquasi near-ring (S-biquasi near-ring) if the near-ring $(T_2, +, \bullet)$ is a S-near-ring.*

We can define the concept of S-sub-biquasi ring, S-sub-biquasi semiring and S-sub-biquasi near-rings and the corresponding ideas as follows.

**DEFINITION 8.3.12:** *Let $(N, +, \bullet)$ be a biquasi ring a proper subset, $(P, +, \bullet)$ of N is said to be a Smarandache sub-biquasi ring (S-sub-biquasi ring) if $P_1 \cup P_2$ is the union of proper subsets; $(P_1, +, \bullet)$ is a S-subring of $N_1$ and $P_2$ is a subnear-ring of $P_2$.*

**DEFINITION 8.3.13:** *Let $(P, +, \bullet)$ be a biquasi semiring. A proper subset $(X, +, \bullet)$ of P is said to be a Smarandache sub-biquasi semiring (S-biquasi near-ring) if $(X_2, +, \bullet)$ is a S-semiring where $X = X_1 \cup X_2$.*

**DEFINITION 8.3.14:** *Let $(S, +, \bullet)$ be a biquasi near-ring. A proper subset $(L, +, \bullet)$ of S is said to be a Smarandache sub-bi quasi near-ring (S-sub-bi quasi near-ring) if $(L_2, +, \bullet)$ is a S-near-ring where $L = L_1 \cup L_2$ and $(L_1, +, \bullet)$ is a semiring.*

Now we define Smarandache bipseudo structures.

**DEFINITION 8.3.15:** *$(Y, +, \bullet)$ is a Smarandache bi pseudo ring (S-bi pseudo ring) if the following conditions are true.*

  i. *$Y = Y_1 \cup Y_2$ are proper subsets of Y.*
  ii. *$(Y_1, +, \bullet)$ is a S-ring.*
  iii. *$(Y_2, +, \bullet)$ is a S-near-ring.*

**DEFINITION 8.3.16:** *Let $(P, +, \bullet)$ be a non-empty set. P is said to be a Smarandache bi pseudo semiring (S-bi pseudo semiring) if the following conditions are true.*

  i. *$P = P_1 \cup P_2$ where $P_1$ and $P_2$ are proper subsets of P.*
  ii. *$(P_1, +, \bullet)$ is a S-ring.*
  iii. *$(P_2, +, \bullet)$ is a S-semiring.*

**DEFINITION 8.3.17:** *Let $(X, +, \bullet)$ be a non-empty set, X is said to be a Smarandache bi pseudo near-ring (S-bi pseudo near-ring) if the following conditions are true.*

  i. *$X = X_1 \cup X_2$ is the union of two proper subsets.*
  ii. *$(X_1, +, \bullet)$ is a S-semiring.*
  iii. *$(X_2, +, \bullet)$ is a S-near-ring.*

**THEOREM 8.3.1:** *Let $(N, +, \bullet)$ be a Smarandache bi pseudo structure then $(N, +, \bullet)$ is a S-bistructure and also S-biquasi structure.*

*Proof:* The proof is left for the reader as the solution is direct by the very definitions.

We define now the notion of bi-ideals.

**DEFINITION 8.3.18:** *Let $(N, +, \bullet)$ be a biquasi ring. A non-empty subset $(P, +, \bullet)$ is called the biquasi bi-ideal of N if $P = P_1 \cup P_2$ and $P_1$ is an ideal of the ring $(N_1, +, \bullet)$ and $P_2$ is an ideal of the near-ring $(N_2, +, \bullet)$. We define $P = P_1 \cup P_2$ to be a Smarandache biquasi bi-ideal (S-biquasi bi-ideal) if $P_1$ is an S-ideal of the ring $(N_1, +, \bullet)$ and $P_2$ is an S-ideal of the near-ring $(N_2, +, \bullet)$.*

On similar lines we define bi-ideals and S-bi-ideals of other biquasi structures.

Further for other bistructures also bi-ideals and S-bi-ideals are defined in an analogous way.

All properties studied/ introduced in case of any bialgebraic structure can also be extended to the bialgebraic quasi structures and bipseudo structures and their Smarandache analogue. Now as a final resort we define Smarandache mixed bialgebraic structures and mixed bialgebraic structures and indicate some interesting results about them.

**DEFINITION 8.3.19:** *Let $(N, +, \bullet)$ be a non-empty set we call N a group-semigroup if the following conditions are true.*

  i. $N = N_1 \cup N_2$; $N_1$ and $N_2$ are proper subsets of N.
  ii. $(N, +)$ is a group.
  iii. $(N, \bullet)$ is a semigroup.

*Example 8.3.8:* Let $(N = Z_{10} \cup Z^o, +, \bullet)$. Clearly N is a group-semigroup for $(Z_{10}, +)$ is a group and $(Z^o, \bullet)$ is a semigroup.

**DEFINITION 8.3.20:** *Let $(P, +, \bullet)$ be a non-empty set. We call P a group-groupoid if the following conditions are true:*

  i. $P = P_1 \cup P_2$ is a union of proper subsets.
  ii. $(P_1, +)$ is a group.
  iii. $(P_2, \bullet)$ is a groupoid.

*Example 8.3.9:* Let $(P = S_3 \cup L_5(2), +, \bullet)$. Clearly $(S_3, \bullet)$ is a group and $(L_5(2), \bullet)$ is a groupoid. Hence P is a group-groupoid.

**DEFINITION 8.3.21:** *Let $(X, +, \bullet)$ be a non-empty set. We say X is a loop-semigroup if the following conditions are satisfied:*

  i. $X = X_1 \cup X_2$ is the proper union of subsets.
  ii. $(X_1, +)$ is a loop.
  iii. $(X_2, \bullet)$ is a semigroup.

*Example 8.3.10:* Let $(P = L_7(3) \cup S(3), +, \bullet)$ is a loop-semigroup.

**DEFINITION 8.3.22:** *Let $(N, +, \bullet)$ be a non-empty set we call N a ring-group if the following conditions are true:*

  i. $N = N_1 \cup N_2$ is the union of proper subsets.
  ii. $(N_1, +, \bullet)$ is a ring.
  iii. $(N_2, \bullet)$ or $(N_2, +)$ is a group.

*Example 8.3.11:* Let $(N = Z \cup S_3, +, \bullet)$ is a ring – group

  i. $(Z, +, \bullet)$ is a ring
  ii. $(S_3, \bullet)$ is a group.

We see N is not a group ring but this set has both flavors. On similar lines one can define semiring group, near-ring group, ring-semigroup, ring-loop and so on.

We define their Smarandache analogues in a very special way.

**DEFINITION 8.3.23:** *Let $(N, +, \bullet)$ be a group-semigroup. We call N a Smarandache group semigroup (S-group-semigroup) if and only if the semigroup is a S-semigroup.*

***Example 8.3.12:*** Let $(N = D_{28} \cup S(7), +, \bullet)$, N is a S-group semigroup.

**DEFINITION 8.3.24:** *Let $(N, +, \bullet)$ be a loop-semigroup, we call N a Smarandache loop-semigroup (S-loop-semigroup) if $(N_1, +)$ is a S-loop and $(N_2, \bullet)$ is a semigroup.*

*If both $(N_1, +)$ and $(N_2, \bullet)$ are S-loops and S-semigroup then we call N a Smarandache strong loop-semigroup (S-strong loop-semigroup).*

**DEFINITION 8.3.25:** *Let $(N, +, \bullet)$ be a ring-group we call N a Smarandache ring-group (S-ring-group) if $(N_1, +, \bullet)$ is a S-ring.*

Thus Smarandache concepts and Smarandache strong concepts can be defined and studied. Now we proceed on to give some applications.

We just recall the definition of a Smarandache semigroup semi automaton.

**DEFINITION 8.3.26:** *Let G be a S-semigroup. A Smarandache S-semigroup semi automaton (S-semigroup semi automaton) in an ordered triple $S(S) = (G, I, \delta)$ where G is a S-semigroup, I is a set of inputs and $\delta: G \times I \to G$ is a transition function. It is easily verified that all group semi automatons are S-S-semigroup semi automatons.*

Thus we see the class of S-S-semigroup semi automatons properly contains the class of group semi-automaton.

***Example 8.3.13:*** Let $G = \{0, 1, 2, \ldots, 5\} = Z_6$ and $I = \{0, 1, 2\}$ be the set of states; define the transition function $\delta : G \times I \to G$ by $\delta(g, i) = g.i \pmod 6$ thus we get the S-semigroup semi-automaton. We give the following state graph representation for this S-S semigroup semi-automaton.

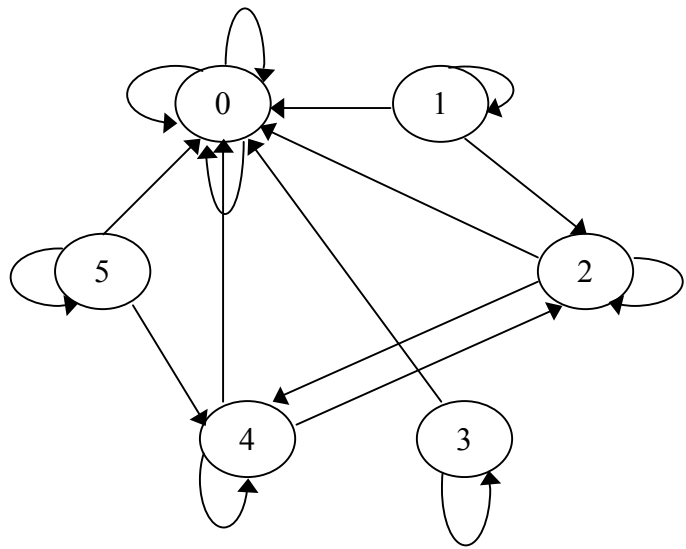

**Figure 8.3.2**

Now we see when we change the function $\delta_1$ for the same S-semigroup, G and $I = \{0, 1, 2\}$ we get a different S-S-semigroup semiautomaton which is given by the following figure; here the transition fuctions. $\delta_1 : G \times I \to G$ is defined by $\delta_1 (g, i) = (g + i) \bmod 6$. We get a nice symmetric S-S-semigroup semiautomaton.

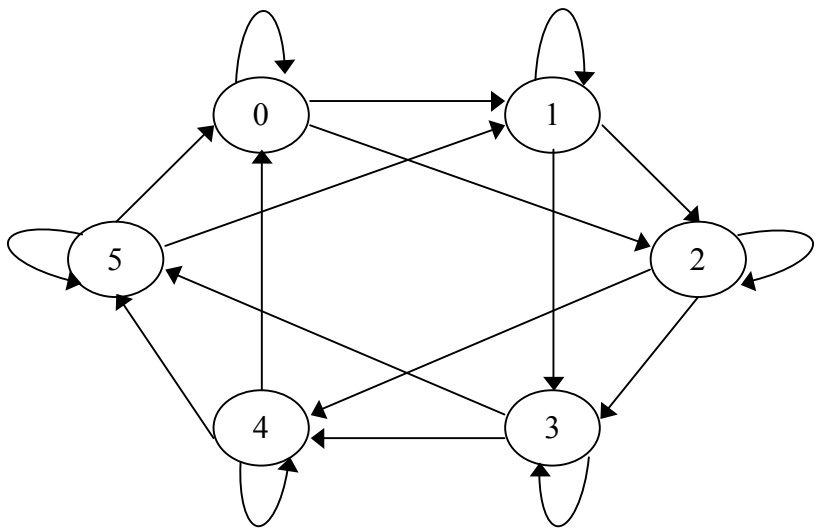

**Figure 8.3.3**

**DEFINITION 8.3.27**: *Let S(S) = (G, I, δ) where G is a S-semigroup I the set of input alphabets and δ the transition function as in case of S-S semigroup semi-automaton. We call a proper subset B ⊂ G with B a S-subsemigroup of G and $I_1$ ⊂ I denoted (B, $I_1$, δ) = S(SB) is called the Smarandache subsemigroup semi-automaton (S-subsemigroup semiautomaton) if δ : B × $I_1$ → B. A S-S-semigroup semi-automaton may fail to have sometimes S-S-subsemigroup subsemi-automaton.*

*Example 8.3.14*: In example 8.3.13 with S = {$Z_6$, I, δ} when we take B = {0, 2, 4} we get a S-S-subsemigroup subsemi-automaton

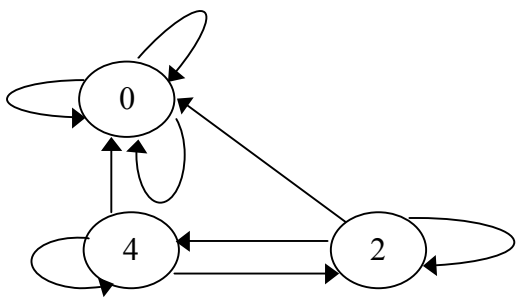

**Figure 8.3.4**

as B = {0, 2, 4} is a S-subsemigroup of $Z_6$.

Now to the best of our knowledge we have not seen the definition of group automaton we define it and then proceed on to get corresponding near-rings.

**DEFINITION 8.3.28**: *A group automaton is an ordered quintuple A = {G, A, B, δ, λ} where G is a group, A is a set of inputs, B is a set of outputs δ : G × A → G and λ : G × A → B is the transition and output functions respectively.*

The Smarandache analogous of this would be as follows.

**DEFINITION 8.3.29**: *A Smarandache S-semigroup automaton (S-S semigroup automaton) is an ordered quintuple $S(A) = (G, A, B, \delta, \lambda)$ where G is a S-semigroup and A set of input alphabets and B set of output alphabets, $\delta$ and $\lambda$ are the transition and output functions respectively.*

*If $\delta$ restricted form $G \times A \to G$ maps a subgroup $T \subset G$ into itself $\delta: T \times A \to T$ then we can say the S-S-semigroup automaton has $(T \subset G, A, B, \delta, \lambda)$ to be a group automaton.*

We call this level I S-S-semigroup automaton. Now we proceed on to define S-S-semigroup semiautomaton and S-S-semigroup automaton of level II.

**DEFINITION 8.3.30**: *Let $S = (G, I, \delta)$ be an ordered triple where G is only a S-semigroup. If $B \subset G$ and B is a group under the operations of G, and if $S^1 = (B, I, \delta)$ is a group semi-automaton with B a group, I is a set of inputs and $\delta: B \times I \to B$ then $S = (G, I, \delta)$ is called the Smarandache S-semigroup semi-automaton of level II (S-S-semigroup semiautomaton II). Similarly we define Smarandache S-semigroup automaton of level II (S-S-semigroup automaton II).*

It is left for the reader to find any relation existing between the two levels of automaton. Now we proceed on to see how these concepts will produce for us the S-near-rings.

**DEFINITION 8.3.31**: *Let $S = (Q, X, \delta)$ be a S-S-semigroup semi-automaton where Q is a S-semigroup under addition '+' i.e. $G \subset Q$ is such that G is a S-semigroup under '+' and G a proper subset of Q/ Let $\delta: G \times X \to G$ is well defined i.e. restriction of $\delta$ to G is well defined so $(G, X, \delta) = C$ becomes a semi-automaton, C is called additive if there is some $x_0 \in X$ with*

$$\delta(q, x) = \delta(q, x_0) + \delta(0, x) \text{ and}$$
$$\delta(q - q', x_0) = \delta(q, x_0) - \delta(q', x_0) \text{ for all } q, q' \in G \text{ and } x \in X.$$

*Then there is some homomorphism*

$$\psi: G \to G \text{ and some map}$$
$$\alpha: X \to G \text{ with } \psi(x_0) = 0 \text{ and } \delta(q, x) = \psi(q) + \alpha(x).$$

*Let $\delta_x$ for a fixed $x \in X$ be the map from $G \to G$; $q \to \delta(q, x)$. Then $\{\delta_x / x \in X\}$ generates a subnear-ring N (P) of the near-ring $(M(G), +, \bullet)$ of all mappings on G; this near-ring N(P) is called as the Smarandache syntactic near-ring (S-syntactic near-ring).*

The varied structural properties are enjoyed by this S-syntactic near-ring for varying S-semigroups. The most interesting feature about these S-syntactic near-ring is that for one S-S-semigroup semi-automaton we have several group semiautomaton depending on the number of valid groups in the S-semigroup. This is the vital benefit in defining S-S-semigroup automaton and S-syntactic near-ring. So for a given S-S-semigroup semi-automaton we can have several S-syntactic near-ring.

Thus the Smarandache notions in this direction has evolved several group automaton for a given S-semigroup. Now when we proceed on to use the bisemigroup and binear-ring structure to define these concepts we get many multipurpose S-S-semi automatons i.e. a single S-S-semigroup automaton can serve several purposes as a finite machine.

**DEFINITION 8.3.32:** *Let $(G, +, \bullet)$ be a S-bisemigroup. A Smarandache S-bisemigroup semi bi-automaton (S-S-bisemigroup semi bi-automaton) is an ordered triple SB (S) = $(G = G_1 \cup G_2, I = I_1 \cup I_2, \delta = \delta_1 \cup \delta_2)$ where G is a S-bisemigroup, $I = I_1 \cup I_2$ is a set of inputs $I_1$ and $I_2$ are proper subsets of I $\delta_i$: $G_i \times I_i \to G_i$ is a transition function i = 1, 2.*

*Example 8.3.15:* $G = Z_3 \cup Z_4$, $I = \{0, 1, 2\} \cup \{a_1, a_2\}$; $\delta = (\delta_1, \delta_2)$ $\delta_1 : Z_3 \times \{2\} \to Z_3$, $\delta_2 : Z_4 \times \{a_1, a_2\} \to Z_4$ given by the following tables:

| $\delta_1$ | 0 | 1 | 2 |
|---|---|---|---|
| 0 | 0 | 0 | 0 |
| 1 | 0 | 1 | 2 |
| 2 | 0 | 2 | 1 |

| $\delta_2$ | 0 | 1 | 2 | 3 |
|---|---|---|---|---|
| $a_1$ | 0 | 2 | 0 | 2 |
| $a_2$ | 1 | 3 | 1 | 3 |

$(G, I, \delta)$ is a S-bisemigroup bisemi-automaton. G is a S-bisemigroup under multiplication.

Thus these two graphs together give the S-bisemigroup bisemi-automaton.

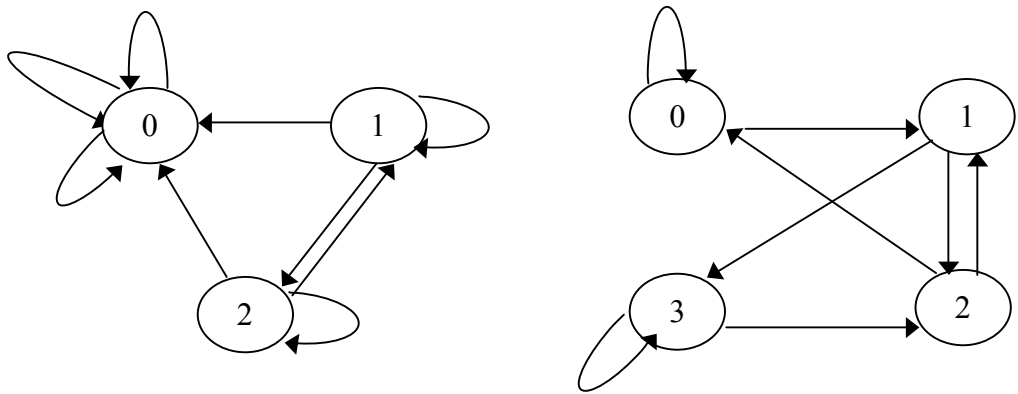

**Figure 8.3.5**

**DEFINITION 8.3.33:** *Let $SB(S) = (G = G_1 \cup G_2, I = I_1 \cup I_2, \delta = (\delta_1, \delta_2))$ be a S-S-semigroup bisemi-automaton. A proper subset $B \subset G$ where $B = B_1 \cup B_2$, $(B_1 \subset G_1$ and $B_2 \subset G_2)$ is called the Smarandache S-sub-bisemigroup bisemi-automaton if $\delta / B_i : B_i \times I_i \to B_i$, i = 1, 2.*

The above example has no S-S-sub-bisemigroup bisemi-automaton.

**DEFINITION 8.3.34:** *A bigroup biautomaton is an ordered quintuple, $A = \{G, A, B, \delta = (\delta_1, \delta_2), \lambda = (\lambda_1, \lambda_2)\}$ where G is a bigroup i.e. $G = G_1 \cup G_2$, A set of inputs, B is set of outputs $\delta_i : G_i \times A \to G_i$, $\lambda_i : G_i \times A \to B$, $i = 1, 2$; is the transition and output functions respectively.*

The Smarandache analogous would be as follows:

**DEFINITION 8.3.35:** *A Smarandache S-bisemigroup biautomaton (S-S-besemigroup biautomaton) is an ordered quintuple $SS(A) = (G = G_1 \cup G_2, A = A_1 \cup A_2, B = B_1 \cup B_2, \delta = (\delta_1, \delta_2), \lambda = (\lambda_1, \lambda_2))$ where G is a S-bisemigroup and A and B are bisets which are input and output alphabets $\delta$ and $\lambda$ are the transition and output functions respectively.*

*If $\delta_i : G_i \times A_i \to G_i$ maps a sub-bigroup $T_i \subset G_i$ into itself $\delta_i : T_i \times A_i \to T_i$ then we can say the SS(A) bi-automaton has $(T \subset G, A, B, \delta, \lambda)$ to be a group bi-automaton.*

We can define level II bi-automaton in a similar way. The task of doing this is left as an exercise to the reader.

**DEFINITION 8.3.36:** *Let $S = (Q, X, \delta)$ be a S-S bisemigroup bisemi-automaton where Q is a S-bisemigroup under '+' i.e. $G \subset Q$ is such that G is a S-bisemigroup under '+' and G is a proper subset of Q. Let $\delta_i : G_i \times X \to G_i$ is well defined that is restriction of $\delta$ to $G_i$ is well defined so that $(G, X, \delta) = C$ becomes a bisemi-automaton. C is called additive if there is some $x_0 \in X$ with $\delta(q, x) = \delta(q, x_0) + \delta(0, x)$ and $\delta(q - q', x) = \delta(q, x) - \delta(q', x)$ for all $q, q' \in G = G_1 \cup G_2$ and $x \in X$.*

*Then there is some bi-homomorphism*

$$\psi : G \to G \text{ and some map}$$
$$\alpha : X \to G \text{ with}$$
$$\psi(x_0) = 0 \text{ and } \delta(q, x) = \psi(q) + \alpha(x).$$

*Let $\delta_x$ for a fixed $x \in X$ be the map from $G \to G$; $q \to \delta(q, x)$; then $\{\delta_x \mid x \in X\}$ generates a sub-binear-ring, N(P) of the binear-ring. $\{M(G) = M(G_1) \cup M(G_2), +, \bullet\}$ of all mappings on $G = G_1 \cup G_2$, this binear-ring N(P) is called the Smarandache bisyntactic binear-ring (S-bisyntactic binear-ring).*

Just as in case of S-semigroups we see in case of S-bisemigroups also we have several bigroup, bisemi-automaton depending on the number of valid bigroups in the S-semigroup. So for a given S-S-bisemigroup bisemi-automaton we can have several S-bi syntactic binear-ring.

Now we see the probable applications of binear-rings to error correcting codes and their Smarandache analogue.

The study of how experiments can be organized systematically so that statistical analysis can be applied in an interesting problem which is carried out by several researchers. In the planning of experiments it often occurs that results are influenced by phenomena outside the control of the experimenter. The introduction of balanced incomplete block design (BIBD) helps in avoiding undesirable influences in the experiment. In general, if we have to test the effect of r different conditions with m possibilities for each condition this leads to a set of r orthogonal latin squares.

A planar near-ring can be used to construct balanced incomplete block designs (BIBD) of high efficiency; we just state how they are used in developing error correcting codes as codes which can correct errors is always desirable than the ones which can only detect errors. In view of this we give the following definition. Just we recall the definition of BIBD for the sake of completeness.

**DEFINITION 8.3.37**: *A balanced incomplete block design (BIBD) with parameters ($v$, $b$, $r$, $k$, $\lambda$) is a pair (P, B) with the following properties.*

    i. *P is a set with $v$ elements.*
    ii. *$B = (B_1, ..., B_b)$ is a subset of p (P) with b elements.*
    iii. *Each $B_i$ has exactly k elements where $k < v$ each unordered pair (p, q) with p, q $\in$ P, p $\neq$ q occurs in exactly $\lambda$ elements in B.*

*The set $B_1, ..., B_b$ are called the blocks of BIBD. Each $a \in P$ occurs in exactly r sets of B. Such a BIBD is also called a ($v$, $b$, $r$, $k$, $\lambda$) configuration or 2 – ($v$, $k$, $\lambda$) tactical configuration or design. The term balance indicates that each pair of elements occurs in exactly the same number of block, the term incomplete means that each block contains less than $v$ - elements. A BIBD is symmetric if $v = b$.*

*The incidence matrix of a ($v$, $b$, $r$, $k$, $\lambda$) configuration is the $v \times b$ matrix $A = (a_{ij})$ where*

$$a_{ij} = \begin{cases} 1 & \text{if } i \in B_j \\ 0 & \text{otherwise} \end{cases},$$

*here i denotes the $i^{th}$ element of the configuration. The following conditions are necessary for the existence of a BIBD with parameters $v$, $b$, $r$, $k$, $\lambda$.*

    1. *$bk = rv$.*
    2. *$r(k – 1) = \lambda (v – 1)$.*
    3. *$b \geq v$.*

*Recall a near-ring N is called planar (or Clay near-ring) if for all equation x o a = x o b + c. (a, b, c $\in$ N, a $\neq$ b) have exactly one solution x $\in$ N.*

***Example 8.3.16***: Let ($Z_5$, '+', $*$) where '+' is the usual '+' and '$*$' is n $*$ 0 = 0, n $*$ 1 = n $*$ 2 = n, n $*$ 3 = n $*$ 4 = 4n for all n $\in$ N. Then 1 $\equiv$ 2 and 3 $\equiv$ 4. N is planar near-ring; the equation x $*$ 2 = x $*$ 3 + 1 with 2 $\neq$ 3 has unique solution x = 3.

A planar near-ring can be used to construct BIBD of high efficiency where by high efficiency 'E' we mean E = $\lambda v / rk$; this E is a number between 0 and 1 and it estimates the quality of any statistical analysis if E $\geq$ 0.75 the quality is good.

Now, how does one construct a BIBD from a planar ring? This is done in the following way.

Let N be a planar ring. Let a ∈ N. Define $g_a$ : N → N by n → n o a and form G = {$g_a$ / a ∈ N}. Call a ∈ N "group forming" if a o N is a subgroup of (N, +). Let us call sets a . N + b (a ∈ N∗, b ∈ N) blocks. Then these blocks together with N as the set of "points" form a tactical configuration with parameters

$$(\upsilon, b, r, k, \lambda) = (\upsilon, \frac{\alpha_1 \upsilon}{|G|} + \alpha_2 \upsilon, \alpha_1 + \alpha_2 |G|, |G|, \lambda)$$

where $\upsilon$ = |N| and $\alpha_1$ ($\alpha_2$) denote the number of orbits of F under the group G \ {0} which consists of entirely of group forming elements. The tactical configuration is a BIBD if and only if either all elements are group forming or just 0 is group forming.

Now how are they used in obtaining error correcting codes. Now using the planar near-ring one can construct a BIBD. By taking either the rows or the columns of the incidence matrix of such a BIBD one can obtain error correcting codes with several nice features. Now instead of using a planar near-ring we can use Smarandache planar near-ring. The main advantage in using the Smarandache planar near-ring is that for one S-planar near-ring we can define more than one BIBD. If there are m-near-field in the S-planar near-ring we can build m BIBD. Thus this may help in even comparison of one BIBD with the another and also give the analysis of the common features.

Thus the S-planar near-ring has more advantages than the usual planar near-rings and hence BIBD's constructed using S-planar near-rings will even prove to be an efficient error correcting codes.

Now we give the binear-ring analog of it.

**DEFINITION 8.3.38:** *A binear-ring N = $N_1 \cup N_2$ is called biplanar if for all equations $x_i$ oa = $x_i$ ob + c (a, b, c ∈ $N_i$ a ≠ b) have exactly one solution $x_i \in N_i$, i =1, 2.*

Now we show how a biplanar binear-ring can be used to construct BIBD of high efficiency where by high efficiency E we mean E = λv/rk this E is a number between 0 and 1 and it estimates the quality of any statistical analysis if E > 0.75 the quality is good.

Here we give the construction of BIBD from a biplanar binear-ring.

Let N be a biplanar binear-ring. Let a ∈ N = $N_1 \cup N_2$. Define $g_a$ : $N_1$ → $N_1$ if a ∈ $N_1$ or $g_a$ : $N_2$→ $N_2$ if a ∈ $N_2$, by

$$n_1 \to n_1 \text{ o a if a} \in N_1$$
$$n_2 \to n_2 \text{ o a if a} \in N_2$$

and from G = $G_1 \cup G_2$ = {$g_a$ | a ∈ $N_1$} ∪ {$g_a$ | a ∈ $N_2$}. Call a ∈ N, bigroup forming if a o $N_1$ is a subgroup or a o $N_2$ is a sub group of ($N_1$, +) and ($N_2$, +) respectively. Let

us call sets $a.N_i + b$ ($a \in N_i$; $b \in N_i$) blocks $i = 1, 2$. Then these blocks together with N as the set of points form a tactical configuration with parameters $(v, b, r, k, \lambda) =$

$$\left( (v_1, \frac{\alpha_1 v_1}{|G_1|} + \alpha_2 v_1, \alpha_2 + \alpha_2|G_1|, |G_1|, \lambda), (v_2, \frac{\alpha_1 v_2}{|G_2|} + \alpha_2 v_2, \alpha_2 + \alpha_2|G_2|, |G_2|, \lambda) \right)$$

where $v_i = |N_i|$ and $\alpha_1$ denotes the number of orbits of F under the group $G_1 \setminus \{0\}$. $\alpha_2$ denotes the number of orbits of F under the group $G_2 \setminus \{0\}$ which consists of entirely of bigroup forming elements. The tactical configuration is a BIBD if and only if either all elements are bigroup forming or just 0 is bigroup forming.

The main advantage in using the Smarandache biplanar binear-ring is that for one S-biplanar binear-ring we define more than one BIBD. If there are m-binear field in the S-biplanar binear-ring we can build m BIBD.

Thus the S-biplanar binear-ring has more advantages than the usual planar binear-rings and hence BIBD's construed using S-biplanar binear-rings will even prove to be an efficient error correcting codes.

**PROBLEMS:**

1. Give an example of a ring loop.
2. Let $(N = Q \cup S(5), +, \bullet)$. Is N a ring semigroup?
3. Can $(N = Z_7 \cup L_7(3), +, \bullet)$ be a ring loop?
4. Find an example of a near-ring- group.
5. What is the structure of $N = Z^0 \cup L_7(3)$?
6. Give an example of a semigroup-group.
7. Give an example of a bi pseudo ring.
8. Is $(N, +, \bullet)$ where $N = (Q \cup S_3, +, \bullet)$ any well defined bi structure?
9. Give an example of a biquasi semiring.
10. When is a biquasi semiring S-biquasi semiring?
11. Illustrate by an example a biquasi near-ring.
12. Give an example of a S-strong loop semigroup.
13. Can a S-strong group groupoid exist? Justify your answer.
14. Build S-bisemi-automaton using S-bisemigroups.
15. Can we build S-cascades using S-algebraic structures?
16. Using $N = (ZS_3 \cup Q; +, \bullet)$ build biplanar binear-ring.
17. Construct S- biplanar binear-rings.

**Chapter 9**

# BISTRUCTURES, BIVECTOR SPACES AND THEIR SMARANDACHE ANALOGUE

This chapter introduces a new algebraic notion called bistructures. Stable and perfect bistructures are defined. The concept of bistructures is different from that of bigroups, bisemigroups etc. The study and introduction of bistructures and its Smarandache analogue is carried out in section one. In section two we introduce bivector spaces and give its Smarandache analogues.

## 9.1 Bistructures and S-Bistructures

In this book we studied two sets and defined operations and obtained bigroups, birings, biloops etc. In case of bistructures we do a different work. We take an arbitrary set and if we can define on it two binary operations such that the set is closed with respect to them then we call it a bistructure. We also give the Smarandache analogue.

**DEFINITION 9.1.1:** *Let X be a non-empty finite or infinite set. If we can define on X two separate distinct binary operations $*$ and $*'$ (such that the binary operators $*$ and $*'$ are closed with respect to X i.e. for every pair of elements $x, y \in X$ we have $x * y \in X$ and $x *' y \in X$). Then we say $\{(X, *), (X, *')\}$ is a bistructure.*

*If $*$ and $*'$ can be found on X such that with respect to both operations it is a same algebraic structure (i.e. groups or semigroups or groupoids,...) then we say $\{(X, *), (X, *')\}$ is a stable bistructure. If for a set no $*$ and $*'$ exists so that X can never be a bistructure then we call such sets adamant sets.*

**Example 9.1.1:** Let $X = \{0, 1, 2, \ldots\}$ be the set of positive integers with zero. Clearly $\{(X, +), (X, \times)\}$ is a stable bistructure as X is a commutative semigroup with identity under both '+' and '×'.

**Example 9.1.2:** $X = \{\text{set of positive rationals excluding zero}\}$. Consider the operation $\{(X, +), (X, \bullet)\}$. Clearly X is not stable bistructure but is a bisubstructure for we can find more than a pair of operations so that it is a bistucture. For $\{(X, \div), (X, \bullet)\}$ where '÷' is the usual division and '•' is the multiplication, X is a bistructure.

Also $\{(X, \div), (X, +)\}$, X is a bistructure. We say X is not stable because X has more than a pair of binary operation under which X is a bistructure. Thus we call these bistructures as complex bistructure.

So the example 9.1.2 is a complex bistructure. Our study is mainly classified under three broad heads:

      i.         Conditions on the set X so that X is a bistructure or not.

ii.   Stable bistructure or complex bistructure.
iii.  Conditions on X so that it is always an adamant set.

**THEOREM 9.1.1:** *The set of positive integers is a stable bistructure.*

*Proof:* Obvious by giving '+' and '×' operations as both $\{(X, +), (X, \times)\}$ are semigroups.

**THEOREM 9.1.2:** *Every stable bistructure is a bistructure.*

*Proof*: Straight-forward and left for the reader to prove.

**DEFINITION 9.1.2:** *Let X be a set if on X, pairs of binary operations can be defined so that (X, *, *'), (X, o, o') is an algebraic structure with two binary operations then X is said to be a strong bistructure. Let X be a set if X is a strong bistructure together with (X, *, *') and (X, o, o') which are same types of algebraic structure say rings or semirings or fields or near-rings then we call X a bistable strong bistructures.*

**THEOREM 9.1.3:** *Every two element set X = {a, b} is a bistructure.*

*Proof*: Define '*' on X by a * b = b, b * a = a, a * a = a and b * b = b then (X, *) is an algebraic structure. Define *' on X by a *' b = a, b *' a = b, a *' a = a and b *' b = b. Then (X, *') is an algebraic structure. So $\{(X, *), (X, *')\}$ is a bistructure.

Now we proceed onto define Smarandache bistructures.

**DEFINITION 9.1.3:** *Let X be a set which is non-empty finite or infinite. If we can find a proper subset $P \subset X$ such that on P we can define two separate binary operations * and *'. Then we say {(X, *), (X, *')} is a Smarandache bistructure (S-bistructure) as {(P, *), (P, *')} is a bistructure. Thus in case of S-bistructures we need not even have *' and * to be completely defined on X.*

*It is sufficient if it is defined on a proper subset. If for the set no such * and *' exists so that X can never be a S-bistructure we call such sets as Smarandache adamant sets (S-adamant sets).*

Thus we see in case of S-bistructures we can have several proper subsets to be bistructures under different operation. Thus this variety of study cannot be done except for Smarandache notions. It is left as an exercise for the reader to prove.

**THEOREM 9.1.4:** *Every S-bistructure in general need not be a bistructure.*

*Example 9.1.3:* Let

$$X = \left\{ \begin{pmatrix} 1 & 0 \\ 0 & 0 \end{pmatrix}, \begin{pmatrix} 0 & 0 \\ 0 & 0 \end{pmatrix}, 1, 0, 1+1=0 \right\}.$$

X is not a bistructure, but X is a S-bistructure for take $P = (0, 1) \subset X$. Clearly $(P, +, \bullet)$ is a bistructure, so X is a S-bistructure. Also take

$$Q = \left\{ \begin{pmatrix} 1 & 0 \\ 0 & 0 \end{pmatrix}, \begin{pmatrix} 0 & 0 \\ 0 & 0 \end{pmatrix} \right\}.$$

It is easily verified (Q, +, •) is a bistructure so X is a S-bistrucure. But X is not at all a bistructure under any operation.

Several interesting results in this direction can be made, a few are recalled here.

**DEFINITION 9.1.4:** *Let (X, +, •) be a bistructure. We say a proper subset P ⊂ X is a sub-bistructure if P itself under the operations of X is a bistructure.*

***Example 9.1.4:*** Let $Z^+$ be the set of integers, $(Z^+, +, \times)$ is a bistructure as $\{(Z^+, +), (Z^+, \times)\}$ are semigroups under '+' and '×'. Take $P = 2Z^+$ now, $\{(2Z^+, +), (2Z^+, \times)\}$ is a sub-bistructure of $(Z^+, +, •)$. Hence the claim.

Depending on the binary operation on the set we can define concepts like normal bi-substructure, ideals etc.

PROBLEMS:

1. Can $Q^+$ be made into a strong bistructure?
2. Can a bistructure be given on the set X = {-1, 1, 0, 3, 5, 7, 9}? Justify your claim.
3. Give an example of a S-bistructure.
4. What is the order of the smallest S-bistructure?
5. Give an example of an adamant set.
6. Illustrate by an example a strong bistructure.
7. Give an example of a bistructure which is not a S-bistructure.

## 9.2 Bivector spaces and S-bivector spaces

In this section we introduce the concept of bivector spaces and S-bivector spaces. The study of bivector spaces started only in 1999 [106]. Here we recall these definitions and extend it to the Smarandache bivector spaces.

**DEFINITION 9.2.1:** *Let $V = V_1 \cup V_2$ where $V_1$ and $V_2$ are two proper subsets of V and $V_1$ and $V_2$ are vector spaces over the same field F that is V is a bigroup, then we say V is a bivector space over the field F.*

*If one of $V_1$ or $V_2$ is of infinite dimension then so is V. If $V_1$ and $V_2$ are of finite dimension so is V; to be more precise if $V_1$ is of dimension n and $V_2$ is of dimension m then we define dimension of the bivector space $V = V_1 \cup V_2$ to be of dimension m + n. Thus there exists only m + n elements which are linearly independent and has the capacity to generate $V = V_1 \cup V_2$.*

*The important fact is that same dimensional bivector spaces are in general not isomorphic.*

***Example 9.2.1:*** Let $V = V_1 \cup V_2$ where $V_1$ and $V_2$ are vector spaces of dimension 4 and 5 respectively defined over rationals where $V_1 = \{(a_{ij})/ a_{ij} \in Q\}$, collection of all $2 \times 2$ matrices with entries from Q. $V_2 = \{$Polynomials of degree less than or equal to $4\}$.

Clearly V is a finite dimensional bivector space of dimension 9. In order to avoid confusion we always follow the following convention very strictly. If $v \in V = V_1 \cup V_2$ then $v \in V_1$ or $v \in V_2$ if $v \in V_1$ then v has a representation of the form $(x_1, x_2, x_3, x_4, 0, 0, 0, 0, 0)$ where $(x_1, x_2, x_3, x_4) \in V_1$ if $v \in V_2$ then $v = (0, 0, 0, 0, y_1, y_2, y_3, y_4, y_5)$ where $(y_1, y_2, y_3, y_4, y_5) \in V_2$.

Thus we follow the notation.

***Notation***: Let $V = V_1 \cup V_2$ be the bivector space over the field F with dimension of V to be m + n where dimension of $V_1$ is m and that of $V_2$ is n. If $v \in V = V_1 \cup V_2$, then $v \in V_1$ or $v \in V_2$ if $v \in V_1$ then $v = (x_1, x_2, …, x_m, 0, 0, …, 0)$ if $v \in V_2$ then $v = (0, 0, …, 0, y_1, y_2, …, y_n)$.

We never add elements of $V_1$ and $V_2$. We keep them separately as no operation may be possible among them. For in example we had $V_1$ to be the set of all $2 \times 2$ matrices with entries from Q where as $V_2$ is the collection of all polynomials of degree less than or equal to 4. So no relation among elements of $V_1$ and $V_2$ is possible. Thus we also show that two bivector spaces of same dimension need not be isomorphic by the following example:

***Example 9.2.2:*** Let $V = V_1 \cup V_2$ and $W = W_1 \cup W_2$ be any two bivector spaces over the field F. Let V be of dimension 8 where $V_1$ is a vector space of dimension 2, say $V_1 = F \times F$ and $V_2$ is a vector space of dimension 6 say all polynomials of degree less than or equal to 5 with coefficients from F. W be a bivector space of dimension 8 where $W_1$ is a vector space of dimension 3 i.e. $W_1 = \{$all polynomials of degree less than or equal to $2\}$ with coefficients from F and $W_2 = F \times F \times F \times F \times F$ a vector space of dimension 5 over F. Thus any vector in V is of the form $(x_1, x_2, 0, 0, 0, …, 0)$ or $(0, 0, y_1, y_2, …, y_6)$ and any vector in W is of the form $(x_1, x_2, x_3, 0, …, 0)$ or $(0, 0, 0, y_1, y_2, …, y_5)$. Hence no isomorphism can be sought between V and W in this set up.

This is one of the marked difference between the vector spaces and bivector spaces. Thus we have the following theorem, the proof of which is left for the reader to prove.

**THEOREM 9.2.1:** *Bivector spaces of same dimension defined over same fields need not in general be isomorphic.*

**THEOREM 9.2.2:** *Let $V = V_1 \cup V_2$ and $W = W_1 \cup W_2$ be any two bivector spaces of same dimension over the same field F. Then V and W are isomorphic as bivector spaces if and only if the vector space $V_1$ is isomorphic to $W_1$ and the vector space $V_2$ is isomorphic to $W_2$, that is dimension of $V_1$ is equal to dimension $W_1$ and the dimension of $V_2$ is equal to dimension $W_2$.*

*Proof:* Straightforward, hence left for the reader to prove.

**THEOREM 9.2.3:** *Let $V = V_1 \cup V_2$ be a bivector space over the field F. W any non empty set of V. $W = W_1 \cup W_2$ is a sub-bivector space of V if and only if $W \cap V_1 = W_1$ and $W \cap V_2 = W_2$ are sub spaces of $V_1$ and $V_2$ respectively.*

*Proof:* Direct; left for the reader to prove.

**DEFINITION 9.2.2:** *Let $V = V_1 \cup V_2$ and $W = W_1 \cup W_2$ be two bivector spaces defined over the field F of dimensions $p = m + n$ and $q = m_1 + n_1$ respectively.*

*We say the map $T: V \to W$ is a bilinear transformation of the bivector spaces if $T = T_1 \cup T_2$ where $T_1 : V_1 \to W_1$ and $T_2 : V_2 \to W_2$ are linear transformations from vector spaces $V_1$ to $W_1$ and $V_2$ to $W_2$ respectively satisfying the following two rules.*

  i. *$T_1$ is always a linear transformation of vector spaces whose first co ordinates are non-zero and $T_2$ is a linear transformation of the vector space whose last co ordinates are non zero.*

  ii. *$T = T_1 \cup T_2$ '$\cup$' is just only a notational convenience.*

  iii. *$T(v) = T_1(v)$ if $v \in V_1$ and $T(v) = T_2(v)$ if $v \in V_2$.*

*Yet another marked difference between bivector spaces and vector spaces are the associated matrix of an operator of bivector spaces which has $m_1 + n_1$ rows and $m + n$ columns where dimension of V is $m + n$ and dimension of W is $m_1 + n_1$ and T is a linear transformation from V to W. If A is the associated matrix of T then.*

$$A = \begin{bmatrix} B_{m_1 \times m} & O_{n_1 \times m} \\ O_{m_1 \times n} & C_{n_1 \times n} \end{bmatrix}$$

*where A is a $(m_1 + n_1) \times (m + n)$ matrix with $m_1 + n_1$ rows and $m + n$ columns. $B_{m_1 \times m}$ is the associated matrix of $T_1 : V_1 \to W_1$ and $C_{n_1 \times n}$ is the associated matrix of $T_2 : V_2 \to W_2$ and $O_{n_1 \times m}$ and $O_{m_1 \times n}$ are non zero matrices.*

***Example 9.2.3:*** Let $V = V_1 \cup V_2$ and $W = W_1 \cup W_2$ be two bivector spaces of dimension 7 and 5 respectively defined over the field F with dimension of $V_1 = 2$, dimension of $V_2 = 5$, dimension of $W_1 = 3$ and dimension of $W_2 = 2$. T be a linear transformation of bivector spaces V and W. The associated matrix of $T = T_1 \cup T_2$ where $T_1 : V_1 \to W_1$ and $T_2 : V_2 \to W_2$ given by

$$A = \begin{bmatrix} 1 & -1 & 2 & 0 & 0 & 0 & 0 & 0 \\ -1 & 3 & 0 & 0 & 0 & 0 & 0 & 0 \\ 0 & 0 & 0 & 2 & 0 & 1 & 0 & 0 \\ 0 & 0 & 0 & 3 & 3 & -1 & 2 & 1 \\ 0 & 0 & 0 & 1 & 0 & 1 & 1 & 2 \end{bmatrix}$$

where the matrix associated with $T_1$ is given by

$$\begin{bmatrix} 1 & -1 & 2 \\ -1 & 3 & 0 \end{bmatrix}$$

and that of $T_2$ is given by

$$\begin{bmatrix} 2 & 0 & 1 & 0 & 0 \\ 3 & 3 & -1 & 0 & 1 \\ 1 & 0 & 1 & 1 & 2 \end{bmatrix}$$

We call $T : V \to W$ a linear operator of both the bivector spaces if both V and W are of same dimension. So the matrix A associated with the linear operator T of the bivector spaces will be a square matrix. Further we demand that the spaces V and W to be only isomorphic bivector spaces. If we want to define eigen bivalues and eigen bivectors associated with T.

The eigen bivector values associated with are the eigen values associated with $T_1$ and $T_2$ separately. Similarly the eigen bivectors are that of the eigen vectors associated with $T_1$ and $T_2$ individually. Thus even if the dimension of the bivector spaces V and W are equal still we may not have eigen bivalues and eigen bivectors associated with them.

*Example 9.2.4:* Let T be a linear operator of the bivector spaces – V and W. $T = T_1 \cup T_2$ where $T_1 : V_1 \to W_1$ dim $V_1$ = dim $W_1$ = 3 and $T_2 : V_2 \to W_2$ where dim $V_2$ = dim $W_2$ = 4. The associated matrix of T is

$$A = \begin{bmatrix} 2 & 0 & -1 & 0 & 0 & 0 & 0 \\ 0 & 1 & 0 & 0 & 0 & 0 & 0 \\ -1 & 0 & 3 & 0 & 0 & 0 & 0 \\ 0 & 0 & 0 & 2 & -1 & 0 & 6 \\ 0 & 0 & 0 & -1 & 0 & 2 & 1 \\ 0 & 0 & 0 & 0 & 2 & -1 & 0 \\ 0 & 0 & 0 & 6 & 1 & 0 & 3 \end{bmatrix}$$

The eigen bivalues and eigen bivectors can be calculated.

**DEFINITION 9.2.3:** *Let T be a linear operator on a bivector space V. We say that T is diagonalizable if $T_1$ and $T_2$ are diagonalizable where $T = T_1 \cup T_2$.*

*The concept of symmetric operator is also obtained in the same way, we say the linear operator $T = T_1 \cup T_2$ on the bivector space $V = V_1 \cup V_2$ is symmetric if both $T_1$ and $T_2$ are symmetric.*

**DEFINITION 9.2.4:** *Let $V = V_1 \cup V_2$ be a bivector space over the field F. We say $\langle , \rangle$ is an inner product on V if $\langle , \rangle = \langle , \rangle_1 \cup \langle , \rangle_2$ where $\langle , \rangle_1$ and $\langle , \rangle_2$ are inner products on the vector spaces $V_1$ and $V_2$ respectively.*

*Note that in $\langle , \rangle = \langle , \rangle_1 \cup \langle , \rangle_2$ the '$\cup$' is just a conventional notation by default.*

**DEFINITION 9.2.5:** *Let $V = V_1 \cup V_2$ be a bivector space on which is defined an inner product $\langle , \rangle$. If $T = T_1 \cup T_2$ is a linear operator on the bivector spaces V we say $T^*$ is an adjoint of T if $\langle T\alpha / \beta \rangle = \langle \alpha / T^* \beta \rangle$ for all $\alpha, \beta \in V$ where $T^* = T_1^* \cup T_2^*$ are $T_1^*$ is the adjoint of $T_1$ and $T_2^*$ is the adjoint of $T_2$.*

The notion of normal and unitary operators on the bivector spaces are defined in an analogous way. T is a unitary operator on the bivector space $V = V_1 \cup V_2$ if and only if $T_1$ and $T_2$ are unitary operators on the vector space $V_1$ and $V_2$.

Similarly T is a normal operator on the bivector space if and only if $T_1$ and $T_2$ are normal operators on $V_1$ and $V_2$ respectively. We can extend all the notions on bivector spaces $V = V_1 \cup V_2$ once those properties are true on $V_1$ and $V_2$.

The primary decomposition theorem and spectral theorem are also true is case of bivector spaces. The only problem with bivector spaces is that even if the dimension of bivector spaces are the same and defined over the same field still they are not isomorphic in general.

Now we are interested in the collection of all linear transformation of the bivector spaces $V = V_1 \cup V_2$ to $W = W_1 \cup W_2$ where V and W are bivector spaces over the same field F.

We denote the collection of linear transformation by B-Hom$_F$(V, W).

**THEOREM 9.2.4:** *Let V and W be any two bivector spaces defined over F. Then B-Hom$_F$(V, W) is a bivector space over F.*

*Proof:* Given $V = V_1 \cup V_2$ and $W = W_1 \cup W_2$ be two bivector spaces defined over the field F. B-Hom$_F$(V, W) = $\{T_1 : V_1 \to W_1\} \cup \{T_2 : V_2 \to W_2\}$ = Hom$_F$($V_1$, $W_1$) $\cup$ Hom$_F$ ($V_2$, $W_2$). So clearly B- Hom$_F$(V,W) is a bivector space as Hom$_F$ ($V_1$, $W_1$) and Hom$_F$ ($V_2$, $W_2$) are vector spaces over F.

**THEOREM 9.2.5:** *Let $V = V_1 \cup V_2$ and $W = W_1 \cup W_2$ be two bivector spaces defined over F of dimension $m + n$ and $m_1 + n_1$ respectively. Then B-Hom$_F$(V,W) is of dimension $mm_1 + nn_1$.*

*Proof*: Obvious by the associated matrices of T.

Thus it is interesting to note unlike in other vector spaces the dimension of $\text{Hom}_F(V, W)$ is mn if dimension of the vector space V is m and that of the vector space W is n. But in case of bivector spaces of dimension $m + n$ and $m_1 + n_1$ the dimension of B-$\text{Hom}_F(V, W)$ is not $(m + n)(m_1 + n_1)$ but $mm_1 + nn_1$, which is yet another marked difference between vector spaces and bivector spaces.

Further even if bivector space V and W are of same dimension but not isomorphic we may not have B-$\text{Hom}_F(V,W)$ to be a bilinear algebra analogous to linear algebra. Thus B-$\text{Hom}_F(V,W)$ will be a bilinear algebra if and only if the bivector spaces V and W are isomorphic as bivector spaces.

Now we proceed on to define the concept of pseudo bivector spaces.

**DEFINITION 9.2.6:** *Let V be an additive group and $B = B_1 \cup B_2$ be a bifield if V is a vector space over both $B_1$ and $B_2$ then we call V a pseudo bivector space.*

*Example 9.2.5:* Let V = R the set of reals, $B = Q(\sqrt{3}) \cup Q(\sqrt{2})$ be the bifield. Clearly R is a pseudo bivector space over B. Also if we take $V_1 = R \times R \times R$ then $V_1$ is also a pseudo bivector space over B.

Now how to define dimension, basis etc of V, where V is a pseudo bivector space.

**DEFINITION 9.2.7:** *Let V be a pseudo bivector space over the bifield $F = F_1 \cup F_2$. A proper subset $P \subset V$ is said to be a pseudo sub-bivector space of V if P is a vector space over $F_1$ and P is a vector space over $F_2$ that is P is a pseudo vector space over F.*

*Example 9.2.6:* Let $V = R \times R \times R$ be a pseudo bivector space over $F = Q(\sqrt{3}) \cup Q(\sqrt{2})$. $P = R \times \{0\} \times \{0\}$ is a pseudo sub-bivector space of V as P is a pseudo bivector space over F.

Interested reader can develop notions in par with bivector spaces with some suitable modifications. Now we proceed on to define Smarandache bivector spaces and give some interesting results about them.

**DEFINITION 9.2.8:** *Let $A = A_1 \cup A_2$ be a K-bivector space. A proper subset X of A is said to be a Smarandache K-bivectorial space (S-K-bivectorial space) if X is a biset and $X = X_1 \cup X_2 \subset A_1 \cup A_2$ where each $X_i \subset A_i$ is S-K-vectorial space.*

**DEFINITION 9.2.9:** *Let A be a K-vectorial bispace. A proper sub-biset X of A is said to be a Smarandache-K-vectorial bi-subspace (S-K-vectorial bi-subspace) of A if X itself is a S-K-vectorial subspace.*

**DEFINITION 9.2.10:** *Let V be a finite dimensional bivector space over a field K. Let $B = B_1 \cup B_2 = \{(x_1,\ldots, x_k, 0 \ldots 0)\} \cup \{(0,0, \ldots, 0, y_1 \ldots y_n)\}$ be a basis of V. We say B is a Smarandache basis (S-basis) of V if B has a proper subset A, $A \subset B$ and $A \neq \phi$, $A \neq B$ such that A generates a bisubspace which is bilinear algebra over K; that is W* *is the sub-bispace generated by A then W must be a K-bi-algebra with the same operations of V.*

**THEOREM 9.2.6:** *Let A be a K bivectorial space. If A has a S-K-vectorial sub-bispace then A is a S-K.vectorial bispace.*

*Proof:* Straightforward by the very definition.

**THEOREM 9.2.7:** *Let V be a bivector space over the field K. If B is a S-basis of V then B is a basis of V.*

*Proof:* Left for the reader to verify.

**DEFINITION 9.2.11:** *Let V be a finite dimensional bivector space over a field K. Let B = $\{v_1, ..., v_n\}$ be a basis of V. If every proper subset of B generates a bilinear algebra over K then we call B a Smarandache strong basis (S-strong basis) for V.*

*Let V be any bivector space over the field K. We say L is a Smarandache finite dimensional bivector space (S-finite dimensional bivector space) over K if every S-basis has only finite number of elements in it.*

All results proved for bivector spaces can be extended by taking the bivector space V = $V_1 \cup V_2$ both $V_1$ and $V_2$ to be S-vector space. Once we have V = $V_1 \cup V_2$ to be a S-bivector space i.e. $V_1$ and $V_2$ are S-vector spaces, we see all properties studied for bivector spaces are easily extendable in case of S-bivector spaces with appropriate modifications.

**PROBLEMS:**

1. Give an example of a bivector space of dimension 3.
2. Illustrate by an example the set of basis for a bivector space V = $V_1 \cup V_2$ of dimension 5 over Q.
3. Show by an example that a bivector space of dimension n where V = $V_1 \cup V_2$ having infinitely many basis over Q but dimension is only n.
4. Let V = R[x] $\cup$ (R × R) be a bivector space over Q. Is V a bivector space over Q? Is V a finite dimensional bivector space? Justify your claim.
5. Let V = Q[x] $\cup$ R be a bivector space over R. Can V be a bivector space over Q? Justify your answer.
6. Is V = Q[x] $\cup$ R[x] a bivector space over R? Substantiate your claim.
7. State and prove spectral theorem for bivector space.
8. Illustrate spectral theorem by an example.
9. State and prove primary decomposition theorem for bivector spaces.
10. Substantiate the proof in problem 9, by a concrete example.
11. Define normal operator on bivector spaces. Illustrate it with examples.
12. Define bilinear algebra.
13. Give an example of a bilinear algebra of dimension 8.
14. Find B-Hom$_F$(V, W) where V = (R × R) $\cup$ R$^1$[x] over R, R$^1$[x] polynomial of degree less than or equal to 3. W = $W_1 \cup W_2$ where $W_1$ = {2 × 2 matrices with entries from R} and $W_2 = R \times R \times R \times R$. What is the dimension of B-$\text{Hom}_R(V,W)$?.

15. Give some interesting properties about pseudo bivector spaces.
16. Let $V = [(a_{ij})]_{n \times n}$ matrices with entries from Q. Is V a pseudo bivector space over the bifield $Q(\sqrt{7}) \times Q(19)$? Justify your claim.
17. Give an example of a S-bivector space.
18. Find a S-bivector space of dimension four.
19. What can be the least dimension of a S-bivector space?

# Chapter 10

# SUGGESTED PROBLEMS

In this section we suggest several problems on bialgebraic structures. As the topic is very new these problem will pave way for research in these notions. There are about 178 problems on these bistructure and their Smarandache analogue. Similar types of problems in every other bialgebraic structure can be coined and solved by the reader.

1. Find necessary and sufficient condition on $G_1$ and $G_2$ so that the order of the sub-bigroups always divides the order of the bigroup. (Here $G = G_1 \cup G_2$ where G is assumed to be a bigroup of finite order).

2. Find necessary and sufficient condition on $G = G_1 \cup G_2$ so that $N(a)$ ($a \in G = G_1 \cap G_2$) is a sub-bigroup of G.

3. Find necessary and sufficient condition on the finite bigroup G so that $C_a = \dfrac{o(G)}{o(N(a))}$ and $o(G) = \sum C_a$ where $a \in G_1$ or $a \in G_2$ or $a \in G_1 \cap G_2$. (Here $G = G_1 \cup G_2$).

4. When is $o(N(a)) = o(G)$ where $G = G_1 \cup G_2$.

5. Do all the classical results that hold good in case of cosets hold good in case of bicosets?

6. Can we say bicoset partitions the bigroups G for any sub-bigroup H of G?

7. Prove every biring in general need not have bi-ideals. Characterize those birings, which has no bi-ideals.

8. Can one prove every finite bidivision rings is a bifield?

9. Obtain some interesting results about group birings.

10. Let RG and SH be any two group birings; where $R = R_1 \cup R_2$ is a biring and $S = S_1 \cup S_2$ is a biring G and H are some groups. Obtain a necessary and sufficient condition so that the group birings are isomorphic; i.e. $RG \cong SH$.

11. Let $R = R_1 \cup R_2$ be the bifield of characteristic 0, if G is a torsion free non abelian group; can the group biring RG have zero divisors. (Solution to this problem will solve the 60 year old zero divisor conjecture for group rings).

12. Find a necessary condition for the biring $R = R_1 \cup R_2$ to be a semiprime biring.

13. Find a necessary and sufficient condition for the group biring to be embeddable in a bidivision ring.

14. Let G be a torsion free group. K any bifield. Is it true that the group biring KG has no proper divisors of zero? Are all units in KG necessarily trival?

15. Find the necessary and sufficient condition for the group biring KG = $K_1G \cup K_2G$ to satisfy a polynomial identity.

16. Introduce the concept of primitive ideals in a group biring and characterize those group birings which are primitive.

17. Find conditions on the group biring to be semisimple.

18. Obtain a necessary and sufficient condition on the group biring RG so that its Jacobson radical J(RG) = (0).

19. Characterize those semigroup birings RS so that

    i. J (RS) = 0.
    ii. J (RS) ≠ 0.

20. Find a sufficient condition for the semigroup biring to be semiprime.

21. Find a necessary and sufficient condition so that the semigroup biring satisfies a standard polynomial identity.

22. Obtain conditions for two semigroup birings RS and $R_1S_1$ where R is a biring and S a semigroup and $R_1$ a biring and $S_1$ a semigroup to be isomorphic.

23. When are the bigroup rings $K_1G$ and $K_2H$ isomorphic, where G = $G_1 \cup G_2$ and H = $H_1 \cup H_2$ are bigroups and $K_1$ and $K_2$ are rings?

24. Characterize those bigroup rings which are

    i. Prime.
    ii. Semiprime.

25. Characterize those bigroup rings, which are embeddable in a bidivision ring.

26. Characterize those bigroup rings, which are semisimple. (or) Find necessary and sufficient condition for the bigroup ring KG where K is a field or a ring and G a bigroup to be semisimple.

27. Find necessary and sufficient condition for the bisemigroup ring to be semisimple.

28. Characterize those bisemigroup rings which are prime.

29. Let KS be a bisemigroup ring; find conditions on the ring K and on the bisemigroup S so that the bisemigroup ring KS is semiprime.

30. Find conditions on the bisemigroup S and on the ring K so that the bi semigroup ring is embeddable in a bidivision ring.

31. Find a necessary and sufficient condition on the bigroup biring to be embeddable in a quad division ring.

32. Let RG and SH be two bigroup birings. Obtain conditions on the birings R and S and on the bigroups G and H so that RG and SH are isomorphic.

33. When is the bigroup biring semiprime?

34. Find conditions on the bigroup and biring so that the bigroup biring is prime.

35. Find conditions on bisemigroup biring to be embeddable in a biquad division ring.

36. When are two bisemigroup brings RS and KP isomorphic? Obtain necessary and sufficient condition for the same.

37. Obtain necessary and sufficient condition for the bisemigroup biring to be

    i. Semiprime.
    ii. Prime.

38. Find some interesting results on quad rings.

39. Does there exist a polynomial quad ring in which every quad ideal is principal?

40. Give an example of a polynomial quad ring in which every polynomial is reducible.

41. Does there exist a polynomial quad ring in which no polynomial is quasi reducible?

42. Let R be a non-associative biring i.e. $R = R_1 \cup R_2$. Is A(R) the associator subring a sub-biring? Find a necessary and sufficient condition for A(R) to be a sub-biring.

43. Let $R = R_1 \cup R_2$ be a non associative biring. Find necessary and sufficient condition on R so that the commutator subring of R is a sub-biring of R.

44. Let $R = R_1 \cup R_2$ be a biring. Find necessary and sufficient condition so that a left (right) quasi regular ideal of R is bi-ideal of R.

45. Obtain a necessary and sufficient condition so that the loop biring RL is embeddable;

    i. In a bidivision ring.
    ii. In a bidomian.

46. When is a loop biring RL semisimple?

47. Find necessary and sufficient condition on the loop L and on the biring $R = R_1 \cup R_2$ so that $J(RL) = W(RL)$.

48. Obtain any interesting result on loop birings.

49. Find the necessary and sufficient condition for the loop birings $R_1L_1$ and $R_2L_2$ to be isomorphic where $R_1$ and $R_2$ are two distinct birings and $L_1$ and $L_2$ are loops.

50. Characterize those loop birings which has every element to be quasi regular.

51. Find necessary and sufficient condition, so that the groupoid biring is embeddable in a non-associative bidivision ring.

52. Prove if G is any finite groupoid and $R = R_1 \cup R_2$ be a bifield of characteristic zero or ring of integers; in the groupoid biring RG we have $J(RG) \subseteq W(RG)$.

53. Find a necessary and sufficient condition for the groupoid biring RG to have $J(RG) = W(RG)$, where $R = R_1 \cup R_2$ is an associative biring.

54. Find a necessary and sufficient condition on the groupoid biring RG so that every element in $RG = R_1G \cup R_2G$ is quasi regular.

55. Does there exist a groupoid biring RG in which no element is quasi regular? ($RG = R_1G \cup R_2G$).

56. Let RG be a groupoid biring. Find the commutator sub-biring L'(RG) of RG. Find conditions on the groupoid G and on the biring $R_1 \cup R_2 = R$ so that in the groupoid biring $RG = R_1G \cup R_2G$ we have

    i. $L'(RG) = RG$.
    ii. $L'(RG) = \{0, 1\}$.
    iii. $L'(RG) \subset RG$.

57. Let $R = R_1 \cup R_2$ be a biring and G a groupoid. Find a necessary and sufficient condition on the groupoid biring $RG = R_1G \cup R_2G$ so that the associator sub-biring A(RG) of RG is such that

    i. $A(RG) = RG$.
    ii. $A(RG) = \{0, 1\}$.
    iii. $A(RG) \subset RG$.

58. Find a necessary and condition for the biloop biring RL to be embeddable in a quad division ring.

59. Find a necessary and sufficient for a biloop biring to be embeddable in a quad field.

60. Find necessary and sufficient condition on the biloop biring RL so that

    i. $J(RL) = (0)$.
    ii. $J(RL) \neq \phi$.

61. Find a necessary and sufficient condition on the biloop biring RL so that

    i. $W(RL) \supseteq J(RL)$.
    ii. $W(RL) = J(RL)$.

62. Find conditions on the biloop L and on the biring R so that the biloop biring RL is a

    i. Moufang quad ring.
    ii. Bol quad ring.
    iii. Bruck quad ring.
    iv. Alternative quad ring.

    Can we say if the biloop is a moufang biloop then the biloop biring will be a Moufang quad ring?

63. Obtain a necessary and sufficient condition on the biloop biring RL to be simple ($L = L_1 \cup L_2$ and $R = R_1 \cup R_2$).

64. Obtain a necessary and sufficient condition on the biloop L and on the biring R so that the biloop biring is semisimple.

65. Let RG be a bigroupoid biring of the bigroupoid $G = G_1 \cup G_2$ over the biring $R = R_1 \cup R_2$. Find a necessary and sufficient condition so that $J(RG) = 0$.

66. Find a necessary and sufficient condition on the bigroupoid biring RG so that $J(RG) = W(RG)$. (R is an associative biring i.e. $R = R_1 \cup R_2$ and $G = G_1 \cup G_2$ is a bigroupoid).

67. Let RG and R'G' be two bigroupoid birings of the bigroupoids G and G' over the birings R and R' respectively, i.e. the birings R and R' are associative commutative birings with unit such that $R = R_1 \cup R_2$ and $R' = R'_1 \cup R'_2$ and G and G' are bigroups such that $G = G_1 \cup G_2$ and $G' = G'_1 \cup G'_2$. The bigroupoid biring $RG = R_1G_1 \cup R_1G_2 \cup R_2G_1 \cup R_2G_2$ and $R'G' = R'_1G_1 \cup R'_2G_1 \cup R'_1G_2 \cup R'_2G_2$. Find necessary and sufficient conditions so that the bigroupoid birings RG and R'G' are isomorphic.

68. Find necessary and sufficient conditions on the bigroupoid $G = G_1 \cup G_2$ and the biring $R = R_1 \cup R_2$ so that in the bigroupoid biring RG we have $J(RG) \subseteq W(RG)$.

69. Find conditions on elements in RG; (where RG is a bigroupoid biring) to be (i) quasi regular in RG (ii) right quasi regular in RG (iii) left quasi regular in RG. (Here if G has elements x such that R is a biring of characteristic 0 then RG has quasi regular elements).

70. Let G be a bigroupoid, R any associative commutative biring with 1. Find a necessary and sufficient condition on the bigroupoid biring to be

    i. Bol quad ring.
    ii. Bruck quad ring.
    iii. Moufang quad ring.
    iv. Alternative quad ring.

71. Does there exists a non associative biring of prime order p?

72. Characterize those S-bigroups which has a unique S-maximal sub-bigroups.

73. Characterize those S-bigroups which has always a S-strong internal direct product.

74. Characterize those S-bigroups G which has only S-internal direct product and never a S-strong internal direct product.

75. Obtain any interesting new results on S-bigroups.

76. Characterize all biloops of finite order, which are Lagrange biloops.

77. Find a necessary and sufficient condition on a finite biloops L so that L is a Cauchy biloops.

78. Characterize those biloops which are p-Sylow biloops.

79. Find a necessary and sufficient condition on a biloop L so that L has no p-Sylow sub-biloops.

80. Find a necessary and sufficient condition on a biloop L so that if L is commutative so is its principal isotope.

81. Obtain conditions on the biloop L so that it is a G-biloop.

82. Obtain an analogous of the edge coloring of the graph $K_{2n}$ using exactly (2n-1) colours by associating with different representation of biloops of order 2n.

83. Characterize those biloops which has non trivial.

    i. Nucleus.
    ii. Moufang center.
    iii. Commutator.
    iv. Associator.

84. Find conditions on biloops L = L$_1$ ∪ L$_2$ so that its representations can match the edge coloring problem.

85. Find a necessary and sufficient condition for a biquasi group (X, +, •) to be a Lagrange bi quasi group.

86. Find a necessary and sufficient condition on the biquasi group (X, +, •) so that X has Cauchy elements.

87. Characterize those biquasi groups which is a Sylow biquasi group.

88. Determine some interesting properties about S-bisemigroups.

89. Does there exists a condition so that every

    i. Bisemiautomaton is a S-bisemiautomaton.
    ii. Biautomaton is a S-biautomaton.

90. Define S-cascades or any other finite machines using S-bisemigroups.

91. Characterize those biloops which are S-biloops.

92. Find some interesting properties of S-biloops.

93. Find conditions on the biloop (L, +, •) so that it has a unique P$_s$(L).

94. Find conditions on the biloop L so that P$_s$(L) = P(L).

95. Obtain a necessary and sufficient conditions on the biloop (L, +, •) so that every S-sub-biloop gives a P$_s$(L) which are distinct.

96. Characterize those biloops which are not S-pseudo associative biloop.

97. Can there be a biloop L in which the associator sub-biloop of L coincides with the S-associator sub-biloop?

98. Give atleast an example of a biloop L in which the commutator sub-biloop of L coincides with the S-commutator sub-biloop.

99. Can there exist a biloop L in which the moufang centre and the S-Moufang centre coincides?

100. Give an example of a biloop, which has Moufang centre but no S-Moufang centre.

101. Find the condition on the biloop (L, +, •) so that the centre coincides with its S-centre (i.e. S-centre is unique and it is equal to the centre of the biloop).

102. Characterize those biquasi loops which are

      i. Lagrange biquasi loop.
     ii. Weakly Lagrange biquasi loop.
    iii. Non Lagrange biquasi loop.

103. Find a necessary and sufficient condition for a biquasi loop to be a Sylow biquasi loop.

104. Characterize those biquasi loops which has

      i. No Cauchy element (other than biquasi loop of prime order).
     ii. Has Cauchy elements.

105. Does there exist a biquasi loop is which every element is a Cauchy element?

106. Find a necessary and sufficient condition for a bigroupoid to be a

      i. S-Lagrange bigroupoid.
     ii. S-p-Sylow bigroupoid.
    iii. S-non Lagrange bigroupoid.

107. Characterize those S-bigroupoids which are S-Cauchy bigroupoids.

108. Construct using Smarandache notions a S-bicascade.

109. Obtain some interesting results about

      i. S-bisemiautomaton.
     ii. S-biautomaton.

110. Characterize those group bisemirings, which are bidivision rings.

111. Obtain some innovative results on problems like isomorphism problem of group bisemirings. To be more precise if $(SG, +, \bullet)$ and $(RH, +, \bullet)$ be any two group bisemirings. When is $SG \cong RH$?

112. Study groupoid bisemirings, SG in general, and study all of its classical properties / problems like

      i. When SG embeddable in a na.semifield?
     ii. G and H be groupoids S and S' be bisemirings; when is $SG \cong S'H$?
    iii. Conditions for SG to be semisimple.
    iv. Can every SG be a prime or a semiprime bisemiring?

113. Study the same problem by replacing the groupoid by the loop L.

114. Does there exist semiprime groupoid bisemiring?

115. Does there exist semiprime loop bisemiring?

116. Find a necessary and sufficient condition for a na bisemiring to be

    i. S-Moufang bisemiring.
    ii. S-Bol bisemiring.
    iii. S-alternative bisemiring.
    iv. S-Jordan bisemiring.

117. Characterize those bisemirings which are

    i. S-strongly commutative.
    ii. S-commutative.
    iii. S-non commutative.

118. Study classical properties on the set of S-bi-ideals of a groupoid bisemiring to be a modular lattice.

119. Classify those bisemirings whose bi-ideals form a modular lattice.

120. Find suitable conditions so that the S-bi-ideals of a non associative bi semiring forms a modular lattice.

121. Characterize those bisemirings whose S-sub-bisemirings does not form a modular lattice.

122. Characterize strict bisemigroups.

123. Find a necessary and sufficient condition for strict bisemigroups to be S-strict bisemigroups.

124. Is it possible to find an upper bound for the number of linearly independent vectors in an bisemivector space?

125. Find a necessary and sufficient condition for a bisemivector space to be a S-bisemivector space.

126. Characterize those S-bisemivector spaces which has only unique basis.

127. Characterize those bi semivector space which has several basis.

128. Does there exist a S-bisemivector space which can have several basis and the number of elements in a basis is also a varying in number.

129. Give an example of a bi pseudo semivector spaces which has several basis.

130. Can there be a S-bipseudo semivector spaces which has a unique basis?

131. Characterize those na birings which has a unique S-commutator biring.

132. Characterize those birings in which the number of S-commutator sub-birings equal the number of S-sub-birings.

133. Obtain some relations between the Jacobson radical and the S-strongly Jacobson radical of a biring.

134. Can any method be formulated to find the number of S-Jacobson radicals of $R = R_1 \cup R_2$?

135. Characterize those group binear-rings which are binear domains.

136. Obtain a necessary and sufficient condition for semigroup binear-rings to be embeddable in a binear domain.

137. Find conditions for bigroup binear-rings to have

    i. Zero divisors.
    ii. Idempotents.
    iii. Ideals.

138. Obtain some interesting properties about quad near-rings.

139. Find conditions for semigroup binear-ring to be Artinian.

140. Define the concept of Marot binear-ring and characterize those binear-rings which are Marot binear-rings.

141. Characterize those group binear-ring for which the bimod (p.q) envelope is a

    i. bisemigroup under '•'.
    ii. bigroup under '•'.

142. Characterize those semigroup bi-rings for which the bimod (p.q) envelope is a

    i. bisemigroup under '•'.
    ii. bigroup under '•'.

143. Find conditions on the bigroup G so that for any near-ring $Z_p$, the mod p envelop $G^*$ is a bigroup under '•'.

144. Study the problem (143) in case of bisemigroup. When is $S^*$ not even closed under '•'?

145. Characterize the conditions under which the bisemigroup $G^*$ is a bisemigroup.

146. Characterize all quad near-rings built using bigroups and binear-rings which are embeddable in non associative quad near domain. (The reader can define quad near domain and quad near fields).

147. Characterize those quad near-rings built using biloops / bigroupoids and binear-rings which cannot be embedded in a quad near field or quad near domain.

148. Find the order of the smallest

    i. quad near-ring.
    ii. quad seminear-ring.
    iii. non associative quad near-ring.
    iv. non associative quad seminear-ring

149. Characterize all biloop near-rings, which are embeddable in a non associative binear domain.

150. Find a necessary and sufficient condition for a group binear-ring to be embeddable in a binear domain or a binear field.

151. Find all biloop near-rings which satisfy

    i. Moufang identity.
    ii. Bol identity.
    iii. Bruck identity.

152. Characterize those non associative quad near-rings which are right alternative and not left alternative.

153. Characterize those non associative quad seminear-rings which are alternative.

154. Find conditions on the quad near-ring so that it is a S-quad near-ring.

155. Obtain a suitable condition on a set X to make it only

    i. a bistructure.
    ii. a strong bistructure.
    iii. a stable bistructure.
    iv. S- bistructure.
    v. perfect bistructure.

156. On a given set of cardinality n, how many S-bistructures can be defined?

157. Obtain any interesting result on bistructures and S-bistructures.

158. Classify binear-rings which are not S-binear-rings.

159. State and prove some interesting results about bivector spaces of finite dimension.

160. Restate Gram-Schmidt orthogonalization process for bivector spaces.

161. Characterize those binear-rings which are right quasi reflexive.

162. Find a necessary and sufficient condition for a binear-ring to be strongly bicommutative.

163. Find condition so that a quad near-ring is a S-quad near-ring.

164. Characterize those binear-rings which has S equiprime bi-ideals.

165. Find a necessary and sufficient condition for a group binear-ring to be embeddable in a binear domain.

166. Find a necessary and sufficient condition for a S-group binear-ring to be a S-binear domain.

167. Obtain a necessary and sufficient condition for a S-semigroup near-ring to be a S-near domain.

168. Does there exist a S-loop near-rings which are non-associative near domains or near-fields?

169. Characterize those S-groupoid seminear-rings which has zero divisors.

170. Can a S-loop seminear-ring be free from zero divisors? Justify or substantiate your claim.

171. When does a group biring NG have $G^*$ its mod p envelope to be

    i. a group?
    ii. a semigroup?

172. For the S-semigroup biring NS find the mod p-envelope to have $S^*$ to be

    i. a group.
    ii. not even closed under an operation.

173. Define S-affine binear-ring in any other way. Obtain some results which are unique about S-affine binear-rings and S-planar binear-rings.

174. What is smallest dimension of a S-bivector space which has no unique basis.

175. State and prove primary decomposition theorem for S-bivector spaces.

176. State and prove Gram-Schmidt or-thogonalization process for S-bivector spaces.

177. Define and obtain some interesting properties about S-linear algebra.

178. Obtain S-diagonalization process for S-bivector spaces.

83. VASANTHA KANDASAMY, W. B., *A note on semigroup rings which are pre p-rings*, Kyungpook Math. J., **34**, 223-225 (1994).

84. VASANTHA KANDASAMY, W. B., *n-ideal near-rings*, J. Math. Res. Expo., **14**, 167-168 (1994).

85. VASANTHA KANDASAMY, W. B., *On strictly right loop rings*, J. Harbin Inst. Sci. Tech., **18**, 116-118 (1994).

86. VASANTHA KANDASAMY, W. B., *Semi idempotents in loop algebra*, J. Inst. Math. Comp. Sci., **7**, 241-244 (1994).

87. VASANTHA KANDASAMY, W. B., *Zero divisors in semi-loop near-rings*, Zeszyty Nauk. Politech. Rzeszowskiej Mat. Fiz., **15**, 79-82 (1994).

88. VASANTHA KANDASAMY, W. B., *A note on the modular loop ring of a finite loop*, Opscula Math., **15**, 109-112 (1995).

89. VASANTHA KANDASAMY, W. B., *Fuzzy subloops of some special loops*, Proc. 26$^{th}$ Iranian Math. Conf., 33-37 (1995).

90. VASANTHA KANDASAMY, W. B., *Group rings which satisfy super ore condition*, Vikram Math. J., **15**, 67-69 (1995).

91. VASANTHA KANDASAMY, W. B., *On the mod p-envelope of $S_n$*, Math. Edu., **29**, 171-173 (1995).

92. VASANTHA KANDASAMY, W. B., *The units of semigroup seminear-rings*, Opscula Math., **15**, 113-114 (1995).

93. VASANTHA KANDASAMY, W. B., *Zero divisors in loop-loop near-ring*, J. Inst. Math. Comp. Sci, **8**, 139-143 (1995).

94. VASANTHA KANDASAMY, W. B., *A note on group rings which are F-rings*, Acta Ciencia Indica, **22**, 251-252 (1996).

95. VASANTHA KANDASAMY, W. B., *I$^*$-rings*, Chinese Quat. J. of Math., **11**, 11-12 (1996).

96. VASANTHA KANDASAMY, W. B., *Idempotents and semi-idempotents in near-rings*, J. Sichuan Univ., **33**, 330-332 (1996).

97. VASANTHA KANDASAMY, W. B., *Idempotents in loop seminear-rings*, Ganit, J. Bangladesh Math. Soc., **16**, 35-39 (1996).

98. VASANTHA KANDASAMY, W. B., *On ordered groupoids and groupoid rings*, J. Math. Comp. Sci., **9**, 145-147 (1996).

# INDEX

# About the Author

Dr. W. B. Vasantha is an Associate Professor in the Department of Mathematics, Indian Institute of Technology Madras, Chennai, where she lives with her husband Dr. K. Kandasamy and daughters Meena and Kama. Her current interests include Smarandache algebraic structures, fuzzy theory, coding/ communication theory. In the past decade she has completed guidance of seven Ph.D. scholars in the different fields of non-associative algebras, algebraic coding theory, transportation theory, fuzzy groups, and applications of fuzzy theory to the problems faced in chemical industries and cement industries. Currently, six Ph.D. scholars are working under her guidance. She has to her credit 241 research papers of which 200 are individually authored. Apart from this she and her students have presented around 262 papers in national and international conferences. She teaches both undergraduate and postgraduate students at IIT and has guided 41 M.Sc. and M.Tech projects. She has worked in collaboration projects with the Indian Space Research Organization and with the Tamil Nadu State AIDS Control Society. She is currently authoring a ten book series on Smarandache Algebraic Structures in collaboration with the American Research Press.

She can be contacted at vasantha@iitm.ac.in
You can visit her on the web at: http://mat.iitm.ac.in/~wbv